\pdfoutput=1
\documentclass[a4paper,reqno]{amsart}
\usepackage[paperheight=11in,paperwidth=8.5in]{geometry}
\usepackage{amsthm,amssymb,amsmath,hyperref}
\usepackage{color, bbm}

\usepackage{tikz}
\usetikzlibrary{decorations.pathreplacing}
\usepackage[utf8]{inputenc}
\usepackage{graphicx,enumerate}
\usepackage{mathtools,bm}
\usepackage{enumitem}

\usepackage{amssymb, amsmath, amsthm, chngcntr, enumerate, mathrsfs, mathtools, pgf, tikz,esint}
\usepackage{graphicx}
\usepackage{dsfont}
\usetikzlibrary{arrows,calc,patterns}
\usepackage[english]{babel}
\usepackage[T1]{fontenc}
\usepackage{hyperref}
\hypersetup{
    colorlinks=true,
    linkcolor=blue,
    filecolor=magenta,
    urlcolor=cyan,
}

\usepackage[normalem]{ulem}
\usepackage{soul}

\numberwithin{equation}{section}
\newtheorem{theorem}{Theorem}[section]

\newtheorem*{theorem*}{Theorem}
\newtheorem{corollary}[theorem]{Corollary}
\newtheorem{observation}[theorem]{Observation}
\newtheorem{obs}[theorem]{Observation}
\newtheorem{lemma}[theorem]{Lemma}
\newtheorem{lem}[theorem]{Lemma}
\newtheorem{proposition}[theorem]{Proposition}

\theoremstyle{definition}

\theoremstyle{definition}
\newtheorem{remark}[theorem]{Remark}

\newtheorem{notation}[theorem]{Notation}
\newtheorem*{notation*}{Notation}

\newtheorem{definition}[theorem]{Definition}
\newtheorem{definition lemma}[theorem]{Definition/Lemma}

\usepackage{mathtools}

\newcommand{\inc}{\mathfrak{i}}
\newcommand{\C}{\mathcal{C}}
\newcommand{\R}{\mathbb{R}}
\newcommand{\Z}{\mathbb{Z}}
\newcommand{\Q}{\mathbb{Q}}
\newcommand{\N}{\mathbb{N}}
\renewcommand{\P}{\mathbb{P}}

\newcommand{\FF}{\mathcal{F}}

\newcommand{\size}{\mathrm{Size\, }}

\newcommand{\energy}{\mathrm{Energy \,}}
\newcommand{\ds}{\displaystyle}
\newcommand{\one}{\mathbf{1}}
\newcommand{\ci}{\tilde{\chi}}
\newcommand{\dist}{\text{dist}}

\newcommand{\ii}{\mathscr}
\newcommand{\mb}{\mathbf}
\newcommand{\tr}{\mathbf}

\newcommand{\THT}{\ensuremath{{\textrm{HTHT}}}}
\newcommand{\BHT}{\ensuremath{\textrm{BHT}}}

\usepackage[framemethod=tikz]{mdframed}
\mdfsetup{linecolor=blue,backgroundcolor=gray!10,roundcorner=10pt,leftmargin=2pt,innerleftmargin=2pt,innerrightmargin=2pt}

\def\beq{\begin{equation}}
\def\eeq{\end{equation}}

\def\beq{\begin{equation}}
\def\eeq{\end{equation}}

\def\l{\lambda}
\def\r{\rho}
\def\L{\Lambda}

\def\r{\mathcal{R}}
\def\g{\gamma}
\def\G{\Gamma}
\def\a{\alpha}
\def\d{\delta}

\def\v{\varphi}
\def\ep{\epsilon}

\def\l{\Lambda}

\def\R{\mathbb{R}}

\def\C{\mathbb{C}}
\def\Z{\mathbb{Z}}

\def\P{\mathbb{P}}

\def\I{\mathcal{I}}

\def\N{\mathbb{N}}

\def\Z{\mathbb{Z}}
\def\beq{\begin{equation}}
\def\eeq{\end{equation}}

\def\beq{\begin{equation}}
\def\eeq{\end{equation}}

\def\r{\rho}
\def\g{\gamma}
\def\G{\Gamma}
\def\a{\alpha}

\def\d{\delta}

\def\ep{\epsilon}

\def\l{\lambda}

\def\pp{P}

\def\R{\mathbb{R}}
\def\C{\mathbb{C}}

\def\Z{\mathbb{Z}}

\def\P{\mathbb{P}}
\def\PP{\mathbb{P}}

\def\N{\mathbb{N}}

\def\Z{\mathbb{Z}}
\def\I{\mathcal{I}}
\def\beq{\begin{equation}}
\def\eeq{\end{equation}}

\def\beq{\begin{equation}}
\def\eeq{\end{equation}}

 \usepackage{todonotes}
\usepackage{xcolor}

\usepackage{bm}

\newcommand{\diam}{\text{ diam}}
\newcommand{\rr}{\mathbb}
\newcommand{\ic}{\mathcal}

\newcommand{\TFC}{\textnormal{{\bf TFC}}}
\newcommand{\sTFC}{\textnormal{{$\mathsf{TFC}$}}}

\newcommand{\bbS}{{\mathbb S}}
\newcommand{\bbA}{{\mathbb A}}
\newcommand{\bbB}{{\mathbb B}}

\newcommand{\bbU}{{\mathbb U}}

\newcommand{\calF}{{\mathcal F}}

\newcommand{\calR}{{\mathcal R}}
\newcommand{\calS}{{\mathcal S}}

\def\el{\mathfrak{l}}
\def\r{\mathfrak{r}}

\newcommand{\supp}{{\textrm{supp\,}}}

\newcommand{\low}{{lo}}
\newcommand{\nonsta}{{hinst}}
\newcommand{\sta}{{hist}}

\makeatletter
\newcommand{\vast}{\bBigg@{4}}
\newcommand{\Vast}{\bBigg@{5}}
\makeatother

\def\Xint#1{\mathchoice
   {\XXint\displaystyle\textstyle{#1}}%
   {\XXint\textstyle\scriptstyle{#1}}%
   {\XXint\scriptstyle\scriptscriptstyle{#1}}%
   {\XXint\scriptscriptstyle\scriptscriptstyle{#1}}%
   \!\int}
\def\XXint#1#2#3{{\setbox0=\hbox{$#1{#2#3}{\int}$}
     \vcenter{\hbox{$#2#3$}}\kern-.5\wd0}}

\def\aver#1{\Xint-_{#1}}

\counterwithout{equation}{section}
\counterwithout{section}{section}

\begin{document}

\title[On a hybrid trilinear Hilbert transform]{On a hybrid trilinear Hilbert transform}

\author[C. Benea]{Cristina Benea}
\address[C. Benea]{Nantes Universit\'e, Laboratoire de Math\'ematiques Jean Leray, Nantes 44322, France}
\email{cristina.benea@univ-nantes.fr}
\author[F. Bernicot]{Fr\'ed\'eric Bernicot}
\address[F. Bernicot]{CNRS - Nantes Universit\'e, Laboratoire de Math\'ematiques Jean Leray, Nantes 44322, France}
\email{frederic.bernicot@univ-nantes.fr}
\author[V. Lie]{Victor Lie}
\address[V. Lie]{Department of Mathematics, Purdue University, IN 47907 USA \& Institute of Mathematics of the Romanian Academy, Bucharest, RO 70700, P.O. Box 1-764, Romania.}
\email{vlie@purdue.edu}


\date{\today}

\begin{abstract}

In this paper we study the boundedness properties of a \emph{hybrid} trilinear Hilbert transform characterized by a \emph{dual regime} which  simultaneously exhibits
 \smallskip
\begin{itemize}
\item \emph{resonant} (zero curvature) features, \emph{i.e.}, modulation invariance properties, and

\item \emph{non-resonant} (non-zero curvature) features.
\end{itemize}
\bigskip

Specifically, given $a\in(0,\infty)\setminus \{1,2\}$ and $\vec{\lambda}:=(\lambda_1,\lambda_2,\lambda_3,\lambda_4)\in (\R\setminus\{0\})^4$ with  $\lambda_1\neq \lambda_2$ and $\lambda_3\neq \lambda_4$, we provide wide-range $L^p$ bounds for the trilinear operator
$$ \THT^a_{\vec{\lambda}} (f_1,f_2,f_3)(x):= p.v. \int_{\mathbb R} f_1(x-\lambda_1 t)\,f_2(x-\lambda_2 t-\lambda_3 t^a)\,f_3(x-\lambda_4 t^a)\,\frac{dt}{t},\quad x \in \R\,.$$
\smallskip

The analysis of this operator reveals a deep interplay between anisotropy, multiscale analysis, and modulation invariance and involves the synthesis of two components:
\smallskip
\begin{itemize}
\item a \emph{low-oscillatory component} which encapsulates an anisotropic Biest operator: this requires a novel approach and is of independent interest;

\item a \emph{high-oscillatory component} which governs the behavior of a hybrid Bilinear Hilbert-Carleson type operator: this relies on a subtle synergy between two properly defined resolutions
    \begin{itemize}
    \item a low-resolution multiscale analysis addressing the cancellation \emph{among} collections of super-Heisenberg quad-tiles; 
    
    \item a high-resolution single scale analysis focusing on the hidden cancellation \emph{within} each super-Heisenberg quad-tile.
    \end{itemize}
\end{itemize}
\medskip

Thus, by combining variable multiresolution, multiscale wave-packet analysis, and Rank II LGC method, our present work establishes a common framework where the diverse methodologies of all of the following naturally converge:
\begin{itemize}
\item the maximal truncated bilinear Hilbert transform, see \cite{Lacey} and \cite{DemeterTaoThiele};

\item the Biest operator studied in \cite{MTTBiest1,MTTBiest2};

\item the hybrid bilinear Hilbert-Carleson operator in \cite{BCa};

\item the non-resonant (curved) trilinear Hilbert transform in \cite{HuLie}.
\end{itemize}
\medskip

Ultimately, this work reveals that the hybrid trilinear transform sits on a critical intermediate scale between the purely non-resonant case treated in \cite{HuLie} and the purely resonant regime---the latter remaining a challenging, major (open) problem in harmonic analysis.

\end{abstract}

\maketitle

\setcounter{tocdepth}{2}
\tableofcontents

\section{Introduction}\label{Intr}

In this paper  we study the $L^p$-boundedness properties of the \emph{hybrid tri-linear Hilbert transform}\footnote{A priori one may consider in all the defining formulas below  $\{f_j\}_{j=1}^3\subset\calS(\R)$. Also, when $t<0$ and $a\not\in\N$ we let $t^{a}:=-\,|t|^{a}$ (same results hold if one chooses instead the convention $t^{a}:=|t|^{a}$).} 
\begin{equation}\label{hybtri}
\THT^a_{\vec{\lambda}} (f_1,f_2,f_3)(x):= p.v. \int_{\mathbb R} f_1(x-\lambda_1 t)\,f_2(x-\lambda_2 t-\lambda_3 t^a)\,f_3(x-\lambda_4 t^a)\,\frac{dt}{t},\quad x \in \R\,,
\end{equation}
where here $a\in(0,\infty)\setminus \{1,2\}$ and $\vec{\lambda}:=(\lambda_1,\lambda_2,\lambda_3,\lambda_4)\in (\R\setminus\{0\})^4$ with  $\lambda_1\neq \lambda_2$ and $\lambda_3\neq \lambda_4$.

The \emph{hybrid} label above captures the two key properties which coexist and describe the nature of \eqref{hybtri}:
\begin{itemize}
\item \textbf{(P1)} The \emph{zero curvature/resonant} feature portrayed by the \emph{modulation invariance} relation
\begin{equation}\label{modinv} 
 |\THT^a_{\vec{\lambda}} (M_{\frac{c}{\lambda_1}}f_1,M_{-\frac{c}{\lambda_2}}f_2, M_{\frac{c \lambda_3}{\lambda_2\lambda_4}}f_3)|=|\THT^a_{\vec{\lambda}} (f_1,f_2,f_3)|\qquad \forall\:c\in\R\,;
\end{equation}

\item \textbf{(P2)} The \emph{nonzero curvature/nonresonant} feature encapsulated within the interaction between the arguments $t$ and $t^a$ of the input functions in  \eqref{hybtri}. 
\end{itemize}

The \emph{motivation} for the present topic is multifaceted and covers the following:

\begin{itemize}
\item Firstly, at the generic and same time conceptual level, this work serves to and is part of the complex program initiated by the authors\footnote{In various collaborative team subsets.} and dedicated to the study of the interplay between zero and non-zero curvature in harmonic analysis. While diverse in aims and realization, this program includes more distantly related topics such as the theme concerning the non-resonant (\cite{HsuL24}) and resonant (\cite{HsuL26}) Carleson-Radon transform, but also directly related predecessors, such as
\begin{itemize}    
\item  the hybrid bilinear Hilbert transform in \cite{GL22};  

\item  the oscillatory Biest operator in \cite{BO};

\item  the non-resonant (\cite{BHL25}) and hybrid (\cite{BCa}) bilinear Hilbert-Carleson operators;

\item  the non-resonant trilinear Hilbert transform in \cite{HuLie}, \cite{HL25}.
    
\end{itemize}

\item Secondly, at a more specific level, this work serves as a natural stepping stone---with the attached intermediate level of complexity---between the purely curved (non-resonant) trilinear Hilbert transform solved in \cite{HuLie} and the purely flat (resonant) trilinear Hilbert transform  which still remains a major open problem in harmonic analysis.\footnote{See Remark \ref{openai1} below.}   

\item Thirdly, there are several interesting connections that one may establish between the main topics of interest herein and other related areas such as ergodic theory,  additive combinatorics and number theory---all these will be succinctly touched upon in Section \ref{Motivation}.
\end{itemize}

Finally, we mention that the relevance of our present study is supported not only by the above motivational component but also by the variety and novelty of the methods required in order to tackle our main results-see Section \ref{Keyideas}.  We conclude our brief introductory discussion with two comments concerning the interference of AI technology with our present work:

\begin{remark}[\textsf{On the OpenAI Preprint Repository}]\label{openai1} After the finalization of the present paper, on October 6, 2026, OpenAI released a repository of over 700 manuscripts 
proposing solutions to numerous long-standing open problems in mathematics. Among these preprints is a claimed proof of the 
boundedness of the resonant trilinear Hilbert transform from 
$L^3(\mathbb{R}) \times L^3(\mathbb{R}) \times L^3(\mathbb{R})$ into $L^1(\mathbb{R})$, 
made available in the official release \href{https://github.com/openai/math/blob/main/preprints/An-L3-bound-for-the-trilinear-Hilbert-transform-October-5-2026/paper.pdf}{repository}. 
However, the proposed proof is written in an idiosyncratic and 
non-standard mathematical exposition, and unlike other selected results 
in the release, it still lacks an official, computer-checked Lean formalization. Furthermore, 
early independent audits using automated error-detection tools like OSVerify 
have flagged potential structural gaps within the manuscript. Pending 
a rigorous and thorough verification by the mathematical community, we shall 
continue to treat the boundedness of this transform as still open.
\end{remark}

\begin{remark}[\textsf{OpenAI and a potential tightening of our main results}]\label{openai2}  \emph{No AI tools were used in the conception of this research or in the drafting of the present manuscript, in keeping with a principle agreed upon by the authors for the current project.} After the work was finalized, one of us checked our theorems with GPT-6 Astra, and its output indicates that the boundedness ranges in our results can likely be extended. We have deliberately not incorporated these machine-generated improvements---see the \emph{AI disclosure} declaration at the end of the paper---so that this paper reflects \emph{exclusively} the authors’ own contribution.
\end{remark}

\subsection{Main results}\label{MR}
\medskip

In this section we focus on providing the necessary context for and the statements of our main results. 

\subsubsection{\underline{The global behavior: the boundedness of the hybrid trilinear Hilbert transform}}

Let us denote by $\bm D$---see Figure 1 below---the open interior of the convex envelope spanned by the following seven extremal points:
\begin{gather*}
A:=(0,\frac{3}{4},\frac{1}{4}),\ B:=(0,\frac{1}{2},\frac{1}{2}),\ C:=(\frac{1}{2},\frac{1}{2} ,0),\ D:=(\frac{1}{2},0,\frac{1}{2}),\\ 
E:=(0,\frac{3}{4},0),\ F:=(0,\frac{1}{2},0),\ G:=(\frac{1}{4},\frac{3}{4},0).
\end{gather*}

\begin{figure}[h]
\includegraphics[scale=0.35]{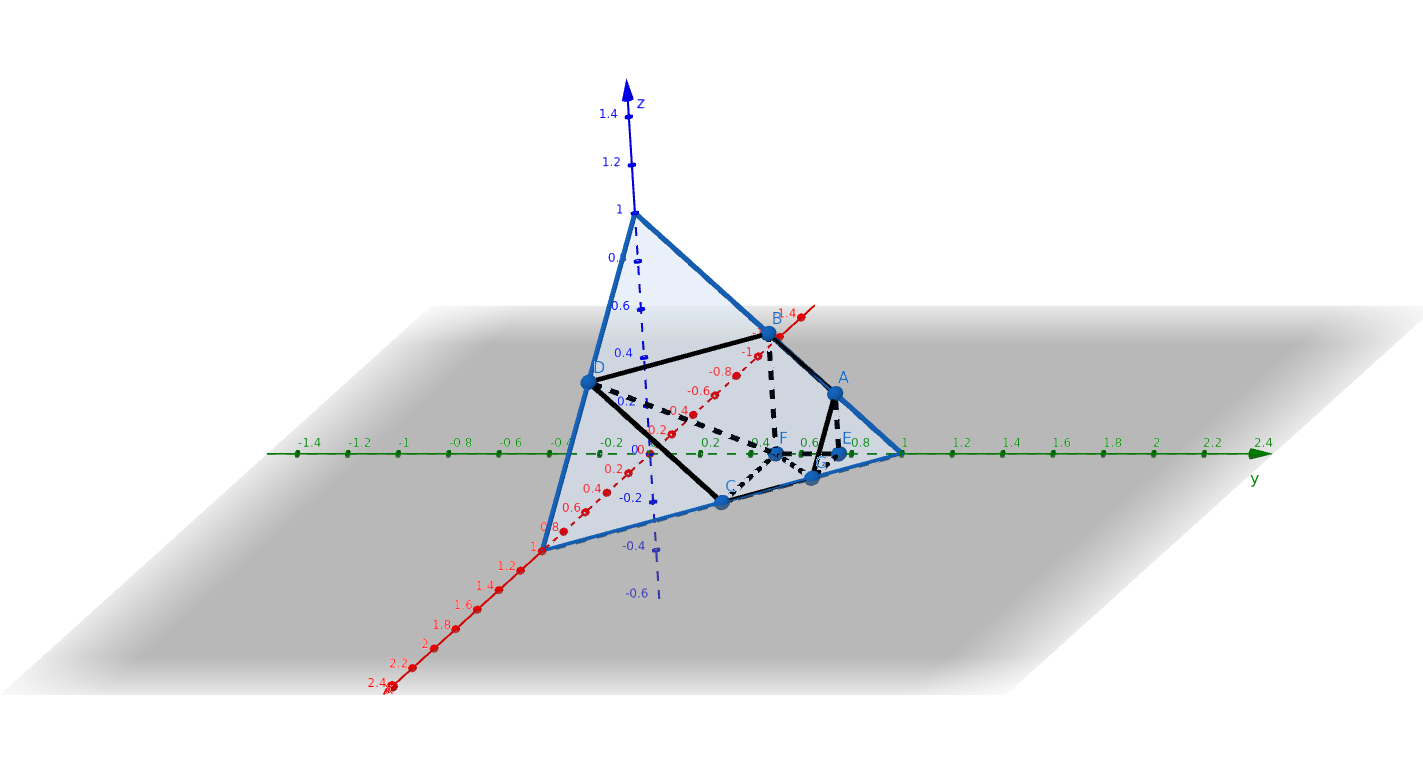}
\caption{The polytope $ \bm D$ is represented in black while its seven extremal points are colored in blue. We observe that five of these points belong to the subplane (in blue) $\{(\alpha_1,\alpha_2,\alpha_3),\ \alpha_i\in(0,1),\ \alpha_1+\alpha_2+\alpha_3=1\}$ arising from the intersection with the region corresponding to the  Banach condition. Indeed, $\bm D$ is obtained as
${\bm D}:=\{\alpha\in {\bm D}' \quad \textrm{and} \quad 0<\alpha_1+\alpha_2+\alpha_3<1\}$ where here ${\bm D}'$ is a suitable wider domain represented in the next figure.}
\end{figure}
Then the following holds:

\begin{theorem} \label{thm:main}  Define the prototype of the hybrid trilinear Hilbert transform in \eqref{hybtri} as
\begin{equation}\label{hththyb}
\THT (f_1,f_2,f_3)(x):= p.v. \int_{\mathbb R} f_1(x-t)\,f_2(x+t+t^3)\,f_3(x-t^3)\,\frac{dt}{t}.
\end{equation}

Let $p_1,p_2,p_3\in (1,\infty)$ be exponents with $(\frac{1}{p_1},\frac{1}{p_2},\frac{1}{p_3})\in {\bm D}$ and set $p\in(1,\infty)$ with
\begin{equation} \frac{1}{p}=\frac{1}{p_1}+\frac{1}{p_2}+\frac{1}{p_3}.\label{eq:holder}
\end{equation}
Then the operator $\THT$ is bounded from $L^{p_1}\times L^{p_2} \times L^{p_3} \rightarrow L^p$.\footnote{Notice in particular that $\THT$ is bounded from $L^{4} \times L^{3} \times L^{4}$ into $L^{6 \over 5}$ (since $(\frac{1}{4}, \frac{1}{3}, \frac{1}{4}) = \frac{1}{2}D + \frac{1}{3} E + \frac{1}{6} F$).} 
\end{theorem}

\medskip
\begin{corollary}\label{extens}   As stated in the introductory discussion, the main result above extends---with the same bounds---to the more general operator defined by \eqref{hybtri}. 
\end{corollary}

\begin{remark} Indeed, the analysis developed in the present paper is only dependent on two key properties:
\begin{itemize}
\item for the \emph{high-frequency component}: the associated bilinear Hilbert-Carleson operator $BC^a$---see the discussion in Sections \ref{Keyideas} and \ref{Motivation}---has to satisfy the hybrid regime $a\in(0,\infty)\setminus \{1,2\}$ treated in \cite{BCa}, thus ruling out the presence of the extra resonances implied by the condition $a\in\{1,2\}$; 

\item for the \emph{low-frequency component}:  the structure of the associated multiplier has to preserve a Biest-like behavior, which, in particular, demands that the $t-$argument of the function  $f_2$, i.e. $\lambda_2 t+\lambda_3 t^a$ to be a non-degenerate combination of the $t-$arguments corresponding to the functions $f_1$ and $f_3$, i.e. $\lambda_1 t$ and $\lambda_4 t^a$. This observation requires a symmetry among the input arguments while simultaneously matching the desired trilinear Hilbert hybrid regime and thus justifies the specific form chosen in \eqref{hybtri}. Consistent with this latter remark, we mention for example that  the 'degenerate' case of $\lambda_3=0$ is not covered by our current analysis.
\end{itemize}
\end{remark}

\subsubsection{\underline{The low frequency behavior: context and preparatives}}   As it turns out, in the main endeavor of understanding the boundedness properties of the hybrid trilinear Hilbert transform one encounters some completely new manifestations of the hybrid behavior which are of independent interest. Indeed, following the expected analysis, the first step is to decompose  $\THT$ into two main blocks: a low frequency and a high frequency component according to the size of the phase of the associated multiplier. While each of these components comes with new challenges, in the present section we focus on the low frequency component: due to the interaction between the modulation invariance and curvature (see \textbf{(P1)} and \textbf{(P2)}), the low frequency component is reshaped into a new class of trilinear anisotropic operators. Thus, an intermediate and necessary step in understanding $\THT$, is represented by the study of some suitable anisotropic trilinear operators involving a double $\BHT$ singularity that can be thought of as an anisotropic version of the Biest operator (see \cite{MTTBiest1,MTTBiest2}).

With these being said, we shift our story from descriptive to concrete, as follows: 

\begin{definition}\label{def:psi-type} Let $\phi\in C_0^{\infty}(\R)$ with $\supp\,\phi\subset[-C, C]$ for some suitable $C>0$.\footnote{The absolute constant $C>0$ is allowed to change from line to line.} We say that 
\begin{itemize}
\item $\phi$ is of \emph{$\Phi$-type} iff  $0\in\supp\,\phi$;

\item $\phi$ is of \emph{$\Psi$-type} iff $\phi(0)=0$. Often, we will tacitly\footnote{This is justified, since elementary reasonings give that any function of \emph{$\Psi$-type} satisfying the initial (weaker) requirement $\phi(0)=0$,  may be decomposed into (a possibly infinite) linear combination of \emph{$\Psi$-type} functions verifying the stronger condition \eqref{hththyb}, with this linear combination discretization having uniform bounded overlapping support property and $l^{\infty}(\Z)$ coefficients.} include in the \emph{$\Psi$-type} definition the stronger, more convenient requirement 
\begin{equation}\label{hththyb}
0\not\in \supp\,\phi\qquad\textrm{and}\qquad \ds \dist(0, \supp (\phi)) \sim \diam (\supp (\phi))\,.
\end{equation}
\end{itemize}

In what follows, we will usually use the notation $\psi$ for $\Psi$-type functions and $\phi$ for $\Phi$-type functions, although these properties should be clear from the context. 
\end{definition}

Next, given $m:\R^3 \to \C$ a suitable symbol, we denote its associated trilinear operator by
\begin{equation} \label{eq:op:anisotropic}
T_{m}(f_1,f_2,f_3)(x) := \iiint_{{\mathbb R}^3}\hat{f_1}(\xi_1) \hat{f_2}(\xi_2) \hat{f_3}(\xi_3) m(\xi_1,\xi_2,\xi_3)  e^{ix(\xi_1+\xi_2+\xi_3)}  \, d\xi_1d\xi_2d\xi_3.
\end{equation}

Now, in the context of our operator $\THT$, the associated symbol $m_{\THT}$ may be decomposed as
\begin{align*}
m_{\THT}(\xi_1,\xi_2,\xi_3) & := 
\sum_{m_1, m_2 \in \Z}\sum_{k \in \Z} \Big( \int_{\R} e^{i t \frac{\xi_1-\xi_2}{2^k}} e^{i t^3 \frac{\xi_3-\xi_2}{2^{3k}}} \rho(t) dt\Big) \psi\big(\frac{\xi_1-\xi_2}{2^{k+m_1}}\big) \psi\big(\frac{\xi_3-\xi_2}{2^{3k+m_2}}\big)  \\
& := \sum_{m_1, m_2 \in \Z}\sum_{k \in \Z} m_{k, m_1, m_2}(\xi_1, \xi_2, \xi_3)\,,
\end{align*}
where the functions $\psi$ above are of $\Psi$-type while also forming a classical Littlewood-Paley partition of unity (see Section \ref{sec:discretization} for more details). Notice that the multiplier $m_{k, m_1, m_2}(\xi_1, \xi_2, \xi_3)$ merges oscillatory integrals with one-scale multipliers associated to the bilinear Hilbert transform. Following the strategy described in the introductory paragraph above, the oscillatory properties of the phase of the multiplier dictate a natural split into two main regimes (with the latter further split into two subcases) which give rise to the following operators:
\begin{itemize}
\item the \emph{low frequency} component: 
\begin{equation} \label{lof}
\ds T^{lo}:=\sum_{m_1, m_2 <0} \ \sum_{k \in \Z} \ \left( \cdots\right)\,;
\end{equation}

\item the \emph{high frequency} component: $$\ds T^{hi}:=\sum_{m_1, m_2\in \Z^2\setminus\Z^2_{-}} \ \sum_{k \in \Z} \ \left( \cdots\right)$$ with this being further decomposed as 
 $$\ds T^{hi}= \ds T^{hinst}\,+\,\ds T^{hist}\qquad \textrm{where}$$
\begin{itemize}
 \item the \emph{high frequency non-stationary} component:
 \begin{equation} \label{hinstf} 
  \ds T^{hinst}:=\sum_{m_1<0 < m_2} \ \sum_{k \in \Z} \ \left( \cdots\right)\ +\sum_{m_2<0 < m_1} \ \sum_{k \in \Z} \ \left( \cdots\right) \ +\sum_{\substack{m_1, m_2 \geq 0 \\ |m_1-m_2|>3}} \ \sum_{k \in \Z} \ \left( \cdots\right)\,;
 \end{equation} 
\item the \emph{high frequency stationary} component:
 \begin{equation} \label{histf} 
 \ds T^{hist}:=\sum_{\substack{m_1, m_2 \geq 0 \\ |m_1-m_2|\leq 3}} \sum_{k \in \Z} \left( \cdots\right)\,.
 \end{equation}
\end{itemize}    
\end{itemize}

In studying these operators, particularly in the low-frequency regime, we need to consider \emph{anisotropic} multipliers of the form
\begin{equation}
\label{eq:form:op:initial:gen}
\sum_{k \in \Z}c_k \, \phi_1\big(\frac{\xi_1-\xi_2}{2^{k}}\big) \phi_2\big(\frac{\xi_3-\xi_2}{2^{3k}}\big), 
\end{equation}
where here the sequence $(c_k)_k\in \ell^{\infty}(\Z)$ and at least one of the functions $\phi_j$ is of $\Psi$-type while the other one is either of $\Psi$-type or of $\Phi$-type. 

We conclude this preparatory section with three key remarks:

\begin{remark}[\textsf{Anisotropic versus isotropic behavior}]\label{an_iso}   In the context of the existing math literature we have the following antithesis:
\begin{itemize}
\item On the one hand, the \emph{isotropic} analogue\footnote{The same conditions apply to the functions $\phi_j$.} of the symbol in \eqref{eq:form:op:initial:gen}, given by
\begin{equation}\label{isocase}
\sum_{k \in \Z}c_k \, \phi_1\big(\frac{\xi_1-\xi_2}{2^{k}}\big) \phi_2\big(\frac{\xi_3-\xi_2}{2^{k}}\big), 
\end{equation}
is a symbol which, together with its derivatives, decays fast away from the line $\{\xi_1=\xi_2=\xi_3\}$. The operators associated with such symbols were studied before, see \emph{e.g.} \cite{Mlinop};

\item On the other hand, the \emph{anisotropic} operators corresponding to \eqref{eq:form:op:initial:gen} are completely new in the context of modulation-invariant operators, with the anisotropy having an important effect that changes altogether the frequency geometry of the symbol.
\end{itemize}
\end{remark}

\begin{remark}[\textsf{Type I versus Type II anisotropic multipliers}]\label{types}  While in the isotropic case above---due to the symmetry of the expression \eqref{isocase}---it is irrelevant which of the functions $\phi_j$ is of $\Psi$-type, in the anisotropic case \eqref{eq:form:op:initial:gen} the situation is more subtle: 
\begin{itemize}
\item We say that $m$ given by \eqref{eq:form:op:initial:gen} is an anisotropic multiplier of \emph{Type I} iff  
\begin{equation}
\label{eq:symp:aniso}
m(\xi_1,\xi_2,\xi_3):=\sum_{k \geq 0} c_k \varphi\big(\frac{\xi_1-\xi_2}{2^{k}}\big) \psi\big(\frac{\xi_3-\xi_2}{2^{3k}}\big)\,+\,\sum_{k \leq 0} c_k \psi\big(\frac{\xi_1-\xi_2}{2^{k}}\big) \varphi\big(\frac{\xi_3-\xi_2}{2^{3k}}\big)
\end{equation}
for $(c_k)_k\in \ell^{\infty}(\Z)$ and $\psi$ of $\Psi$-type. Equivalently, $m$ is of \emph{Type I} if the $\Psi$-type function acts on the larger frequency scale. 

It is worth mentioning that in this situation the treatment of the corresponding $T_m$ is simpler in nature and has as a consequence the boundedness of  $T_m$ for a larger range than the one corresponding to Type II case discussed below---see Theorem \ref{thm:proof:main:1} and Remark \ref{rem:42}.
\smallskip

\item Conversely, we say that $m$ given by \eqref{eq:form:op:initial:gen} is an anisotropic multiplier of \emph{Type II}, if it is \emph{not} of Type I. 

Notice that via a standard summation procedure, any $m$ of Type II may be rewritten\footnote{The precise meaning of this will be revealed in Section \ref{sec:discretization1}.}---up to symbols of Type I or other well behaved symbols as for example \eqref{isocase}---as
\begin{equation}\label{typeII}
m(\xi_1, \xi_2, \xi_3)=\sum_{k \geq 0} c_k \, \Big( \sum_{j=0}^{3k} \psi\big(\frac{\xi_3-\xi_2}{2^{j}}\big) \Big) \psi\big(\frac{\xi_1-\xi_2}{2^{k}}\big)+ \sum_{k \leq 0} \, c_k \, \Big( \sum_{j=k}^{0} \psi\big(\frac{\xi_1-\xi_2}{2^{j}}\big) \Big)\psi\big(\frac{\xi_3-\xi_2}{2^{3k}}\big)\,,
\end{equation}
where here, as usual by now,  $(c_k)_k\in \ell^{\infty}(\Z)$ and $\psi$ is of $\Psi$-type.
\end{itemize}
\end{remark}

\begin{remark}[\textsf{Anisotropic Biest operators}]\label{anbies} Guided by \eqref{eq:symp:aniso} and \eqref{typeII} we define the following class of multipliers:
\begin{equation}\label{anbiest}
m(\xi_1, \xi_2, \xi_3)=\sum_{k \geq 0} c_k \, \Big( \sum_{j=0}^{3k} d_{j} \psi\big(\frac{\xi_3-\xi_2}{2^{j}}\big) \Big) \psi\big(\frac{\xi_1-\xi_2}{2^{k}}\big)+ \sum_{k \leq 0} \, c_k \, \Big( \sum_{j=k}^{0} d_{j} \psi\big(\frac{\xi_1-\xi_2}{2^{j}}\big) \Big)\psi\big(\frac{\xi_3-\xi_2}{2^{3k}}\big)\,,
\end{equation}
where here $(c_k)_k,\, (d_k)_k\in \ell^{\infty}(\Z)$ and $\psi$ is of $\Psi$-type.

It is not hard to see that this new class naturally extends\footnote{Strictly speaking, in each of the formulas \eqref{eq:symp:aniso}--\eqref{anbiest} the functions $\psi, \phi$ are allowed to change form line to line. Moreover, one may also substitute in the formulas above $\psi\big(\frac{\xi_3-\xi_2}{2^{j}}\big)$ by $j-$dependent functions $\psi_j\big(\frac{\xi_3-\xi_2}{2^{j}}\big)$  as long as suitable support overlapping and uniform smoothness control are imposed on the family $\{\psi_j\}_{j\in\Z}$. This would allow a much more flexible transitioning from \eqref{eq:symp:aniso} to \eqref{anbiest} and more general statements absorbing the coefficients $\{d_j\}_j$ into the properties of the family  $\{\psi_j\}_{j\in\Z}$  (and similarly for $\{c_k\}_{k\in\Z}$). However, for the main goal of our paper non on these are required.}  both Type I and Type II multiplier classes, and, as a consequence, any good control on the multiplier \eqref{anbiest} will immediately imply a similar good control on the multiplier represented by \eqref{eq:form:op:initial:gen}.

Given its above relevance for our present study, for $m$ obeying \eqref{anbiest}, we refer to the associated $T_m$ as an \emph{anisotropic Biest} operator whose boundedness properties will be provided in Theorem \ref{thm:main-ani:biest} below.

We mention here that the above name is justified by the realization that \eqref{anbiest} may be interpreted as the \emph{anisotropic} analogue of the \emph{isotropic} symbol
\begin{equation}
\label{eq:middle:biest}
\sum_{k \in \Z} c_k \Big(\sum_{j \in \Z \atop j \leq k} d_j \,  \psi\big(\frac{\xi_1-\xi_2}{2^{j}}\big) \Big) \psi\big(\frac{\xi_3-\xi_2}{2^{k}}\big).
\end{equation}
Indeed, while strictly speaking, the Biest operator studied originally in \cite{MTTBiest2} has symbol $\ds \one_{\xi_1<\xi_2<\xi_3}$, the latter work has been generalized in \cite{Joeun-IteratedTrilFourierIntg} to a context which essentially matches \eqref{eq:middle:biest}.
\end{remark}

\subsubsection{\underline{The low frequency behavior: the main results}}

Let us denote by ${\bm D}'$---see Figure 2 below---the open interior of the convex envelope spanned by the following seventeen extremal points:
\begin{gather*}
 (0, \frac{1}{2}, 0),\ 
(\frac{1}{2}, 0, \frac{1}{2}),\  
(\frac{3}{4}, \frac{1}{8}, \frac{3}{4}),\ 
(\frac{1}{2}, \frac{3}{4}, \frac{1}{2}), \ 
(\frac{3}{4}, \frac{1}{2}, \frac{1}{4}), \  
(\frac{1}{2}, \frac{3}{4}, 0),\\ 
(\frac{3}{4}, \frac{1}{2}, \frac{3}{4}),\  
(\frac{5}{8}, \frac{5}{8}, 0),\ 
(\frac{3}{4}, \frac{1}{8}, \frac{5}{8}),\ 
(\frac{1}{2}, \frac{1}{2}, 0),\
(\frac{5}{8}, \frac{1}{8}, \frac{3}{4}),\
(0, \frac{3}{4}, 0),\\
(0,\frac{3}{4}, \frac{1}{2}),\
(0, \frac{5}{8}, \frac{5}{8}),\ 
(\frac{2}{3}, 0, \frac{2}{3}),\ 
(\frac{1}{4}, \frac{1}{2}, \frac{3}{4}),\ 
(0, \frac{1}{2}, \frac{1}{2}).
\end{gather*}

\begin{figure}[h]
\includegraphics[scale=0.35]{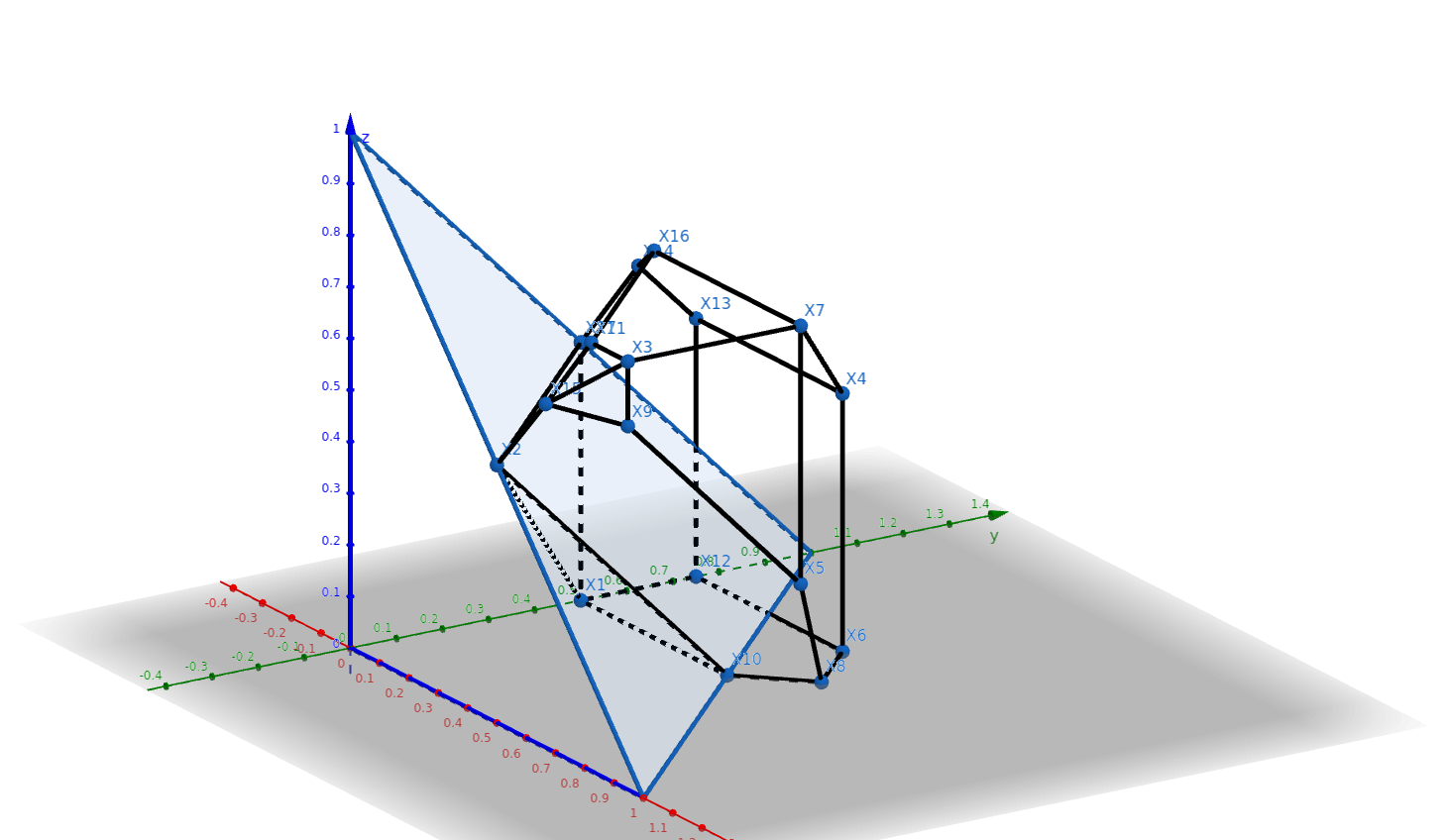}
\caption{The drawing representing the polytope ${\bm D}'$ with its seventeen extremal points. The plane $\{\alpha=(\alpha_1,\alpha_2,\alpha_3),\ \alpha_1+\alpha_2+\alpha_3=1\}$ in displayed in blue. Notice that only the points $X_1=F,X_{12}=E$ are strictly inside the Banach range.}
\end{figure}

Then, the following holds:
\begin{theorem}[\textsf{Model for the low frequency component of $\THT$}]
\label{thm:ani}
Let $(c_k)_{k \in \Z}$ be a bounded sequence of complex numbers and $m$ a symbol as in  \eqref{eq:form:op:initial:gen} with at least one of the functions $\phi_1$ and $\phi_2$ of $\Psi$-type. Then, $T_m$ is bounded from $L^{p_1} \times L^{p_2} \times L^{p_3} \to L^p$ whenever $(\frac{1}{p_1},\frac{1}{p_2},\frac{1}{p_3})\in{\bm D}'$ and $\frac{1}{p_1}+\frac{1}{p_2}+\frac{1}{p_3}=\frac{1}{p}$. Moreover, we have that
 $$ \|T_m\|_{L^{p_1}\times L^{p_2} \times L^{p_3} \to L^p} \lesssim_{p_1,p_2,p_3,\phi_1,\phi_2} \|(c_k)_{k}\|_{\ell^\infty(\Z)}.$$
\end{theorem}

Next, we notice that symbols of type \eqref{eq:form:op:initial:gen} are morally equivalent to $\ds \mathbbm{m}((\xi_1-\xi_2)^3, \xi_3-\xi_2)$, where $\mathbbm{m}\in C^\infty({\mathbb R}^2 \setminus\{(0,0)\})$ is a $2$-dimensional Mikhlin symbol obeying for sufficiently many multi-indices
\begin{equation} |\partial_{\eta_1}^{\alpha_1} \partial_{\eta_2}^{\alpha_2} \mathbbm{m}(\eta_1,\eta_2)|\lesssim_{\alpha_1,\alpha_2} |(\eta_1,\eta_2)|^{-(\alpha_1+\alpha_2)}\,. \label{eq:Mikhlin} 
\end{equation}

As a consequence, the following may be proved with similar methods as the ones used for the resolution of Theorem \ref{thm:ani}:

\begin{theorem} [\textsf{Anisotropic Mikhlin multiplier}] \label{thm:anisotropic:Mikhlin}
Let $\mathbbm{m}$ be a 2-dimensional Mikhlin symbol, satisfying \eqref{eq:Mikhlin} and let $m$ be defined by $\ds m(\xi_1, \xi_2, \xi_3):=\mathbbm{m}((\xi_1-\xi_2)^3, \xi_3-\xi_2) $. Then, the trilinear operator $T_m$ defined by \eqref{eq:op:anisotropic} is bounded\footnote{The range in Theorem \ref{thm:anisotropic:Mikhlin} consists of the polytope ${\bm D}$ to which we add the face given by the interior of the points $A,B,C,D,G$ embedded into $\{\alpha:=(\alpha_1,\alpha_2,\alpha_3),\ \alpha_1+\alpha_2+\alpha_3=1\}$. For example, $T_m$ is bounded from $L^4\times L^2 \times L^4$ into $L^1$.}
 from $L^{p_1} \times L^{p_2} \times L^{p_3} \to L^p$ provided that $(\frac{1}{p_1},\frac{1}{p_2},\frac{1}{p_3})\in{\bm D}'$ and $\frac{1}{p_1}+\frac{1}{p_2}+\frac{1}{p_3}=\frac{1}{p}$. 
\end{theorem}
%
In effect, as a consequence of the line of thought developed throughout the Remarks \ref{an_iso}---\ref{anbies}, we have that the above two theorems are direct consequences of the following:

\begin{theorem}  [\textsf{Anisotropic Biest operator}] \label{thm:main-ani:biest}
The anisotropic Biest operator $T_{anBiest}:=T_{m}$ having as symbol the multiplier $m$ defined by \eqref{anbiest} is a bounded trilinear operator from $L^{p_1} \times L^{p_2} \times L^{p_3} \to L^p$ whenever $(\frac{1}{p_1},\frac{1}{p_2},\frac{1}{p_3})\in{\bm D}'$ and $\frac{1}{p_1}+\frac{1}{p_2}+\frac{1}{p_3}=\frac{1}{p}$. Moreover, we have 
$$\|T_{anBiest}\|_{L^{p_1}\times L^{p_2} \times L^{p_3} \to L^p} \lesssim_{p_1,p_2,p_3} \, \|(c_j)_{j}\|_{\ell^\infty} \, \|(d_k)_{k}\|_{\ell^\infty}.$$
\end{theorem}

We conclude this section with the following 

\begin{remark}[\textsf{Open problem: the nature of structure for well-behaved Marcinkiewicz multipliers}] Inspecting the statements above, we notice that all the main operators involve particular cases of trilinear symbols having the form $m(\xi_1,\xi_2,\xi_3)=\mathsf{m}(\xi_1-\xi_2,\xi_3-\xi_2)$ with $\mathsf{m}$ satisfying the Marcinkiewicz condition
\begin{equation}
 \big|\partial_{\eta_1}^{n_1}\partial_{\eta_2}^{n_2} \mathsf{m}(\eta_1,\eta_2)\big|\lesssim |\eta_1|^{-n_1} |\eta_2|^{-n_2}. \label{eq:Marcin}
 \end{equation}
Once at this point, we record the presence of the following dichotomy:
\begin{itemize}
\item On the one hand, we have that condition \eqref{eq:Marcin} is not sufficient in order to guarantee the boundedness of the corresponding trilinear operator---see \cite[Section 7]{GrafakosKalton} for examples of unbounded bilinear operators corresponding to our trilinear operators by letting $f_2\equiv 1$; 
 
\item On the other hand, all the results in the present section suggest that by adding into the picture a multiplicative structure paired with an anisotropic character one is able to obtain non trivial boundedness results. 
 \end{itemize} 
 As a consequence of the above, it is natural to investigate necessary and sufficient conditions for establishing the boundedness of operators associated with Marcinkiewicz symbols relative to suitable, natural extra-structures imposed on these symbols.
\end{remark}

\subsection{Outline of the main ideas behind the proof and structure of the paper}\label{Keyideas}
\medskip

As already evoked in the brief introductory paragraph, see \textbf{(P1)}---\textbf{(P2)}, the key relevance of our $\THT$ operator is its \emph{hybrid} nature meant to capture the interplay between 
\begin{itemize}
\item \textbf{(\underline{P1})} \emph{zero/resonant curvature} features;

\item \textbf{(\underline{P2})} \emph{non-zero/non-resonant curvature} features.
\end{itemize}

This subtle balance between the two competing factors propagate throughout the study of our main operator and is reflected into the  analysis of the two key components of $\THT$---see the decomposition \eqref{dectht} in our preliminary Section \ref{sec:discretization1}:
\begin{itemize}
\item the \emph{low oscillatory} component ($T^{\low}$): in this instance, its associated multiplier enjoys both  
\begin{itemize}
\item \textbf{(a)} translation invariance derived from \textbf{(\underline{P1})}, and

\item \textbf{(b)} anisotropy derived from \textbf{(\underline{P2})}.
\end{itemize}
\smallskip

\item the \emph{high oscillatory}\footnote{As it turns out, the non-stationary subcomponent $T^{\nonsta}$ can be treated in a similar spirit with the low frequency component due to the fast decay encoded in the associated symbol; the latter is a direct manifestation of the non-stationary character, \emph{i.e.}, the lack of stationary points for the phase of the multiplier.} component ($T^{\sta}$): in this instance, focussing on the stationary subcomponent, the multiplier of the operator satisfies on top of \textbf{(a)} and \textbf{(b)} above also condition
\begin{itemize}
\item \textbf{(c)} a highly oscillatory twisted phase which is a consequence of \textbf{(\underline{P2})}.  
\end{itemize}     
    
\end{itemize}

Now, in accordance with the above binary division, our paper splits into two main parts:
\medskip

\subsubsection{\underline{Part I. The boundedness of the anisotropic Biest operator. Treatment of $T^{\low}$ and $T^{\nonsta}$}.}$\qquad$
\smallskip

The joint features \textbf{(a)} and \textbf{(b)} above create a new type of object, referred to as an anisotropic Biest operator---see Remark \ref{anbies} and Theorem \ref{thm:main-ani:biest}---which is of independent interest. The behavior of this new object as captured in Theorem \ref{thm:main-ani:biest} controls the behavior of both the low frequency component (Theorem \ref{thm:low}) and of the high frequency non-stationary component (Theorem \ref{thm:non-sta}) and further implies the corresponding statements in Theorem \ref{thm:ani} and Theorem \ref{thm:anisotropic:Mikhlin}. All these results  are addressed in Part I of our paper covering Sections \ref{sec:anis:model:red}--\ref{sec:low}, as follows:  
\begin{itemize}
\item Section \ref{sec:anis:model:red} serves as a warm-up stage focussing on several, rather standard reductions:
\begin{itemize}
\item In Section \ref{subsec:anis:model:red} we show that Theorem \ref{thm:anisotropic:Mikhlin} addressing the behavior of anisotropic Mikhlin multipliers is easily derived from Theorem \ref{thm:ani} which in turn is a consequence of Theorem \ref{thm:main-ani:biest};

\item In Section \ref{typeI1} we provide an elementary treatment to a special family of Type I multipliers;

\item In Section \ref{reduct}, we show that the control over Type I and Type II multipliers as well as of Theorem \ref{thm:main-ani:biest} may be reduced to Theorem \ref{thm:proof:main:1} and Theorem \ref{thm:proof:main:2}. 
\end{itemize}

\item Section \ref{sec:proof:m:+} provides the resolution of Theorem \ref{thm:proof:main:1} concerned with the boundedness behavior of (generic) Type I multipliers. In direct correspondence with properties \textbf{(a)} and \textbf{(b)} listed at the beginning of our section, the analysis of our Type I multiplier involves:
\begin{itemize}    
\item \textbf{(a1)} a wave-packet decomposition which is consistent with, and further extends the $\BHT$ time-frequency analysis; 

\item \textbf{(b1)} involvement of time-frequency regions which have area larger than one, hence the necessity of considering super-Heisenberg tiles. 
\end{itemize}
The juxtaposition of these two features create a "\emph{dilated four-linear $\BHT$-type analysis}", which, as expected, comes with several challenges:
\begin{itemize} 
\item firstly, one has to consider quadrilinear forms (see \eqref{eq:model:almost:tiles}) associated with $4-$super-Heisenberg tiles $P=(P_1, P_2, P_3, P_4)$ ranging through a rank 1 family $\P$  (see \eqref{eq:def:collection:P}) which capture the geometric description in Figure 3, Section \ref{sec:red:op:anis:local:L2};  

\item secondly, within each $P_j$, we apply a suitable Heisenberg tile partition and consider the corresponding adapted wave-packet family: this naturally points toward an $L^2$ based energy quantity associated with each $P_j$ (see \eqref{grup}) which then propagates throughout the key estimates---see \emph{e.g.} the size estimate in Proposition \ref{prop:size1}. In particular, this approach is responsible for the local $L^2$ range restriction for the boundedness range in Theorem \ref{thm:proof:main:1}. Thus, any attempt of extending the range requires some new refinements/ideas which are the topic of the more general result treated in the following section.  
\end{itemize} 

\item  Section \ref{sec:proof:m:2} addresses Theorem \ref{thm:proof:main:2} which covers the generic anisotropic Biest multipliers and thus subsumes both Type I and Type II multipliers. As expected, while the analysis performed in this case is far more delicate that the one in the previous case, it is still consistent with  \textbf{(a)} and \textbf{(b)} in the introductory paragraph, and thus requires:
\begin{itemize}    
\item \textbf{(a2)} a modulation invariant wave-packet decomposition, which this time, unlike \textbf{(a1)}, it involves two rank-one families of Heisenberg tritiles, thus  encapsulating a Biest-type scenario; 

\item \textbf{(b2)} unlike the classical, isotopic Biest scenario, our present setting involves a \emph{distinct scale} "composition/itteration" of two BHT-type model operators.
\end{itemize}     
The juxtaposition of these two features create an "\emph{anisotropic Biest-type analysis}", which, again, comes with several challenges on its own:
\begin{itemize} 
\item firstly, inspecting the derived discretized model operator in \eqref{eq:def:form:order} we notice that its anisotropic character enforces a lack of scale symmetry in the associated Gabor coefficients \eqref{eq:coeffA}, with the third defining sequence $ a_{Q_3}^{(3)}$ involving a joint $(f_3, f_4)$ behavior. This latter fact, propagates to the size estimates in Proposition \ref{prop:size:2} and becomes a first difficulty as well as a first reason for the final restrictive boundedness range;   

\item secondly, due to the same anisotropic character, there are new difficulties arising for a good control over the energy (see \eqref{eq:change:order:form:P} and \eqref{energydef}). Indeed, while in the classical (isotropic) Biest scenario, the symmetry between the two variables $\xi_1,\xi_3$ in the isotropic analogue of \eqref{eq:tril:aniso:op:bi:est} allows the authors in \cite{MTTBiest1,MTTBiest2} to obtain via a certain trick a good energy control, in our current setting there is no natural substitute. Thus, in order to tackle this issue, we first introduce a suitable maximally truncated operator\footnote{A very rudimentary version such an operator was also involved, in different context, in \cite[Section 10.1]{HuLie}.} (see \eqref{eq:def:H*}) whose control---expressed in \eqref{eq:bddnedd:max:bessel:function}---is obtained via a well known maximal inequality of Bourgain (\cite{Bou87-erg}). We mention here that a similar reasoning was employed for treating the maximal truncated bilinear Hilbert transform (\cite{Lacey}, \cite{DemeterTaoThiele}). With this done, we need to introduce a new notion of (tile) {\it incidence}---see \eqref{eq:def:incidence}---which controls via \eqref{eq:dualization:energy} and \eqref{eq:crucial:en:est:2} the energy, and, in turn, is controlled by the multiplicity counting function of a forest---see \eqref{eq:upper:bound:incidence}. Finally, employing the nontrivial loss in the energy estimate offered by Proposition \ref{prop:energy}, a variational argument involving the size of the incidence provides us with the key control in Propositions \ref{prop:algo:2}. Notice however, that the above mentioned nontrivial loss propagates throughout these estimates thus preventing us from recovering the full boundedness range of the classical Biest operator.
\end{itemize} 

\item Finally, in Section \ref{sec:low} we provide the treatment of $T^{\low}$ and  $T^{\nonsta}$ corresponding to Theorem \ref{thm:low} and Theorem \ref{thm:non-sta}, respectively. Indeed, as it turns out, the boundedness behavior of both these components is derived via some standard reasonings from the corresponding properties of the anisotropic Biest operator stated in Theorem \ref{thm:main-ani:biest} which have been understood\footnote{Recall that Theorem \ref{thm:main-ani:biest} was reduced in  Section \ref{sec:anis:model:red} to Theorem \ref{thm:proof:main:1} and Theorem \ref{thm:proof:main:2}.}  in Sections \ref{sec:proof:m:+} and \ref{sec:proof:m:2}.
\end{itemize}

\subsubsection{\underline{Part II. Treatment of the high frequency stationary component \bf{$T^{\sta}$}}.}$\qquad$
\medskip

The second part of our paper is devoted to the proof of Theorem \ref{thm:main} which concerns the boundedness properties of
$$T^{\sta}\approx\ds \sum_{m \geq 0} \sum_{k \in\Z} T_{k,m}\,$$
where here  $T_{k,m}$ is  associated with the multiplier
\beq\label{multkmas}
\mathfrak{m}_{k,m}(\xi_1,\xi_2,\xi_3):= \left(\int_{\R}e^{-i\,\frac{\xi_1-\xi_2}{2^k}\,t}\,e^{-i\,\frac{\xi_3-\xi_2}{2^{3k}}\,t^3}\,\rho(t)\,dt\right)\,
\psi\left(\frac{\xi_1-\xi_2}{2^{k+m}}\right)\,\psi\left(\frac{\xi_3-\xi_2}{2^{3k+m}}\right),
\eeq
and $(k,m)\in \Z\times\N$ represent
\begin{itemize}
\item $k$ is the physical scale capturing the magnitude in the $t$ integrand;

\item $m$ is the frequency scale capturing the magnitude of the phase multiplier. 
\end{itemize}
The entire fabric of the proof of Theorem \ref{thm:main} relies now on the interplay between two types of interactions:
\begin{itemize} 
\item $(\mathcal{I} 1)$ the interaction between physical scale $k$ and frequency scale $m$;

\item $(\mathcal{I} 2)$ the interaction between the zero and nonzero curvature features \textbf{(\underline{P1})} and \textbf{(\underline{P2})}, respectively, which transfer to the interaction between points \textbf{(a)} on the one hand, and, points \textbf{(b)} and \textbf{(c)} on the other hand.    
\end{itemize}

Following these, and consistent with the \textbf{(\underline{P1})}--\textbf{(\underline{P2})} interaction, the approach in Part II is divided into a preamble and two main subcomponents:
\begin{itemize} 
\smallskip
\item \underline{\textsf{Preamble:}} This is the content of Section \ref{sec:discretization}, which encompasses a preliminary discussion harmonizing $(\mathcal{I} 1)$ and $(\mathcal{I} 2)$:
    
In Section \ref{subsec:prell} we introduce the quadrilinear form $\L_{k,m}(f_1,f_2,f_3,f_4):=\int_{\R} T_{k,m}(f_1,f_2,f_3)\,f_4,$  and lay out an outline for approaching Theorem \ref{thm:main}---see Steps 1 to 3 and the key Theorem \ref{thm:case1} therein.

In Section \ref{lres} we introduce a preliminary \emph{low resolution discretization}, which, for each pair $(k,m)$ expresses $\L_{k,m}=\sum_{j \in \Z}\sum_{\el\in\Z} \sum_{\r\sim 2^{2k}} \L_{k,m}^{\el,\r, j}$ with each component $\L_{k,m}^{\el,\r, j}$ being adapted to a (four-)super-Heisenberg tile $P=P(k,m,\el,\r, j)=(P_1,P_2,P_3,P_4)$ such that $P_j$ determines the time-frequency localization of the input $f_j$. Moreover, for each fixed $m$ the collection of all  $\{P(k,m,\el,\r, j)\}_{k,\el,\r,j}$ forms a rank one family of (four super-Heisenberg) tiles. 

In Section \ref{subsec:63} we introduce a preliminary \emph{high resolution discretization} of each $\L_{k,m}^{\el,\r, j}$  into a superposition of wave-packets interacting with each of the input functions---see Proposition \ref{prop:GFF}. This discretized model is achieved via Rank I LGC method, it implies a partitioning of each super-Heisenberg tile $P_j$ into a suitable family of Heisenberg tiles,  and is \emph{not absolutely summable}.\footnote{In the sense explained in Remark \ref{summability}.}
\smallskip

\item \underline{\textsf{Part II.1. A high-resolution single scale analysis:}}  This is an expression of \textbf{(\underline{P2})}, since the focus is on exploiting the cancellation hidden in the non-zero curvature feature of the multiplier \eqref{multkmas}. Concretely, building on the intuition gained in Section  \ref{subsec:63}, our main goal is to obtain an exponential decay (smoothing estimate) of the form $2^{-\delta\min\{2|k|,m\}}$ ($\delta>0$ absolute constant) for some suitable $L^2$ based norm applied to $\L_{k,m}^{\el,\r, j}$---see \eqref{F1} part of Observation \ref{Reduct}. Consistent with the  $(\mathcal{I} 1)$ interaction this decay estimate is analyzed in four distinct regimes: $k\geq \frac{m}{4}$  in Section \ref{sec:onescale_k_large}, $0\leq k\leq \frac{m}{8}$ in Section \ref{Diagkpsmall124}, $\frac{m}{8}<k<\frac{m}{4}$ in Section \ref{Diagkinterm} and finally $k<0$ (touched upon only very briefly due to the symmetric character of $\THT$) in Section \ref{Seckneg}. The crux of all this single-scale analysis relies on the Rank II LGC method: via the so-called joint (maximal) Fourier coefficients, this method converts the non-absolutely summable Rank I LGC discretized model into an absolutely summable one. Moreover, a key ingredient in achieving our above stated goal is the reduction of \eqref{F1} to the analysis of the hybrid Bilinear Carleson operator $BC^a$ in \cite{BCa} by further developing the methods introduced in \cite{HuLie}.
\smallskip      

\item \underline{\textsf{Part II.2. A low-resolution multiscale analysis:}}   This is an expression of \textbf{(\underline{P1})} addressing the modulation invariance character of $\THT$ and having as its main focus the propagation of the smoothing/decay estimate obtained at Part II.1. from \emph{one} scale to \emph{arbitrarily many} scales. Thus, the key phenomenon here is capturing the cancellation \emph{in between} super-Heisenberg tiles $P$ belonging to a rank one family of tiles as opposed to the cancellation \emph{within} a super Heisenberg tile which was the topic of Part II.1. This is precisely what explains the passage from the high (single scale) to the low (multiscale) resolution, underlying the multiscale modulation invariance analysis we are to perform on $T^{\sta}$. Indeed, noticing the (double) translation invariance property structure within \eqref{multkmas} relative to the $\xi_2-\xi_1$ or $\xi_3-\xi_2$ arguments it is natural to expect that our analysis will share many common points with the ones employed for the bilinear Hilbert transform and the Biest operator. All these are covered in Sections \ref{sec:bilinear:analysis}, \ref{sec:bilinear:analysis2} and \ref{sec:bilinear:analysis3} mirroring the $(k,m)$ interaction regimes in Part II.1.
\end{itemize}
\smallskip

\subsection{A brief historical background, context and motivation}\label{Motivation} In this final introductory section we present a succinct overview on the historical development, context for our present work as well as its underlying motivation.

\subsubsection{Historical background} Given $n\in\N$ and $\vec{\g}:=(\g_1,\ldots, \g_n)$ with $\{\g_j\}_{j=1}^n$ real planar smooth curves obeying some suitable nondegeneracy conditions, we define the generic \textsf{$n$-linear \emph{$\vec{\g}$-Hilbert transform}} by
\beq\label{nhilb}
H_{n,\vec{\g}}(f_1,\ldots,f_n)(x):= \textrm{p.v.}\int_{\R}\prod_{j=1}^n f_j(x-\g_j(t))\,\frac{dt}{t}
\eeq
where here all functions $\{f_j\}_{j=1}^{n}$ may be apriori considered within the Schwarz class.
\smallskip

In this paper we study the boundedness behavior of \eqref{nhilb} under two key prescribed conditions:
\begin{itemize}
\item \underline{\emph{$n=3$ multilinearity}}, thus moving beyond the well understood bilinear setting, and, crucially, 

\item the presence of a \underline{\emph{hybrid regime}}, that is,  when \eqref{nhilb} has both 

\begin{itemize}
\item \emph{resonance} features, or equivalently, modulation invariance properties,  and

\item \emph{non-resonance} features, or equivalently, non-zero curvature properties.
\end{itemize}

\end{itemize}

While the previous sections together with a subsequent discussion in the present section clarify and offer more context to the meaning, relevance  and motivation for considering the specific problem of the trilinear Hilbert transform in the hybrid regime, for now, we keep our exposition concise and only focus on a brief historical overview:

\begin{itemize}
\item \textbf{(I)} The origin of the generic program concerning \eqref{nhilb} traces back to the classical work of M. Riesz (\cite{MR28}) in the $n=1$ case as well as 
to the seminal work of M. Lacey and C. Thiele (\cite{lt1}, \cite{lt2}) on the bilinear Hilbert transform ($n=2$) in the \emph{resonant} regime, \emph{i.e.}, the case $\g_j(t)=\lambda_j t$ with $\l_1\not=\l_2\in \R$. The $n=2$ \emph{non-resonant} regime was first treated in various settings with diverse, heterogenous  methods, see \cite{Li13}, \cite{Lie15}, \cite{LX16}, \cite{Lie18}, \cite{GL20}, \cite{DR22}. Later, this was integrated within a dedicated, homogeneous program concerning the wide applicability range of the so-called (Rank I) LGC-methodology---see \cite{Lie19} and \cite{GL22}, with the latter also encapsulating hybrid regimes corresponding to the situation
\beq\label{bilhilb}
 H_{2,\vec{\g}} (f_1,f_2)(x):= p.v. \int_{\mathbb R} f_1(x-\lambda_1 t)\,f_2(x-\lambda_2 t-\lambda_3 t^a)\,\frac{dt}{t},\quad x \in \R\,,
\eeq
where here $\l_1\not=\l_2\in\R$, $\l_3\in\R$ and $a\in(0,\infty)\setminus\{1,2\}$. 

Thus, the $n=1$ and $n=2$ cases are by now essentially\footnote{Indeed, there are only few instances left open, among which we mention the question on the maximal boundedness range for the resonant bilinear Hilbert transform.} settled.
\medskip

\item \textbf{(II)} In the process of transitioning towards the $n=3$ case, a new, intermediate, but central class of objects was introduced initially in \cite{BCa} and later in \cite{BHL25}, corresponding to two distinct regimes---see \textbf{(II.1)} and \textbf{(II.2)} below. Indeed, the theory\footnote{For a more general framework as well as a further motivational discussion please see the Introduction in \cite{BHL25}.} of the so-called Bilinear Hilbert-Carleson operators concerns the study of the boundedness properties for   
\beq\label{nhilBHCarl}
BC_{\vec{\G}}(f_1,f_{2})(x):= \textrm{p.v.}\int_{\R}f_1(x-\underline{\g}_1(t))\,f_2(x-\underline{\g}_2(t))\,\frac{e^{i\,\g(x,t)}}{t}\,dt\,,
\eeq
where here $\vec{\G}(x,t):=(\underline{\g}_1(t),\underline{\g}_{2}(t), \g(x,t))$ with each $\underline{\g}_j$ a suitable piecewise smooth planar curve and  $\g(x,t)$ a piecewise smooth function in $t$ and only measurable in $x$.

The importance of considering these type of objects is two-folded:
\begin{itemize}
\item On the one hand, the analysis of \eqref{nhilBHCarl} reveals a genuine new manifestation, which, confirms the higher complexity degree of these objects relative to their bilinear Hilbert analogues: indeed, unlike the previous harmonic analysis literature, the standard wave-packet discretized model associated with \eqref{nhilBHCarl} is \emph{lacking the absolute summability property even within a single spatial scale}\footnote{Please see Remark \ref{summability} for the precise meaning of this.}. This challenge is overcome in \cite{BCa} via a new method  which introduced the concept of joint (maximal) Fourier coefficients. Relying on the latter, the authors therein regrouped the information carried by the conditionally summable wave packet discretized model into an absolutely summable (correlative) time-frequency model. 

\item On the other hand, an even deeper analysis uncovers deep connections between different regimes of Bilinear Hilbert-Carleson operators and their (suitably defined) trilinear Hilbert analogues with all these connections being detailed at item \textbf{(III)} below.  
\end{itemize}

Touching now briefly upon both  the chronology and the regime analysis we have:
\begin{itemize}
\smallskip 
 
\item \textbf{(II.1)} The prototypical \emph{hybrid regime} corresponding in \eqref{nhilBHCarl} to $\underline{\g}_1(t)=t$, $\underline{\g}_2(t)=-t$ and $\g(x,t)=\lambda(x) t^a$ with $\lambda(\cdot)$ real measurable and $a\in (0,\infty)\setminus\{1,2\}$ was resolved in \cite{BCa}. Here the hybrid reference is justified by the fact that the resulting bilinear operator---let's call it for simplicity $BC_{H}$--obeys
    \begin{itemize}
    \item a resonant (\textit{i.e.}, modulation invariance) feature inherited from the resonant bilinear Hilbert transform, concretely 
 \[|BC_{H} (M_c f_1, M_c f_2)| = |BC_{H} (f_1, f_2)|;\]
    
    \item a non-resonant (\textit{i.e.}, non-zero curvature) feature, that is, the kernel representation of $BC_{H}$ involves a maximal (generalized) modulation in the presence of curvature. 
    \end{itemize}
\smallskip    

\item \textbf{(II.2)} The prototypical \emph{non-resonant regime}  corresponding in \eqref{nhilBHCarl} to $\underline{\g}_1(t))=t$, $\underline{\g}_2(t)=t^2$ and $\g(x,t)=\lambda(x) t^a$ with $\lambda(\cdot)$ and $a\in (0,\infty)\setminus\{1,2\}$ was resolved in \cite{BHL25}. Notice that in this situation our operator has the more amenable feature of lacking modulation invariance properties but this simplification is counterbalanced by the anisotropic character of the arguments within the input functions.
\end{itemize}
\medskip

\item \textbf{(III)} Finally, the approach of the case $n=3$ followed the natural complexity hierarchy dictated by the non-zero to zero curvature spectrum:
\begin{itemize}
\smallskip

\item \textbf{(III.1)} The \emph{purely non-resonant regime} was settled in \cite{HuLie}, \cite{HL25} where the authors proved that the curved trilinear Hilbert transform defined by taking in \eqref{nhilb} $\g_j(t)=t^{\a_j}$ with $\{\a_j\}_{j=1}^3\subset (0,\infty)$ \emph{pairwise distinct}, is bounded from $L^{p_1}(\R)\times L^{p_2}(\R)\times L^{p_3}(\R)\,\rightarrow\,L^r(\R)$ within the range\footnote{Assuming wlog that $\{\a_j\}_{j=1}^3$ is strictly increasing, one can also add to the range the value $p_2=\infty$.} $\frac{1}{p_1}+\frac{1}{p_2}+\frac{1}{p_3}=\frac{1}{r}$, $1<p_j<\infty$ and $r>\frac{1}{2}$.

The resolution of this regime relied on developing the so called Rank II-LGC method which originated in the insight provided by the (maximal) joint Fourier coefficient used for treating the hybrid Bilinear Hilbert Carleson operator. Moreover, a key fundamental step for proving the boundedness of the curved trilinear Hilbert transform was the control of the associated non-resonant Bilinear Carleson operator obtained by taking in \eqref{nhilBHCarl} $\underline{\g}_j(t)=t^{\a_j}$ for $j\in\{1,2\}$ and $\g(x,t)=\lambda(x) t^{\a_3}$. 
\smallskip

\item \textbf{(III.2)} The \emph{hybrid regime} is settled in the present paper by allowing in \eqref{nhilb} both modulation invariance (zero curvature) features and purely non-zero curvature features. Specifically, we show that given $a\in(0,\infty)\setminus \{1,2\}$ and $\vec{\lambda}:=(\lambda_1,\lambda_2,\lambda_3,\lambda_4)\in (\R\setminus\{0\})^4$ with  $\lambda_1\neq \lambda_2$ and $\lambda_3\neq \lambda_4$, the trilinear operator
\beq\label{hybHilbtransf}    
 \THT^a_{\vec{\lambda}} (f_1,f_2,f_3)(x):= p.v. \int_{\mathbb R} f_1(x-\lambda_1 t)\,f_2(x-\lambda_2 t-\lambda_3 t^a)\,f_3(x-\lambda_4 t^a)\,\frac{dt}{t},\quad x \in \R
 \eeq
obeys a wide-range of $L^p$ bounds. The modulation invariance features are encoded in the relation 
$$ |\THT^a_{\vec{\lambda}} (M_{\frac{c}{\lambda_1}}f_1,M_{-\frac{c}{\lambda_2}}f_2, M_{\frac{c \lambda_3}{\lambda_2\lambda_4}}f_3)|=|\THT^a_{\vec{\lambda}} (f_1,f_2,f_3)|\qquad \forall\:c\in\R\,,$$
while the nonzero curvature features are encapsulated within the interaction between the arguments $t$ and $t^a$ appearing as inputs within  \eqref{hybHilbtransf}. 

This object extends the behavior of \eqref{bilhilb} while also encapsulating many features of the purely curved trilinear Hilbert transform discussed at the previous item. The boundedness of \eqref{hybHilbtransf} involves a vast array of tools, among which, yet again, a crucial ingredient is offered by the hybrid Bilinear Hilbert Carleson operator discussed ai item \textbf{(II.1)}.
\end{itemize}
\end{itemize}
\smallskip

\subsubsection{Context and rationale} As already mentioned at the very beginning of our Introduction, this work is part of a vast research program aim to better understand the interplay between zero and  non-zero curvature in the harmonic analysis realm. For more details on all these, we refer the interested reader to the extended introductions in \cite{BCa}, \cite{GL22}, \cite{HuLie} and \cite{HsuL26}.

In what follows, we particularize \eqref{hybHilbtransf} to its essential representative 
\beq\label{hthtrep}
\THT (f_1,f_2,f_3)(x):= p.v. \int_{\mathbb R} f_1(x-t)\,f_2(x+t+t^3)\,f_3(x-t^3)\,\frac{dt}{t}\,,
\eeq
and start our discussion by a comparison with its closest, already understood, relatives
\begin{itemize}
\item{The bilinear Hilbert-Carleson operator $BC^a$}, $a\in(0,\infty)\setminus\{1,2\}$, (\cite{BCa}), given by
 $$ BC^a(f_1,f_2)(x) := \sup_{\lambda>0} \Big| p.v. \int f_1(x-t) f_2(x+t) e^{i\lambda t^a} \frac{dt}{t}\Big|$$
 \item{The purely curved trilinear Hilbert transform},  (\cite{HuLie}),  given by
$$ H_{3,(1,2,3)}(f_1,f_2,f_3)(x):=p.v. \int_{\mathbb R} f_1(x-t)f_2(x-t^2)f_3(x-t^3) \frac{dt}{t}.$$
\end{itemize}

We first observe that both these operators involve a non-zero curvature feature: for $BC^3$ via the curve $(t,t^3)$, and, for $H_{3,(1,2,3)}$ via the curve $(t,t^2,t^3)$. Moreover, $BC^3$ also enjoys a zero curvature feature in the form of a $\BHT$ modulation invariance. Building on these, the $\THT$ operator serves as a natural next hierarchical object in the form of a \emph{hybrid trilinear} operator enjoying the following properties: 
\begin{itemize}
\item {\textsf{Translation invariance}:} For $\tau_{x_0}$ the translation of $x_0$, one has for every $x_0\in{\mathbb R}$ and $f_i\in{\mathcal S}({\mathbb R})$
$$ \THT(\tau_{x_0}f_1,\tau_{x_0}f_2, \tau_{x_0}f_3)=\tau_{x_0} \THT(f_1,f_2,f_3);$$  
\item {\textsf{Modulation invariance}:} For $M_{\xi_0}$ the modulation by $\xi_0$, one has for every $\xi_0\in{\mathbb R}$ and $f_i\in{\mathcal S}({\mathbb R})$
$$ \THT(M_{\xi_0} f_1,M_{\xi_0} f_2, M_{\xi_0} f_3)=M_{3\xi_0} \THT(f_1,f_2,f_3).$$  
\item {\textsf{Asymptotic scale symmetry}:}\footnote{This symmetry is not required in our proof but, however, it is a matter of convenience since it allows a symmetric treatment of the case $k<0$ relative to the case $k\geq 0$.} The behavior for $|t|\to 0$ is the same as the one for $|t| \to \infty$ by exchanging the roles of $f_1$ and $f_3$.
\end{itemize}

With these in mind, we move to the following comparative analysis:

\begin{itemize}
\smallskip

\item {\textsf{[Asymptotic behavior and relationship to $BC^a$]:}} Isolating the asymptotic behavior of $t$ in \eqref{hthtrep}, we have
\begin{itemize}
\smallskip
\item on the one hand, when $|t|\ll 1$ then $t+t^3 \sim t$ and so, at the heuristic level, $\THT$ is expected to behave similar to $ p.v. \int_{|t|<1} f_1(x-t)\,f_2(x+t)\,f_3(x-t^3)\,\frac{dt}{t}$;
  
\item on the other hand, when $|t|\gg 1$ then $t+t^3 \sim t^3$ and so, at the heuristic level, $\THT$ is expected to behave similar to $ p.v. \int_{|t|>1} f_1(x-t)\,f_2(x+t^3)\,f_3(x-t^3)\,\frac{dt}{t}$.
\end{itemize}
Now both the trilinear expressions written explicitly above are the direct trilinear counterpart of the bilinear maximal operator $BC^a$ for $a\in\{\frac{1}{3}, 3\}$. 

Consistent with the overarching philosophy discussed in Section 1.2 of \cite{BHL25} stating that there exits a deep connection between the $n$-linear Hilbert transform \eqref{nhilb} and the associated  $(n-1)$-linear Carleson operator, our study of $\THT$ is directly linked to $BC^{a}$, with an explicit reference appearing in Proposition \ref{sgscale}.
\smallskip

\item {\textsf{[Complexity hierarchy]:}} The hybrid transform $\THT$ serves as a natural intermediate stepping stone between the purely non-zero curvature setting represented by $H_{3,(1,2,3)}$ and the purely zero curvature setting embodied by the classical trilinear Hilbert transform given by
$$ (f_1,f_2,f_3) \mapsto p.v. \int_{\mathbb R} f_1(x+t)\,f_2(x+2t)\,f_3(x+3t)\,\frac{dt}{t}\,,$$
with the latter representing a major topic of interest in harmonic analysis---for more on this see Remark \ref{openai1} as well as a detailed account for the relevance of $THT$ in Section 1.3.1, \cite{HuLie}.
\smallskip

\item  {\textsf{[$\BHT$ and Biest behavior incorporated]:}}  Inspecting \eqref{hthtrep}, we observe a twisted double $\BHT$-structure interacting via curvature: indeed, on the one hand we have the $(-t,t)$ structure associated with the input pair $(f_1,f_2)$, and, on the other hand, we have the $(t^3,-t^3)$ structure associated with the input pair $(f_2,f_3)$. This type of double $\BHT$-structure is reminiscent of the {\it Biest structure}
\begin{equation}\label{Bies}
\iiint_{\xi_1 < \xi_2 < \xi_3} \hat{f}_1(\xi_1) \hat{f}_2(\xi_2) \hat{f}_3(\xi_3) e^{2\pi i x (\xi_1 + \xi_2 + \xi_3)} \, d\xi_1 \, d\xi_2 \, d\xi_3=
\text{p.v.} \iint_{\mathbb{R}^2} f_1(x-t) \, f_2(x+t+s) \, f_3(x-s) \, \frac{dt}{t} \, \frac{ds}{s} 
\end{equation} 
    with its classical version above studied in \cite{MTTBiest1,MTTBiest2}, and a further, recent oscillatory version completed in \cite{BO}. However, in contrast with these works, our $\THT$ operator adds two extra layers of difficulty:
 \begin{itemize} 
\smallskip 
 \item \emph{it is significantly more singular}: the integral kernel in \eqref{hthtrep} is $\int_{\R} \left(\ldots\right)\,\frac{dt}{t}$ as opposed to the additional smoothing effect of $\iint_{\R^2} \left(\ldots\right)\,\frac{dt}{t}\,\frac{ds}{s}$ in \eqref{Bies};
 \smallskip  
\item \emph{it has an anisotropic character}: recalling the antithetical discussion when comparing \eqref{anbiest} with \eqref{eq:middle:biest}, the anisotropy makes the geometry of the time-frequency analysis significantly more involved---see the extended discussion in Section \ref{Keyideas} and in particular the necessity of involving both a low and high resolution analysis corresponding to super-Heisenberg and Heisenberg tiles, respectively.
 \end{itemize}  
 \end{itemize}

Thus, in a nutshell, the operator $\THT$ is the first example of a trilinear operator having a level of complexity which is compatible with all of the following features:
\begin{itemize}
\item a  \emph{hybrid character} combining both zero and non-zero curvature features;

\item the \emph{same degree of kernel singularity} as the trilinear Hilbert transform;

\item an \emph{anisotropic nature} while obeying a modulation invariance structure;

\item \emph{subsumes} the behaviors of $\BHT$, Biest, $BC^a$ and, in part, the purely curved trilinear Hilbert transform $ H_{3,(1,2,3)}$.
\end{itemize}
\medskip

\subsubsection{Further motivation and connections---an outline} Beyond its immediate relevance to harmonic analysis, the framework developed in this paper is conceptually motivated by certain related phenomena in fields such as ergodic theory,  additive combinatorics and number theory. While we do not establish a formal dictionary between these settings, the structural parallels suggest that our methods might offer a useful heuristic perspective for studying the problems summarized in what follows. However, a rigorous extension of the ideas in our present paper to the broader settings briefly\footnote{For a more detailed overview, we refer the reader to the extended introductions of \cite{HuLie}.} described  below lies beyond the scope of the present work.

To keep our exposition focused, in what follows, we state the relevant related problems only briefly, bypassing their historical development and associated bibliographical references. We refer the reader to Section 1.3.3  of \cite{HuLie} for a more thorough discussion. These being said, we have:

\begin{itemize}
\item \textbf{(A)} \textsf{Ergodic theory setting:} Let $(X, \mu, T)$ be a $(\sigma-)$finite measure preserving system with $T:\:X\,\rightarrow\,X$ an invertible bimeasurable map. For any $k,\,N\in\N$, any polynomials\footnote{Here $\Z[n]$ designates the space of all polynomials $P(n)$ in one indeterminate $n$ with
integer coefficients.} $P_1,\ldots, P_k\in\Z[n]$ and any measurable functions $f_1,\ldots, f_k\in L^{0}(X)$ define the associated non-conventional polynomial multiple ergodic average as
\begin{eqnarray}\label{Nonconvmultergodav}
A_{N}^{P_1,\ldots,\,P_k}(f_1,\ldots,\,f_k)(x):=\frac{1}{N} \sum_{n=1}^N f_1(T^{P_1(n)}x)\ldots f_k(T^{P_k(n)}x)\,.
\end{eqnarray}
The \emph{main goal} here concerns the almost everywhere convergence of these averages, with the latter playing an important role in detecting recurrence point patterns.
\smallskip

\item \textbf{(B)} \textsf{Number theory/additive combinatorics setting:} Let as before $k\in \N$ and $P_1,\ldots,\,P_k$ polynomials (with no constant terms) in $K[y]$ where here one may take $K$ as $\Z$, a finite field $F_q$ or $\R$. The \emph{main goal} in this setting is to identify qualitative and/or quantitative conditions on a set $S\subset K$ that guarantee the existence of and/or lower bounds for polynomial progressions $x,\,x+P_1(y),\ldots, x+P_k(y)$ lying within the set $S$. 
\end{itemize}

Both the above topics have a rich history, and, to our knowledge, both these problems are open for the case $k=3$ in the situation in which the polynomials are \emph{allowed to be linearly dependent} (in a nontrivial manner) hence with \emph{no distinct degree requirement}.

We conclude our introductory section by noticing that our main results in Theorem \ref{thm:main} and Corollary \ref{extens} represent the singular integral Euclidean analogue of the two topics above precisely in \emph{the $k=3$ case with polynomial linear dependence allowed}.

\subsection*{Acknowledgments}

The third author was supported by NSF grant DMS-2400238 and by the Simons Travel grant MPS-TSM-00008072. He also wishes to express his gratitude to the University of Nantes for providing a nice research environment during his several visits to collaborate with his coauthors. 

This project has been in progress over several years, and the authors would like to warmly thank both the Centre International de Rencontres Mathématiques ({\it CIRM}) in Marseille and the {\it Centro di Ricerca Matematica Ennio De Giorgi} in Pisa. Each center hosted the authors for a productive ``{\it Research in Pairs}'' week, offering exceptional conditions that greatly facilitated our work on this project.




\section{A preliminary discretization of $\THT$: the main three components}
\label{sec:discretization1}

In this section we provide an outline for a preliminary discretization of $\THT$ and briefly explain the strategy to be developed in order to prove the main result of our paper, that is, Theorem \ref{thm:main}. Before this though, we make a short detour for setting the context for some of the frequently used notation.

\subsection{Preparatives: notations}

For an easy reference, we list a few notations that will appear throughout the current work.

Let $N\gg 1$ be a large enough natural number.\footnote{Throughout the paper $N$ is allowed to vary from line to line.} For an interval $I\subset \mathbb{R}$ with center $c_I$, $\one_I$ denotes its characteristic function, $\chi_I$ stands for the non-local version, defined by
$$\chi_I(x):=\Big( 1+\frac{\dist(x,I)}{|I|}\Big)^{-N},$$
and $\ci_{I}$ its smooth version
$$\ci_I(x):=\Big( 1+\frac{|x-c_I|^2}{|I|^2}\Big)^{-N/2}\,.$$

Given $j,k\in\Z$ we set $J_j^k:= [j 2^{-k},(j+1) 2^{-k}]$ and let $J^k:=J_1^k$. To describe the refined physical intervals taking into account the dilation parameter $m$ (or $3m$), given $a, m >0$ fixed and $q\in\Z$ we further define $\ds I_q^{am,k}:=[q \,2^{-\frac{am}{2}-k} ,(q+1) 2^{-\frac{am}{2}-k}]$. Notice that $J_j^k=I_j^{0,k}$. 

We often use $\approx$ to mean equality up to a superposition\footnote{Series with fast decaying coefficients.} of similar terms; this will be used especially when we perform certain decompositions.

For two positive quantities $A$ and $B$, we write $A \simeq B$ whenever there exists $C>0$ an universal constant such that $C^{-1} \cdot A \leq B \leq C \cdot A$.

For an integer variable appearing in a summation, we write $\ell \sim a$ to say that the variable $\ell$ takes about $a$ integer values within an interval of the form $[a(1+b), a(2+b)]$, where $b$ is a fixed parameter whose value is clear from the context, although we do not always make it explicit. In most of the cases $b=0$.

For $p$ a Lebesgue exponent with $0<p\leq \infty$, its dual exponent is denoted $p'$. Notice that $p'$ can be negative if $p \leq 1$, which corresponds to quasi-Banach $L^p$ spaces.

\subsection{Preliminary discretization}\label{preldiscret}
Isolating the kernel of the operator $\THT$, we perform a Whitney decomposition of the kernel relative to the $0$ and $\infty$ singularities:
\begin{equation}\label{tdiscret}
\frac{1}{t}=\sum_{k\in\Z} \rho_k(t)=:\sum_{k\in\Z} 2^k\rho(2^{k} t),
\end{equation}
where here $\rho\in C^\infty_0({\mathbb R})$ is a smooth function compactly supported away from the origin and having the cancellation property $\int_{\mathbb R} \rho(t) dt=0$. 

With these, we have
\begin{equation}\label{op-de}
 \THT (f_1,f_2,f_3)=\sum_{k\in\Z} T_k(f_1,f_2,f_3)\,,
\end{equation}
where
\begin{equation}\label{op1}
T_k(f_1,f_2,f_3)(x):=\int_{\R} f_1(x-t)\,f_2(x+t+t^3)\,f_3(x-t^3)\, 2^k \rho(2^{k} t)\,dt\,.
\end{equation}
Dualizing the operators in \eqref{op-de}--\eqref{op1} we have the quadrilinear form relations
\begin{align}
\L(f_1,f_2,f_3,f_4)&:=\langle \THT(f_1,f_2,f_3),f_4\rangle \nonumber \\
& =\L_{+}(f_1,f_2,f_3,f_4)\,+\L_{-}(f_1,f_2,f_3,f_4)\,,
\label{fm}
\end{align}
where
\begin{equation}\label{fmd}
\L_{\pm}(f_1,f_2,f_3,f_4):=\sum_{k\in\Z_{\pm}} \L_k(f_1,f_2,f_3,f_4)\,,
\end{equation}
with
\begin{equation}\label{fm1}
\L_k(f_1,f_2,f_3,f_4):=\int_{\R}\int_{\R} f_1(x-t)\,f_2(x+t+t^3)\,f_3(x-t^3)\,f_4(x)\, 2^{k}\rho(2^{k} t)\,dt\,dx\,.
\end{equation}
Further on, by expressing \eqref{fm1} in a multiplier form, we get
\begin{equation}\label{fm2}
\L_k(f_1,f_2,f_3,f_4)=\iiint_{\R^3} \hat{f}_1(\xi_1)\,\hat{f}_2(\xi_2)\,\hat{f}_3(\xi_3)\,\hat{f}_4(-\xi_1-\xi_2-\xi_3)\, \mathfrak{m}_{k}(\xi_1,\xi_2,\xi_3)\,d\xi_1\,d\xi_2\,d\xi_3\,,
\end{equation}
with
\beq\label{multk}
\mathfrak{m}_{k}(\xi_1,\xi_2,\xi_3):= \int_{\R}e^{-i\,\frac{\xi_1-\xi_2}{2^k}\,t}\,e^{-i\,\frac{\xi_3-\xi_2}{2^{3k}}\,t^3}\,\rho(t)\,dt.
\eeq
Following the natural steps, we are going to separate our analysis according to the stationary/non-stationary behavior of the phase associated to the multiplier above, that is
\beq\label{ph}
\lambda(t):=-\frac{\xi_1-\xi_2}{2^k}t - \frac{\xi_3-\xi_2}{2^{3k}}t^3.
\eeq
We thus employ a Littlewood-Paley decomposition for a suitable smooth function $\psi\in C^\infty_0([-4, -1] \cup [1,4])$
$$ 1 = \sum_{m_1,m_2 \in{\mathbb Z}} \psi\left(\frac{\xi_1-\xi_2}{2^{k+m_1}}\right)\,\psi\left(\frac{\xi_3-\xi_2}{2^{3k+m_2}}\right),$$
and plug the latter into the integral formula of the symbol $\mathfrak{m}_k$:
\begin{align}
\mathfrak{m}_{k}(\xi_1,\xi_2,\xi_3) & = \sum_{m_1,m_2\in{\mathbb Z}} \int_{\R} e^{-i\,\frac{\xi_1-\xi_2}{2^k}\,t}\,e^{-i\,\frac{\xi_3-\xi_2}{2^{3k}}\,t^3} \psi\left(\frac{\xi_1-\xi_2}{2^{k+m_1}}\right)\,\psi\left(\frac{\xi_3-\xi_2}{2^{3k+m_2}}\right) \,\rho(t)\,dt \nonumber\\
& := \sum_{m_1,m_2 \in {\mathbb Z}} \mathfrak{m}_{k,m_1,m_2}(\xi_1,\xi_2,\xi_3).\label{multm}
\end{align}

Now, restricted to each symbol $\mathfrak{m}_{k,m_1,m_2}$, we have:
\begin{itemize}
\item $\lambda$ has a critical point within the support of $\rho$ iff $|m_1-m_2|\leq 3$;
\item if $|m_1-m_2|> 3$ then $\lambda$ has a magnitude of order $\max\{2^{m_1},2^{m_2}\}$.
\end{itemize}

Directed by the above behavior we are invited to consider the following decomposition:
\begin{align*}
\mathfrak{m}_{k} & = \Big(\sum_{m_1,m_2 < 0} \mathfrak{m}_{k,m_1,m_2}\Big) \\
& + \Big(\sum_{m_1< 0 \leq m_2} \mathfrak{m}_{k,m_1,m_2} + \sum_{m_2< 0 \leq m_1} \mathfrak{m}_{k,m_1,m_2} + \sum_{m_1,m_2\geq 0 \atop |m_1-m_2|> 3} \mathfrak{m}_{k,m_1,m_2}\Big) \\
& + \Big(\sum_{m_1,m_2\geq 0 \atop |m_1-m_2|\leq 3} \mathfrak{m}_{k,m_1,m_2}\Big) \\
 & := \mathfrak{m}_{k}^{\low} + \mathfrak{m}_{k}^{\nonsta} + \mathfrak{m}_k^{\sta}\,.
\end{align*}

As expected, with the obvious correspondences, we associate to these symbols the operators $T_{k,m_1,m_2}$, $T_k^{\low}$, $T_k^{\sta}$, $T_k^{\nonsta}$ and further define
\begin{equation}
\label{eq:def:low:nonsta:sta}
T^{\low}:=\sum_{k \in \Z}T_k^{\low}, \qquad T^{\sta}:=\sum_{k \in \Z}T_k^{\sta} \quad \text{and} \quad  T^{\nonsta}:=\sum_{k \in \Z}T_k^{\nonsta}.
\end{equation}
Thus, consistent with the description given in \eqref{lof}--\eqref{histf}, we have achieved the desired decomposition of our main operator into three main components: low frequency, high frequency non-stationary, and, high frequency stationary:
\begin{equation}\label{dectht}
 \THT=\ds T^{\low}+ T^{hi}=T^{\low}+ T^{\nonsta}+T^{\sta}\,.
\end{equation}

\subsection{The three main components: statements}

With all these settled, our main result, Theorem \ref{thm:main}, is a direct consequence of the following:

\begin{theorem} [\textsf{Low frequency}]\label{thm:low} If $1<p_1,p_2,p_3< \infty$ with $(\frac{1}{p_1},\frac{1}{p_2},\frac{1}{p_3})\in{\bm D}$ and $p$ given by \eqref{eq:holder} then 
$$ \|T^{\low}\|_{L^{p_1} \times L^{p_2} \times L^{p_3} \rightarrow L^p} \lesssim 1.$$
\end{theorem}

\begin{theorem} [\textsf{High frequency non-stationary}] \label{thm:non-sta} If $1<p_1,p_2,p_3< \infty$ with $(\frac{1}{p_1},\frac{1}{p_2},\frac{1}{p_3})\in{\bm D}$ and $p$ given by \eqref{eq:holder} then 
$$ \| T^{\nonsta}\|_{L^{p_1} \times L^{p_2} \times L^{p_3} \rightarrow L^p} \lesssim 1.$$
\end{theorem}

\begin{theorem}  [\textsf{High frequency stationary}] \label{thm:sta} If $1<p_1,p_2,p_3 < \infty$ with $0<\frac{1}{p}:=\frac{1}{p_1}+\frac{1}{p_2}+\frac{1}{p_3}<1$ then
$$ \| T^{\sta}\|_{L^{p_1} \times L^{p_2} \times L^{p_3} \rightarrow L^p} \lesssim 1.$$
\end{theorem}

We end this section with the following comment:

\begin{itemize}
\item  Theorems \ref{thm:low} and \ref{thm:non-sta} rely on Theorem \ref{thm:main-ani:biest}. All these results will be proved in Part I of our paper. 

\item Theorem \ref{thm:sta} encapsulates the interaction between high oscillation\footnote{Associated with the phase of the multiplier.} and time-frequency analysis and involves a dual multiresolution analysis:
\begin{itemize}  
\item a high multiresolution analysis addressing the modulation invariance feature of $T^{\sta}$ and involving a BHT-type analysis of suitable defined super-Heisenberg tiles;

\item a low multiresolution analysis addressing the non-zero curvature feature of $T^{\sta}$ and capturing the cancellation within each single super-Heisenberg tile via Rank II LGC method. 
\end{itemize}
    The proof of Theorem \ref{thm:sta} will be delivered in Part II of our paper.
\end{itemize}

\part[I]{THE BOUNDEDNESS OF THE ANISOTROPIC BIEST OPERATOR AND DERIVED COROLLARIES. TREATMENT OF $T^{\low}$ AND $T^{\nonsta}$}
\label{partI}
In this first part, our main focus will be the proof of Theorem \ref{thm:main-ani:biest}, a result which is of independent interest. From this, we will then deduce Theorem \ref{thm:ani}, Theorem \ref{thm:anisotropic:Mikhlin}, Theorem \ref{thm:low} and Theorem \ref{thm:non-sta}.

\section{Preliminary reductions}
\label{sec:anis:model:red}

In this section we have the following goals:
\begin{itemize}
\item In Section \ref{subsec:anis:model:red} we show that via a suitable discretization applied to Mikhlin symbols, Theorem \ref{thm:anisotropic:Mikhlin} follows from Theorem \ref{thm:ani} which in turn, follows from Theorem \ref{thm:main-ani:biest}.
\item In Section \ref{typeI1} we provide an elementary treatment to a special family of Type I multipliers---recall \eqref{eq:symp:aniso}---for which both functions involved in the summand are of $\Psi$-type; 
\item Finally, in Section \ref{reduct}, we reduce the treatment of both the generic Type I and Type II multipliers as well as of Theorem \ref{thm:main-ani:biest} to Theorem \ref{thm:proof:main:1} and Theorem \ref{thm:proof:main:2}. 
\end{itemize}

\subsection{Reduction of Theorem \ref{thm:anisotropic:Mikhlin} and Theorem \ref{thm:ani} to Theorem \ref{thm:main-ani:biest}.}
\label{subsec:anis:model:red}

Recalling \eqref{eq:Mikhlin}, we know---see \emph{e.g.} \cite[Chapter 2.13]{multilinear_harmonicII}---that a classical Mikhlin symbol ${\mathbbm{m}}(\eta_1, \eta_2)$ can be written as
\[
\sum_{k \in \Z} {\mathbbm{m}}(\eta_1, \eta_2) \phi_1 \Big(  \frac{\eta_1}{2^k} \Big) \phi_2 \Big(  \frac{\eta_2}{2^k} \Big),
\]
where at least one of the functions $\phi_1$ and $\phi_2$ is of $\Psi$-type, supported away from $0$, and the other function is either of $\Psi$-type or $\Phi$-type, supported around the origin. We note that in the case when $\mathbbm{m}(\eta_1, \eta_2)$ is supported away from the line $\eta_2=0$ the function $\phi_2$ is of $\Psi$-type, supported away from $0$. 

Implementing this decomposition, $\ds {\mathbbm{m}}(\eta_1^3, \eta_2)$ can be written as $\ds \sum_{k \in \Z} {\mathbbm{m}}(\eta_1^3, \eta_2) \phi_1 \Big(  \frac{\eta_1^3}{2^k} \Big) \phi_2 \Big(  \frac{\eta_2}{2^k} \Big)$. If we set $\tilde \phi_1(\eta)=\phi_1(\eta^3)$ (thus $\tilde \phi_1$ is smooth and of the same type as $\phi_1$), this becomes morally equivalent to
\[
\sum_{\kappa \in \Z} {\mathbbm{m}}(\eta_1^3, \eta_2) \tilde{\phi}_1 \Big(  \frac{\eta_1}{2^{\kappa}} \Big) \phi_2 \Big(  \frac{\eta_2}{2^{3 \kappa}} \Big).
\]
Next, we decompose each term in the sum into 2-variable Fourier series on rectangles of size $2^{\kappa} \times 2^{3 \kappa}$ in order to represent the expression above as a sum of terms of the form
\begin{align}
\label{eq:symbol:aniso:Four}
\sum_{j_1, j_2 \in \Z} \sum_{\kappa \in \Z} \, c^{j_1, j_2}_{\kappa} \, \tilde{\phi}_1 \Big(  \frac{\eta_1}{2^{\kappa}} \Big)  e^{i \frac{\eta_1}{2^{\kappa}} j_1} \phi_2 \Big(  \frac{\eta_2}{2^{3 \kappa}} \Big) e^{i \frac{\eta_2}{2^{3\kappa}} j_2}.
\end{align}
As a consequence of the Mikhlin property satisfied by $\mathbbm{m}$, we notice that the Fourier coefficients $c^{j_1, j_2}_{\kappa}$ decay uniformly in $\kappa$ and that, for any $M \geq 1$, one has
\begin{equation*}
\big|c^{j_1, j_2}_\kappa\big| \lesssim_M (1+|j_1|+|j_2|)^{-M}.
\end{equation*}
Our goal is to reduce matters to the situation when we have a $\Psi$-type function associated to the dominating frequency. We look for example at the situation in which $\tilde \phi_1$ is of $\Psi$-type, $\phi_2$ is of $\Phi$-type and $\kappa \geq 0$ so that we need to deal with
\begin{align}
\label{eq:symbol:aniso:Four:+}
\sum_{j_1, j_2 \in \Z} \sum_{\kappa \geq 0} \, c^{j_1, j_2}_{\kappa} \, \tilde{\phi}_1 \Big(  \frac{\eta_1}{2^{\kappa}} \Big)  e^{i \frac{\eta_1}{2^{\kappa}} j_1} \phi_2 \Big(  \frac{\eta_2}{2^{3 \kappa}} \Big) e^{i \frac{\eta_2}{2^{3\kappa}} j_2}.
\end{align}
We observe that 
\[
\phi_2 \Big(  \frac{\eta_2}{2^{3 \kappa}} \Big) e^{i \frac{\eta_2}{2^{3\kappa}} j_2}= \phi_2 \Big(  \frac{\eta_2}{2^{3 \kappa}} \Big) \big( \, e^{i \frac{\eta_2}{2^{3\kappa}} j_2}- 1\big) + \phi_2 \Big(  \frac{\eta_2}{2^{3 \kappa}} \Big).
\] 
Now, the first term is of $\Psi$-type, and for the second term we express the $\Phi$-type function as a sum of $\Psi$-type functions:
\[
\phi_2 \Big(  \frac{\eta_2}{2^{3 \kappa}} \Big)= \phi_2 \Big(  \frac{\eta_2}{2^{\kappa}} \Big) + \sum_{j=k}^{3k-1}\left( \phi_2 \Big(  \frac{\eta_2}{2^{j+1}}  \Big)- \phi_2 \Big(  \frac{\eta_2}{2^{j}} \Big)\right):= \phi_2 \Big(  \frac{\eta_2}{2^{\kappa}} \Big) + \sum_{j=k}^{ 3k-1} \psi_2 \Big(  \frac{\eta_2}{2^{j}} \Big).
\]
With this, the expression in \eqref{eq:symbol:aniso:Four:+} becomes
\begin{align*}
&\sum_{j_1, j_2 \in \Z} \sum_{\kappa \geq 0} \, (1+|j_2|)^{N} c^{j_1, j_2}_{\kappa} \, \tilde{\phi}_1 \Big(  \frac{\eta_1}{2^{\kappa}} \Big)  e^{i \frac{\eta_1}{2^{\kappa}} j_1} (1+|j_2|)^{-N} \phi_2 \Big(  \frac{\eta_2}{2^{3 \kappa}} \Big) \big( e^{i \frac{\eta_2}{2^{3\kappa}} j_2}-1 \big) \\
+ & \sum_{j_1, j_2 \in \Z} \sum_{\kappa \geq 0} \, c^{j_1, j_2}_{\kappa} \, \tilde{\phi}_1 \Big(  \frac{\eta_1}{2^{\kappa}} \Big)  e^{i \frac{\eta_1}{2^{\kappa}} j_1}  \phi_2 \Big(  \frac{\eta_2}{2^{ \kappa}} \Big)  \\
+ &  \sum_{j_1, j_2 \in \Z} \sum_{\kappa \geq 0} \, c^{j_1, j_2}_{\kappa} \, \tilde{\phi}_1 \Big(\frac{\eta_1}{2^{\kappa}} \Big)  e^{i \frac{\eta_1}{2^{\kappa}} j_1} \Big(  \sum_{j=\kappa}^{ 3 \kappa-1} \psi_2 \Big(  \frac{\eta_2}{2^{j}} \Big)\Big). 
\end{align*}
The first term corresponds (modulo superpositions) to the situation in which two $\Psi$-type functions with uniform bounded derivatives appear. The second term is an admissible isotropic symbol and for the last term we interchange the summation in $\kappa$ and $j$ to get
\begin{align*}
\sum_{j_1, j_2 \in \Z} \,  \sum_{j \geq 0}  \psi_2 \Big(  \frac{\eta_2}{2^{j}} \Big) \Big( \sum_{\kappa>j/3}^j \, c^{j_1, j_2}_{\kappa} \, \tilde{\phi}_1 \Big(  \frac{\eta_1}{2^{\kappa}}\Big)  \, e^{i \frac{\eta_1}{2^{\kappa}} j_1} \Big).
\end{align*}
Invoking the boundedness of the Biest operator \cite{MTTBiest2}, the above is reduced to the study of 
\begin{align*}
\sum_{j_1, j_2 \in \Z} \,  \sum_{j \geq 0}  \psi_2 \Big(  \frac{\eta_2}{2^{j}} \Big) \Big( \sum_{0 \leq \kappa \leq {j \over 3}} \, c^{j_1, j_2}_{\kappa} \, \tilde{\phi}_1 \Big(  \frac{\eta_1}{2^{\kappa}}\Big)  \, e^{i \frac{\eta_1}{2^{\kappa}} j_1} \Big),
\end{align*}
so that we find ourselves in the situation in which the $\Psi$-type function $\psi_2$ is associated to the dominating scale. A similar argument allows us to deal with the summation over $\kappa \leq 0$.

Thus $\mathbbm{m}(\eta_1^3, \eta_2)$ can be expressed as a superpositions of terms which belong to one of the following four categories, in the increasing order of complexity:
\begin{itemize}[leftmargin=2em]
\item (A) The $2D$ isotropic case, with at least one of the functions is of $\Psi$-type:
$$\ds \sum_{\kappa \in \Z}  \, c_{\kappa} \, \phi_1 \Big(  \frac{\eta_1}{2^{\kappa}} \Big) \, \phi_2 \Big(  \frac{\eta_2}{2^{\kappa}} \Big)\,.$$
\item (B) The $2D$ jointly vanishing Type I multiplier case, when both $\psi_1$  and $\psi_2$ are of $\Psi$-type:
$$\ds \sum_{\kappa \in \Z}  \, c_{\kappa} \, \psi_1 \Big(  \frac{\eta_1}{2^{\kappa}} \Big) \, \psi_2 \Big(  \frac{\eta_2}{2^{3 \kappa}} \Big)\,.$$
\item (C) The $2D$ generic Type I multiplier case\footnote{In this situation, only the function associated with the largest frequency is assumed to be of $\Psi$-type.}:
$$\ds \sum_{\kappa \geq 0} \, c_{\kappa} \, \varphi_1 \Big(  \frac{\eta_1}{2^{\kappa}} \Big) \, \psi_2 \Big(  \frac{\eta_2}{2^{3 \kappa}} \Big)\,+\,\ds \sum_{\kappa \leq 0} \, c_{\kappa} \, \psi_1 \Big(  \frac{\eta_1}{2^{\kappa}} \Big) \, \varphi_2 \Big(  \frac{\eta_2}{2^{3 \kappa}} \Big)\,.$$
\item (D) The $2D$ generic Type II multiplier case: 
 $$\ds \sum_{j \geq 0} \, \Big[  \sum_{j/3 \leq \kappa \leq j}  \, c_{\kappa} \, \psi_1 \Big(  \frac{\eta_1}{2^{\kappa}}\Big)   \Big] \,  \psi_2 \Big(  \frac{\eta_2}{2^{j}} \Big)\,+\,\ds \sum_{j \leq 0} \,  \psi_1 \Big(  \frac{\eta_1}{2^{j}} \Big) \, \Big[  \sum_{j \leq \kappa \leq j/3}  \, c_{\kappa} \, \psi_2 \Big(  \frac{\eta_2}{2^{3 \kappa}}\Big)   \Big]\,.$$ 
\end{itemize}

Applying all of the above to the setting described within the statement of Theorem \ref{thm:anisotropic:Mikhlin}, we notice that understanding the symbol $\mathbbm{m}((\xi_1-\xi_2)^3, \xi_3-\xi_2)$ reduces to studying one of the following four cases corresponding to the ones above:
\begin{itemize}[leftmargin=2em]
\item (A1) The isotropic case (with at least one of the functions is of $\Psi$-type):
$$\ds \sum_{\kappa \in \Z}  \, c_{\kappa} \, \phi_1 \Big(  \frac{\xi_1-\xi_2}{2^{\kappa}} \Big) \, \phi_2 \Big(  \frac{\xi_3-\xi_2}{2^{\kappa}} \Big)\,.$$
\item (B1) The jointly vanishing Type I multiplier case (with both functions $\psi_1$  and $\psi_2$ of $\Psi$-type):
$$\ds \sum_{\kappa \in \Z}  \, c_{\kappa} \, \psi_1 \Big(  \frac{\xi_1-\xi_2}{2^{\kappa}} \Big) \, \psi_2 \Big(  \frac{\xi_3-\xi_2}{2^{3 \kappa}} \Big)\,.$$ 
\item (C1) The generic Type I multiplier case:
 $$\ds \sum_{\kappa \geq 0} \, c_{\kappa} \, \varphi_1 \Big(  \frac{\xi_1-\xi_2}{2^{\kappa}} \Big) \, \psi_2 \Big(  \frac{\xi_3-\xi_2}{2^{3 \kappa}} \Big) \,+\,\ds \sum_{\kappa \leq 0} \, c_{\kappa} \, \psi_1 \Big(  \frac{\xi_1-\xi_2}{2^{\kappa}} \Big) \, \varphi_2 \Big(  \frac{\xi_3-\xi_2}{2^{3 \kappa}} \Big)\,. $$
\item (D1) The generic Type II multiplier case:
$$\ds \sum_{j \geq 0} \, \Big[  \sum_{j/3 \leq \kappa \leq j}  \, c_{\kappa} \, \psi_1 \Big(  \frac{\xi_1-\xi_2}{2^{\kappa}}\Big)   \Big] \,  \psi_2 \Big(  \frac{\xi_3-\xi_2}{2^{j}} \Big)\,+\,\ds \sum_{j \leq 0} \,  \psi_1 \Big(  \frac{\xi_1-\xi_2}{2^{j}} \Big) \, \Big[  \sum_{j \leq \kappa \leq j/3}  \, c_{\kappa} \, \psi_2 \Big(  \frac{\xi_3-\xi_2}{2^{3 \kappa}}\Big)   \Big]\,.$$ 
\end{itemize}
We now have that Theorem \ref{thm:anisotropic:Mikhlin} may be reduced to Theorem \ref{thm:ani} which further follows from Theorem \ref{thm:main-ani:biest}, based on the following steps:
\begin{itemize}
\item (A2) The symbol in $(A1)$ is singular with respect to the subspace $\Gamma':= (1,1,1)\cdot{\mathbb R}$ and hence it gives rise to an isotropic trilinear operator such as the one studied in \cite{Mlinop}. The latter operator satisfies $L^{p_1}\times L^{p_2}\times L^{p_3}\to L^p$ estimates, where $1<p_1, p_2, p_3, p \leq \infty$, $\ds \frac{1}{p_1}+\frac{1}{p_2}+\frac{1}{p_3}=\frac{1}{p}$ and $0<\frac{1}{p} < \frac{5}{2}$.

\item (B2) The operator associated to the symbol in $(B1)$ is bounded inside the Banach range, thanks to Proposition \ref{prop:facile} treated in the following section.

\item (C2) The operator associated to the symbol in $(C1)$ is bounded inside the local $L^2$ range as a consequence of Theorem \ref{thm:proof:main:1}. 
Moreover---see Remark \ref{rem:42}---the boundedness range for both operators in (B2) and (C2) can be extended via an interpolation procedure relative to  ${\bm D}'$ with the latter range derived from employing Theorem \ref{thm:main-ani:biest}.

\item (D2) Finally, the operator with symbol given by $(D1)$ is treated by Theorem \ref{thm:main-ani:biest}, modulo some classical isotropic \emph{biest}-type terms which are bounded in a wider range \cite{Joeun-IteratedTrilFourierIntg, MTTBiest2}.
\end{itemize}

\subsection{Jointly vanishing Type I multipliers}\label{typeI1}


In this section we present a special case of Theorem \ref{thm:ani}, corresponding to the situation when both  functions $\phi_1$ and $\phi_2$ are of $\Psi$-type. This result is going to be particularly useful in Section \ref{sec:low}, which is why we state it in the form needed therein. As it turns out, for the special case treated here we can prove bounds covering the full Banach range.

\begin{proposition}[\textsf{Jointly vanishing Type I multipliers}]
\label{prop:facile}
Let $(c^{m_1, m_2}_k)_{k \in \Z}$ be a bounded sequence of complex numbers with  
$$\ds C:= \sup_{m_1, m_2 \in \Z}\|(c^{m_1, m_2}_k)_{k \in \Z}\|_{\ell^\infty} <+\infty\,.$$ 
Let $\phi_1$ and $\phi_2$ be two $\Psi$-type functions with uniform bounded derivatives.\footnote{One can see from the proof that it is enough to have uniform bounds for the first 4 derivatives of $\phi_1$ and $\phi_2$, as this guarantees enough decay of the Fourier coefficients appearing in the decompositions we perform.} For $m_1, m_2$ integer parameters, we let $T_m$ be the trilinear operator associated to the symbol
\begin{align*}
m(\xi_1, \xi_2, \xi_3):= \sum_{k \in \Z} c^{m_1, m_2}_k \, \phi_1\Big( \frac{\xi_1-\xi_2}{2^{k+m_1}} \Big) \, \phi_2\Big( \frac{\xi_3-\xi_2}{2^{3k+m_2}} \Big) .
\end{align*}
Then for all\footnote{In fact, we cannot have $p_1=p_2=\infty$, nor $p_3=p'=\infty$.} Lebesgue exponents $1<p_1, p_2, p_3, p' \leq \infty$ with $\frac{1}{p_1}+\frac{1}{p_2}+\frac{1}{p_3}=\frac{1}{p}$, $T_m$ is bounded from $L^{p_1} \times L^{p_2} \times L^{p_3}$ into $L^p$, with an operator norm bounded by $C$.
\begin{proof}
As in \cite{MTTBiest2}, the strategy depends whether $|\xi_1-\xi_2| \leq |\xi_3-\xi_2|$ or  $|\xi_3-\xi_2| \leq |\xi_1-\xi_2|$. Due to symmetry, we only analyze the first situation, which arises when $\ds k \geq \frac{m_1-m_2}{2}$. Then the idea is to replace $\ds \phi_2\Big( \frac{\xi_3-\xi_2}{2^{3k+m_2}} \Big)$ by $\phi_2\Big( \frac{\xi_3-\frac{\xi_1+\xi_2}{2}}{2^{3k+m_2}} \Big)$, up to an error term that needs to be estimated. Indeed, we notice that
\begin{align*}
\phi_2\Big( \frac{\xi_3-\xi_2}{2^{3k+m_2}} \Big)= \phi_2\Big( \frac{\xi_3-\frac{\xi_1+\xi_2}{2}}{2^{3k+m_2}} \Big) + \frac{\xi_1-\xi_2}{2^{3k+m_2+1}} \, \int_0^1 \phi_2' \Big( \frac{\xi_3-\xi_2}{2^{3k+m_2}}  \, \tau + \frac{\xi_3-\frac{\xi_1+\xi_2}{2}}{2^{3k+m_2}} (1-\tau) \Big) d\tau.
\end{align*}

Thus 
\begin{align*}
\sum_{k \geq \frac{m_1-m_2}{2}} c^{m_1, m_2}_k \, \phi_1\Big( \frac{\xi_1-\xi_2}{2^{k+m_1}} \Big) \,  \phi_2\Big( \frac{\xi_3-\xi_2}{2^{3k+m_2}} \Big) = \sum_{k \geq \frac{m_1-m_2}{2}} c^{m_1, m_2}_k \, \phi_1\Big( \frac{\xi_1-\xi_2}{2^{k+m_1}} \Big) \,  \phi_2\Big( \frac{\xi_3-\frac{\xi_1+\xi_2}{2}}{2^{3k+m_2}} \Big) &\\
+ \sum_{k \geq \frac{m_1-m_2}{2}} \frac{c^{m_1, m_2}_k}{2^{2k+m_2-m_1}} \, \frac{\xi_1-\xi_2}{2^{k+m_1}} \,\phi_1\Big( \frac{\xi_1-\xi_2}{2^{k+m_1}} \Big) \, \int_0^1 \phi_2' \Big( \frac{\xi_3-\xi_2}{2^{3k+m_2}}  \, \tau + \frac{\xi_3-\frac{\xi_1+\xi_2}{2}}{2^{3k+m_2}} (1-\tau) \Big) d\tau.&
\end{align*}

The error term (the second one) is easily bounded: first, due to the extra decay, one can fix the scale $k$. Then, in the Banach range, when $p\geq 1$, one can also fix the parameter $\tau \in [0, 1]$ since one can interchange the integration in $x$ and in $\tau$. The resulting operator becomes
\begin{align*}
\iint_{\R^2} f_1\big( x-s_1- \frac{1-\tau}{2} s_2  \big) \, f_2\big( x+s_1- \frac{1+\tau}{2} s_2  \big) \, f_3\big( x- s_2 \big) \Phi^{k+m_1}_1(s_1) \Phi^{3k+m_2}_2(s_2) d s_1 ds_2,
\end{align*}  
with $\Phi^{k+m_1}_1$ and $\Phi^{3k+m_2}_2$ $L^1$-normalized functions, and it is not difficult to see that it is bounded within the Banach range.

For the main term, we consider the 4-linear form associated to the trilinear symbol, which becomes
\begin{align}
\label{eq:est:square:func:Prop:facile}
\sum_{k \geq \frac{m_1-m_2}{2}} c^{m_1, m_2}_k \, \int_{\xi_1+\ldots+\xi_4=0} \prod_{j=1}^4 \hat f_j(\xi_j)  \phi_1\Big( \frac{\xi_1-\xi_2}{2^{k+m_1}} \Big) \,  \phi_2\Big( \frac{\xi_3-\frac{\xi_1+\xi_2}{2}}{2^{3k+m_2}} \Big) d \xi_1 \ldots d \xi_4 &\\
 = \int_{\R} \sum_{k \geq \frac{m_1-m_2}{2}}    \Big(  c^{m_1, m_2}_k  \iint_{\R^2}  \hat{f_1}(\xi_1) \hat{f_2}(\xi_2) e^{ix(\xi_1+\xi_2)}  \phi_1\Big( \frac{\xi_1-\xi_2}{2^{k+m_1}} \Big)d\xi_1d\xi_2 \Big) &\\
 \cdot  \Big(  \iint_{\R^2}  \hat{f_3}(\xi_3) \hat{f_4}(\xi_4) e^{ix(\xi_3+\xi_4)}  \phi_2\Big( \frac{3\xi_3+\xi_4}{2^{3k+m_2+1}} \Big)d\xi_3 \, d\xi_4 \Big) \, dx.& \nonumber
\end{align}

Applying Cauchy-Schwarz, we decouple the information in $(f_1, f_2)$ and $(f_3, f_4)$ by making appear two bilinear BHT-type square functions. Once the information decoupled, we use Rademacher functions $(\epsilon_k)_k$ in order to linearize the square functions and reduce ourselves to bilinear operators associated to symbols
\begin{equation}
\label{eq:symbol:lin:Rad}
\sum_{k \geq \frac{m_1-m_2}{2}}  \epsilon_k \, c^{m_1, m_2}_k \, \phi_1\Big( \frac{\xi_1-\xi_2}{2^{k+m_1}} \Big) \qquad \text{and } \qquad \sum_{k \geq \frac{m_1-m_2}{2}}  \epsilon_k \, \phi_2\Big( \frac{3\xi_3+\xi_4}{2^{3k+m_2}} \Big),
\end{equation}
respectively. Since
\[
\eta \mapsto \sum_{k \geq \frac{m_1-m_2}{2}}  \epsilon_k \, c^{m_1, m_2}_k \, \phi_1\Big( \frac{\eta}{2^{k+m_1}} \Big) \qquad \text{and } \qquad \eta \mapsto \sum_{k \geq \frac{m_1-m_2}{2}}  \epsilon_k \, \phi_2\Big( \frac{\eta}{2^{3k+m_2}} \Big)
\]
are one-variable Mikhlin symbols due to the hypothesis that the smooth functions $\phi_1, \phi_2$ are of cancellative $\Psi$-type, the expressions in \eqref{eq:symbol:lin:Rad} satisfy the Mikhlin condition \eqref{eq:Mikhlin} and we deduce the boundedness of these bilinear square functions in the same range as that of the bilinear Hilbert transform, independently of the parameters $m_1, m_2$. However, since we used duality in order to make appear the two square functions, the range of boundedness of the operator $T_m$ is narrowed to the Banach range.
\end{proof}
\end{proposition}

\subsection{Control on Type I, Type II multipliers and Theorem \ref{thm:main-ani:biest} as derivations of Theorem \ref{thm:proof:main:1} and Theorem \ref{thm:proof:main:2}}\label{reduct}

We start our discussion with a result stating the local-$L^2$ boundedness of operators associated with (generic) Type I multipliers. The proof of this result will be provided in Section \ref{sec:proof:m:+}.  

\begin{theorem}
\label{thm:proof:main:1}
Let $T_m$ be a trilinear operator with a generic Type I symbol 
$$m(\xi_1,\xi_2,\xi_3)= \sum_{k \geq 0} c_k \phi\big(\frac{\xi_1-\xi_2}{2^{k}}\big) \psi\big(\frac{\xi_3-\xi_2}{2^{3k}}\big).$$
Then $T_m$ is bounded from $L^{p_1}\times L^{p_2}\times L^{p_3}$ to $L^p$, where 
$2<p_1, p_2, p_3, p' \leq \infty$ and $\ds \frac{1}{p_1}+\frac{1}{p_2}+\frac{1}{p_3}=\frac{1}{p}$. Moreover
$$ \|T_{m}\|_{L^{p_1}\times L^{p_2} \times L^{p_3} \to L^p} \lesssim \|(c_k)_k\|_{\ell^\infty}.$$
\end{theorem}

\begin{remark} \label{rem:42} Note that the Type I symbols are a particular examples of anisotropic Biest multipliers studied in Theorem \ref{thm:proof:main:2} below. Consequently, the boundedness range in Theorem \ref{thm:proof:main:1} can be extended via an interpolation procedure with the range obtained in Theorem \ref{thm:proof:main:2}.\footnote{A similar strategy will be employed later on in Section \ref{sec:bilinear:analysis} in order to deal with the multiscale analysis of $T^{\sta}$.} This improves the behavior of Type I multipliers relative to the ones discussed in Theorem \ref{thm:anisotropic:Mikhlin}. 
\end{remark}

In order to state our next result, we need to introduce some notation: we denote by ${\bm D}''$ the open interior of the convex hull of the following fifteen extremal points:
\begin{gather*}
(\frac{1}{2},0,0),\ (0, \frac{1}{2}, 0),\ (\frac{3}{4}, \frac{1}{2}, 0),\
(0, \frac{2}{3}, \frac{2}{3}),\
(\frac{1}{2}, \frac{3}{4}, 0),\
(\frac{3}{4}, 0, 0), \\
(\frac{3}{4}, 0, \frac{2}{3}), \ 
(0, \frac{3}{4}, 0),\ 
(\frac{3}{4}, \frac{1}{2}, 1),\ 
(\frac{2}{3}, 0, \frac{2}{3}),\ 
(0, \frac{3}{4}, \frac{2}{3}),\ 
(\frac{1}{2}, 0, \frac{1}{2}), \\
(0, \frac{1}{2}, \frac{1}{2}),\ 
(\frac{1}{2}, \frac{3}{4}, 1),\ 
(\frac{1}{2}, \frac{1}{2}, 1).
\end{gather*}

We are now ready to state the following

\begin{theorem}
\label{thm:proof:main:2}
Let $(c_j)_{j \in \Z}$ and $(d_k)_{k \in \Z}$ be two bounded sequences of complex numbers and let $T_m$ be a trilinear operator with symbol $m(\xi_1,\xi_2,\xi_3)$ given by
\begin{equation}\label{s1}
\sum_{k \geq 0} d_k \Big( \sum_{j=k}^{3k} c_{j} \psi\big(\frac{\xi_1-\xi_2}{2^{j}}\big) \Big) \psi\big(\frac{\xi_3-\xi_2}{2^{3k}}\big)
\end{equation}
or
\begin{equation}\label{s2}
\sum_{k \geq 0} d_k \Big(\sum_{j=0}^{k} c_{j} \psi\big(\frac{\xi_1-\xi_2}{2^{j}}\big)  \Big) \psi\big(\frac{\xi_3-\xi_2}{2^{3k}}\big).
\end{equation}
Then $T_m$ is bounded from $L^{p_1}\times L^{p_2}\times L^{p_3}$ to $L^p$, where $\ds \frac{1}{p_1}+\frac{1}{p_2}+\frac{1}{p_3}=\frac{1}{p}$ and $(\frac{1}{p_1},\frac{1}{p_2},\frac{1}{p_3})\in{\bm D}''$. Moreover
$$ \|T_{m}\|_{L^{p_1}\times L^{p_2} \times L^{p_3} \to L^p} \lesssim \|(c_j)_{j}\|_{\ell^\infty} \, \|(d_k)_{k}\|_{\ell^\infty}.$$
\end{theorem}

We conclude this section with the following

\begin{remark} \label{remaniso} i) We recall that the equivalence in terms of the boundedness behavior between the two symbols in \eqref{s1} and \eqref{s2}
is guaranteed by the boundedness of their sum. The latter represents essentially the \emph{isotropic biest} symbol in \eqref{eq:middle:biest} whose boundedness behavior covers a wider range. 

Consequently, in Section \ref{sec:proof:m:2}, we will only address Theorem \ref{thm:proof:main:2} for symbols defined by \eqref{s2}.

\noindent ii) Notice that Theorem \ref{thm:main-ani:biest} is a consequence of Theorem \ref{thm:proof:main:2} and its symmetric counterpart obtained by exchanging the role of $\xi_1,\xi_3$ when $k\leq 0$. Indeed, we have that
$$ {\bm D}'=\{\alpha=(\alpha_1,\alpha_2,\alpha_3)\,|\,\alpha\in{\bm D}'' \quad \textrm{and} \quad (\alpha_3,\alpha_2,\alpha_1) \in{\bm D}''\}.$$
\end{remark}

\section{Proof of Theorem \ref{thm:proof:main:1}: The anisotropic case involving a single collection of super-Heisenberg tiles}
\label{sec:proof:m:+}

In this section, we focus on proving Theorem \ref{thm:proof:main:1} which requires a suitable control on the trilinear operator 
\begin{equation}
\label{eq:tri:initial:1}
T(f_1, f_2, f_3)(x)= \sum_{k \geq 0} c_k \int_{\R^3} \hat f_1(\xi_1) \hat f_2(\xi_2) \hat f_3(\xi_3) \phi\big(  \frac{\xi_1-\xi_2}{2^k}  \big)  \psi\big( \frac{\xi_3-\xi_2}{2^{3k}}   \big) e^{i x (\xi_1+\xi_2+\xi_3)} d \xi_1 d \xi_2 d \xi_3,
\end{equation}
with  $\{c_k\}_{k\in\Z}\in\ell^{\infty}(\Z)$, $\psi$ of $\Psi$-type and $\phi$ of $\Phi$-type (recall Definition \ref{def:psi-type}).\footnote{Here we may take $\phi$ with $\supp\,\phi\subseteq [-4, 4]$.}

\subsection{Time-frequency terminology}
\label{subsec:terminology}

For the reader's convenience, in this section we recall some classical time-frequency definitions and notations that are useful in performing the multi-scale analysis.

Thus, $\ic D$ will denote the collection of dyadic intervals and $\ic D^{\sigma}$ the collection of $\sigma$-shifted dyadic intervals with $\ds \sigma \in \big\{ 0, \frac{1}{3}, \frac{2}{3}  \big\}$ (see for example \cite[Section 4]{MTTBiest2}).

\begin{definition}
A \emph{Heisenberg tile} is a rectangle $I \times \omega$ of area $1$, with $I$ a dyadic interval and $\omega$ a (possibly) shifted dyadic interval. More generally, a tile is a rectangle $I \times \omega$, not necessarily of area one\footnote{If a tile $P$ has the associated area greater than one we may choose to label $P$ as a super-Heisenberg tile.}, with $I \in \ic D$ and $\omega \in \ic  D^{\sigma}$.
A $4$-tile $P=(P_1, P_2, P_3, P_4)$ is a $4$-tuple of tiles having the same area. In most cases, the $P_j$ share the same spatial intervals (denoted by $I_P$, if that is the case).
\end{definition}

\begin{definition}
A \emph{wave packet} associated to a Heisenberg tile $I \times \omega$ is an $L^2$-normalized function $\phi$ (sometimes denoted $\phi_{I \times \omega}$) so that $\phi$ and its derivatives decay\footnote{We say $\phi$ is \emph{adapted} to $I$.} fast away from $I$ and $\hat \phi$ is supported on $\omega$.
\end{definition}

In Sections \ref{sec:onescale_k_large} and \ref{Diagkpsmall124} we perform a meticulous analysis allowing to obtain a suitable key single scale decay estimate. In order to simplify the notation we use the following convention:
\begin{notation}
\label{not:norms:TF:reg}
For a time frequency region $R=I_R \times \omega_R$, we consider a finite collection $\mathsf W_R$ of mutually disjoint Heisenberg tiles contained in $I_R \times \omega_R$ . If $f$ is a Schwartz function, we set 
\[
F:=\sum_{I \times \omega \in \mathsf W_R} \langle f, \phi_{I \times \omega} \rangle \,  \phi_{I \times \omega},
\]
and we use the notation 
\begin{equation}
\label{eq:notation:time:freq}
\|f\|_{L^2(I_R \times \omega_R)}= \|f^{\omega_R}\|_{L^2(I_R)}:= \big( \sum_{I \times \omega \in \mathsf W_R} | \langle f, \phi_{I \times \omega} \rangle |^2   \big)^{1 \over 2}.
\end{equation}

While the collection $\mathsf W_R$ will not be specified in such a notation, its definition will be clear from the context. Sometimes it will be convenient to only emphasize the spatial interval, case in which we will write
\[
\|f\|_{\ii L^2(I_R)}:= \|f\|_{L^2(I_R \times \omega_R)}.
\]
In most cases, $\mathsf W_R$ will consist of a union of Heisenberg tiles with spatial intervals having the same length $\ell_1$. Moreover, for an intermediate dyadic interval $J \subset I_R$ with $\ell_1\leq |J| \leq |I_R|$, we use the following convention:
\begin{equation}
\label{eq:notation:time:freq:tails}
\|f\|^2_{\tilde{\ii L}^2(J)}:= \sum_{\# \in \Z} (1+|\#|)^{-M}  \sum_{\substack{I \times \omega \in \mathsf W_R \\ I \subset J_{\#} \cap I_R}} | \langle f, \phi_{I \times \omega} \rangle |^2,
\end{equation}
where here $ J_{\#}$ designates the $\#$-shift of $J$.

Notice that $\ds \|f\|_{L^2(I_R \times \omega_R)}$ is equivalent to the $L^2$ norm of the function $F$, which may be represented as the time-frequency localization of $f$ to the region $I_R \times \omega_R$.
\end{notation}

Given $\P$ a collection of $4$-tiles, we highlight sub-collections having more structure: 
\begin{definition}
For any $j \in \{ 1,2,3,4\}$, we say that $\tr{T} \subset \P $ is a $j$-tree with top datum $(I_{\tr{T}}, \xi_{\tr{T}}=(\xi_{\tr{T}}^{(1)}, \xi_{\tr{T}}^{(2)}, \xi_{\tr{T}}^{(2)}, \xi_{\tr{T}}^{(4)}))$ provided that 
\[
P \in \tr{T} \Leftrightarrow I_P \subseteq I_{\tr{T}} \text{    and    } \xi_{\tr{T}}^{(j)} \in 3 \omega_{P_j}.
\]
In this case, we say that $\tr{T}$ is \emph{$j$-lacunary} provided that for any $P \in \tr{T}$, $\ds \xi_{\tr{T}}^{(j)} \in 3 \omega_{P_j} \setminus  \omega_{P_j}$, and we say that $\tr{T}$ is \emph{$j$-overlapping} provided that for any $P \in \tr{T}$, $\ds \xi_{\tr{T}}^{(j)} \in \omega_{P_j}$.

If $\tr{T}_1, \ldots, \tr{T}_N$ is a finite sequence of $j$-lacunary trees, we say that this sequence forms a \emph{chain of strongly $j$-disjoint trees} provided
\begin{enumerate}[leftmargin=2em]
\item[-] the trees are disjoint: $P_j \neq P'_j$ for any $P \in \tr{T}_{i_1}, P' \in \tr{T}_{i_2}$ with $i_1 \neq i_2$;
\item[-] whenever $P \in \tr{T}_{i_1}, P' \in \tr{T}_{i_2}$ with $i_1 \neq i_2$ are so that $\ds 2 \omega_{P_j} \cap 2 \omega_{P'_j} \neq \emptyset$, then if $\ds |\omega_{P_j}|<|\omega_{P'_j}|$, one has $I_{P'} \cap I_{\tr{T}_{i_1}} = \emptyset$ and vice-versa;
\item[-] whenever $P \in \tr{T}_{i_1}, P' \in \tr{T}_{i_2}$ with $i_1 < i_2$ are such that $\ds 2 \omega_{P_j} \cap 2 \omega_{P'_j} \neq \emptyset$, then if $\ds |\omega_{P_j}|=|\omega_{P'_j}|$, one has $I_{P'} \cap I_{\tr{T}_{i_1}} = \emptyset$.
\end{enumerate}
\end{definition}

We also define the notion of ``sizes'' for coefficients indexed by the collection of tiles.
\begin{definition}[\textsf{Size}]
\label{def:size:bis}
Given $\P' \subseteq \P$ a collection of $4$-tiles, any $j \in \{ 1, 2, 3, 4  \}$ and coefficients $\mathfrak{a}_j:=(a_{P_j})_{P\in \P}$ we define the associated size as
\begin{equation}
\label{eqdef:size:j}
\size_{\P'}^j(\mathfrak{a}_j) := \sup_{\tr{T} \subset \P'}\  \frac{1}{|I_{\tr{T}}|^\frac{1}{2}} \big( \sum_{P \in \tr{T}} |a_{P_j}|^2   \big)^{1 \over 2}
\end{equation}
where the supremum in $\tr{T}$ varies over $j$-lacunary trees in $\P'$.
\end{definition}
We recall the self-improvement property on the size, due to John-Nirenberg (see Lemma 6.4 in \cite{MTTBiest2}]): 
\begin{equation}
\size_{\P'}^j(\mathfrak{a}_j) \simeq \sup_{\substack{\tr{T} \subset \P' \\ \tr{T}\, j\text{-lac tree}}}\  \frac{1}{|I_{\tr{T}}|} \big\| \big( \sum_{P \in \tr{T}} |a_{P_j}|^2 \frac{\chi_{I_P}}{|I_P|}  \big)^{1 \over 2} \big\|_{L^{1,\infty}(I_{\tr{T}})}.
\label{eq:JN}
\end{equation}

\subsection{Reduction to a symmetric $4$-(sub)linear discrete form}
\label{sec:red:op:anis:local:L2}
We proceed with the study of $T$ by considering more precise Fourier localizations for the functions involved. Our approach consists of firstly, localizing all the information at the dominant frequency scale $2^{3k}$, and subsequently, using spatial localizations at the dominant spatial scale $2^{-k}$. We are thus led to considering time-frequency regions of area $ 2^{2k} \gg 1$, above the Heisenberg principle regime.

\begin{figure}[h!]
\centering
\includegraphics[scale=0.7]{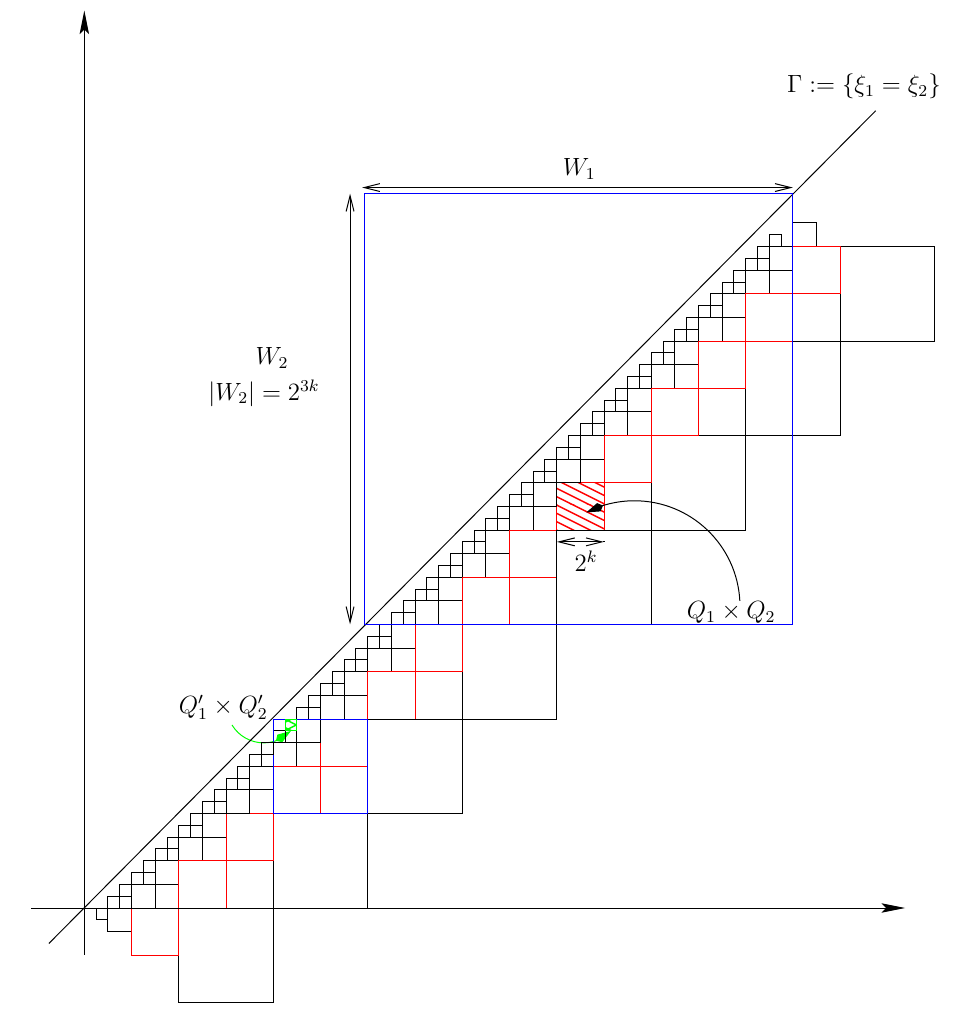}
\caption{Examples of a frequency tiles $Q_1\times Q_2$ in red (resp. $Q_1'\times Q_2'$ in green) and its enlarged tile $W_1 \times W_2$ in blue with $W_1=W_2$. }
\end{figure}

More precisely, for every $k \geq 0$ we consider $\mathcal Q_k$ to be the covering of $\ds \supp \phi \Big( \frac{\xi_1-\xi_2}{2^k}\Big)$ (which can be either the region $\ds \{ (\xi_1, \xi_2) : |\xi_1-\xi_2| \leq 2^{k} \}$ or $\ds \{ (\xi_1, \xi_2) : |\xi_1-\xi_2| \simeq 2^{k} \}$, depending whether $\phi$ is a $\Phi$- or a $\Psi$-type function) by squares $Q_1 \times Q_2$ of sidelength $\simeq 2^k$. And we let $\ic W_k$ to be a covering of the region $\ds \{ (\xi_2, \xi_3) : |\xi_2-\xi_3| \simeq 2^{3k}\}$ by squares $W_2 \times W_3$ of sidelength $\sim 2^{3k}$. The strategy, as illustrated in \cite{multilinear_harmonicII}, is to localize $\ds \phi\big(  \frac{\xi_1-\xi_2}{2^k}  \big)$ onto squares $Q_1\times Q_2 \in \ic Q_k$ and $\ds \psi\big( \frac{\xi_3-\xi_2}{2^{3k}}\big)$ onto squares $W_2 \times W_3 \in \ic W_{3k}$ in order to tensorize the symbol via Fourier series. Using partitions of unity associated to these coverings, we need to deal with the symbol
\begin{align*}
\sum_{k \geq 0} \sum_{Q_1\times Q_2 \in \ic Q_k} \sum_{W_2 \times W_3 \in \ic W_{3k}}c_k \phi\big(  \frac{\xi_1-\xi_2}{2^k}  \big) \, \phi_{Q_1 \times Q_2}(\xi_1, \xi_2) \, \psi\big(  \frac{\xi_3-\xi_2}{2^{3k}}  \big) \, \phi_{W_2 \times W_3}(\xi_2, \xi_3).   
\end{align*}
The nontrivial contribution in the above sum comes under the condition that $\ds Q_2 \cap W_2 \neq \emptyset$, where $|Q_2|=2^k, |W_2|=2^{3k}$ with $k \geq 0$. Under the assumption that the collections $\mathcal Q_k$ and $\mathcal W_k$ are sparse (in the sense of \cite[Chapter 6.1]{multilinear_harmonicII}), if $\ds Q_2 \cap W_2 \neq \emptyset$, we can say that $W_2$ is the unique interval of length $2^{3k}$ belonging to $\ic W_{3k}$ which intersects $Q_2 \in \ic Q_k$. We write this $\ds \hat Q_2=W_2$.

\begin{notation} \label{not:hat}
\label{not:dilation} If $\mathcal G_1$ and $\mathcal G_2$ are two sparse (potentially shifted) dyadic grids and $\omega \in \mathcal G_1$ and $W \in \mathcal G_2$ are two (shifted) dyadic intervals,  we write
\[
\hat \omega=W
\]
to say that $\omega \cap W \neq \emptyset$ and
\begin{equation}
|\omega|= \begin{cases}
 |W|^{1 \over 3}, &\quad \text{if   } |\omega| \geq 1 \\
 |W|^3  &\quad \text{if   } |\omega| < 1.
\end{cases}
\end{equation}
We note that the frequency intervals will be in the first situation and that different shifted dyadic grids might be involved, whilst the spatial interval will be in the second situation, both belonging to the classical dyadic grid.
\end{notation}
As we carry on the time-frequency analysis, we need to study the associated $4$-linear form $\Lambda$ for which it is important to obtain frequency information on the dualizing function. This is why, to any square $\ds W_2\times W_3 \in \ic W_{3k}$ we associate the cube $\bar W := W_1 \times W_2 \times W_3 \times W_4$ of sidelength $\simeq 2^{3k}$, where
\[
W_1=W_2,\qquad  \text{and} \qquad   W_4 \subset -(W_1+W_2+W_3).
\]
The frequency representation of $\Lambda$ involves an integral over the manifold $\ds \{ \xi_1+\ldots+\xi_4=0 \}$; so if we know that $\xi_2 \in Q_2$ with $\hat Q_2=W_2$ and $\xi_3 \in W_3$, then $\xi_1 \in W_1$ because $\hat Q_1=W_2=W_1$ and $\xi_4 \in -(W_1+W_2+W_3)$. So $W_4$ is one of the finitely many intervals of length comparable to $2^{3k}$ contained in $-(W_1+W_2+W_3)$.

Already the fact that $\ds W_2 \times W_3 \in \ic W_{3k}$ implies that $W_2$ is a shift of $W_3$ (as in the study of the bilinear Hilbert transform). From the above we deduce that knowing $W_{i}$ for $i\in \{ 1, 2, 3, 4 \}$ results in knowing $\bar W$ completely, since the other $W_j$ are fixed shifts of $W_{i}$.

We define
\begin{equation}
\label{eq:def:cubes:anis:bht}
\overline{\ic W}_{3k}:= \big\{ \bar W= W_1 \times W_2 \times W_3 \times W_4 :  W_2 \times W_3 \in  \ic W_{3k} \big\} \qquad \text{and} \qquad \overline{\ic W}:= \cup_{k \geq 0} \overline{\ic W}_{3k}.
\end{equation}

After implementing the Fourier series decomposition, the 4-linear form $\Lambda$ associated to the operator in \eqref{eq:tri:initial:1} becomes a superposition of 4-linear forms of type
\begin{align*}
\sum_{k \geq 0} c_k\, \sum_{W_2 \times W_3 \in \ic W_{3k}}  \sum_{\substack{Q_1\times Q_2 \in \ic Q_k \\ \hat Q_2=W_2}} \int_{\xi_1+\ldots+\xi_4=0} \prod_{j=1}^4 \hat f_j(\xi_j) \phi_{Q_1}(\xi_1) \phi_{Q_2}(\xi_2) \phi_{W_3}(\xi_3) \phi_{W_3}(\xi_4) d \xi.
\end{align*}

Next, we perform a wave packet decomposition for the frequency-localized functions: if $Q_1 \times Q_2 \in \ic Q_k$, then we write for example
\begin{equation}
\label{eq:Gabor:dec:f_1}
\hat f_1(\xi_1) \phi_{Q_1}(\xi_1)= \sum_{ \substack{I: |I| = 2^{-k}}} \langle f_1, \phi_{I \times Q_1} \rangle  \hat{\phi}_{I \times Q_1}(\xi_1),
\end{equation}
where the sum ranges over the collection of dyadic intervals of length $2^{-k}$ and $\phi_{I \times Q_1}$ is an $L^2$-normalized wave packet associated to the time-frequency region of area 1 given by $I \times Q_1$, supported in frequency onto $Q_1$. Notice that for $f_1$ and $f_2$ the wave packets are associated to regions of sizes $2^{-k} \times 2^k$, while for $f_3$ and $f_4$ they are associated to time-frequency regions of sizes $2^{-3k} \times 2^{3k}$. In frequency the dominant scale is $2^{3k}$, and at the same time the dominant space scale is $2^{-k}$.

In this way the $4$-linear form $\Lambda(f_1, f_2, f_3, f_4)=\langle T(f_1,f_2,f_3),f_4\rangle$ becomes a superposition of terms, which correspond to \emph{model quadrilinear forms}, of the type
\begin{align}
\label{eq:model:almost:tiles}
\sum_{k \geq 0} \sum_{W_2 \times W_3 \in \ic W_{3k}}  \sum_{\substack{Q_1\times Q_2 \in \ic Q_k \\ \hat Q_2=W_2}}\  \sum_{\substack{I : |I|=2^{-k}}} \ \sum_{\substack{  J \subset I \\ |J|=2^{-3k}}} |I|^{-1} \langle f_1, \phi_{I \times Q_1} \rangle \, \langle f_2, \phi_{I \times Q_2} \rangle \, \langle f_3, \phi_{J \times W_3} \rangle \, \langle f_4, \phi_{J \times W_4} \rangle.
\end{align}
The coefficients $(c_k)_{k \in \Z}$ appearing in \eqref{eq:tri:initial:1} are assimilated into the wave packets, which becomes morally equivalent to assuming them to be $\equiv \, 1$.

Next, we do a Cauchy-Schwartz in the $J$ parameter and in the $Q_2$ parameter, recalling that $Q_1$ is unequivocally determined by $Q_2$. Thus $\ds \Lambda(f_1, f_2, f_3, f_4)$ is bounded above by
\begin{align}
\label{eq:bound:above:Lambda}
\sum_{k \geq 0} \sum_{\bar W \in \overline{\ic W}_{3k}} \sum_{\substack{I : |I|=2^{-k}}} |I|^{-1} \, \Big( \sum_{\hat Q_1 = W_1}  \big|\langle f_1, \phi_{I \times Q_1} \rangle\big|^2 \Big)^{1 \over 2} \, \Big( \sum_{\hat Q_2=W_2}  \big|\langle f_2, \phi_{I \times Q_2} \rangle\big|^2 \Big)^{1 \over 2} &\\  \Big(  \sum_{\substack{ \hat J = I}}  \big|\langle f_3, \phi_{J \times W_3} \rangle\big|^2 \Big)^{1 \over 2} \, \Big(  \sum_{\substack{ \hat J=I}}  \big|\langle f_4, \phi_{J \times W_4} \rangle\big|^2 \Big)^{1 \over 2}. & \nonumber
\end{align}

Now we switch to a more standard time-frequency notation: a generic cube $\overline W \in \overline{\ic W}$, which has sidelength $2^{3k}$, and a spatial interval $I$ of length $2^{-k}$ will produce a $4$-tile $P=(P_1, P_2, P_3, P_4)$, where for any $1 \leq j \leq 4$, 
\[
P_j:=I_P \times \omega_{P_j}=I \times W_j. 
\]
That is, a tile $P$ consists of four rectangles, each of area $2^{2k}$, all having the same spatial interval $I_P=I$ of length $2^{-k}$, with the frequency information yielding a cube of sidelength $2^{3k}$, $\overline{W} \in \overline{\ic W}_k$.

From here on, the collection of all possible tiles $P=(P_1, P_2, P_3, P_4)$ will be denoted $\P$:
\begin{equation}
\label{eq:def:collection:P}
\P:=\big\{ (P_1, P_2, P_3, P_4) : P_j=I \times W_j \text{ for all }1 \leq j \leq 4, |I|^3=|W_j|^{-1}, W_1\times W_2 \times W_3 \times W_4 \in \overline{\ic W}\big\}.
\end{equation}
As usual, via standard arguments, we may assume that $\P$ is finite. In this new terminology \eqref{eq:bound:above:Lambda} becomes
\begin{align}
\sum_{P \in \P} |I_P|^{-1} & \big( \sum_{\omega: \hat \omega=\omega_{P_1}} \big| \langle f_1, \phi_{I_P \times \omega} \rangle \big|^2 \big)^\frac{1}{2} \big( \sum_{\omega: \hat \omega=\omega_{P_2}} \big| \langle f_2, \phi_{I_P \times \omega} \rangle \big|^2 \big)^\frac{1}{2} \nonumber \\
& \cdot \big( \sum_{I: \hat I = I_P}  \big| \langle f_3, \phi_{I \times \omega_{P_3}} \rangle \big|^2 \big)^\frac{1}{2} \big( \sum_{I: \hat I = I_P}  \big| \langle f_4, \phi_{I \times \omega_{P_4}} \rangle \big|^2 \big)^\frac{1}{2}.\label{eq:form:quadri-sublin:Lambda^{(1)}}
\end{align}

Later on we will prefer to exploit {\it space separation} between Heisenberg boxes, which is why we use the following result:
\begin{lemma}
\label{lemma:wide:to:thin}
With the above notation, we have for any $M>0$
\begin{equation}
\label{eq:wide:thing:decay}
\big( \sum_{\omega: \hat \omega=\omega_{P_1}} \big| \langle f_1, \phi_{I_P \times \omega} \rangle \big|^2 \big)^\frac{1}{2} \lesssim_{M} \sum_{\# \in \Z} (1+ |\#|)^{-M } \big( \sum_{I: \hat I = I_{P, \#}}  \big| \langle f_1, \phi_{I \times \omega_{P_1}} \rangle \big|^2 \big)^\frac{1}{2}, 
\end{equation}
where $I_{P, \#}$ denotes the interval having the same length as $I_{P}$, but which is shifted $\#$ units.\footnote{If $\ds I_P=[2^{-k}n, 2^{-k}(n+1)]$, then $\ds I_{P, \#}:=[2^{-k}(n+\#), 2^{-k}(n+\#+1)]$.}
\end{lemma}

\begin{remark}
\label{remark:equiv:l2:tiles}
Notice that a reverse inequality of the same type holds, so that modulo tails, we do not lose any information:
\begin{equation*}
\big( \sum_{I: \hat I = I_{P}}  \big| \langle f_1, \phi_{I \times \omega_{P_1}} \rangle \big|^2 \big)^\frac{1}{2} \lesssim_{M} \sum_{\# \in \Z} (1+ |\#|)^{-M } \big( \sum_{\omega: \hat \omega=\omega_{P_1}^{\pm}} \big| \langle f_1, \phi_{I_{P,\#} \times \omega} \rangle \big|^2  \big)^\frac{1}{2},
\end{equation*}
with $\omega_{P_1}^{\pm}=\omega_{P_1}+\epsilon|\omega_{P_1}|$ for $\epsilon\in\{-1,0,1\}$. This also remains true if we choose a different way of cutting the time-frequency tiles : for any $0 \leq \ell \leq 2k$, 
\begin{equation*}
\big( \sum_{I: \hat I = I_{P}}  \big| \langle f_1, \phi_{I \times \omega_{P_1}} \rangle \big|^2 \big)^\frac{1}{2} \lesssim_{M} \sum_{\# \in \Z} (1+ |\#|)^{-M } \big(
\sum_{\substack{\omega:\omega \subset 3 \omega_{P_1} \\ |\omega|= 2^{- \ell}|\omega_{P_1}| }} \sum_{\substack{I : \,I \subset I_{P, \#} \\ |I|\cdot |\omega|=1}} \big| \langle f_1, \phi_{I \times \omega} \rangle \big|^2  \big)^\frac{1}{2}.
\end{equation*}

\end{remark}

\begin{proof}[Proof of Lemma \ref{lemma:wide:to:thin}]
First, since all the $\omega$ appearing in the left-hand side of the above equation are contained inside $\omega_{P_1}$, we can see that
\[
\langle f_1, \phi_{I_P \times \omega} \rangle= \langle \Pi_{\omega_{P_1}} f_1, \phi_{I_P \times \omega} \rangle,
\]
where $\Pi_{\omega_{P_1}} f_1$ represents a smooth projection of $f_1$ onto the interval $\omega_{P_1}$. Next, we use the fact that
\[
\big( \sum_{\omega: \hat \omega=\omega_{P_1}} \big| \langle f_1, \phi_{I_P \times \omega} \rangle \big|^2 \big)^\frac{1}{2}= \big( \sum_{\omega: \hat \omega=\omega_{P_1}} \big| \langle \Pi_{\omega_{P_1}} f_1, \phi_{I_P \times \omega} \rangle \big|^2 \big)^\frac{1}{2} \lesssim \| \Pi_{\omega_{P_1}} f_1 \cdot \chi_{I_P}  \|_2
\]
and continue with a Fourier series decomposition of $\Pi_{\omega_{P_1}} f_1$, which naturally localizes in space on intervals of length $|\omega_{P_1}|^{-1}=|I_P|^{3}$. More precisely, as in \eqref{eq:Gabor:dec:f_1}, we have
\[
\Pi_{\omega_{P_1}} f_1(x)=\sum_{J : |J|= |\omega_{P_1}|^{-1}} \langle f_1, \phi_{J \times \omega_{P_1}} \rangle \phi_{J \times \omega_{P_1}}(x).
\]
These can be further re-organized as $\ds \Pi_{\omega_{P_1}} f_1(x)=\sum_{\# \in \Z }\sum_{J : \hat J = I_{P, \#}} \langle f_1, \phi_{J \times \omega_{P_1}} \rangle \phi_{J \times \omega_{P_1}}(x)$ and thanks to the triangle inequality, we obtain 
\begin{align*}
\| \Pi_{\omega_{P_1}} f_1 \cdot \chi_{I_P}  \|_2 & = \| \big(\sum_{\# \in \Z }\sum_{J : \hat J = I_{P, \#}} \langle f_1, \phi_{J \times \omega_{P_1}} \rangle \phi_{J \times \omega_{P_1}} \big) \cdot \chi_{I_P}  \|_2 \\
&\lesssim \sum_{\# \in \Z} \| \big( \sum_{J : \hat J = I_{P, \#}} \langle f_1, \phi_{J \times \omega_{P_1}} \rangle \phi_{J \times \omega_{P_1}} \big) \cdot \chi_{I_P}  \|_2. 
\end{align*}
By expanding the square, we notice that 
\begin{align*}
\| \big( \sum_{J : \hat J = I_{P, \#}} \langle f_1, \phi_{J \times \omega_{P_1}} \rangle \phi_{J \times \omega_{P_1}} \cdot \chi_{I_P}  \|_2^2 \lesssim   \sum_{J : \hat J = I_{P, \#}}  \big|\langle f_1, \phi_{J \times \omega_{P_1}} \rangle \big|^2  (1+|\#|)^{-M} \sum_{J' : \hat J' = I_{P, \#}} \Big( 1+\frac{\dist(J, J')}{|J|}  \Big)^{-M_1},
\end{align*}
thus concluding that \eqref{eq:wide:thing:decay} holds.
\end{proof}

Applying this lemma we can bound the expression in \eqref{eq:form:quadri-sublin:Lambda^{(1)}} by a superposition of terms of the form
\begin{align}\label{grup}
\Lambda_{\P}(f_1, f_2, f_3, f_4):&=  \sum_{P \in \P} |I_P|^{-1} \big( \sum_{I: \hat I = I_{P}}  \big| \langle f_1, \phi_{I \times \omega_{P_1}} \rangle \big|^2 \big)^\frac{1}{2} \big( \sum_{I: \hat I = I_{P}}  \big| \langle f_2, \phi_{I \times \omega_{P_2}} \rangle \big|^2 \big)^\frac{1}{2} \\
& \qquad \qquad \quad \cdot \big( \sum_{I: \hat I = I_P}  \big| \langle f_3, \phi_{I \times \omega_{P_3}} \rangle \big|^2 \big)^\frac{1}{2} \big( \sum_{I: \hat I = I_P}  \big| \langle f_4, \phi_{I \times \omega_{P_4}} \rangle \big|^2 \big)^\frac{1}{2} \\
&:= \sum_{P \in \P} |I_P|^{-1} f_1(P_1) \cdot f_2(P_2) \cdot f_3(P_3) \cdot f_4(P_4).
\end{align}
When we compare the above expression to the quadrilinear form associated to the operator in \eqref{eq:tri:initial:1}, we notice that a lot of information has been lost : on one hand, all the functions are localized in $L^2$. And on the other hand, while in \eqref{eq:tri:initial:1} the joint frequency information for $f_1$ and $f_2$ was supported on $\{ |\xi_1-\xi_2| \leq 2^k \}$  or $\{ |\xi_1-\xi_2| \sim 2^k \}$ (depending on the nature of the function $\phi$), now the joint frequency information is supported on $\{ |\xi_1-\xi_2| \leq 2^{3k} \}$.

\subsection{Boundedness of the quadri(sub)linear form $\Lambda_{\P}(f_1, f_2, f_3, f_4)$}~\\
\label{subsec:bddness:4-lin:main}
Here we perform a detailed time-frequency analysis of the $4$-(sub)linear form $\Lambda_{\P}(f_1, f_2, f_3, f_4)$. We note that our collection of tiles $\P$ is not a classical \emph{rank 1} collection as in \cite{Mlinop} (due to the anisotropy), but it enjoys somehow similar properties: for any $P \in \P$, if we know one of the tiles $\{P_1= P_2,P_3,P_4\}$ then we can uniquely determine the other components. So although $\Lambda_{\P}$ bears a resemblance to other operators that have already appeared in the literature, it does not directly fall under the scope of existing models, in particular since our  tiles have area $\gg 1$.

For $j\in\{1,2,3,4\}$ and $\P$ a collection of $4$-tiles, we set
$$ \size_{\P'}^j(f_j):=\size_{\P'}^j(\mathfrak{f_j}) \qquad \textrm{with} \qquad \mathfrak{f_j}:=(f_j(P_j))_{P\in\P}.$$

\begin{lemma}[\textsf{Size-energy discretization}] 
\label{lemma:decomposition:lemma:4:lin:form:no:curv} Let $\P' \subset \P$ be a collection of $4$-tiles, let $j \in \{1, 2, 3, 4\}$ and $\lambda >0$. If $\size_{\P'}^j(f_j) \leq \lambda$, then there exists a decomposition $\ds \P':= \P'_{(1)} \sqcup \P'_{(2)}$ so that 
\[
\size_{\P'_{(1)}}^j(f_j) \leq { \lambda \over 2}
\]
and $\P'_{(2)}$ can be written as a union of $j$-trees (lacunary or overlapping), $\ds \P'_{(2)} = \bigcup_{\tr{T} \in \ii F_{\lambda}} \tr{T}$ satisfying 
\begin{equation}
\label{eq:dec:imp:space}
\sum_{\tr{T} \in \ii F_{\lambda}} |I_{\tr{T}}| \lesssim \lambda^{-2} \| f_j \|_2^2.
\end{equation}  
\end{lemma}

\begin{proof} The proof is in most parts standard and we follow for example the presentation in \cite[Ch. 6]{multilinear_harmonicII}. We initialize the collection $\ii F_{\lambda}:=\emptyset$; we proceed by selecting $\tr{T}_1$, a maximal $j$-lacunary (upwards) tree having the property that 
\begin{equation}
\label{eq:stop:time:cond}
\frac{1}{|I_{\tr{T}_1}|} \big( \sum_{P \in \tr{T}_1} f_j(P_j)^2   \big)^{1 \over 2} \geq \frac{\lambda^2}{4}
\end{equation}
and so that $\xi_{\tr{T}_1}^{(j)}$ is maximal and add it to $\ii F_{\lambda}$.

Then we also add to $\ii F_{\lambda}$ the tree $\tr{T}_1'$, the $j$-overlapping tree having the same top datum. In order to ensure enough separation between the trees, we also remove the $j$-lacunary and $j$-overlapping trees having datum $(I_{\tr{T}_1, \#}, \xi_{\tr{T}_1})$, where $\# \in \{-1 , 1\}$ -- that is we remove the trees having the same frequency information as $\tr{T}$, with the space information contained in the left-translate or right-translate of $I_{\tr{T}}$. We denote these $\tr{T}_{1, \#}$ and $\tr{T}'_{1, \#}$, where $\# \in \{-1 , 1\}$. We add them to $\ii F_{\lambda}$, and remove them from the collection $\P'$. Afterwards we begin again this operation.

Since the collection of tiles is finite, this produces a finite chain of strongly $j$-disjoint trees $\tr{T}_1, \ldots, \tr{T}_N$, each coming with five associated trees. Thus we are left with proving 
\[
\sum_{\ell=1}^N |I_{\tr{T}_{\ell}}| \lesssim \lambda^{-2} \| f_j \|_2^2,
\]
which follows via Bessel-type arguments from the estimate
\begin{equation}
\label{eq:bessel:tree:top}
\Big\| \sum_{\ell=1}^N \sum_{P \in \tr{T}_{\ell}} \sum_{I : \hat I = I_P}  \langle f_j, \phi_{I \times \omega_{P_j}} \rangle  \phi_{I \times \omega_{P_j}}     \Big\|_2 \lesssim \lambda \cdot \Big( \sum_{\ell=1}^N |I_{\tr{T}_{\ell}}| \Big)^{1 \over 2}.
\end{equation}
In the classical case the tiles have area 1 and have wave packets associated to them, so one can use the almost orthogonality between wave packets in order to deduce the desired estimate. Here our regions have area $\gg 1$ and we need to assess the almost orthogonality between \textit{sums} of wave packets, which makes the proof more technical. More precisely, in the standard case the expression $\ds \sum_{I : \hat I = I_P}  \langle f_j, \phi_{I \times \omega_{P_j}} \rangle  \phi_{I \times \omega_{P_j}}$ which appears above is replaced by $\langle f_j, \phi_{I_P \times \omega_{P_j}} \rangle \phi_{I_P \times \omega_{P_j}}$.

Squaring the left-hand side of \eqref{eq:bessel:tree:top}, we get
\begin{equation}
\label{eq:bessel:3:cases}
\sum_{\ell, \ell'=1}^N \sum_{P \in \tr{T}_{\ell}} \sum_{P' \in \tr{T}_{\ell'}}  \sum_{I: \hat I=I_P} \sum_{I': \hat {I'}=I_{P'}}  \langle f_j, \phi_{I \times \omega_{P_j}} \rangle \langle f_j, \phi_{I' \times \omega_{P'_j}} \rangle \langle \phi_{I \times \omega_{P_j}}, \phi_{I' \times \omega_{P'_j}} \rangle.
\end{equation}
Since the wave packets are compactly supported in frequency, the term $\ds \langle \phi_{I \times \omega_{P_j}}, \phi_{I' \times \omega_{P'_j}} \rangle$ is non-zero if and only if $\ds \omega_{P_j} \cap \omega_{P'_j} \neq \emptyset$, i.e. if $\ds \omega_{P_j} = \omega_{P'_j}$ or if $\ds \omega_{P_j} \subsetneq \omega_{P'_j}$ (and symmetrically if $\ds \omega_{P'_j} \subseteq \omega_{P_j}$).

By symmetry, it suffices to bound the terms
\begin{align*}
& \sum_{\ell=1}^N \sum_{P \in \tr{T}_{\ell}}  \sum_{\substack{\ell'=1}}^N     \sum_{\substack{P' \in \tr{T}_{\ell'} \\  \omega_{P_j}= \omega_{P'_j}}}  \sum_{I: \hat I=I_P} \sum_{I': \hat {I'}=I_{P'}}   \langle f_j, \phi_{I \times \omega_{P_j}} \rangle \langle f_j, \phi_{I' \times \omega_{P'_j}} \rangle \langle \phi_{I \times \omega_{P_j}}, \phi_{I' \times \omega_{P'_j}} \rangle \\
+ & \sum_{\ell=1}^N \sum_{P \in \tr{T}_{\ell}}  \sum_{\substack{\ell'=1}}^N     \sum_{\substack{P' \in \tr{T}_{\ell'} \\  \omega_{P_j} \subsetneq \omega_{P'_j}}}   \sum_{I: \hat I=I_P} \sum_{I': \hat {I'}=I_{P'}} \langle f_j, \phi_{I \times \omega_{P_j}} \rangle \langle f_j, \phi_{I' \times \omega_{P'_j}} \rangle \langle \phi_{I \times \omega_{P_j}}, \phi_{I' \times \omega_{P'_j}} \rangle :=I+II.
\end{align*}

For the first term, we use the trivial inequality $|a\cdot b| \leq a^2+b^2$ and symmetry to get
\begin{align*}
|(I)| &\lesssim \sum_{\ell=1}^N \sum_{P \in \tr{T}_{\ell}}  \sum_{I: \hat I=I_P} \big| \langle f_j, \phi_{I \times \omega_{P_j}} \rangle \big|^2  \Big(   \sum_{\substack{\ell'=1}}^N \sum_{\substack{P' \in \tr{T}_{\ell'} \\  \omega_{P_j}= \omega_{P'_j}}} \sum_{I': \hat {I'}=I_{P'}}  \big| \langle \phi_{I \times \omega_{P_j}}, \phi_{I' \times \omega_{P'_j}} \rangle \big|  \Big) .
\end{align*}
For a fixed tile $P$, all the other tiles $P'$ so that $\omega_{P_j}= \omega_{P'_j}$ (be it in the same tree as $P$ or in another tree of the chain of strongly $j$-disjoint trees) have mutually disjoint spatial intervals, so that, using the size constraint on the collection of tiles, we obtain $|(I)|\lesssim \sum_{\ell=1}^N |I_{\tr{T}_{\ell}}| \cdot  \lambda^2$.

For $(II)$, we apply Cauchy-Schwarz :
\begin{align*}
|(II)|&\lesssim \sum_{\ell=1}^N \sum_{P \in \tr{T}_{\ell}} \sum_{\substack{\ell'=1 \\ \ell' \neq \ell}}^N  \sum_{\substack{P' \in \tr{T}_{\ell'} \\  \omega_{P_j} \subsetneq \omega_{P'_j}}} \big(  \sum_{I: \hat I=I_P}   \big|  \langle f_j, \phi_{I \times \omega_{P_j}} \rangle \big|^2    \big)^{1 \over 2} \big(  \sum_{I': \hat{I'}=I_{P'}}   \big|  \langle f_j, \phi_{I' \times \omega_{P'_j}} \rangle \big|^2    \big)^{1 \over 2} \\
&\qquad \qquad \qquad \cdot   \big( \sum_{I: \hat I=I_P} \sum_{I': \hat {I'}=I_{P'}} \big| \langle \phi_{I \times \omega_{P_j}}, \phi_{I' \times \omega_{P'_j}} \rangle    \big|^2   \big)^{1 \over 2} \\
&\lesssim  \sum_{\ell=1}^N \sum_{P \in \tr{T}_{\ell}} \sum_{\substack{\ell'=1 \\ \ell' \neq \ell}}^N  \sum_{\substack{P' \in \tr{T}_{\ell'} \\  \omega_{P_j} \subsetneq \omega_{P'_j}}} \lambda |I_P|^{1 \over 2} \cdot \lambda |I_{P'}|^{1 \over 2}  \big( \sum_{I: \hat I=I_P} \sum_{I': \hat {I'}=I_{P'}} \big| \langle \phi_{I \times \omega_{P_j}}, \phi_{I' \times \omega_{P'_j}} \rangle    \big|^2   \big)^{1 \over 2}. 
\end{align*}

Now we fix $\ell$ and $P \in \tr{T}_{\ell}$ and we would like to estimate
\begin{align*}
\sum_{\substack{\ell'=1 \\ \ell' \neq \ell}}^N  \sum_{\substack{P' \in \tr{T}_{\ell'} \\  \omega_{P_j} \subsetneq \omega_{P'_j}}} |I_{P'}|^{1 \over 2} \Big( \sum_{I: \hat I=I_P} \sum_{I': \hat {I'}=I_{P'}} \big| \langle \phi_{I \times \omega_{P_j}}, \phi_{I' \times \omega_{P'_j}} \rangle    \big|^2   \Big)^{1 \over 2}.
\end{align*}
We note that the trees appearing in the sum above all have similar frequency information, but due to the selection algorithm the $I_{P'}$ are all mutually disjoint.

For a fixed $\ell' \neq \ell$, and for any $P' \in \tr{T}_{\ell'}$ so that $\omega_{P_j} \subsetneq \omega_{P'_j}$, we use the fact that $|I_P|\gg |I_{P'}|$ (which implies for all $I, I'$ appearing in the summation that $|I|\gg |I'|$) and the standard estimate 
\[
\big| \langle \phi_{I \times \omega_{P_j}}, \phi_{I' \times \omega_{P'_j}} \rangle    \big| \lesssim \frac{|I'|^{1 \over 2}}{|I|^{1 \over 2}} \Big( 1+ \frac{\dist(I, I')}{|I|}   \Big)^{-M}.
\]
Since we made sure that the trees are well separated in space , we have in fact that
\begin{align*}
 \sum_{I: \hat I=I_P} \sum_{I': \hat {I'}=I_{P'}} \big| \langle \phi_{I \times \omega_{P_j}}, \phi_{I' \times \omega_{P'_j}} \rangle    \big|^2  \lesssim  \sum_{I: \hat I=I_P} \sum_{I': \hat {I'}=I_{P'}} \frac{|I'|}{|I|} \Big(\frac{\dist(I, I')}{|I|}   \Big)^{-2M}.
\end{align*}
We would like to show that 
\begin{equation}
\label{eq:est:trees:energy:wanted}
 \sum_{I: \hat I=I_P} \sum_{I': \hat {I'}=I_{P'}} \big| \langle \phi_{I \times \omega_{P_j}}, \phi_{I' \times \omega_{P'_j}} \rangle    \big|^2   \lesssim \frac{|I_{P'}|}{|I_P|} \Big(\frac{\dist(I_P, I_{P'})}{|I_P|}   \Big)^{-2 \tilde M}.
\end{equation}
Given the definitions of $I$ and $I'$, we have 
\begin{align}
\label{eq:sum:en:1}
 \sum_{I: \hat I=I_P} \sum_{I': \hat {I'}=I_{P'}} \big| \langle \phi_{I \times \omega_{P_j}}, \phi_{I' \times \omega_{P'_j}} \rangle    \big|^2  \lesssim   \frac{|I_{P'}|^3}{|I_P|^3}  \sum_{I: \hat I=I_P} \sum_{I': \hat {I'}=I_{P'}} \Big(\frac{\dist(I, I')}{|I|}   \Big)^{-2M}.
\end{align}
We fix $I'$ and sum over $I$: 
\[
\sum_{I: \hat I=I_P} \Big(\frac{\dist(I, I')}{|I|}   \Big)^{-2M}= \sum_{I: \hat I=I_P} \Big(\frac{\dist(I, I')}{|I_P|^3}   \Big)^{-2M} \lesssim \Big(\frac{\dist(I_P, I')}{|I_P|^3}   \Big)^{-M}= \Big(\frac{\dist(I_P, I')}{|I_P|}   \Big)^{-M} |I_P|^{2M}.
\]

Now we use the fact that $I_P$ and $I_{P'}$ are well separated and
\[
\Big(\frac{\dist(I_P, I')}{|I_P|}   \Big)^{-M} \lesssim  \frac{1}{|I'|} \int_{I'} \Big( \frac{\dist(I_P, x)}{|I_P|}   \Big)^{-M} dx= \frac{1}{|I_{P'}|^3} \int_{I'} \Big( \frac{\dist(I_P, x)}{|I_P|} \Big)^{-M} dx.
\]
Summing in $I'$ we get 
\begin{align*}
\sum_{I: \hat I=I_P} \sum_{I': \hat {I'}=I_{P'}} \Big(\frac{\dist(I, I')}{|I|}   \Big)^{-2M} & \lesssim \frac{|I_P|^{2M}}{|I_{P'}|^3} \int_{I_{P'}} \Big( \frac{\dist(I_P, x)}{|I_P|} \Big)^{-M} dx \lesssim \frac{|I_P|^{2M}}{|I_{P'}|^2} \Big( \frac{\dist(I_P, I_{P'})}{|I_P|} \Big)^{-M}.
\end{align*}

Returning to the estimate \eqref{eq:sum:en:1}, we obtain 
\begin{align}
\label{eq:small:to:big}
 \sum_{I: \hat I=I_P} \sum_{I': \hat {I'}=I_{P'}} \big| \langle \phi_{I \times \omega_{P_j}}, \phi_{I' \times \omega_{P'_j}} \rangle    \big|^2   & \lesssim   \frac{|I_{P'}|^3}{|I_P|^3}  \frac{|I_P|^{2M}}{|I_{P'}|^2} \Big( \frac{\dist(I_P, I_{P'})}{|I_P|} \Big)^{-M} \\
 &\lesssim   \frac{|I_{P'}|}{|I_P|}  \Big( \frac{\dist(I_P, I_{P'})}{|I_P|} \Big)^{-M} \cdot |I_P|^{2M-2}, \nonumber
\end{align}
which implies the desired \eqref{eq:est:trees:energy:wanted} since $|I_P|=2^{-k} \leq 1$.

Once we have \eqref{eq:est:trees:energy:wanted}, the argument follows the same ideas as in the classical case, see for example \cite[Section 6.5]{multilinear_harmonicII}:
\begin{align*}
|(III)| & \lesssim \lambda^2  \sum_{\ell=1}^N \sum_{P \in \tr{T}_{\ell}} \sum_{\substack{\ell'=1 \\ \ell' \neq \ell}}^N  \sum_{\substack{P' \in \tr{T}_{\ell'} \\  \omega_{P_j} \subsetneq \omega_{P'_j}}}  |I_P|^{1 \over 2} \cdot  |I_{P'}|^{1 \over 2}  \frac{|I_{P'}|^{1 \over 2}}{|I_P|^{1 \over 2}} \Big(\frac{\dist(I_P, I_{P'})}{|I_P|}   \Big)^{- \tilde M}\\
&\lesssim \lambda^2 \sum_{\ell=1}^N \sum_{P \in \tr{T}_{\ell}}  \int_{(3 I_{\mb T})^c} \Big( \frac{\dist(I_P, x)}{|I_P|}   \Big)^{- \tilde M} dx \lesssim \lambda^2  \big( \sum_{\ell=1}^N  |I_{\tr{T}_{\ell}}|  \big).
\end{align*}

\end{proof}

Once we have the decomposition in Lemma \ref{lemma:decomposition:lemma:4:lin:form:no:curv}, we apply it iteratively to obtain the following:
\begin{proposition}
\label{prop:decompo:full}
Let $\P$ a finite collection of $4$-tiles as in \eqref{eq:def:collection:P} and $n_0 \in \Z$ so that $\ds \max_{1 \leq j \leq 4 } \Big(  \frac{\size_{\P}^j(f_j)}{\|  f_j \|_2}   \Big) \simeq  2^{-n_0} $. Then one can decompose $\P$ as 
\begin{equation}
\label{eq:coll:dec}
\P= \bigcup_{n\geq n_0} \bigcup_{j \in \{ 1, 2, 3, 4  \}} \bigcup_{\tr{T} \in \ii T_n^{(j)}} \tr{T}:= \bigcup_{n\geq n_0}  \P_n,
\end{equation} 
where for each $n \geq n_0$ we have, for $\P_n \neq \emptyset$, $\ds \size_{\P_n}^j(f_j) \leq \min \big( 2^{-n} \|f_j\|_2, \size_{\P}^j(f_j) \big)$ and each  non-empty $\ii T_n^{(j)}$ is a collection of trees satisfying 
\[
\sum_{\tr{T} \in \ii T_n^{(j)}}|I_{\tr{T}}| \lesssim 2^{2n}.
\] 
\end{proposition}

As in the case of the bilinear Hilbert transform (the cubes in $\overline{\ic W}$ from \eqref{eq:def:cubes:anis:bht} have very similar properties), each tree contains at least two lacunary directions, which means that for any tree $\tr{T} \subset \P$
\begin{equation}
\label{eq:tree:est}
\big| \Lambda_{\tr{T}}(f_1, f_2, f_3, f_4) \big| \lesssim \prod_{j=1}^4 \size_{\tr{T}}^j (f_j)   \cdot |I_{\tr{T}}|.
\end{equation}

Now that we have decomposed the collection $\P$ into unions of trees according to the levels of the sizes, we follow the presentation in \cite[Lemma 6.12]{multilinear_harmonicII} to get 
\begin{proposition}\label{prop:final:est:0} 
If $\P$ is a finite collection of $4$-tiles as in \eqref{eq:def:collection:P}, then 
\begin{align*}
\big| \Lambda_{\P}(f_1, f_2, f_3, f_4)  \big| \lesssim & \big( \size_{\P}^1(f_1)\big)^{\theta_1}  \big( \size_{\P}^2(f_2)\big)^{\theta_2}  \big( \size_{\P}^3(f_3)\big)^{\theta_3} \big( \size_{\P}^4(f_4)\big)^{\theta_4} \\
&\cdot \|f_1\|_2^{1- \theta_1} \|f_2\|_2^{1- \theta_2} \|f_3\|_2^{1- \theta_3} \|f_4\|_2^{1- \theta_4},
\end{align*}
for any $0 < \theta_1, \theta_2, \theta_3, \theta_4 <1$ with $\theta_1+\theta_2+\theta_3+\theta_4=1$.

Moreover, if $I_0$ is a fixed spatial interval and $\P(I_0):=\{P\in\P, I_P\subset I_0\}$ denotes the collection of $4$-tiles in $\P$ with spatial support in $I_0$, then 
\begin{align}
\label{eq:local:est}
\big| \Lambda_{\P(I_0)}(f_1, f_2, f_3, f_4)  \big|  \lesssim & \big( \size_{\P(I_0)}^1(f_1)\big)^{\theta_1}  \big( \size_{\P(I_0)}^2(f_2)\big)^{\theta_2}  \big( \size_{\P(I_0)}^3(f_3)\big)^{\theta_3} \big( \size_{\P(I_0)}^4(f_4)\big)^{\theta_4} \\
&\cdot \|f_1 \cdot \chi_{I_0}\|_2^{1- \theta_1} \|f_2 \cdot \chi_{I_0}\|_2^{1- \theta_2} \|f_3 \cdot \chi_{I_0}\|_2^{1- \theta_3} \|f_4 \cdot \chi_{I_0}\|_2^{1- \theta_4}, \nonumber
\end{align}
for any $0 < \theta_1, \theta_2, \theta_3, \theta_4 <1$ with $\theta_1+\theta_2+\theta_3+\theta_4=1$.
\end{proposition}

If classically the sizes correspond to localized Calder\'on-Zygmund operators, this is not the case anymore in our current context since we employ collections of super-Heiseberg tiles whose areas are greater than one. Due to that, the sizes only have an $L^2$ character:
\begin{proposition} \label{prop:size1} For $\P$ a finite collection of $4$-tiles as in \eqref{eq:def:collection:P} and $j\in\{1,2,3,4\}$, one has
$$ \size_{\P}^j(f_j) \lesssim \sup_{P\in \P} \left(\frac{1}{|I_P|}\int |f_j(x)|^2 \chi_{I_P}(x) dx \right)^{1/2}.$$
\end{proposition}

\begin{remark} 
One could invoke John-Nirenberg in order to obtain $L^1$ sizes as in \cite[Lemma 6.11]{multilinear_harmonicII}; however, that would compel us to work with the smaller tiles (of Heisenberg type), which seems incompatible with the anisotropic character of our operator, for technical reasons.
\end{remark}

Combining the above results, we can finally conclude:

\begin{proof}[Proof of Theorem \ref{thm:proof:main:1}]
By implementing the local estimate \eqref{eq:local:est} of Proposition \ref{prop:final:est:0} above in the helicoidal method machinery of \cite[Section 3.3]{BM-sparse} -- see Theorem 12 therein -- we obtain
\begin{align*}
\big| \Lambda_{\P}(f_1, f_2, f_3, f_4)  \big| \lesssim \int_{\R} M_{2}(f_1)(x) \, M_{2}(f_2)(x) \, M_{2}(f_3)(x)\, M_{2}(f_4)(x) dx,
\end{align*}
where $M$ is the Hardy-Littlewood maximal operator and $M_2(f):=M(|f|^2)^{1/2}$. This implies the boundedness in the local-$L^2$ range of the operator $T_m$, thus concluding the proof of Theorem \ref{thm:proof:main:1}. 
\end{proof}
\section{Proof of Theorem \ref{thm:proof:main:2}: The anisotropic Biest case involving two collections of Heisenberg tiles}
\label{sec:proof:m:2}
In this section we are providing the proof of Theorem \ref{thm:proof:main:2} which concerns the behavior of the trilinear operators $T_m$ associated to the symbol
\begin{align}
\label{eq:tril:aniso:op:bi:est}
m(\xi_1,\xi_2,\xi_3)=\sum_{k \geq 0} d_k \, \Big( \sum_{0 \leq j \leq k} c_j \, \psi\big( \frac{\xi_1-\xi_2}{2^{j}}  \big) \Big)  \,\psi\big( \frac{\xi_3-\xi_2}{2^{3 k}}  \big).
\end{align}
 
\subsection{Reductions}
As in \cite{MTTBiest2}, \cite{flag_paraproducts_Revista} and in \cite{Joeun-IteratedTrilFourierIntg}, we first want to express $T_m$ as a local composition of bilinear operators. This is achieved by replacing the $\xi_3-\xi_2$ variable in \eqref{eq:tril:aniso:op:bi:est} by $\xi_3-\frac{\xi_1+\xi_2}{2}$, modulo error terms. Indeed, we notice that 
\begin{align}
\psi\big( \frac{\xi_3-\xi_2}{2^{3k}}\big)- \psi\big( \frac{\xi_3-\frac{\xi_1+\xi_2}{2}}{2^{3k}}\big)= \frac{\xi_1-\xi_2}{2 \cdot 2^{3k}} \int_0^1 \psi'\big(\frac{\xi_3-\frac{\xi_1+\xi_2}{2}+\tau \frac{\xi_1-\xi_2}{2}}{2^{3k}}\big) d \tau,
\end{align}
so that $m(\xi_1, \xi_2, \xi_3 )$ in \eqref{eq:tril:aniso:op:bi:est} becomes
\begin{align}
\label{eq:tril:form:Taylor}
& \sum_{k \geq 0} d_k \, \Big( \sum_{0 \leq j \leq k} c_j \, \psi\big( \frac{\xi_1-\xi_2}{2^{j}}  \big) \Big) \psi\big( \frac{\xi_3-\frac{\xi_1+\xi_2}{2}}{2^{3k}} \big) \\
&+\sum_{k \geq 0} 2^{-2k} \, d_k \, \Big( \sum_{0 \leq j \leq k}  c_j \, 2^{j-k} \psi\big( \frac{\xi_1-\xi_2}{2^{j}}\big)  \cdot  \frac{\xi_1-\xi_2}{2^{j}}\Big) \int_0^1 \psi'\Big(\frac{\xi_3-\frac{1-\tau}{2} \xi_1 -\frac{1+\tau}{2} \xi_2}{2^{3k}}\Big) d \tau \nonumber \\
&:= \mathfrak{m}(\xi_1, \xi_2, \xi_3)+\mathfrak{r}(\xi_1, \xi_2, \xi_3). \nonumber
\end{align}

The above decomposition naturally results in $T_m=T_{\mathfrak{m}}+T_{\mathfrak{r}}$, so the boundedness of $T_m$ from Theorem \ref{thm:proof:main:2} is a consequence of the following two statements.
\begin{proposition}
\label{prop:control:remainder}
The operator $T_{\mathfrak{r}}$ is bounded from $L^{p_1}\times L^{p_2} \times L^{p_3}$ into $L^p$ whenever $1\leq p_1,p_2\leq \infty$, $1<p_3\leq \infty$ with $\ds \frac{1}{p_1}+\frac{1}{p_2}+\frac{1}{p_3}=\frac{1}{p}$ and moreover
\[
\big\| T_{\mathfrak{r}} \big\|_{L^{p_1}\times L^{p_2} \times L^{p_3} \to L^{p}} \lesssim \|d_k\|_{\ell^\infty} \, \|c_j\|_{\ell^\infty}.
\]
\end{proposition}

\begin{theorem} \label{thm:model:biest:ani} The operator $T_{\mathfrak{m}}$ is bounded from $L^{p_1}\times L^{p_2} \times L^{p_3}$ into $L^p$ whenever  $\ds \frac{1}{p_1}+\frac{1}{p_2}+\frac{1}{p_3}=\frac{1}{p}$ with $(\frac{1}{p_1},\frac{1}{p_2},\frac{1}{p_3})\in{\bm D}''$ and moreover
\[
\big\| T_{\mathfrak{m}} \big\|_{L^{p_1}\times L^{p_2} \times L^{p_3} \to L^{p}} \lesssim \|d_k\|_{\ell^\infty} \, \|c_j\|_{\ell^\infty}.
\]
\end{theorem}

We first prove Proposition \ref{prop:control:remainder}, which, thanks to the extra $2^{-2k}$ decay appearing in \eqref{eq:tril:form:Taylor}, reduces to the analysis of a one-scale operator similar to the one discussed in Proposition \ref{prop:Tk}. The proof of Theorem \ref{thm:model:biest:ani} is more involved and will be carried out afterwards, in the remainder of the section. 

\begin{proof}[Proof of Proposition \ref{prop:control:remainder}]
In the remainder term $\mathfrak{r}(\xi_1, \xi_2, \xi_3)$ the extra $2^{-2k}$ decay allows us to fix both $k \geq 0$ and $j\in[0,k]$ in the summation. More precisely, we first focus on the fixed-scale trilinear operator, denoted $T_{\mathfrak{r}_{k,j}}$, with symbol\footnote{Here $\tilde \psi(\eta):=\eta \cdot \psi(\eta)$ and $\Phi_j(\xi)= \Phi(\xi/2^j)$.}
\[
\mathfrak{r}_{k,j}(\xi_1, \xi_2, \xi_3) := \tilde\psi_j\big(\xi_1-\xi_2\big) \cdot \int_0^1 \psi'_{3k}\Big(\xi_3-\frac{1-\tau}{2} \xi_1 -\frac{1+\tau}{2} \xi_2\Big) d \tau.
\]
In space, $T_{\mathfrak{r}_{k,j}}$ can be represented as 
\begin{align*}
T_{\mathfrak{r}_{k,j}}(f_1, f_2, f_3)(x) = \int_0^1 \iint f_1(x-u+t\frac{1-\tau}{2}) f_2(x+u+t \frac{1+\tau}{2}) f_3(x-t) \check{\psi'}_{3k}(t)\check{\tilde \psi}_{j}(u) \, dtdu d\tau.
\end{align*}
Since the functions $\check{\psi'}_{3k}$ and $\check{\tilde \psi}_{j}$ are $L^1$-normalized, the boundedness of $T_{\mathfrak{r}_{k,j}}$ in the Banach range (when $p\geq 1$) is a consequence of Minkowski's inequality.

For the quasi-Banach situation (when $p<1$), we notice that $u$ is at scale $2^{-j}$ and $t$ at scale $2^{-3k}\leq 2^{-j}$, so we do the same reasoning as for the one-scale operator in Proposition \ref{prop:Tk}: we localize everything to intervals of length $2^{-j}$. For every interval $J$ of length $2^{-j}$ and $p<1$
\begin{align*}
 \|T_{\mathfrak{r}_{k,j}}(f_1, f_2, f_3)\|_{L^p(J)} & \leq |J|^{\frac{1}{p}-1} \|T_{\mathfrak{r}_{k,j}}(f_1, f_2, f_3)\|_{L^1(J)} \\
& \lesssim |J|^{\frac{1}{p}-1} \int_0^1 \iiint |f_1(x-u+t\frac{1-\tau}{2})| |f_2(x+u+t \frac{1+\tau}{2})| \\
& \qquad \qquad \qquad \qquad \cdot|f_3(x-t)| |\check{\psi'}_{3k}(t)| |\check{\tilde \psi}_{j}(u)| {\bf 1}_J(x) dx dt du d\tau.
\end{align*}
Using the localization from $\check{\psi'}_{3k}(t)$ and $\check{\tilde \psi}_{j}(u)$, one may control the expressions appearing above using
$$ |f_1(x-u+t\frac{1-\tau}{2})| |\check{\psi'}_{3k}(t)|^{1 \over 3}|\check{\tilde \psi}_{j}(u)|^{1 \over 2}{\bf 1}_{J}(x) \lesssim |J|^{-{1\over 3}} 2^k |f_1(x-u+t\frac{1-\tau}{2})| \chi_{J}(x-u+t\frac{1-\tau}{2})$$ 
and similarly for $f_2$ and $f_3$. Then the change of variables 
$$(x,u,t) \rightarrow (y_1,y_2,y_3):=(x-u+t\frac{1-\tau}{2},x+u+t\frac{1+\tau}{2},x-t)$$
has Jacobian determinant independent of $\tau\in(0,1)$ and so we deduce that $\|T_{\mathfrak{r}_{k,j}}(f_1, f_2, f_3)(x)\|_{L^p(J)}$ is bounded above by
\begin{align*}
 |J|^{\frac{1}{p}-2}2^{3k} \int_0^1\Big( \iiint |f_1(y_1)\chi_{J}(y_1)| |f_2(y_2)\chi_{J}(y_2)| |f_3(y_3)\chi_{J}(y_3)|\chi_{[0,2^{3k}]}(y_1+y_2-2y_3) dy_1dy_2dy_3 \Big) d\tau.
\end{align*}
Then use first H\"older inequality in $y_3$ and then in $y_1,y_2$, that give
\begin{align*}
\|T_{\mathfrak{r}_{k,j}}(f_1, f_2, f_3)(x)\|_{L^p(J)} & \lesssim |J|^{\frac{1}{p}-2}2^{3k} \|f_1\chi_{J}\|_{p_1} |J|^{1 \over p_1'} \|f_2\chi_{J}\|_{p_2} |J|^{1 \over p_2'} \|f_3\chi_{J}\|_{p_3} 2^{-{3k \over p_3'}} \\
& \lesssim |J|^{1 \over p_3} 2^{3k \over p_3} \|f_1\chi_{J}\|_{p_1} \|f_2\chi_{J}\|_{p_2} \|f_3\chi_{J}\|_{p_3}.
\end{align*}
So by the H\"older scaling on the exponents, it can be summed over the intervals $J$ of length $2^{-j}$ and we get (uniformly in $j,k$)
\begin{align*}
\|T_{\mathfrak{r}_{k,j}}(f_1, f_2, f_3)\|_{L^p} \lesssim 2^{3k-j \over p_3} \|f_1\|_{L^{p_1}} \|f_2\|_{L^{p_2}} \|f_3\|_{L^{p_3}}.
\end{align*}
Then, since
$$ T_{\mathfrak{r}}=\sum_{k\geq 0} 2^{-2k}d_k \sum_{j=0}^k  2^{j-k} c_j T_{\mathfrak{r}_{k,j}}$$
with uniformly bounded coefficients $(d_k)_k$ and $(c_j)_j$, we obtain (by subadditivity with $p_3>1$) that $T_{\mathfrak{r}}$ is also bounded as stated.
\end{proof}

Now, we turn to the proof of Theorem \ref{thm:model:biest:ani}. Proceeding as in Section \ref{sec:red:op:anis:local:L2} (along the lines of \cite{MTTBiest2} and \cite{multilinear_harmonicII}), the quadrilinear form associated to $T_{\mathfrak{m}}$ can be written as a superposition of quadrilinear forms of the type
\begin{align}
\label{eq:def:4:lin:model:coeff}
\Lambda_{\rr Q, \rr P}(f_1, f_2, f_3, f_4):=
\sum_{P \in \rr P} d_P |I_P|^{-\frac{1}{2}} \Big \langle \sum_{\substack{ Q \in \rr Q \\ \widehat{\omega}_{Q_3} \subseteq \omega_{P_1} }} c_Q |I_Q|^{-\frac{1}{2}} \langle f_1, \phi_{Q_1}^1 \rangle  \langle f_2, \phi_{Q_2}^2 \rangle \phi_{Q_3}^3 , \phi_{P_1}^1   \Big \rangle \langle  f_3, \phi_{P_2}^2 \rangle \langle  f_4, \phi_{P_3}^3 \rangle, 
\end{align}
where $\rr P$ and $\rr Q$ are both rank-1 families\footnote{We say that a family of $n$-tiles is of \emph{rank 1} if each $n$-tile $P=(P_1, \ldots, P_n)$ belonging to the collection is determined by any of its components $P_j$ with $1 \leq j \leq n$. See also \cite[Def. 6.3]{multilinear_harmonicII}.} of tri-tiles, $(c_Q)_{Q \in \rr Q}, (d_P)_{P \in \rr P}$ are uniformly bounded sequences of complex numbers and $\phi^j_{P_j},\phi^j_{Q_j}$ are wave packets associated to the corresponding tile. For simplicity, the wave packets will be denoted by $\phi_{P_j},\phi_{Q_j}$. Since the coefficients $(c_Q)_{Q \in \rr Q}, (d_P)_{P \in \rr P}$ can be incorporated into the wave packets (for example $c_Q$ in $\phi_{Q_1}$ and $d_P$ in $\phi_{P_3}$), we can consider them to be equal to 1.

We recall that a tritile $P=(P_1, P_2, P_3)$ is defined by $\ds P_j=I_P \times \omega_{P_j}$ with the spatial $I_P$ a dyadic interval, and the frequency intervals $\omega_{P_j}$ belonging to fixed shifted dyadic collections. Moreover, the $\ds \omega_{P_1} \times \omega_{P_2} \times \omega_{P_3}$ are Whitney cubes with respect to the singularity $\{\eta_1=\eta_2\}$ and intersect the plane $\{\eta_1+\eta_2+\eta_3=0\}$. A similar statement holds for the tritiles in $\Q$.

Theorem \ref{thm:model:biest:ani} is a direct consequence of the following result:
\begin{theorem} \label{thm:model:biest:ani:2}
Let $\rr P$, $\rr Q$ be two finite, rank-1 collections of tritiles. Then the quadrilinear form  $\Lambda_{\Q,\P}$ is bounded from $L^{p_1}\times L^{p_2} \times L^{p_3} \times L^{p_4'}$ into $\C$ for all Lebesgue exponents $p_1, p_2, p_3, p_4$ with $\ds \frac{1}{p_1}+\frac{1}{p_2}+\frac{1}{p_3}=\frac{1}{p_4}$ and $(\frac{1}{p_1},\frac{1}{p_2},\frac{1}{p_3})\in{\bm D}''$. Moreover, the norm of the quadrilinear form is independent of the collections $\P$ and $\Q$.
\end{theorem}

\begin{remark}
As we study the quadrilinear form $\Lambda_{\Q,\P}$, we can obtain restricted weak-type estimates within a wider range than the one stated in Theorem \ref{thm:model:biest:ani:2} above. Some of these can only be used for expressing the boundedness of certain the adjoints of $T_{\mathfrak m}$ inside a part of the quasi-Banach regime, so they are not useful for bounding $T_{\mathfrak m}$ itself. In Subsection \ref{subsec:range} we make explicit the full range of restricted weak-type estimates that can be obtained with our techniques -- this will be given by a collection of 25 extremal points (with some redundancies). If we add the constraints that $1<p_1, p_2, p_3 \leq \infty$, we obtain precisely ${\bm D}''$, as stated in the theorem above.

\end{remark}

By changing the order of summation, we write $\Lambda_{\rr Q, \rr P}(f_1,f_2,f_3,f_4)$ as 
{\fontsize{10}{10}
\begin{equation}
\label{eq:def:form:order}
 \Lambda_{\rr Q, \rr P}(f_1,f_2,f_3,f_4) = \sum_{Q\in\Q} |I_Q|^{-1/2} \langle f_1,\phi_{Q_1}\rangle \langle f_2,\phi_{Q_2}\rangle \langle \sum_{P\in\P \atop \widehat\omega_{Q_3} \subset \omega_{P_1}} |I_P|^{-1/2} \langle f_3, \phi_{P_2}\rangle \langle f_4, \phi_{P_3}\rangle \phi_{P_1}, \phi_{Q_3} \rangle.
\end{equation}}

Thus, it can be represented as $\Lambda_{\rr Q, \rr P}(f_1,f_2,f_3,f_4) = \sum\limits_{Q\in\Q} |I_Q|^{-1/2} \cdot a_{Q_1}^{(1)} \cdot a_{Q_2}^{(2)} \cdot a_{Q_3}^{(3)}$, with the three families of coefficients defined by 
\begin{align}
\label{eq:coeffA}
\begin{split}
& a_{Q_1}^{(1)}:=\langle f_1,\phi_{Q_1}\rangle, \qquad \qquad  a_{Q_2}^{(2)}:=\langle f_2,\phi_{Q_2}\rangle, \\
& a_{Q_3}^{(3)}:=\langle F_{Q_3}, \phi_{Q_3} \rangle, \qquad \textrm{where} \qquad F_{Q_3}:= \sum_{P\in\P \atop \widehat{\omega}_{Q_3} \subset \omega_{P_1}} |I_P|^{-1/2} \langle f_3, \phi_{P_2}\rangle \langle f_4, \phi_{P_3}\rangle \phi_{P_1}.
\end{split}
\end{align}
As we did before in Section \ref{sec:proof:m:+}, we aim to decompose $\rr Q$ into subcollections where we have more information - the trees. This is done according to how much information concentrates along them, i.e. depending on the sizes associated to each of the families of coefficients $\ds \big( a_{Q_j}^{(j)}  \big)_{Q \in \rr Q}$.

\begin{proposition} \label{prop:algo:1}
Let $\rr Q$ be a finite rank-1 collection of tritiles, $i=1, 2$ and let $s_i \in \Z$ be defined by $\ds  2^{s_i} \simeq \size_{\rr Q}^i (\mathfrak{a}^{(i)})$, with $\mathfrak{a}^i=(a_{Q_i}^{(i)})_{Q \in \rr Q}$. Then we have the decomposition $ \ds \Q := \bigsqcup_{m_i\leq s_i} \Q_{m_i}^{(i)}$ where each nonempty $\Q_{m_i}^{(i)}$ is a union of trees $\ds \Q_{m_i}^{(i)}:= \bigcup_{\tr{T}\in {\FF_{m_i}^{(i)}}} \tr{T}$ satisfying
\[
\size_{\Q_{m_i}^{(i)}}^i (\mathfrak{a}^{(i)}) \leq 2^{m_i} \qquad \text{and} \qquad \sum_{\tr{T}\in \FF_{m_i}^{(i)}} |I_{\tr{T}}|\lesssim 2^{-2m_i} \|f_i\|_2^2.
\]
\end{proposition}

When we want to attain a similar decomposition according to the sizes of the third family of coefficients, the procedure becomes more complicated, in particular because the expression $F_{Q_3}$ in \eqref{eq:coeffA} depends as well on the tile $Q_3$. $F_{Q_3}$ is in fact a model operator associated to the rank-1 collection of tritiles $\rr P$, and which interacts with the $Q_3$ tile. While we can still perform the decomposition with respect to the magnitude of $\ds \size_{\rr Q}^3 \big( (a_{Q_3}^{(3)})_Q \big)$, the real difficulty is estimating how many trees are associated to each level set.

We obtain the following result:
\begin{proposition} \label{prop:algo:2} [\textsf{Size-energy discretization}]
Let $\Q$ and $\P$ be two finite, rank-1 collections of tritiles. Let $(a_{Q_3}^{(3)})_{Q \in \Q}$ be a sequence as in \eqref{eq:coeffA} and $s_3 \in \Z$ such that $\ds  2^{s_3} \simeq \size_{\rr Q}^3 (\mathfrak{a}^{(3)})$, with $\mathfrak{a}^{(3)}=(a_{Q_3}^{(3)})_Q$. Then we can decompose $\Q$ as $\ds \Q := \bigsqcup_{m_3\leq s_3} \Q_{m_3}$, where each non-empty $\Q_{m_3}$ is a union of trees $\ds \Q_{m_3}:= \bigcup_{\tr{T}\in {\FF_{m_3}}} T$ satisfying $\ds \size_{\rr Q_{m_3}}^3 (\mathfrak{a}^{(3)}) \leq 2^{m_3}$ and\footnote{The implicit constant in \eqref{eq:est:nb:trees:a_Q_3} depends on $\epsilon$, but we omit this from the notation.} 
\begin{align}
\label{eq:est:nb:trees:a_Q_3}
\sum_{\tr{T}\in \FF_{m_3}} |I_{\tr{T}}|& \lesssim 2^{-(4+\epsilon)m_3} \|f_3\|_2^{2\beta} \, \|f_4\|_2^{2(1-\beta)} \cdot \size_{\P}(f_3)^{(1-\beta)(4+\epsilon)} \, \size_{2,\Q}(f_3)^{\beta(2+\epsilon)} \\
 & \quad \cdot \size_{\P}(f_4)^{\beta (4+\epsilon)} \, \size_{2,\Q}(f_4)^{(1-\beta)(2+\epsilon)}, \nonumber
\end{align}
for any arbitrary $\beta\in(0,1)$ and any $\epsilon>0$.
\end{proposition}

In the statement, we use\footnote{Of course, $\size_{\P}(f_i)$ and $\size_{\P, 2}(f_i)$ are defined in a similar way.} the following \emph{sizes} : for $i=1,..,4$
$$ \size_{\Q}(f_i):= \sup_{Q\in\Q} \big(\frac{1}{|I_Q|} \int |f_i(x)| \chi_{I_Q}(x) dx\big)$$
and its $L^2$-version
$$ \size_{2,\Q}(f_i):= \sup_{Q\in\Q} \big(\frac{1}{|I_Q|} \int |f_i(x)|^2 \chi_{I_Q}(x) dx\big)^{1/2}.$$

We recall that for $i=1,2,3$ and $\mathfrak{a}:=(a^{(i)}_{Q_i})_{Q\in \rr Q}$ a sequence of coefficients then
\[
 \size_{\rr Q}^i (\mathfrak{a}) := \sup_{ \tr{T} \text{ $i$-lacunary tree} \atop \tr{T}\subset \Q} \big(\frac{1}{|I_{\tr{T}}|}  \sum_{Q \in \tr{T}} |a_{Q_i}^{(i)}|^2  \big)^{1 \over 2}.
\]
In the case of our current coefficients, given by \eqref{eq:coeffA}, we have the following estimates. 

\begin{proposition}[\textsf{Size estimate}] \label{prop:size:2}
Let $\rr P, \rr Q$ be two rank-1 collections of tritiles and $\big(a_{Q_j}^{(j)}\big)_{Q \in \Q}$ the families of coefficients given by \eqref{eq:coeffA}, for $1 \leq j \leq 3$. Then we have
\begin{itemize}
\item if $i=1,2$
$$ \size_{\rr Q} (\mathfrak{a}^{(i)}) \lesssim \size_{\rr Q}(f_i)=\sup_{Q\in \rr Q} \frac{1}{|I_Q|} \int |f_i(x)| \chi_{I_Q}(x) dx.$$
\item if $i=3$, we have two estimates
$$ \size_{\rr Q} (\mathfrak{a}^{(3)}) \lesssim \size_{\rr P}(f_3) \cdot \size_{\rr P}(f_4)= \big(\sup_{P\in \rr P} \frac{1}{|I_P|} \int |f_3(x)| \chi_{I_P}(x) dx\big) \cdot \big(\sup_{P\in \rr P} \frac{1}{|I_P|} \int |f_4(x)| \chi_{I_P}(x) dx\big),$$
and also for any $v_3,v_4\in (2,\infty)$ with $\frac{1}{v_3}+\frac{1}{v_4}=\frac{1}{2}$
$$ \size_{\rr Q} (\mathfrak{a}^{(3)}) \lesssim \size_{\rr Q,v_3,v_4}(f_3,f_4) := \sup_{Q\in \rr Q} \big( \frac{1}{|I_Q|} \int |f_3(x)|^{v_3} \chi_{I_Q}(x) dx\big)^{1/v_3} \cdot \big(\frac{1}{|I_Q|} \int |f_4(x)|^{v_4} \chi_{I_Q}(x) dx\big)^{1/v_4}.$$
\end{itemize}
\end{proposition}

These crucial estimates will be proved later in the current section, we are first going to explain how we conclude the proof of Theorem \ref{thm:model:biest:ani:2}.

\subsection{Boundedness of the $4$-linear form}
\label{subsec:range}

We put together all the previous results in order to achieve the boundedness of $\Lambda_{\rr Q, \rr P}$, as presented in Theorem \ref{thm:model:biest:ani:2}. 

\begin{proof}[Proof of Theorem \ref{thm:model:biest:ani:2}]
For $i=1,2$, we set\footnote{The quantities $S_j$, $1 \leq j \leq 3$ depend on the collections $\P$ and $\Q$, but we remove this dependency from the notation in order to simplify the presentation.}
$$ S_i:=\size_{\Q}(f_i)$$
and we denote 
$$S_3= \min\{\size_{\rr Q,v_3,v_4}(f_3,f_4),\size_{\P}(f_3)\cdot \size_{\P}(f_4)\}.$$

So by using the previous statements, we deduce that for coefficients $\theta_1,\theta_2, \theta_3\in (0,1)$ with $\theta_1+\theta_2+\theta_3=1$, we can bound the $4$-linear form from \eqref{eq:def:form:order} in the following way
\begin{align*}
\big| \Lambda_{\Q, \P}(f_1, f_2,  f_3, f_4) \big|  \lesssim \sum_{m_1: 2^{m_1} \leq S_1} \sum_{m_2: 2^{m_2} \leq S_2} \sum_{m_3: 2^{m_3} \leq S_3} 2^{m_1+m_2+m_3} \sum_{\tr{T} \in \FF_{m_1} \cap \FF_{m_2} \cap \FF_{m_3}} |I_{\tr{T}}| \qquad & \\
  \lesssim \sum_{m_1: 2^{m_1} \leq S_1} \sum_{m_2: 2^{m_2} \leq S_2} \sum_{m_3: 2^{m_3} \leq S_3} 2^{m_1+m_2+m_3} \big(\sum_{\tr{T} \in \FF_{m_1} } |I_{\tr{T}}|\big)^{\theta_1} \cdot \big(\sum_{\tr{T}\in \FF_{m_2} } |I_{\tr{T}}|\big)^{\theta_2} \cdot \big(\sum_{\tr{T}\in \FF_{m_3}} |I_{\tr{T}}|\big)^{\theta_3}&.
\end{align*}
So by choosing $2\tilde \alpha_3:=4-4\beta$ and $2\tilde \alpha_4:=4\beta$ with $\beta$ from Proposition \ref{prop:algo:2}, one has
\begin{align}
\label{eq:cond:1}
\begin{cases}
& \theta_1+\theta_2+\theta_3=1, \qquad 0 < \theta_1, \theta_2, \theta_3 <1 \\
&\tilde \alpha_3+\tilde \alpha_4=2, \qquad 0 < \tilde \alpha_3, \tilde \alpha_4<2
\end{cases}
\end{align}
and for $\epsilon>0$ also from Proposition \ref{prop:algo:2}
\begin{align*}
\Big| \Lambda_{\Q, \P}(f_1, f_2, f_3, f_4) \Big| \lesssim  \sum_{m_1: 2^{m_1} \leq S_1} \sum_{m_2: 2^{m_2} \leq S_2} \sum_{m_3: 2^{m_3} \leq S_3}  2^{m_1} 2^{m_2} 2^{m_3} \big(  2^{-2m_1} \|f_1\|_2^2  \big)^{\theta_1} \big(  2^{-2m_2} \|f_2\|_2^2  \big)^{\theta_2}  \\
 \quad \cdot \Big( 2^{-(4+\epsilon)m_3} 
 \size_{\P}(f_3)^{\tilde \alpha_3 (2+\frac{\epsilon}{2})} \size_{2,\Q}(f_3)^{\tilde \alpha_4 \cdot (1+\frac{\epsilon}{2})} \cdot \size_{\P}(f_4)^{\tilde \alpha_4 (2+\frac{\epsilon}{2})} \size_{2,\Q}(f_4)^{\tilde \alpha_3 \cdot (1+\frac{\epsilon}{2})} \cdot
  \|f_3\|_2^{\tilde \alpha_4} \|f_4\|_2^{\tilde \alpha_3} \Big)^{\theta_3}. 
\end{align*}
Consequently, if we have in addition
\begin{equation}
\label{eq:cond:2}
0<\theta_1<\frac{1}{2}, 0<\theta_2<\frac{1}{2} \text{     and     } 0<\theta_3<\frac{1}{4+\epsilon}:=\frac{1}{4}-\tilde \epsilon,
\end{equation}
then the geometric series above are summable and we have
\begin{align}
\label{eq:est:4Lin:form:size:En}
\Big| \Lambda_{\Q, \P}(f_1, f_2, f_3, f_4) \Big| & \lesssim  S_1^{(1-2 \theta_1)} S_2^{(1-2 \theta_2)} S_3^{(1-(4+\epsilon) \theta_3)} \cdot
\size_{\P}(f_3)^{\tilde \alpha_3 (2+\frac{\epsilon}{2})\theta_3} \size_{2,\Q}(f_3)^{\tilde \alpha_4\cdot (1+\frac{\epsilon}{2}) \cdot\theta_3} \\
& \quad \cdot \size_{\P}(f_4)^{\tilde \alpha_4 (2+\frac{\epsilon}{2})\theta_3} \size_{2,\Q}(f_4)^{\tilde \alpha_3\cdot (1+\frac{\epsilon}{2}) \cdot\theta_3}  \|f_1\|_2^{2 \theta_1}  \|f_2\|_2^{2 \theta_2}  \|f_3\|_2^{\tilde \alpha_4 \theta_3} \|f_4\|_2^{ \tilde \alpha_3 \theta_3}. \nonumber
\end{align}

Now, to prove boundedness in the stated range, we only need to prove restricted weak-type estimates close to each of the extremal points in the range, the strong-type result within the convex hull following via multilinear interpolation,

So we fix $E_1, E_2, E_3, E_4$ sets of finite measure. In proving restricted weak-type estimates for the quadrilinear form, the bad\footnote{For more on multilinear interpolation, see for example \cite{multilinear_harmonicII}.} index can take any value in $\{1,2,3,4\}$, so the analysis and the construction of the exceptional subset will depend upon this.

\begin{itemize}[leftmargin=1em]
\item[$\bullet$] If the bad index is the 4th one, we set
\[
\Omega:= \bigcup_{j=1}^4 \Big\lbrace x: M(\one_{E_j})(x) > C \frac{|E_j|}{|E_4|}   \Big\rbrace,
\]
and $\ds E_4':=E_4 \setminus \Omega$ (which is a major subset of $E_4$ for $C$ a large enough numerical constant). Let $f_1,f_2,f_3,f_4$ be functions with $ |f_i|\leq {\bf 1}_{E_i}$ for $i\in\{1,2,3\}$ and $f_4$ with $|f_4|\leq {\bf 1}_{E_4'}$ and we look to control $\Lambda_{\Q,\P}(f_1,f_2,f_3,f_4)$.

For $d,d'\geq 0$, we define the subcollections $\Q_d,\P_{d'}$ as the subcollections of tri-tiles with
\begin{equation}
\label{eq:assumption:distance}
1 + \frac{\dist(I_Q, \Omega^c)}{|I_Q|} \simeq 2^d \qquad \textrm{and} \qquad 1 + \frac{\dist(I_P, \Omega^c)}{|I_P|} \simeq 2^{d'}.
\end{equation}
So $ \Q:=\bigcup_{d\geq 0} \Q_d$ and $\P:=\bigcup_{d'\geq 0} \P_{d'}$. For every $d,d'\geq 0$ fixed, we apply \eqref{eq:est:4Lin:form:size:En} to $\Q_d,\P_{d'}$.
Then we have
\[
S_1 \lesssim 2^d \frac{|E_1|}{|E_4|}, \quad S_2 \lesssim 2^d \frac{|E_2|}{|E_4|}, 
\]
and\footnote{$M$ will always be an integer large enough which could vary from a line to another one.}
\[
S_3 \lesssim \min\Big( \Big( 2^d \frac{|E_3|}{|E_4|} \Big)^{1 \over v_3} 2^{-Md}, \Big( 2^{d'} \frac{|E_3|}{|E_4|} \Big) 2^{-Md'}  \Big).
\]

All these plugged into \eqref{eq:est:4Lin:form:size:En} imply that\footnote{$\delta,\delta'$ will be positive parameters which could vary from line to line.} 
\begin{align*}
\Big| \Lambda_{\Q_d, \P_{d'}}(f_1, f_2, f_3, f_4) \Big| \lesssim  2^{-d \delta}  2^{-d' \delta'} &\Big( \frac{|E_1|}{|E_4|} \Big)^{1-2 \theta_1} \Big( \frac{|E_2|}{|E_4|} \Big)^{1-2 \theta_2}
\left( \Big(\frac{|E_3|}{|E_4|} \Big)^{{\gamma \over v_3}+(1-\gamma)}  \right)^{(1-(4+\epsilon) \theta_3)} \\ & \qquad \Big( \frac{|E_3|}{|E_4|} \Big)^{\tilde \alpha_3\theta_3(2+\frac{\epsilon}{2}) + \tilde \alpha_4 \theta_3(\frac{1}{2}+\frac{\epsilon}{4})} \cdot |E_1|^{\theta_1} |E_2|^{\theta_2}|E_3|^{\frac{\tilde \alpha_4}{2} \theta_3} |E_4|^{\frac{\tilde \alpha_3}{2}\theta_3} \\
\lesssim 2^{-d \delta}  2^{-d' \delta'} & |E_1|^{1-\theta_1} |E_2|^{1-\theta_2} |E_3|^{\beta_3} |E_4|^{-\theta_3-\beta_3}
\end{align*}
with
$$ \beta_3:=\big(1-\frac{\gamma}{v_3'}\big)\big(1-(4+\epsilon)\theta_3\big)+\theta_3(2+\tilde \alpha_3 + O(\epsilon)).$$
After summation over $d,d'\geq 0$, we obtain
\begin{align*}
\big| \Lambda_{\Q, \P}(f_1, f_2, f_3, f_4) \big| & \leq \sum_{d,d'\geq 0} \big| \Lambda_{\Q_d, \P_{d'}}(f_1, f_2, f_3, f_4) \big| \\
& \lesssim  |E_1|^{1-\theta_1} |E_2|^{1-\theta_2} |E_3|^{\beta_3} |E_4|^{-\theta_3-\beta_3}.
\end{align*}

Notice that \eqref{eq:cond:1} and \eqref{eq:cond:2} in fact imply the constraints
\begin{equation}
\label{eq:cond:3}
\begin{gathered}
0 < \theta_1,  \theta_2<\frac{1}{2}, \quad 0< \theta_3< \frac{1}{4}-\tilde \epsilon, \quad \theta_1+\theta_2+\theta_3=1,
  \\
  0<\tilde \alpha_3,\tilde \alpha_4 <2, \quad \tilde \alpha_3+\tilde \alpha_4=2.
 \end{gathered}
\end{equation}
So extremal points for the exponents will be given by the following situations:\footnote{Here we use the notation $\frac{1}{p} \sim \alpha$ to say that $\frac{1}{p}$ is extremely close to $\alpha$, while still remaining in the expected range: for Banach Lebesgue spaces, $0 \leq \frac{1}{p} \leq 1$, and $-\infty < \frac{1}{p} \leq 0$ for H\"older conjugates corresponding to quasi-Banach Lebesgue spaces.}
\begin{itemize}[leftmargin=1.6em]
\item[-] if $\theta_1 \sim \frac{1}{4}$, $\theta_2 \sim \frac{1}{2}$ and $\theta_3 \sim \frac{1}{4}$, which means $\epsilon,\tilde \epsilon\sim 0$ and $\beta_3\sim \frac{1}{4}(2+\tilde \alpha_3)$, hence
\[
\frac{1}{p_1} \sim \frac{3}{4},\quad \frac{1}{p_2} \sim \frac{1}{2},\quad \frac{1}{p_3} \sim \frac{2+\tilde \alpha_3}{4},\quad \frac{1}{p_4'} \sim -\frac{3+\tilde \alpha_3}{4}.
\]
We set aside (taking $\tilde \alpha_3\sim 0$ and $\tilde \alpha_3\sim 2$ respectively)
\[
A_1=\Big( \frac{3}{4}, \frac{1}{2}, \frac{1}{2}, - \frac{3}{4}  \Big), \qquad A_2=\Big( \frac{3}{4}, \frac{1}{2},1, - \frac{5}{4}  \Big).
\]
\item[-] if $\theta_2 \sim \frac{1}{4}$, $\theta_1 \sim \frac{1}{2}$ and $\theta_3 \sim \frac{1}{4}$, which means
\[
\frac{1}{p_1} \sim \frac{1}{2},\quad \frac{1}{p_2} \sim \frac{3}{4},\quad \frac{1}{p_3} \sim \frac{2+\tilde \alpha_3}{4},\quad \frac{1}{p_4'} \sim -\frac{3+\tilde \alpha_3}{4}.
\]
We set aside (taking $\tilde \alpha_3\sim 0$ and $\tilde \alpha_3\sim 2$ respectively)
\[
A_3=\Big( \frac{1}{2}, \frac{3}{4}, \frac{1}{2}, - \frac{3}{4}  \Big), \qquad A_4=\Big( \frac{1}{2}, \frac{3}{4}, 1, - \frac{5}{4}  \Big).
\]
\item[-] finally, if $\theta_3 \sim 0$ and $\epsilon\sim 0$, then $\theta_1 \sim \frac{1}{2}$ and $\theta_2 \sim \frac{1}{2}$, which means $\beta_3\sim 1-\frac{\gamma}{v_3'}$
\[
\frac{1}{p_1} \sim \frac{1}{2}, \quad \frac{1}{p_2} \sim \frac{1}{2}, \quad \frac{1}{p_3} \sim 1-\frac{\gamma}{v_3'}, \quad \frac{1}{p_4'} \sim -1+\frac{\gamma}{v_3'}.
\]
We get (with $\frac{\gamma}{v_3'}\sim 0$ or $\frac{\gamma}{v_3'}\sim 1$ respectively)
\[
A_5=\Big( \frac{1}{2}, \frac{1}{2}, 1, - 1  \Big), \qquad A_6=\Big( \frac{1}{2}, \frac{1}{2}, 0, 0  \Big).
\]
\end{itemize}

\item[$\bullet$] When the bad index is the third one, we exchange the roles of $E_3$ and $E_4$; more precisely, we set
\[
\Omega:= \bigcup_{j=1}^4 \Big\lbrace x: M(\one_{E_j})(x) > C \frac{|E_j|}{|E_3|}   \Big\rbrace,
\]
and under the assumptions \eqref{eq:assumption:distance} we have 
\begin{align*}
S_1 \lesssim 2^d \frac{|E_1|}{|E_3|}, \quad S_2 \lesssim 2^d \frac{|E_2|}{|E_3|}, S_3 \lesssim \min \left( \Big( 2^d\frac{|E_4|}{|E_3|} \Big)^{1 \over v_4} 2^{-Md'}, \Big( 2^{d'} \frac{|E_4|}{|E_3|} \Big) 2^{-Md'}   \right). 
\end{align*} 
Hence (we only inverse the roles of $E_3,E_4$ and $v_3,v_4$, $\tilde \alpha_3,\tilde \alpha_4$ which are symmetric)
\begin{align*}
\Big| \Lambda_{\P, \Q}(f_1, f_2, f_3, f_4) \Big| \lesssim  
|E_1|^{1-\theta_1} |E_2|^{1-\theta_2} |E_3|^{-\theta_3-\beta_4} |E_4|^{\beta_4}
\end{align*}
with 
$$ \beta_4:=\big(1-\frac{\gamma}{v_4'}\big)\big(1-(4+\epsilon)\theta_3\big)+\theta_3(2+\tilde \alpha_4+O(\epsilon)).$$
As previously, such estimates produce $6$ extremal points
\begin{align*}
A_7=\Big( \frac{3}{4}, \frac{1}{2},  - \frac{3}{4} , \frac{1}{2} \Big), \qquad A_8=\Big( \frac{3}{4}, \frac{1}{2},  - \frac{5}{4}, 1 \Big), \quad A_9=\Big( \frac{1}{2}, \frac{3}{4},  - \frac{3}{4} , \frac{1}{2} \Big),\\
A_{10}=\Big( \frac{1}{2}, \frac{3}{4},  - \frac{5}{4}, 1 \Big), \qquad A_{11}=\Big( \frac{1}{2}, \frac{1}{2},  - 1, 1 \Big), \quad A_{12}=\Big( \frac{1}{2}, \frac{1}{2},  0 , 0 \Big).
\end{align*}

\item[$\bullet$] When the bad index is the first one, we set
\[
\Omega:= \bigcup_{j=1}^4 \Big\lbrace x: M(\one_{E_j})(x) > C \frac{|E_j|}{|E_1|}   \Big\rbrace.
\]
We do not have any decay concerning sizes associated to the collection $\P$ (hence we cannot use any decomposition in parameter $d'$), so all we can use is that
\begin{align*}
S_1 \lesssim 2^{-Md}, \quad S_2 \lesssim 2^d \frac{|E_2|}{|E_1|}, \quad S_3 \lesssim \min \left(2^{d \over 2} \Big( \frac{|E_3|}{|E_1|} \Big)^{1 \over v_3} \Big( \frac{|E_4|}{|E_1|} \Big)^{1 \over v_4},1 \right)
\end{align*} 
and
\[
\size_{\P} f_3, \size_{\P} f_4 \lesssim 1.
\]
Overall that yields for any $\gamma\in[0,1]$
\begin{align*}
\Big| \Lambda_{\Q_d, \P}(f_1, f_2, f_3, f_4) \Big| \lesssim  2^{-d \delta} |E_1|^{\beta_1} |E_2|^{1-\theta_2} |E_3|^{\beta_3}   |E_4|^{\beta_4},
\end{align*}
with
$$ \beta_1:=\theta_2-2\theta_3-\frac{\gamma}{2} + 2 \gamma \theta_3, \quad \beta_3:=\frac{\gamma}{v_3}+\theta_3\big(\tilde \alpha_4-\frac{4\gamma}{v_3}\big)$$
and
$$ \beta_4:=\frac{\gamma}{v_4}+\theta_3\big(\tilde \alpha_3 -\frac{4\gamma}{v_4}\big).$$ 
So we can sum over $d\geq 0$ and that gives us the following extremal points of boundedness:
\begin{itemize}[leftmargin=1.6em]
\item[-] if $\theta_1 \sim \frac{1}{4}$ -- so that $\theta_2 \sim \frac{1}{2}$ and $\theta_3 \sim \frac{1}{4}$ (for $\epsilon\sim 0$) -- we get
\[
\frac{1}{p_1} =\beta_1 \sim 0,\quad \frac{1}{p_2} \sim \frac{1}{2},\quad \frac{1}{p_3}=\beta_3 \sim \frac{\tilde \alpha_4}{4},\quad \frac{1}{p_4'} \sim \frac{\tilde \alpha_3}{4}.
\]
This produces, for $(\tilde \alpha_3,\tilde \alpha_4) \sim (0,2)$ or $\sim (2,0)$ respectively, 
\[
A_{13}=\Big(0, \frac{1}{2}, \frac{1}{2}, 0  \Big), \qquad A_{14}=\Big( 0, \frac{1}{2}, 0 , \frac{1}{2} \Big).
\]

\item[-] if $\theta_2 \sim \frac{1}{4}$, and so that $\theta_1 \sim \frac{1}{2}$ and $\theta_3 \sim \frac{1}{4}$ (for $\epsilon\sim 0$), we get
\[
\frac{1}{p_1} \sim -\frac{1}{4}, \quad \frac{1}{p_2} \sim \frac{3}{4}, \quad \frac{1}{p_3} \sim \frac{\tilde \alpha_4}{4},\quad \frac{1}{p_4'} \sim \frac{\tilde \alpha_3}{4}.
\]
This produces, for $(\tilde \alpha_3,\tilde \alpha_4) \sim (0,2)$ or $\sim (2,0)$ respectively, 
\[
A_{15}=\Big( -\frac{1}{4}, \frac{3}{4}, \frac{1}{2}, 0 \Big), \qquad A_{16}=\Big( -\frac{1}{4}, \frac{3}{4}, 0 , \frac{1}{2} \Big).
\]

\item[-] if $\theta_3\sim 0$ with $\theta_1 \sim \frac{1}{2}$ and $\theta_2 \sim \frac{1}{2}$ (for $\epsilon\sim 0$), we get
\[
\frac{1}{p_1} \sim \frac{1-\gamma}{2}, \quad \frac{1}{p_2} \sim \frac{1}{2}, \quad \frac{1}{p_3} \sim \frac{\gamma}{v_3},\quad \frac{1}{p_4'} \sim \frac{\gamma}{v_4}.
\]
This produces, for $(\tilde \alpha_3,\tilde \alpha_4) \sim (0,2)$ or $\sim (2,0)$ respectively and $\gamma\in\{0,1\}$ the following extremal points 
\[
A_{17}=\Big( \frac{1}{2}, \frac{1}{2}, 0, 0 \Big), \qquad A_{18}=\Big( 0, \frac{1}{2}, 0 , \frac{1}{2}\Big) \qquad A_{19}=\Big( 0, \frac{1}{2}, \frac{1}{2} , 0 \Big).
\]

\end{itemize}

\item[$\bullet$] Finally, if the bad index is the second one, we get (by symmetry between $\theta_1,\theta_2$)
\begin{align*}
A_{20}=\Big(\frac{1}{2},0,\frac{1}{2},0  \Big), \quad A_{21}=\Big( \frac{1}{2},0, 0, \frac{1}{2} \Big), \quad A_{22}=\Big(\frac{3}{4},  -\frac{1}{4}, \frac{1}{2},0  \Big), \\
A_{23}=\Big(\frac{3}{4}, -\frac{1}{4} , 0 ,\frac{1}{2} \Big),  \quad A_{24}=\Big(\frac{1}{2}, 0 , 0 , \frac{1}{2}  \Big), \quad A_{25}=\Big(\frac{1}{2},  0, \frac{1}{2}  , 0  \Big).
\end{align*}
\end{itemize}

At the end, we obtain boundedness of the quadrilinear form $\Lambda_{\Q, \P}$ in the range given by the convex hull of the $A_i$'s (we notice that $A_{6}=A_{12}=A_{17}$, $A_{13}=A_{19}$, $A_{14}=A_{18}$, $A_{20}=A_{25}$, $A_{21}=A_{24}$) and by the use of Sage Software\footnote{The Sage Software was used solely for computing the extremal points of ${\bm D}''$.}) we obtain boundedness in ${\bm D}''$ given by the intersection with $[0,1]^3$.

Lastly, we remark that the estimates that we obtain do depend implicitly on the distance from a Lebesgue quadruple $(\frac{1}{p_1}, \frac{1}{p_2}, \frac{1}{p_3}, \frac{1}{p_4'})$ (with $\frac{1}{p_1}+\frac{1}{p_2}+\frac{1}{p_3}+\frac{1}{p_4'}=1$) to an extremal point $A_j$, with $1 \leq j \leq 25$. This is expected and the reason why we only claim boundedness in the interior of ${\bm D}''$.
\end{proof}

We are now left with proving the \emph{size estimates} in Proposition \ref{prop:size:2} and the decomposition algorithm in Proposition \ref{prop:algo:2}. 

\subsection{Size control}
\label{sec:control:size:a_{Q_3}}

The first part of Proposition \ref{prop:size:2} is standard and corresponds to Lemma 2.13 and 6.13 in \cite{multilinear_harmonicII}. It remains to check the second part.

\begin{proof}[Proof of Proposition \ref{prop:size:2}]
If $i=1, 2$, the estimate for $\size_{\rr Q}^i (\mathfrak{a}^{(i)})$, is standard, we only focus on the case $i=3$. According to the definition, we have
\[
 \size_{\rr Q}^3 (\mathfrak{a}^{(3)}) = \sup_{ \tr{T} \text{ $3$-lacunary tree} \atop  \tr{T} \subset \Q} \big(\frac{1}{|I_{\tr{T}}|} \sum_{Q \in \tr{T}} |a_{Q_3}^{(3)}|^2  \big)^{1 \over 2}.
\]

So let $\tr{T}$ be a 3-lacunary tree in $\rr Q$ with top datum $(I_{\tr{T}}, \xi_{\tr{T}})$. We would like to prove that 
\begin{equation}
\label{eq:ineq:ell:2a}
\big( \frac{1}{|I_{\tr{T}}|} \sum_{Q \in \tr{T}} |a_{Q_3}^{(3)}|^2  \big)^{1 \over 2} \lesssim \big(\frac{1}{|I_{\tr{T}}|} \int |f_3(x)|^{v_3} \chi_{I_{\tr{T}}}(x) dx\big)^{1/v_3} \cdot  \big(\frac{1}{|I_{\tr{T}}|} \int |f_4(x)|^{v_4} \chi_{I_{\tr{T}}}(x) dx\big)^{1/v_4},
\end{equation}
which implies the second estimate on $\size_{\rr Q}^3 (\mathfrak{a}^{(3)})$ in Proposition \ref{prop:size:2}. The functions $f_3$ and $f_4$ are associated to the collection $\P$ of tritiles, so we need so understand the interactions between tiles in $\P$ and the tree $\tr{T}$. If $Q \in \tr{T}$ and $P \in \P$ is so that its contribution in \eqref{eq:coeffA} is nonzero, then we necessarily have $\hat \omega_{Q_3} \subset \omega_{P_1}$, which further implies that $\xi_{\tr T} \in 3 \omega_{P_1}$ and $|I_P| \ll |I_Q| \leq |I_{\tr T}|$.

This observation motivates the following definition: for any $\# \in \Z$, we recall that $I_{\tr{T}, \#}$ represents the $\#$-shift of $I_{\tr{T}}$ and we set
\begin{equation}
\label{eq:def:tiles:P:as:tree}
\rr P_{\tr{T}, \#}:= \Big\lbrace  P \in \rr P : I_P \subseteq I_{\tr{T}, \#} \text{  and there exists  } Q' \in \tr{T} \text{   so that   } \hat \omega_{Q'_3} \subseteq \omega_{P_1}  \Big\rbrace.
\end{equation}
Observe that every $\rr P_{\tr{T}, \#}$ has the structure of an 1-overlapping tree, which is lacunary in the second and third directions. As expected, the main contribution will come from the term corresponding to $\#=0$.

Indeed, we first notice that for any $Q \in \tr{T}$, 
\begin{align*}
a_{Q_3}^{(3)}=\sum_{\# \in \Z} \Big\langle    \sum_{P\in\P_{\tr{T}, \#} \atop \widehat{\omega}_{Q_3} \subset \omega_{P_1}} |I_P|^{-1/2} \langle f_3, \phi_{P_2}\rangle \langle f_4, \phi_{P_3}\rangle \phi_{P_1} , \phi_{Q_3} \Big\rangle:= \sum_{\# \in \Z} a_{Q_3}^{(3), \#}.
\end{align*}
Then by triangle inequality we have $ \big(\sum\limits_{Q \in \tr{T}} |a_{Q_3}^{(3)}|^2  \big)^{1 \over 2} \leq \sum\limits_{\# \in \Z} \big(\sum\limits_{Q \in \tr{T}} |a_{Q_3}^{(3), \#}|^2  \big)^{1 \over 2}$. 

For $\#=0$, we set $P_{\tr{T}}=P_{\tr{T},0}$. Next we dualize and make explicit $\big( \sum_{Q \in \tr{T}} |a_{Q_3}^{(3), 0}|^2  \big)^{1 \over 2}$: for $\ds  \big( \epsilon_{Q_3} \big)_{Q_3 \in \tr{T}}$ of $\ell^2$-norm equal to $1$, we have
\begin{align*}
\big( \sum_{Q \in \tr{T}} |a_{Q_3}^{(3),0}|^2  \big)^{1 \over 2}& =  \sum_{Q \in \tr{T}} \epsilon_{Q_3} \, a_{Q_3}^{(3), 0} = \sum_{Q \in \tr{T}} \epsilon_{Q_3} \sum_{\substack{ P \in \rr P_{\tr{T}} \\ \hat \omega_{Q_3} \subseteq \omega_{P_1} }} |I_P|^{-\frac{1}{2}}  \langle  f_3, \phi_{P_2} \rangle \langle  f_4, \phi_{P_3} \rangle  \langle \phi_{P_1} , \phi_{Q_3} \rangle. 
\end{align*}

Now we switch the order of summation, and use that $\rr P_{\tr{T}}$ is ($2$ and $3$)-lacunary: 
\begin{align*}
\big( \sum_{Q \in \tr{T}} |a_{Q_3}^{(3),0}|^2  \big)^{1 \over 2} & \lesssim \sum_{P \in \rr P_{\tr{T}}} |I_P|^{- \frac{1}{2}} \big|\langle f_3, \phi_{P_2} \rangle \big| \big|\langle f_4, \phi_{P_3} \rangle \big|
\big| \langle  \sum_{\substack{Q \in \tr{T} \\ \hat \omega_{Q_3} \subseteq \omega_{P_1} }} \epsilon_{Q_3}  \phi_{Q_3} , \phi_{P_1} \rangle    \big| \\
& \lesssim \int_{\R} \Big(\sum_{P \in \rr P_{\tr{T}}}  \big|\langle f_3, \phi_{P_2} \rangle \big|^2 \frac{\one_{I_P}(x)}{|I_P|}  \Big)^{1 \over 2} \Big(\sum_{P \in \rr P_{\tr{T}}}  \big|\langle f_4, \phi_{P_3} \rangle \big|^2 \frac{\one_{I_P}(x)}{|I_P|}  \Big)^{1 \over 2} \\
&  \qquad \quad \cdot \sup_{P \in \rr P_{\tr{T}}} \one_{I_P}(x) \big| \langle  \sum_{\substack{Q \in \tr{T} \\ \hat \omega_{Q_3} \subseteq \omega_{P_1} }} \epsilon_{Q_3}  \phi_{Q_3} , \frac{1}{|I_P|^{1 \over 2}} \phi_{P_1} \rangle    \big|   dx.
\end{align*}
If $P \in \rr P_{\tr{T}}$ is fixed, then we can re-write $\ds \sum_{\substack{Q \in \tr{T} \\ \hat \omega_{Q_3} \subseteq \omega_{P_1} }} \epsilon_{Q_3}  \phi_{Q_3} $ as $\Pi_{\omega_{P, \tr{T}}} \big(  \sum_{Q \in \tr{T} }  \epsilon_{Q_3} \phi_{Q_3} \big)(x),$
where the interval $\omega_{P, \tr{T}}$ is given by $[\xi_{\tr{T}}-3|\omega_{P_1}|^\frac{1}{3}, \xi_{\tr{T}}+ 3|\omega_{P_1}|^\frac{1}{3}]$. Hence we get for any $x \in \R$,
\[
\sup_{P \in \rr P_{\tr{T}}} \one_{I_P}(x) \big| \langle  \sum_{\substack{Q \in \tr{T} \\ \hat \omega_{Q_3} \subseteq \omega_{P_1} }} \epsilon_{Q_3}  \phi_{Q_3} , \frac{1}{|I_P|^{1 \over 2}} \phi_{P_1}^1 \rangle    \big| \lesssim M \circ M \big(   \sum_{Q \in \tr{T} }  \epsilon_{Q_3} \phi_{Q_3}  \big)(x),
\]
with $M$ being the Hardy-Littlewood maximal operator. H\"older's inequality for $\frac{1}{v_3}+\frac{1}{v_4}=\frac{1}{2}$ yields
\begin{align*}
\big( \sum_{Q \in \tr{T}} |a_{Q_3}^{(3),0}|^2  \big)^{1 \over 2} & \lesssim \Big\| \Big(\sum_{P \in \rr P_{\tr{T}}}  \big|\langle f_3, \phi_{P_2} \rangle \big|^2 \frac{\one_{I_P}}{|I_P|}  \Big)^{1 \over 2}  \Big\|_{v_3} \Big\| \Big(\sum_{P \in \rr P_{\tr{T}}}  \big|\langle f_4, \phi_{P_3} \rangle \big|^2 \frac{\one_{I_P}}{|I_P|}  \Big)^{1 \over 2}  \Big\|_{v_4}  \cdot \| M \circ M \big(   \sum_{Q \in \tr{T} }  \epsilon_{Q_3} \phi_{Q_3} \big) \|_2 \\
& \lesssim \Big\| \Big(\sum_{P \in \rr P_{\tr{T}}}  \big|\langle f_3, \phi_{P_2} \rangle \big|^2 \frac{\one_{I_P}}{|I_P|}  \Big)^{1 \over 2}  \Big\|_{v_3} \Big\| \Big(\sum_{P \in \rr P_{\tr{T}}}  \big|\langle f_4, \phi_{P_3} \rangle \big|^2 \frac{\one_{I_P}}{|I_P|}  \Big)^{1 \over 2}\Big\|_{v_4}. 
\end{align*}
Invoking the $L^{v_j}$-boundedness of the square function, we obtain 
\begin{equation}
\label{eq:what:we:wanted:loc1}
\big(\frac{1}{|I_{\tr T}|} \sum_{Q \in \tr{T}} |a_{Q_3}^{(3),0}|^2  \big)^{1 \over 2} \lesssim \big(\frac{1}{|I_{\tr{T}}|} \int |f_3(x)|^{v_3} \chi_{I_{\tr{T}}}(x) dx\big)^{1/v_3} \cdot  \big(\frac{1}{|I_{\tr{T}}|} \int |f_4(x)|^{v_4} \chi_{I_{\tr{T}}}(x) dx\big)^{1/v_4},
\end{equation}
which is precisely \eqref{eq:ineq:ell:2a}, with $(a_{Q_3}^{(3)})_{Q \in \tr T}$ replaced by $(a_{Q_3}^{(3),0})_{Q \in \tr T}$. 

On the other hand, since $\P_{\tr T}$ is 2-and 3-lacunary, we also have 
\begin{equation}
\label{eq:what:we:wanted:loc2}
\big(\frac{1}{|I_{\tr T}|} \sum_{Q \in \tr{T}} |a_{Q_3}^{(3),0}|^2  \big)^{1 \over 2} \lesssim \size_{\P_{\tr T}} (f_3) \cdot \size_{\P_{\tr T}} (f_4) \lesssim  \size_{\P} (f_3) \cdot \size_{\P} (f_4).
\end{equation}
Similar estimates remains true for $|\#|\leq 2$.

If $|\#|$ is big, we can gain from the decay of the wave packets. First, we notice that for any $Q \in \tr{T}$, 
\begin{align*}
|a_{Q_3}^{(3), \#}| \lesssim  & \int_{\R }  \Big(\sum_{P \in \rr P_{\tr{T}, \#}}  \big|\langle f_3, \phi_{P_2} \rangle \big|^2 \frac{\one_{I_P}(x)}{|I_P|}  \Big)^{1 \over 2}   \cdot \Big(\sum_{P \in \rr P_{\tr{T}, \#}}  \big|\langle f_4, \phi_{P_3} \rangle \big|^2 \frac{\one_{I_P}(x)}{|I_P|}  \Big)^{1 \over 2} dx   \cdot \sup_{P \in \rr P_{\tr{T}}^{\#}} \frac{|\langle  \phi_{P_1}, \phi_{Q_3} \rangle| }{|I_P|^{1 \over 2}}\\
&\lesssim |I_Q|^{- \frac{1}{2}} \Big( \frac{\dist(I_Q, I_{\tr{T}, \#})}{|I_Q|} \Big)^{-M} \Big\| \Big(\sum_{P \in \rr P_{\tr{T}, \#}}  \big|\langle f_3, \phi_{P_2} \rangle \big|^2 \frac{\one_{I_P}}{|I_P|}  \Big)^{1 \over 2}  \Big\|_{q_3} \Big\| \Big(\sum_{P \in \rr P_{\tr{T}, \#}}  \big|\langle f_4, \phi_{P_3} \rangle \big|^2 \frac{\one_{I_P}}{|I_P|}  \Big)^{1 \over 2}\Big\|_{q_4}, 
\end{align*}
for any $q_3, q_4 >1$ with $\frac{1}{q_3}+\frac{1}{q_4}=1$.

Since 
\[
\sum_{Q \in \tr T} |I_Q|^{-1} \Big( \frac{\dist(I_Q, I_{\tr{T}, \#})}{|I_Q|} \Big)^{-2M} \lesssim |\#|^{-2M} |I_{\tr T}|^{-1},
\]
we obtain after summing over $Q \in \tr T$
\begin{align*}
\big( \frac{1}{|I_{\tr{T}}|} \sum_{Q \in \tr{T}} |a_{Q_3}^{(3), \#}|^2  \big)^{1 \over 2} \lesssim  |\#|^{-M} \Big\| \Big(\frac{1}{|I_{\tr T, \#}|}\sum_{P \in \rr P_{\tr{T}, \#}}  \big|\langle f_3, \phi_{P_2} \rangle \big|^2 \frac{\one_{I_P}}{|I_P|}  \Big)^{1 \over 2}  \Big\|_{q_3} \Big\| \Big(\frac{1}{|I_{\tr T, \#}|} \sum_{P \in \rr P_{\tr{T}, \#}}  \big|\langle f_4, \phi_{P_3} \rangle \big|^2 \frac{\one_{I_P}}{|I_P|}  \Big)^{1 \over 2}\Big\|_{q_4}, 
\end{align*}
for any $q_3, q_4 >1$ with $\frac{1}{q_3}+\frac{1}{q_4}=1$.

On one hand, the boundedness of the square function and H\"older implies \eqref{eq:ineq:ell:2a} with $(a_{Q_3}^{(3)})_{Q \in \tr T}$ replaced by $(a_{Q_3}^{(3),\#})_{Q \in \tr T}$, and with some extra decay in $|\#|$. So if we now sum over $\# \in \Z$ we can deduce that \eqref{eq:ineq:ell:2a} holds.

On the other hand, if above we invoke John-Nirenberg (we recall that $\rr P_{\tr{T}, \#}$ is a 2-ans 3-lacunary tree in $\P$), we obtain that 
\[
\big(\frac{1}{|I_{\tr T}|} \sum_{Q \in \tr{T}} |a_{Q_3}^{(3),\#}|^2  \big)^{1 \over 2} \lesssim  2^{-10 |\#|}\size_{\P} (f_3) \cdot \size_{\P} (f_4).
\]
Summing again over $\# \in \Z$, we deduce that 
\[
\big(\frac{1}{|I_{\tr T}|} \sum_{Q \in \tr{T}} |a_{Q_3}^{(3)}|^2  \big)^{1 \over 2} \lesssim \size_{\P_{\tr T}} (f_3) \cdot \size_{\P_{\tr T}} (f_4),
\]
which implies the first estimate on $\size_{\rr Q}^3 (\mathfrak{a}^{(3)})$ in Proposition \ref{prop:size:2}.

\end{proof}

\subsection{Energy control}
\label{sec:control:energy:a_{Q_3}}
We are left with proving Proposition \ref{prop:algo:2}. As in Section \ref{subsec:bddness:4-lin:main}, the general statement follows by iterating a \emph{decomposition result} similar to Lemma \ref{lemma:decomposition:lemma:4:lin:form:no:curv}: we consider $\Q' \subset \Q$ a generic finite collection of tiles so that $\size_{\rr Q'}^3 (\mathfrak{a}^{(3)}) \leq 2^{m_3+1} \leq  \size_{\rr Q}^3 (\mathfrak{a}^{(3)})$ and we start selecting, according to their frequency position, the lacunary trees satisfying
\begin{equation}
\label{eq:tree:selection:cond}
2^{2m_3} \leq \frac{1}{|I_{\tr{T}}|} \sum_{Q\in \tr{T}} |a_{Q_3}^{(3)}|^2 \leq 2^{2(m_3+1)}.
\end{equation}
With each lacunary tree we also remove an overlapping tree having the same top datum, in order to ensure a good separation between the selected trees. Indeed, after this process comes to an end, we have obtained a collection of $3$-lacunary, strongly 3-disjoint trees $\FF_{m_3}$, each coming with an overlapping companion. We are left with assessing $\sum_{\tr{T}\in \FF_{m_3}}|I_{\tr{T}}|$, which is the most involved estimate in this Part I -- this will be the goal of the remaining part of this subsection.

\medskip

In what follows we consider the collections $\P$ and $\Q$ to be fixed and we are going to work solely with the collection of \emph{tiles} $\ds \big\{  Q_3 , Q \in \Q  \big \}$. Since our initial family of tritiles is of rank-1, there is a one-to-one correspondence between the two collections and for simplicity we denote them both $\Q$. From here on we follow the terminology in \cite{DemeterTaoThiele} used in the study of maximal multilinear operators, as our proof of Proposition \ref{prop:algo:2} uses certain maximal Bessel inequalities which were initially introduced in \cite{Lacey} for bounding the bilinear maximal function. We call \emph{forest} a collection of lacunary, mutually strongly disjoint trees. Given a forest $\mathfrak{F}$, the associated \emph{counting function} is defined by 
\begin{equation}
\label{eq:def:counting:func}
 {\mathcal N}_{\mathfrak F} := \sum_{\tr{T} \in \mathfrak F} {\bf 1}_{I_{\tr{T}}}.
\end{equation}
We also use the notation
\begin{equation}
\label{eq:tiles:forest}
\Q_{\mathfrak F} := \bigcup_{\tr T \in \mathfrak F} \tr T
\end{equation}
for the collection of tiles (or tritiles) associated to a forest $\mathfrak F$ in $\Q$.

Proposition \ref{prop:algo:2} follows from the following:
\begin{proposition}
\label{prop:control:en:F_m3}
Let $m_3 \in \Z$, $\epsilon >0$ and let $\FF_{m_3}$ be a forest of tiles in $\Q$. Let $(a_{Q_3}^{(3)})_{Q_3 \in \Q}$ be the sequence defined in \eqref{eq:coeffA}, depending in particular on the rank-1 collection of tritiles $\P$. Assume that 
\begin{equation}
\label{eq:cond:selected:trees}
\text{ for each $\tr T \in \FF_{m_3}$, } \qquad 2^{2m_3} \leq \frac{1}{|I_{\tr{T}}|} \sum_{Q\in \tr{T}} |a_{Q_3}^{(3)}|^2 \leq 2^{2(m_3+1)}
\end{equation}
\begin{equation}
\label{eq:cond:selected:subtrees}
\text{ and for every sub-tree $\tilde{\tr{T}}\subset \tr{T}\in \FF_{m_3}$,\quad $\ds \frac{1}{|I_{\tilde{\tr{T}}}|} \sum_{Q\in \tilde{\tr{T}}} |a_{Q_3}^{(3)}|^2 \leq 2^{2(m_3+1)}$}.
\end{equation}
Then we have, for any $\beta \in (0, 1)$ and any $\epsilon>0$, 
\begin{align}
\| {\mathcal N}_{\FF_{m_3}}  \|_1=\sum_{\tr{T}\in \FF_{m_3}} |I_{\tr{T}}|& \lesssim 2^{-(4+\epsilon)m_3} \|f_3\|_2^{2\beta} \|f_4\|_2^{2(1-\beta)} \cdot \size_{\P}(f_3)^{(1-\beta)(4+\epsilon)} \size_{2,\Q}(f_3)^{\beta (2+\epsilon)} \nonumber \\
 & \quad \cdot \size_{\P}(f_4)^{\beta (4+\epsilon)} \size_{2,\Q}(f_4)^{(1-\beta)(2+\epsilon)}. \label{eq:est:nb:trees:a_{Q_3}}
\end{align}
\end{proposition}

In order to prove the estimate \eqref{eq:est:nb:trees:a_{Q_3}}, we use \emph{good-$\lambda$} estimates for the counting function associated to a forest. More precisely, we use Corollary 9.7 in \cite{DemeterTaoThiele} which states that it is enough to control, for any subforest $\mathfrak F' \subset \mathfrak F$, the $L^1$ norm of $\ic N_{\mathfrak F'}$ and the $BMO$ norm of $\mathfrak F'$ by a certain power of the $L^\infty$ norm of $\ic N_{\mathfrak F'}$. 

Given a forest $\mathfrak{F}$, we recall that its $BMO$ norm is defined by 
\begin{equation}
\label{eq:def:BMO:norm:forests}
\| \mathfrak{F}  \|_{\textrm{BMO}} := \sup_{J \in \ic D}  \frac{1}{|J|} \sum_{\substack{\tr{T} \in \mathfrak{F} \\ I_{\tr{T}} \subseteq J}} |I_{\tr{T}}|.
\end{equation}

Via Corollary 9.7 of \cite{DemeterTaoThiele}, the estimate \eqref{eq:est:nb:trees:a_{Q_3}} is a consequence of
\begin{proposition}
\label{prop:L1:BMO:est:forest}
Under the hypotheses of Proposition \ref{prop:control:en:F_m3}, we let $\ic F$ be a subforest of $\ic F_{m_3}$ satisfying conditions \eqref{eq:cond:selected:trees} and \eqref{eq:cond:selected:subtrees}. Let $\epsilon >0$ and let $\alpha_3, \alpha_4 \in (0, 1)$ with $\alpha_3+\alpha_4=1$. Then 
\begin{equation}
\label{eq:L1:est:forest}
\|\ic N_{\ic F}\|_1 \lesssim 2^{-2 \, m_3} \, \|\ic N_{\ic F}\|_\infty^{\frac{1}{2}+\epsilon} \,  \|f_3\|_2^{2 \, \alpha_3} \, \|f_4\|_2^{2 \, \alpha_4} \, \size_{\P}(f_3)^{2 - 2 \alpha_3} \,  \size_{\P}(f_4)^{2 - 2 \alpha_4}
\end{equation}
and
\begin{equation}
\label{eq:BMO:est:forest}
\| \ic F \|_{\textrm{BMO}} \lesssim 2^{-2 \, m_3} \, \|\ic N_{\ic F}\|_\infty^{\frac{1}{2}+\epsilon} \,  \big( \size_{\Q, 2} (f_3) \big)^{2 \alpha_3} \, \big( \size_{\Q, 2} (f_4) \big)^{2 \alpha_4} \, \size_{\P}(f_3)^{2 - 2 \alpha_3} \,  \size_{\P}(f_4)^{2 - 2 \alpha_4}.
\end{equation}
\end{proposition}

\begin{proof}[Proof of Proposition \ref{prop:control:en:F_m3}]

We are going to apply Corollary 9.7 of \cite{DemeterTaoThiele} with $\mu>0$ such that $\frac{1}{\mu}=\frac{1}{2}+2\epsilon$ and define $\epsilon'$ such that $(4+\epsilon')\cdot (1-4\epsilon)=4$. Then we have
$$ 2\frac{\mu}{\mu-1}=4+\epsilon' \qquad \textrm{and} \qquad \frac{2}{\mu-1}=2+\epsilon'.$$
So Corollary 9.7 of \cite{DemeterTaoThiele} and Proposition \ref{prop:L1:BMO:est:forest} imply that with $\beta=\alpha_3$ (and so $\alpha_4=1-\beta$)
\begin{align*}
\| {\mathcal N}_{\FF_{m_3}}  \|_1 & \lesssim 2^{-(4+\epsilon')m_3} \|f_3\|_2^{2\beta} \|f_4\|_2^{2(1-\beta)} \cdot \size_{\P}(f_3)^{(1-\beta)(4+\epsilon')} \size_{2,\Q}(f_3)^{\beta (2+\epsilon')} \nonumber \\
 & \quad \cdot \size_{\P}(f_4)^{\beta (4+\epsilon')} \size_{2,\Q}(f_4)^{(1-\beta)(2+\epsilon')}. 
\end{align*}
This is exactly \eqref{eq:est:nb:trees:a_{Q_3}} with $\epsilon'$ instead of $\epsilon$. Since $\epsilon$ can be  arbitrarily small and so $\epsilon'\sim 16\epsilon$ as well, that concludes the proof.
\end{proof}

We are left with proving Proposition \ref{prop:L1:BMO:est:forest}. The estimates \eqref{eq:L1:est:forest} and \eqref{eq:BMO:est:forest} are quite technical and require a certain preparation of the collection of tiles. We recall a few results which will be useful later on.
 
First, we will need to avoid a pile-up of the tree tops, for which we invoke a version of Lemma 3.14 in \cite{Lacey}. We adapt this latter lemma to our context using collections of trees in $\Q$ since two disjoint trees could have the same spatial interval.
\begin{lemma}
\label{lemma:Lacey:control:tops}
Let $\mathfrak F$ be a family of trees in $\Q$ and let $\ds \mathcal N_{\mathfrak F}(x):= \sum_{\tr T \in \mathfrak F} \one_{I_{\tr T}}(x)$. Then one can split $\mathfrak F:=\mathfrak F^{(1)} \cup \mathfrak F^{(2)}$ so that 
\begin{equation}
\label{eq:cond:subset:grid:control}
\sum_{\tr T \in \mathfrak F^{(1)}} \big(  M \one_{I_{\tr T}} \big) ^2 (x) \lesssim \|\ic N_{\mathfrak F}\|_{\infty}
\end{equation}
and $\ds \sum_{\tr T \in \mathfrak F^{(2)}}|I_{\tr T}| \leq \frac{1}{100} \sum_{\tr T \in \mathfrak F}|I_{\tr T}|$.
\end{lemma}
\begin{proof}
We set $K:=\| \mathcal N_{\mathfrak F} \|_\infty$ and we notice that we can decompose $\ds \mathfrak F = \bigsqcup_{j=1}^K \mathfrak F_j$, where each subforest $\mathfrak F_j$ contains trees with mutually disjoint spatial intervals and thus has a counting function satisfying 
$$\ds \sum_{\tr T \in \mathfrak F_j} \one_{I_{\tr T}}(x) \leq 1\,.$$ 
We now use \cite[Lemma 3.14]{Lacey} for each collection $\{ I_{\tr T} : \tr T \in \mathfrak F_j\}$ in order to obtain the corresponding decomposition of the subforests : $\mathfrak F_j= \mathfrak F_j^{(1)} \cup \mathfrak F_j^{(2)}$ with
\[
\sum_{\tr T \in \mathfrak F_j^{(1)}} \big( M \one_{I_{\tr T}} \big)^2(x) \lesssim 1
\]
and $\ds \big| \bigcup_{\tr T \in  \mathfrak F_j^{(2)}} I_{\tr T} \big|= \sum_{\tr T \in  \mathfrak F_j^{(2)}} |I_{\tr T}| \leq \frac{1}{1000}  \sum_{\tr T \in  \mathfrak F_j}|I_{\tr T}|$. With these we set
\[
\mathfrak F^{(1)}:=\bigcup_{j=1}^K \mathfrak F^{(1)}_j \qquad \text{and} \qquad \mathfrak F^{(2)}=\bigcup_{j=1}^K \mathfrak F^{(2)}_j
\]
and notice that the two collections satisfy the desired conditions.
\end{proof}

As mentioned earlier, we will also need a certain maximal Bessel inequality associated to a forest $\mathfrak F$, which only holds if we remove an exceptional set of tiles in $\Q_{\mathfrak F}$. 
Indeed, Theorem 10.1 in \cite{DemeterTaoThiele} produces the following result via dualization:
\begin{theorem}(\cite[Theorem 10.1]{DemeterTaoThiele})
\label{thm:max:Bessel}
Let $\tilde \Q$ be a collection of tiles and $\mathfrak F$ a forest of tiles in $\tilde \Q$. Let $\mu>0$ be a positive parameter. There exists an exceptional set $\Q_* \subset \tilde \Q_{\mathfrak F}$ (depending on the parameter $\mu$) such that
\begin{equation}
\label{eq:smallness:exc}
\big| \bigcup_{Q \in \Q_*}I_Q \big| \leq \frac{1}{10} \frac{\|\ic N_{\mathfrak F}\|_1}{\|\ic N_{\mathfrak F}\|_\infty}
\end{equation}
and such that the operator $H^*_{ \tilde\Q_{\mathfrak F}} : \ell^2( \tilde \Q_{\mathfrak F} \setminus \Q_*) \to L^2(\R)$, defined by
\begin{equation}
H^*_{ \tilde \Q_{\mathfrak F}} \big( \big( c_{Q_3} \big)_{Q \in \tilde \Q_{\mathfrak F} \setminus \Q_* } \big):= \sup_{\kappa \in \Z} \Big| \sum_{T \in \mathfrak F} \sum_{Q \in \tilde\Q_{\mathfrak F} \setminus \Q_*} c_{Q_3} \phi_{Q_3}(x) \cdot \one_{\{ |I_{Q_3}|>2^\kappa \}}(x)  \Big|
\end{equation}
is bounded, satisfying
\begin{equation}
\big\|  H^*_{ \tilde \Q_{\mathfrak F}} \Big(\big( c_{Q_3} \big)_{Q \in \tilde \Q_{\mathfrak F} \setminus \Q_* }\Big)\big\|_2 \lesssim \|\ic N_{\mathfrak F}\|_\infty^{\frac{1}{\mu}}  \, \| \big( c_{Q_3} \big)_{Q \in \tilde\Q_{\mathfrak F} \setminus \Q_* }  \|_{\ell^2}.
\end{equation}
\end{theorem}

With this we can begin proving our result
\begin{proof}[Proof of Proposition \ref{prop:L1:BMO:est:forest}]
We recall that we intend to prove the estimates \eqref{eq:L1:est:forest} and \eqref{eq:BMO:est:forest}, which concern $\|\ic N_{\ic F}\|_1$ and $\| \ic F  \|_{BMO}$. First we apply Lemma \ref{lemma:Lacey:control:tops} to the collection $\ds \mathfrak I:= \{ I_{\tr T} : \tr T \in \ic F  \}$. This will naturally split $\ds \ic F= \ic F^{(1)} \sqcup \ic F^{(2)}$, and since
\[
\big\|  \ic N_{\ic F}  \big\|_1 \simeq \big\|  \ic N_{\ic F^{(1)}}  \big\|_1 
\]
it will be enough to prove \eqref{eq:L1:est:forest} with the left hand side replaced by $\|\ic N_{\ic F^{(1)}}\|_1$. Notice that the trees in $\ic F^{(1)}$ still satisfy conditions \eqref{eq:cond:selected:trees} and \eqref{eq:cond:selected:subtrees}, but we also have 
\begin{equation}
\label{eq:top:control}
\sum_{\tr T \in \ic F^{(1)}} \big( M \one_{I_{\tr T}} \big)^2(x) \lesssim \|\ic N_{\ic F} \|_\infty,
\end{equation}
a condition that will be important in Proposition \ref{prop:energy} below.

Next we apply Theorem \ref{thm:max:Bessel} to the forest $\ic F^{(1)}\subset \Q$ to produce an exceptional set of tiles $\ds \Q_* \subset \Q_{\ic F^{(1)}}$. We set $\ds  \Omega:= \bigcup_{Q\in \Q_*} I_Q$, and we claim that 
\begin{equation} \label{eq:omega}
\sum_{Q \in \Q_{\ic F^{(1)}} \atop I_Q \subset \Omega} |a_{Q_3}^{(3)}|^2 \leq \frac{1}{2} 2^{2m_3} \|{\mathcal N}_{\FF^{(1)}} \|_1.
\end{equation}
For any $\tr T \in \ic F^{(1)}$, we consider the collection of maximal dyadic intervals in $\{I_Q \subset \Omega,\ Q\in \tr{T}\}$. For any interval $I$ in this collection we have, due to the subtree condition \eqref{eq:cond:selected:subtrees}, that
$$ \sum_{Q\in \tr{T} \atop I_Q \subset I} |a_{Q_3}^{(3)}|^2 \leq |I \cap I_{\tr{T}}| \cdot 2^{2m_3+2}.$$
By summing this information for all such maximal intervals $I$, we get
$$ \sum_{Q\in \tr{T} \atop I_Q \subset \Omega} |a_{Q_3}^{(3)}|^2 \leq |\Omega \cap I_{\tr{T}}| \cdot 2^{2m_3+2}.$$
Now we sum over $\tr T \in \ic F^{(1)}$:
\begin{align*}
 \sum_{\tr{T}\in \FF^{(1)}} \sum_{Q\in \tr{T} \atop I_Q \subset \Omega} |a_{Q_3}^{(3)}|^2 & \leq \sum_{\tr{T}\in \FF^{(1)}} |\Omega \cap I_{\tr{T}}| \cdot 2^{2m_3+2} \leq 2^{2m_3+2} \|{\mathcal N}_{\FF^{(1)}} \|_\infty |\Omega| \\
 & \leq \frac{1}{2} 2^{2m_3} \|{\mathcal N}_{\FF^{(1)}} \|_1,
 \end{align*}
where we have used the smallness of $|\Omega|$ as given by condition \eqref{eq:smallness:exc}.

Since the trees in $\FF^{(1)}$ satisfy condition \eqref{eq:cond:selected:trees}, we have 
\begin{align*}
\big\|  \ic N_{\ic F}  \big\|_1 \simeq \big\|  \ic N_{\ic F^{(1)}}  \big\|_1 & \leq 2^{-2 \, m_3} \sum_{ \tr T \in \ic F^{(1)}} \sum_{Q \in \tr T} |a_{Q_3}|^2 \\
& \leq 2^{-2 \, m_3} \sum_{ \tr T \in \ic F^{(1)}} \sum_{Q \in \tr T  \atop I_Q \subset \Omega} |a_{Q_3}|^2 + ^{-2 \, m_3} \sum_{ \tr T \in \ic F^{(1)}} \sum_{Q \in \tr T \setminus \Q_*} |a_{Q_3}|^2.
\end{align*}
Due to \eqref{eq:omega}, we deduce that
\begin{align*}
\big\|  \ic N_{\ic F}  \big\|_1 \simeq \big\|  \ic N_{\ic F^{(1)}}  \big\|_1 \simeq 2^{-2 \, m_3} \sum_{ \tr T \in \ic F^{(1)}} \sum_{Q \in \tr T \setminus \Q_*} |a_{Q_3}|^2.
\end{align*}

Notice that \eqref{eq:L1:est:forest} follows from proving 
\begin{align}
\label{eq:est:coeff:F^1}
\Big( \sum_{ \tr T \in \ic F^{(1)}} \sum_{Q \in \tr T \setminus \Q_*} |a_{Q_3}|^2 \Big)^{1 \over 2} \lesssim  \|\ic N_{\ic F}\|_\infty^{\frac{1}{4}+{\epsilon \over 2}} \,  \|f_3\|_2^{ \alpha_3} \, \|f_4\|_2^{ \alpha_4} \, \size_{\P}(f_3)^{1 - \alpha_3} \,  \size_{\P}(f_4)^{1 -\alpha_4}.
\end{align}

For simplicity, we introduce a notation which allows us to reshuffle the collection $\ic F^{(1)}$: $\ic F'$ is the collection of trees defined by 
\[
\ic F':= \big\{ \tr T \setminus \Q_* : \tr T \in \ic F^{(1)}   \big\}.
\]
Here we want to preserve the structure of the forest, so $\tr T \setminus \Q_* $ remains a tree with the same top datum as $\tr T$, although it might be empty. In particular, we have $\|\ic N_{\ic F'} \|_1 = \big\|  \ic N_{\ic F^{(1)}}  \big\|_1 \eqsim \big\|  \ic N_{\ic F}  \big\|_1$.

Notice that the trees in $\ic F'$ do not necessarily verify \eqref{eq:cond:selected:trees}, but they do verify \eqref{eq:cond:selected:subtrees}. The left hand of \eqref{eq:est:coeff:F^1} becomes $\ds \Big( \sum_{ \tr T \in \ic F'} \sum_{Q \in \tr T } |a_{Q_3}|^2 \Big)^{1 \over 2}$, which we can easily dualize: there exists a sequence $\ds \big( c_{Q_3} \big)_{Q \in \Q_{\ic F'}}$ with $\ds \| \big( c_{Q_3} \big)_{Q \in \Q_{\ic F'}} \|_{\ell^2}=1$ so that 
\[
\big( \sum_{ \tr T \in \ic F'} \sum_{Q \in \tr T } |a_{Q_3}|^2 \big)^{1 \over 2}=  \sum_{ \tr T \in \ic F'} \sum_{Q \in \tr T } c_{Q_3} \, a_{Q_3}= \sum_{ \tr T \in \ic F'} \sum_{Q \in \tr T } c_{Q_3} \, \big\langle \sum_{P\in\P \atop \hat \omega_{Q_3} \subset \omega_{P_1}} |I_P|^{-1/2} \langle f_3, \phi_{P_2}\rangle \langle f_4, \phi_{P_3}\rangle \phi_{P_1}, \phi_{Q_3} \big\rangle.
\]

By changing the order of summation, the term above can be written as
\begin{align}
\label{eq:change:order:form:P}
\sum_{P\in\P} |I_P|^{-{1\over 2}}  \langle \sum_{\tr{T}\in \FF'} \sum_{Q\in \tr{T} \atop \hat \omega_{Q_3} \subset \omega_{P_1}} c_{Q_3} \phi_{Q_3},\phi_{P_1} \rangle \, \langle f_3,\phi_{P_2}\rangle \, \langle f_4,\phi_{P_3}\rangle := \sum_{P\in\P} |I_P|^{-{1\over 2}} \cdot b_{P_1}^{(1)} \cdot b_{P_2}^{(2)} \cdot b_{P_3}^{(3)}.
\end{align}

Before moving on, we notice that Theorem \ref{thm:max:Bessel} yields in particular a bound for the function
\begin{equation}
\label{eq:def:H*}
H^*_{\Q_{\ic F'}}\big(  \big( c_{Q_3} \big)_{Q \in \Q_{\ic F'} }  \big) (x) :=\sup_{\kappa \in \Z/3} \Big| \sum_{T \in \FF '} \sum_{Q \in \Q_{\FF'}} c_{Q_3} \phi_{Q_3}(x) \cdot \one_{\{ |\omega_{Q_3}|<2^\kappa \}}(x)  \Big|,
\end{equation}
where $\mu=\frac{1}{\epsilon}$, for $\epsilon>0$ as in the hypotheses of Proposition \ref{prop:L1:BMO:est:forest}.

Indeed we get that 
\begin{equation}
\label{eq:bddnedd:max:bessel:function}
\big\|  H^*_{\Q_{\ic F'}}\Big(  \big( c_{Q_3} \big)_{Q \in \Q_{\ic F'} }  \Big) \big\|_2 \lesssim \|\ic N_{\ic F^{(1)}}\|_{\infty}^{\epsilon} \, \| \big( c_{Q_3} \big)_{Q \in \Q_{\ic F'}} \|_{\ell^2}. 
\end{equation}

We now return to \eqref{eq:change:order:form:P}. Since $\P$ is a rank-1 collection, we can invoke right away the abstract \emph{size-energy estimate} of \cite[Proposition 6.5]{MTTBiest2}. We recall that the \emph{energy} of a sequence  $\mathfrak{b}^{(j)}:=(b_{P_j}^{(j)})_{P\in \P}$ is defined by 
\begin{equation}
\label{energydef}
\energy_{\P}^j(\mathfrak{b}):= \sup_{n} \sup_{{\mathbb F}} 2^n \big(\sum_{\tr{T}' \in {\mathbb F}} |I_{\tr{T}'}|\big)^{1/2}
\end{equation}
where the supremum is taken over $n \in \Z$ and $\rr F$ ranges over all collections of $j$-disjoint trees $\tr{T}'\subset \P$ satisfying
\begin{itemize}
\item for every $\tr{T}' \in{\mathbb F}$, $\ds  2^{2n} |I_{\tr{T}'}|\leq \sum_{P\in \tr{T}'} |b_{P_j}|^2$; 
\item for every subtree $\tilde{\tr{T}}' \subset \tr{T}'\in {\mathcal F}$, we have $\ds \sum_{P\in \tilde{\tr{T}}'} |b_{P_j}|^2 \leq 4\cdot 2^{2n} |I_{\tilde{\tr{T}}'}|$.
\end{itemize}

Applying this result, we obtain, for any $\alpha_3,\alpha_4\in(0,1)$ with $\alpha_3+\alpha_4=1$,
\begin{equation}
\label{eq:*sizeenergyy}
2^{m_3} \,\|\ic N_{\ic F} \|_1^{\frac{1}{2}} \lesssim \energy_{\P}^1 (\mathfrak{b}^{(1)}) \cdot 
\Big( \size_{\P}^2 (\mathfrak{b}^{(2)})\Big)^{1-\alpha_3}\,\Big(\energy_{\P}^2 (\mathfrak{b}^{(2)})\Big)^{\alpha_3} \cdot \Big(\size_{\P}^3 (\mathfrak{b}^{(3)})\Big)^{1-\alpha_4}\Big(\energy_{\P}^3(\mathfrak{b}^{(3)})\Big)^{\alpha_4}.
\end{equation}
We already know (see, for example,  \cite[Lemma 6.7 and 6.8]{MTTBiest2}) that, for $j=2, 3$
\begin{align*}
\size_{\P}^j (\mathfrak{b}^{(j)}) \lesssim \size_{\P}(f_{j+1}), \qquad \energy_{\P}^j (\mathfrak{b}^{(j)}) \lesssim \| f_{j+1} \|_2.
\end{align*}

We are left with estimating $\energy_{\P}^1 (\mathfrak{b}^{(1)})$, which is more involved given that 
\[ b_{P_1}= \langle \sum_{\tr{T}\in \FF'} \sum_{Q\in \tr{T} \atop \hat \omega_{Q_3} \subset \omega_{P_1}} c_Q \phi_{Q_3},\phi_{P_1} \rangle.
\]
This step will result from Proposition \ref{prop:energy} below, which shows that $\energy_{\P}^1(\mathfrak{b}^{(1)}) \lesssim \|{\mathcal N}_{\FF}\|_\infty^{{1\over 4}+\epsilon}$, for $\epsilon>0$.\footnote{Again, the implicit constant does depend on $\epsilon$.}

With this, we can now conclude the proof of \eqref{eq:L1:est:forest}. Indeed, by applying \eqref{eq:*sizeenergyy}, the previous size-energy estimates and \eqref{eq:energy:P}, we have
$$ 2^{m_3} \,\|\ic N_{\ic F} \|_1^{\frac{1}{2}}\lesssim \|{\mathcal N}_{\FF}\|_\infty^{{1\over 4}+\epsilon} \cdot \size_{\P}(f_3)^{1-\alpha_3} \|f_3\|_2^{\alpha_3} \cdot \size_{\P}(f_4)^{1-\alpha_4}\|f_4\|_2^{\alpha_4},$$
which proves \eqref{eq:L1:est:forest}.

The estimate \eqref{eq:BMO:est:forest} follows in a similar manner: we fix $J$ a dyadic interval and we resume the proof of estimate \eqref{eq:L1:est:forest} with the difference that now we only consider the trees in $\ic F$ with spatial intervals contained inside $J$. Everything follows the same reasoning, except that we have a localized energy estimate:
\[
\energy_{\P}^1(\mathfrak{b}^{(1)}) \lesssim \sup_{P \in \P} \Big( 1+\frac{\dist(I_P, J)}{|I_P|} \Big)^{-100} \|{\mathcal N}_{\FF}\|_\infty^{{1\over 4}+\epsilon},
\]
which can be tracked in the proof of Lemma \ref{eq:lemma:loss}, as a consequence of \eqref{est:decay:en}.
\end{proof}

We are left with estimating $\energy_{\P}^1(\mathfrak{b}^{(1)})$, where $\P$ is a rank-1 collection of tritiles -- which is done in the following proposition:

\begin{proposition} \label{prop:energy} Let $\Q$, $\FF, \FF^{(1)}, \ic F'$, $\epsilon>0$ be as in the proof Proposition \ref{prop:L1:BMO:est:forest} so that, in particular, \eqref{eq:bddnedd:max:bessel:function} and \eqref{eq:top:control} hold. Let $\P$ be a rank-1 collection of tritiles and let the sequence $\mathfrak{b}^{(1)}$ be defined by \eqref{eq:change:order:form:P}. Then
\begin{equation}
\label{eq:energy:P}
\energy_{\P}^1(\mathfrak{b}^{(1)}) \lesssim \|{\mathcal N}_{\FF}\|_\infty^{{1\over 4}+\epsilon}.
\end{equation}
\end{proposition}

To prove it, we need to introduce the notion of \textit{incidence} of a tile in $\rr P$, with respect to the forest $\ic F' \subset \rr Q$. For this, we consider structures in $\P$ which are relevant with respect to the forest $\ic F' \subset \Q$. First we set\footnote{With the above notation, notice that if $I_P \cap 2I_{\tr T} \neq \emptyset$, then $I_P \subseteq 3I_{\tr T}$ since $|I_P|\ll |I_{\tr T}|$. } 
\begin{equation}
\label{eq:def:tiles:tree:P_T}
\rr P_{\tr{T}}=\rr P_{\tr{T}}^{0}:= \Big\lbrace  P \in \rr P : I_P \subseteq  3 I_{\tr{T}}  \text{  and there exists  } Q' \in \tr{T} \text{   so that   } \hat \omega_{Q'_3} \subseteq \omega_{P_1}  \Big\rbrace
\end{equation}
and for $\ell \geq 1$,
\begin{equation}
\label{eq:def:tiles:P:as:tree:ell}
\rr P_{\tr{T}}^{\ell}:= \Big\lbrace  P \in \rr P : I_P \cap 2^{\ell}I_{\tr T}=\emptyset, \quad I_P \cap 2^{\ell+1}I_{\tr T}\neq \emptyset  \text{  and there exists  } Q' \in \tr{T} \text{   so that   } \hat \omega_{Q'_3} \subseteq \omega_{P_1}  \Big\rbrace.
\end{equation}
Notice that every $\P_{\tr T}^{\ell}$, with $\ell \geq 0$, is an $1$-overlapping tree, so it is lacunary in the second and third directions.

\begin{definition}[\textsf{Incidence}]
\label{def:incidence}
 Let $\FF'$ be the forest appearing in the proof of Proposition \ref{prop:L1:BMO:est:forest} and let $\P$ be a rank-1 collection of tritiles. For any $P \in \P$, we define the \emph{incidence of $P$ with respect to the forest $\ic F'$} to be 
\begin{equation}
\label{eq:def:incidence}
\inc(P):= \sum_{\ell\geq 0} 2^{-3\ell} \cdot \sharp \{\tr{T}\in \FF',\ P\in \P_{\tr{T}}^\ell\}.
\end{equation} 
\end{definition} 
 
This notion is motivated by the estimate in Lemma \ref{eq:lemma:loss} below.

In order to work with $\energy_{\P}^1(\mathfrak{b}^{(1)})$ and use the particular form of the coefficients $b_{P_1}^{(1)}$, we dualize the energy by making appeal to \cite[Lemma 6.3]{MTTBiest2}. Thus there exists $\rr F$ a collection of strongly $1$-disjoint trees in $\rr P$ and complex coefficients $\ds (d_{P_1})_{P \in \P_{\rr F}}$ so that $\ds \| (d_{P_1})_{P \in \P_{\rr F}}  \|_{\ell^2}=1$ and
\begin{equation}
\label{eq:dualization:energy}
\energy_{\P}^1(\mathfrak{b}^{(1)}) \simeq \sum_{P\in\P_{\rr F}} d_{P_1} b_{P_1}^{(1)} = \sum_{\tr{T}' \in \rr F } \sum_{P \in \tr{T}'} \sum_{\tr{T} \in \FF'} \sum_{\substack{ Q \in \tr{T} \\ \hat \omega_{Q_3} \subseteq \omega_{P_1}}} c_{Q_3} d_{P_1} \langle  \phi_{P_1}, \phi_{Q_3}  \rangle.
\end{equation}

This quantity can be estimated in the following way:
\begin{lemma}
\label{eq:lemma:loss}
Let $\P, \Q$ be rank-1 collections of tritiles, $\P' \subset \P$ a subcollection and $\ic F' \subset \Q$ a forest. Let $\rr F$ be a collection of strongly $1$-disjoint trees in $\P'$. Let $\ds (c_{Q_3})_{Q \in \Q}$ and $\ds (d_{P_1})_{P \in \P}$ be two sequences of coefficients. Then 
\begin{align}
\label{eq:crucial:en:est:2}
& \Big|  \sum_{\tr{T}' \in \rr F} \sum_{P \in \tr{T}'} \sum_{\tr{T} \in \FF'} \sum_{\substack{ Q \in \tr{T} \\ \hat \omega_{Q_3} \subseteq \omega_{P_1}}} c_{Q_3} d_{P_1} \langle  \phi_{P_1}, \phi_{Q_3}  \rangle   \Big|   \\
&  \lesssim \sum_{\ell \geq 0} 2^{-10 \ell} \Big( \sum_{\tr{T} \in \FF'} \sum_{Q \in \tr{T}}|c_{Q_3}|^2 \Big)^{1 \over 2}  \Big( \sum_{\tr{T}' \in \rr F} \sum_{P \in \tr{T}'} |d_{P_1}|^2 \Big)^{1 \over 2} \cdot \big( \sup_{P \in \tr{T}' \in \rr F} \# \{ \tr{T} \in \FF' : P\in\P_{\tr{T}}^\ell\} \big)^{1 \over 2}.   \nonumber
\end{align}
\end{lemma}

\begin{proof} Notice that in the expression above a couple of tiles $Q \in \tr T \in  \FF'$ and $P \in \P'_{\rr F}$ have a non-zero contribution if and only if $\hat \omega_{Q_3} \subset \omega_{P_1}$. That implies in particular that $\xi_{\tr T} \in 3 \omega_{P_1}$ and $|I_P|\leq |I_Q|\leq |I_{\tr{T}}|$. 

Recalling the definitions \eqref{eq:def:tiles:tree:P_T} and \eqref{eq:def:tiles:P:as:tree:ell}, we observe that all the tri-tiles $P \in \P$ which have a non-zero contribution in \eqref{eq:crucial:en:est:2} belong to some $\P_{\tr T}^{\ell}$, for some $\tr T \in \FF'$ and $\ell \geq 0$.

We denote by $A$ the left hand side of \eqref{eq:crucial:en:est:2}. With the above notation and after changing the order of summation this becomes
\begin{align*}
A=\big|\sum_{\ell\geq 0} \sum_{\tr{T}' \in \rr F} \sum_{\tr{T} \in \FF'} \sum_{P \in \tr{T}' \cap \P_{\tr T}^{\ell}} d_{P_1}  \sum_{\substack{ Q \in \tr{T} \\ \hat{\omega}_{Q_3} \subset \omega_{P_1}}} c_{Q_3} \langle \phi_{Q_3} , \phi_{P_1}  \rangle \big|.
\end{align*}
In order to track the (expected) decay in $\ell$, we notice that
$$\langle \phi_{Q_3}, \phi_{P_1} \rangle = \langle \tilde\phi_{Q_3}^{\tr T}, \phi_{P_1} \,\ci_{I_{\tr{T}}}\rangle,$$
where $\ds \ci_{I_{\tr T}}(x):= (1+\frac{|x-c(I_T)|^2}{{|I_{\tr T}|^2}})^{-N}$ and $\tilde\phi_{Q_3}^{\tr T}(x):= \phi_{Q_3}(x) \ci^{-1}_{I_{\tr T}}(x)$ is defined in frequency by
\[
\ii F \big( \tilde\phi_{Q_3}^{\tr T}\big)(\xi) = q_{I_{\tr T}}(\partial_{\xi})(\hat \phi_{Q_3})(\xi).
\]
Observe that $\ds q_{I_{\tr T}}(\partial_{\xi})$ is a differential polynomial with coefficients depending on $I_{\tr T}$ so in particular $\tilde\phi_{Q_3}^{\tr T}$ is still supported in frequency on $\omega_{Q_3}$.

Consequently, $A$ can be written as
\begin{align*}
\big|\sum_{\ell\geq 0} \sum_{\tr{T}' \in \rr F} \sum_{\tr{T} \in \FF'} \sum_{P \in \tr{T}' \cap \P_{\tr T}^{\ell}} d_{P_1}  \langle \Pi_{W_{\tr T, P}} \big( \sum_{   \substack{ Q \in \tr{T} \\ \hat{\omega}_{Q_3} \subset \omega_{P_1}}} c_{Q_3}  \tilde\phi_{Q_3}^{\tr T} \big), \phi_{P_1} \,\ci_{I_{\tr{T}}} \rangle \big|,
\end{align*}
where $\Pi_{W_{\tr T, P}}$ represents a smooth frequency projection operator associated to the interval $ W_{\tr T, P}:= [\xi_{\tr{T}}-|\omega_{P_1}|^{1 \over 3}, \xi_{\tr{T}}+ |\omega_{P_1}|^{1 \over 3}]$.

Then $A$ is bounded above by
\begin{align}
A& \lesssim \sum_{\ell\geq 0} \sum_{\tr{T}' \in \rr F} \sum_{\tr{T} \in \FF'} \sum_{P \in \tr{T}' \cap \P_{\tr{T}}^\ell} |d_{P_1}| \frac{1}{|I_P|^{1 \over 2}} \int_{\R} M\big(   \sum_{ Q \in \tr{T}} c_{Q_3} \phi_{Q_3} \ci_{I_{\tr{T}}}^{-1} \big)(x) \cdot \chi_{I_P}(x) \ci_{I_{\tr{T}}}(x) dx \label{est:decay:en} \\
&\lesssim \sum_{\ell\geq 0} 2^{-10 \ell} \sum_{\tr{T}' \in \rr F} \sum_{\tr{T} \in \FF'}  \sum_{P \in \tr{T}' \cap \P_{\tr{T}}^{\ell}} |d_{P_1}| |I_P|^{1 \over 2} \frac{1}{|I_P|} \int_{\R}  M\big(\sum_{ Q \in \tr{T}} c_{Q_3} \phi_{Q_3} \ci_{I_{\tr{T}}}^{-1}\big)(x) \cdot \chi_{I_P}(x) dx. \nonumber
\end{align}

Now we use the fact that the trees in the forest $\rr F$ are $1$-disjoint: if they share information in frequency then they must be disjoint in space. Hence all the tiles in $\P_{\tr{T}}^\ell$ have mutually disjoint space intervals since they all contain $\xi_{\tr{T}}$ in frequency. 

After Cauchy-Schwarz in $\tr{T}' \in \rr F, \tr{T} \in \FF'$, and $P \in \tr{T}'$, we bound the above estimate by
\begin{align*}
A & \lesssim \sum_{\ell\geq 0} 2^{-10\ell} \Big( \sum_{\tr{T}' \in \rr F } \sum_{\tr{T} \in \FF'} \sum_{P \in \tr{T}' \cap \P_{\tr{T}}^\ell} |d_{P_1}|^2   \Big)^{1 \over 2} \cdot  \\
& \qquad \qquad \Big( \sum_{\tr{T}' \in \rr F} \sum_{\tr{T} \in \FF'} \sum_{P \in \tr{T}' \cap \P_{\tr{T}}^\ell} |I_P|   \big| \inf_{y \in I_P} \big( M \circ M \big(\sum_{ Q \in \tr{T}} c_{Q_3} \phi_{Q_3} \ci_{I_{\tr{T}}}^{-1}\big) (y)\big|^2 \Big)^{1 \over 2} \\
&\lesssim \sum_{\ell\geq 0} 2^{-10 \ell} \Big( \sum_{\tr{T}' \in \rr F} \sum_{P \in \tr{T}'} |d_{P_1}|^2  \cdot \big( \sup_{P \in \tr{T}' \in \rr F} \# \{ \tr{T} \in \FF' : P\in\P_{\tr{T}}^\ell\} \big)  \Big)^{1 \over 2} \cdot \\
& \qquad \qquad  \Big( \sum_{\tr{T} \in \FF'} \|\sum_{Q \in \tr{T}} c_{Q_3} \phi_{Q_3} \ci_{I_{\tr{T}}}^{-1}\|_2^2  \Big)^{1 \over 2}.
\end{align*}
This proves our Lemma since the $L^2$ norm in the display above satisfies, due to usual orthogonality arguments in a fixed tree,
\[
 \|\sum_{Q \in \tr{T}} c_{Q_3} \phi_{Q_3} \ci_{I_{\tr{T}}}^{-1}\|_2 \lesssim (\sum_{Q \in \tr T} |c_{Q_3}|^2)^{1 \over 2}.
\]
\end{proof}

We remark that the last term in the right hand side of \eqref{eq:crucial:en:est:2} resembles $\ds \sup_{P \in \P_{\rr F} \cap \P'}  \mathfrak{i}(P)$ which was the motivational factor for introducing the notion of incidence.

\begin{lemma}
Let $\ic F'$ be a forest in $\Q$ and $\P$ a collection of tritiles. Then for every $P \in \P$, 
\begin{equation}
\label{eq:upper:bound:incidence}
\mathfrak{i}(P) \lesssim \| \sum_{\tr{T} \in \FF'} \big( M({\bf 1}_{I_{\tr{T}}})\big)^2 \|_{\infty} 
\end{equation}

\begin{proof}
We notice that whenever $P \in \P_{\tr T}^{\ell}$, for $\ell \geq 0$ and $\tr T$ a tree, we have
 \[
  M({\bf 1}_{I_{\tr{T}}})(x) \cdot \one_{I_P}(x) \simeq 2^{-\ell} \one_{I_P}(x)
 \]
and so
\[
 \sharp \{\tr{T} \in \FF',\ P\in \P_{\tr{T}}^\ell\} \cdot  \one_{I_P}(x) \lesssim 2^{2 \, \ell} \sum_{\tr T \in \FF'}  \big(M{\bf 1}_{I_{\tr{T}}}\big)^2(x) \cdot  \one_{I_P}(x).
\]
This implies the desired estimate
\[
\inc(P)=\sum_{\ell\geq 0} 2^{-3\ell} \cdot \sharp \{\tr{T}\in \FF',\ P\in \P_{\tr{T}}^\ell\} \lesssim \| \sum_{\tr{T} \in \FF'} \big( M({\bf 1}_{I_{\tr{T}}})\big)^2 \|_{\infty}.
\]
\end{proof}
\end{lemma}

We are now ready to prove the \emph{energy estimate}:

\begin{proof}[Proof of Proposition \ref{prop:energy}]
We begin by dualizing the energy: as in \eqref{eq:dualization:energy}, there exists $\rr F$ a collection of strongly 1-disjoint trees and $\ds (d_{P_1})_{P \in \P_{\rr F}}$ a sequence of complex numbers with  $\ds \| (d_{P_1})_{P \in \P_{\rr F}}  \|_{\ell^2}=1$ so that 
\begin{equation*}
\energy_{\P}^1(\mathfrak{b}^{(1)}) \simeq \sum_{P\in\P_{\rr F}} d_{P_1} b_{P_1}^{(1)} = \sum_{\tr{T}' \in \rr F } \sum_{P \in \tr{T}'} \sum_{\tr{T} \in \FF'} \sum_{\substack{ Q \in \ic F' \\ \hat{\omega}_{Q_3} \subseteq \omega_{P_1}}} c_{Q_3} d_{P_1} \langle  \phi_{P_1}, \phi_{Q_3}  \rangle.
\end{equation*}

Next, we split the collection $\P_{\rr F}$ in two parts, according to the incidence of the tiles: for $\delta\in(0,1)$ to be chosen later, we set
$$\P_1:=\{P\in \P_{\rr F},\ \inc(P)\leq \|{\mathcal N}_{\FF'}\|_\infty^\delta\} \qquad \textrm{and} \qquad \P_2:=\{P\in \P_{\rr F},\ \|{\mathcal N}_{\FF'}\|_\infty^\delta< \inc(P)\}.$$ 

As a consequence of Lemma \ref{eq:lemma:loss}, we have 
\begin{equation}
\label{eq:prop:energy1}
\Big| \sum_{\tr{T}' \in \rr F} \sum_{P\in \tr{T}' \atop P\in \P_1} \sum_{\tr{T}\in\FF'} \sum_{Q\in \tr{T} \atop \hat{\omega}_{Q_3} \subset \omega_{P_1}} c_{Q_3} d_{P_1} \langle \phi_{Q_3}, \phi_{P_1} \rangle \Big| \lesssim \|{\mathcal N}_{\FF'}\|_\infty^{\delta/2}.
\end{equation}

Next we would like to estimate
\begin{equation}
\label{eq:prop:energy:2}
\sum_{\tr{T}'\in\rr F} \sum_{P\in \tr{T}' \atop P\in \P_2} \sum_{\tr{T}\in\FF'} \sum_{Q\in \tr{T} \atop \hat{\omega}_{Q_3} \subset \omega_{P_1}} c_{Q_3} d_{P_1} \langle \phi_{Q_3}, \phi_{P_1} \rangle.
\end{equation}
Since the incidence of the tiles in $\P_2$ is large, the term above is more difficult to control.

By definition of the incidence $\inc(P)$, we have for arbitrary nonnegative quantities $x_P$ that
\begin{equation}
\label{eq:conseq:incid}
\sum_{\tr{T}'\in\rr F} \sum_{P\in \tr{T}' \cap \P_2} x_P = \sum_{\ell \geq 0} \sum_{\tilde {\tr{T}}\in\FF'} \sum_{P\in \P_{\tilde {\tr{T}}}^\ell \cap \P_2} 2^{-3\ell} \frac{x_P}{\inc(P)}.
\end{equation}
We would like to apply the identity above to the expression
\begin{equation*}
x_P:= \Big|d_{P_1} \langle \phi_{P_1}, \sum_{\tr{T}\in\FF'} \sum_{Q\in \tr{T} \atop \hat{\omega}_{Q_3} \subset \omega_{P_1}} c_{Q_3} \phi_{Q_3} \rangle\Big|.
\end{equation*}

Notice that
\begin{align*}
x_P  = \Big| d_{P_1} \langle \phi_{P_1}, \sum_{\tr{T}\in\FF'} \sum_{Q\in \tr{T} \atop |\omega_{Q_3}| \leq |\omega_{P_1}|^{1/3}} c_{Q_3} \phi_{Q_3} \rangle \Big|  \leq |d_{P_1}| |I_P|^{1/2} \big(\frac{1}{|I_P|} \int H^*(y)  \chi_{I_P}(y) dy \big),
\end{align*}
with $H^*$ given by \eqref{eq:def:H*}:
$$ H^*(x):=  H^*_{\Q_{\ic F'}}\big(  \big( c_{Q_3} \big)_{Q \in \Q_{\ic F'} }  \big) (x)  =\sup_{\kappa \in \Z/3} |\sum_{\tr{T}\in\FF'} \sum_{Q\in \tr{T} \atop |\omega_{Q_3}|\leq 2^{\kappa}} c_{Q_3} \phi_{Q_3}(x)|$$ 
and satisfying $\ds \| H^* \|_2 \lesssim \|\mathcal N_{\ic F'}\|_{\infty}^{\epsilon}$ -- see \eqref{eq:bddnedd:max:bessel:function} (coming from Theorem \ref{thm:max:Bessel}).

Combining all of the above and applying Cauchy-Schwarz, we get
\begin{align*}
&\big|  \sum_{\tr{T}'\in \rr F} \sum_{P\in \tr{T}' \atop P\in \P_2} \sum_{\tr{T}\in\FF'} \sum_{Q\in \tr{T} \atop \hat{\omega}_{Q_3} \subset \omega_{P_1}} c_{Q_3} d_{P_1} \langle \phi_{Q_3}, \phi_{P_1} \rangle \big| \\
& \qquad \lesssim \sum_{\ell\geq 0} \sum_{\tr{T}\in\FF'} \sum_{P\in \P_{\tr{T}}^\ell \atop P\in \P_2} 2^{-3\ell} \frac{|d_{P_1}|}{\inc(P)} |I_P|^{1/2} \big(\frac{1}{|I_P|} \int H^*(y) \chi_{I_P}(y) dy\big) \\
& \qquad \lesssim  \Big(\sum_{\ell\geq 0} \sum_{\tr{T}\in\FF'} \sum_{P\in \P_{\tr{T}}^\ell} 2^{-3\ell} \frac{|d_{P_1}|^2}{\inc(P)}\Big)^{1/2} \Big(\sum_{\ell\geq 0} \sum_{\tr{T}\in\FF'} \sum_{P\in \P_{\tr{T}}^\ell \atop P\in\P_2} 2^{-3\ell} \frac{|I_P|}{\inc(P)} \big(\frac{1}{|I_P|} \int H^*(y) \chi_{I_P}(y) dy\big)^2 \Big)^{1/2}.
 \end{align*}
 
Due to \eqref{eq:conseq:incid}, the first quantity is controlled by 
 \begin{align*}
 \Big(\sum_{\ell\geq 0} \sum_{\tr{T}\in\FF'} \sum_{P\in \P_{\tr{T}}^\ell \cap \P_2} 2^{-3\ell} \frac{|d_{P_1}|^2}{\inc(P)}\Big)^{1/2} & = \Big(\sum_{P\in \P_2} |d_{P_1}|^2\Big)^{1/2} \lesssim 1.
 \end{align*}
 For the second quantity, we involve the maximal function and the definition of $\P_2$ to get
\begin{align*}
& \Big(\sum_{\ell\geq 0} \sum_{\tr{T}\in\FF'} \sum_{P\in \P_{\tr{T}}^\ell \cap \P_2} 2^{-3\ell} \frac{|I_P|}{\inc(P)} \big(\frac{1}{|I_P|} \int H^*(y) \chi_{I_P}(y) dy\big)^2 \Big)^{1/2} \\
& \qquad \qquad \lesssim \|{\mathcal N}_{\FF'}\|_\infty^{-\delta/2} \Big(\sum_{\ell\geq 0} \sum_{\tr{T}\in\FF'} \sum_{P\in \P_{\tr{T}}^\ell} 2^{-3\ell} \int_{I_P} M[ H^*](x)^2 dx \Big)^{1/2}.
\end{align*}
Notice that for every $\tr{T}$ fixed, the tiles $P\in \P_{\tr{T}}^\ell$ are overlapping in frequency ($\omega_{P_1}$ contains $\xi_{\tr{T}}$) and so the intervals $(I_P)_{P\in \P_{\tr{T}}^\ell}$ are mutually disjoint because $\rr F$ is a collection of strongly $1$-disjoint trees. Moreover, all these intervals are contained in $2^\ell I_{\tr{T}}$, so
\begin{align*}
& \Big(\sum_{\ell\geq 0} \sum_{\tr{T}\in\FF'} \sum_{P\in \P_{\tr{T}}^\ell \cap \P_2} 2^{-3\ell} \frac{|I_P|}{\inc(P)} \big(\frac{1}{|I_P|} \int H^*(y) \chi_{I_P}(y) dy\big)^2 \Big)^{1/2} \\
& \qquad \qquad \lesssim\|{\mathcal N}_{\FF'}\|_\infty^{-\delta/2} \Big(\sum_{\ell \geq 0} \sum_{\tr{T}\in\FF'} 2^{-3 \ell} \int M[ H^*](x)^2 \chi_{2^\ell I_{\tr{T}}}(x) dx \Big)^{1/2}. 
\end{align*}
On the other hand, since $\FF^{(1)}$ and $\FF'$ have the same top datum, we have by \eqref{eq:top:control}
\[
\sum_{\tr{T}\in\FF'} \chi_{2^\ell I_{\tr{T}}}(x) \lesssim 2^{2\ell} \sum_{\tr{T}\in \FF'} [M({\bf 1}_{I_{\tr{T}}})(x)]^2 \lesssim
2^{2\ell} \|{\mathcal N}_{\FF}\|_\infty.
\]
Summing over $\ell \geq 0$ and using the $L^2$-control of $H^*$ we conclude
\begin{equation}
\label{eq:est:energy:2}
\big| \sum_{\tr{T}'\in\rr F} \sum_{P\in \tr{T}' \cap \P_2} \sum_{\tr{T}\in\FF'} \sum_{Q\in \tr{T}\atop \hat{\omega}_{Q_3} \subset \omega_{P_1}} c_{Q_3} d_{P_1} \langle \phi_{Q_3}, \phi_{P_1} \rangle \big|\lesssim_\epsilon \|{\mathcal N}_{\FF}\|_\infty^{(1-\delta)/2+\epsilon}.
\end{equation}

Given \eqref{eq:prop:energy1} and \eqref{eq:est:energy:2}, if we set $\delta=\frac{1}{2}+\epsilon$, we obtain \eqref{eq:energy:P}.
\end{proof}

\section{Treatment of the low frequency and high frequency non-stationary components}
\label{sec:low}
In this section, we provide the proofs of Theorems \ref{thm:low} and \ref{thm:non-sta}, carefully tracking the range of boundedness of our operators and detailing the way in which the previous multilinear time-frequency analysis is implemented for each of them. Indeed, as we are going to explain, the low oscillatory and high oscillatory non-stationary components of $\THT$ can be essentially reduced to the model operators studied in the last three sections.

\subsection{A preamble}

This section collects preliminary facts and expository remarks. Though not strictly necessary for understanding the $T^{\low}$ or $T^{\nonsta}$ components, these results provide helpful context for the subsequent work.

Recalling the decomposition and notation introduced in Section \ref{preldiscret}, we record here that each single-scale $T_k,T_{k,m_1,m_2}$ operators are uniformly bounded in $k,m_1,m_2$: 

\begin{proposition} \label{prop:Tk} Let $p_1,p_2,p_3\in[1,\infty]$ be fixed exponents and set $\frac{1}{p}:=\frac{1}{p_1}+\frac{1}{p_2}+\frac{1}{p_3}$. If $p\geq \frac{1}{2}$ then $T_k$ and $T_{k,m_1,m_2}$ are bounded from $L^{p_1} \times L^{p_2} \times L^{p_3} \to L^p$, uniformly in $k,m_1,m_2 \in {\mathbb Z}$.
\end{proposition}

\begin{proof} If $p\geq 1$, then the result is very easy since by Minkowski inequality, the $L^p$-norm can be put inside/under the integral along $t$. So there is nothing to prove in the Banach range.

Hence we now focus on the quasi-Banach range $p<1$. Let us first consider the case $k\geq 0$. In such a case, we write $T_k(f_1,f_2,f_3)(x) = \tilde T_k(f_1^k,f_2^k,f_3^k)(2^k x)$
with 
$$\tilde T_k (F_1,F_2,F_3)(x) := \int_{\R} F_1(x-t)\,F_2(x+t+2^{-2k}t^3)\,F_3(x-2^{-2k} t^3)\, \rho(t)\,dt$$
and $f_i^k:=f_i(2^{-k} \cdot)$. Because of the H\"older scaling $\frac{1}{p}=\frac{1}{p_1}+\frac{1}{p_2}+\frac{1}{p_3}$, we have
$$ \|T_k\|_{L^{p_1} \times L^{p_2} \times L^{p_3} \to L^p} = \|\tilde T_k\|_{L^{p_1} \times L^{p_2} \times L^{p_3} \to L^p}.$$

So let us now check and prove the boundedness of this new operator $\tilde T_k$. We decompose the real line ${\mathbb R}$ into intervals $I$ of length one and so
$$ \|\tilde T_k (F_1,F_2,F_3)\|_p = \big( \sum_{I} \|\tilde T_k (F_1,F_2,F_3)\|_{L^p(I)}^p \big)^{1/p}.$$
Now, since $p<1$, we have for every interval $I$ of length $|I|=1$
\begin{align*}
\|\tilde T_k (F_1,F_2,F_3)\|_{L^p(I)} & \leq \| \tilde T_k (F_1,F_2,F_3)\|_{L^1(I)} \\
& \leq \iint_{I \times [1/2,2]} |F_1(x-t)| |F_2(x+t+2^{-2k}t^3)| |F_3(x-2^{-2k} t^3)| \, dxdt.
\end{align*}
By performing a change of variables, this last quantity is bounded by $\|F_1\|_{L^\infty(4I)}\|F_2\|_{L^1(4I)}\|F_3\|_{L^1(4I)}$ (because $k\geq 0$ so the localisation of $x\in I$ translates into a similar localization for the three functions) and all the other permutations (putting the $L^\infty$ norms on $F_2$ or $F_3$). So by interpolation it is bounded also by 
$$ \|F_1\|_{L^{p_1}(4I)}\|F_2\|_{L^{p_2}(4I)}\|F_3\|_{L^{p_3}(4I)}$$
because $p\geq \frac{1}{2}$ and $|I|=1$.
We then obtain the desired boundedness, since by H\"older, we have
\begin{align*}
\|\tilde T_k (F_1,F_2,F_3)\|_p & \lesssim \big( \sum_{I} \|F_1\|_{L^{p_1}(4I)}^p \|F_2\|_{L^{p_2}(4I)}^p\|F_3\|_{L^{p_3}(4I)}^p \big)^{1/p} \\
& \lesssim \big( \sum_{I} \|F_1\|_{L^{p_1}(4I)}^{p_1}\big)^{1/p_1} \big( \sum_{I} \|F_2\|_{L^{p_2}(4I)}^{p_2} \big)^{1/p_2} \big( \sum_{I} \|F_3\|_{L^{p_3}(4I)}^{p_3} \big)^{1/p_3} \\
& \lesssim \|F_1\|_{L^{p_1}} \|F_2\|_{L^{p_2}}\|F_3\|_{L^{p_3}},
\end{align*}
as desired.

For $k\leq 0$, we follow a very similar argument (exchanging the role of the scale $2^k$ by the scale $2^{3k}$ which is now the smallest), with the corresponding operator
$$\bar T_k (F_1,F_2,F_3)(x) := \int_{\R} F_1(x-2^{2k}t)\,F_2(x+2^{2k}t+t^3)\,F_3(x-t^3)\, \rho(t)\,dt.$$
The analogous result for $T_{k, m_1,m_2}$ follows the same steps and hence will be omitted here.
\end{proof}

Notice that the boundedness range for the single scale operator is significantly wider than the one corresponding to the full operator $\THT$ in Theorem \ref{thm:main}. This phenomenon is a consequence of the difficulty of controlling the multi-scale interactions in the hybrid scenario, and, interestingly enough, it sharply contrasts with the purely curved analogue treated in \cite{HuLie} and \cite{HL25} in which the boundedness range of the full operator matches---up to end-points---the  boundedness range of the corresponding single scale operator.

\subsection{The low oscillation component $T^{\low}$}

\label{subsec:low}

\begin{proof}[Proof of Theorem \ref{thm:low}]
When there is essentially no oscillation, i.e. when $m_1,m_2\leq 0$, the phase appearing in the definition of the symbol $\mathfrak{m}_{k,m_1,m_2}(\xi_1, \xi_2, \xi_3)$ morally behaves like $\ds \max \big( 2^{m_1}, 2^{m_2} \big)$. In order to take advantage of this decay, we use Taylor expansions for the complex exponentials:
\begin{align*}
&\Big(  \int_{\R}  e^{-i\,\frac{\xi_1-\xi_2}{2^k}\,t -i\,\frac{\xi_3-\xi_2}{2^{3k}}\,t^3}  \,\rho(t)\,dt  \Big) \psi\left(\frac{\xi_1-\xi_2}{2^{k+m_1}}\right)\,\psi\left(\frac{\xi_3-\xi_2}{2^{3k+m_2}}\right)\\
=&\sum_{j_1, j_2 \geq 0 } \frac{1}{{j_1}!} \frac{1}{{j_2}!} \Big( \int_{\R}  t^{j_1} \, t^{3 \, j_2} \,\rho(t)\,  dt \Big)  \Big(  -i\frac{\xi_1-\xi_2}{2^k}\Big)^{j_1}  \Big(-i\frac{\xi_3-\xi_2}{2^{3k}}\Big)^{j_2} \psi\left(\frac{\xi_1-\xi_2}{2^{k+m_1}}\right)\,\psi\left(\frac{\xi_3-\xi_2}{2^{3k+m_2}}\right) \\
=&:\sum_{j_1, j_2 \geq 0 } \frac{1}{{j_1}!} \frac{1}{{j_2}!} c(j_1, j_2) \Big(  \frac{\xi_1-\xi_2}{2^k}\Big)^{j_1}  \Big(\frac{\xi_3-\xi_2}{2^{3k}}\Big)^{j_2} \psi\left(\frac{\xi_1-\xi_2}{2^{k+m_1}}\right)\,\psi\left(\frac{\xi_3-\xi_2}{2^{3k+m_2}}\right),
\end{align*}
with coefficients $c(j_1,j_2)$ satisfying the following properties:
\begin{itemize}
\item $c(0,0)=\int \rho(t) dt=0$ by the cancellation property on $\rho$;
\item $\rho$ being supported on $[-4,4]$, $|c(j_1,j_2)| \lesssim 4^{j_1+3j_2}$.
\end{itemize}
In what follows, we split our analysis into several cases:

\medskip
\noindent \underline{Case $j_1=j_2=0$.} The corresponding term is vanishing because $c(0,0)=\int_{\R} \rho(t) dt=0$.

\medskip
\noindent \underline{Case $j_1,j_2\geq 1$.}  In this case, we have a certain decay both in parameters $2^{m_1}$ and $2^{m_2}$ (recall that $m_1,m_2\leq 0$) -- the associated symbol is
\[
\sum_{k \in \Z} \sum_{m_1, m_2<0} 2^{m_1\, j_1} 2^{m_2 \, j_2}  \Big(  \frac{\xi_1-\xi_2}{2^{k+m_1}}\Big)^{j_1} \psi\left(\frac{\xi_1-\xi_2}{2^{k+m_1}}\right)\,  \Big(\frac{\xi_3-\xi_2}{2^{3k+m_2}}\Big)^{j_2}  \psi\left(\frac{\xi_3-\xi_2}{2^{3k+m_2}}\right).
\]
If we fix $m_1, m_2<0$, it is enough to study and obtain bounds independent of $m_1$ and $m_2$ for operators associated to the symbol
\begin{equation}
\sum_{k \in \Z} \Big(  \frac{\xi_1-\xi_2}{2^{k+m_1}}\Big)^{j_1} \psi\left(\frac{\xi_1-\xi_2}{2^{k+m_1}}\right)\,  \Big(\frac{\xi_3-\xi_2}{2^{3k+m_2}}\Big)^{j_2}  \psi\left(\frac{\xi_3-\xi_2}{2^{3k+m_2}}\right),
\label{eq:symbolm1m2}
\end{equation}
which takes the form
\begin{equation}
\label{eq:underlying:symbol}
\sum_{k \in \Z}  \tilde \psi_{(j_1)}\left(\frac{\xi_1-\xi_2}{2^{k+m_1}}\right)\, \tilde \psi_{(j_2)}\left(\frac{\xi_3-\xi_2}{2^{3k+m_2}}\right).
\end{equation}
Clearly, both functions $ \tilde \psi_{(j_1)}$ and $ \tilde \psi_{(j_2)}$ are of $\Psi$-type; therefore we can invoke Proposition \ref{prop:facile} in order to deduce the boundedness of the operator within the Banach range.

So operators corresponding to symbols given by \eqref{eq:symbolm1m2} are bounded on $L^{p_1}\times L^{p_2}\times L^{p_3}\to L^p$ for $1<p_1,p_2,p_3,p < \infty $, uniformly with respect to $m_1,m_2$. This concludes the boundedness of the corresponding component in $T^{\low}$.

\medskip
\noindent \underline{Case $j_1=0$ and $j_2\geq 1$.} We have to deal with 
\[
\sum_{j_2 \geq 1} \frac{1}{{j_2}!} c(0,j_2)\sum_{k \in \Z} \sum_{m_1, m_2<0}\,  2^{j_2m_2}  \, \psi\left(\frac{\xi_1-\xi_2}{2^{k+m_1}}\right)  \Big(\frac{\xi_3-\xi_2}{2^{3k+m_2}}\Big)^{j_2}  \psi\left(\frac{\xi_3-\xi_2}{2^{3k+m_2}}\right). 
\]
Since we have enough decay in $j_2 \geq 1$ and in $m_2 \leq -1$, we consider these parameters fixed.

The summation in $m_1$ is split into two regimes: $m_1\in [m_2/3,-1]$ and $m_1 \leq m_2/3$. In this way, we need to consider operators associated to the symbol
\begin{align*}
&\sum_{m_2/3 \leq m_1 \leq -1} \sum_{k \in \Z} \psi\left(\frac{\xi_1-\xi_2}{2^{k+m_1}}\right)  \Big(\frac{\xi_3-\xi_2}{2^{3k+m_2}}\Big)^{j_2}  \psi\left(\frac{\xi_3-\xi_2}{2^{3k+m_2}}\right) \\
+ & \sum_{k \in \Z} \varphi\left(\frac{\xi_1-\xi_2}{2^{k+\frac{m_2}{3}}}\right)  \Big(\frac{\xi_3-\xi_2}{2^{3k+m_2}}\Big)^{j_2}  \psi\Big(\frac{\xi_3-\xi_2}{2^{3k+m_2}}\Big):=I+II 
\end{align*}
Each term appearing in the expression $(I)$, corresponding to $m_1\in [m_2/3,0]$, can be handled by invoking again Proposition \ref{prop:facile} thus deducing that the associated operators are bounded inside the Banach range. So summing in $m_1\in[m_2/3,0]$ gives a bound with a growth of order $(1+|m_2|)$ which is acceptable and summable over $m_2\leq 0$ due to the decay $2^{j_2 m_2}$ with $j_2\geq 1$.
For the second term $(II)$, after a change of variables\footnote{Here we slightly abuse the language. Strictly speaking, this change of variable assumes that $\frac{m_2}{3}\in\Z$. In general one has to involve a superposition reasoning depending on the congruence classes modulo $3$.}  in the summation $k+\frac{m_2}{3} \to k$, we can invoke Theorem \ref{thm:main-ani:biest} (and Theorem \ref{thm:proof:main:1}) to deduce that the corresponding component is bounded, uniformly in $m_2$. Thanks to the extra decay coming from $2^{j_2 m_2}$ (with $j_2\geq 1$ and $m_2 \leq -1$), we can sum in $m_2$ and deduce that the part of $T^{\low}$ for $j_1=0<j_2$ is bounded from $L^{p_1}\times L^{p_2}\times L^{p_3}\to L^p$ for $(\frac{1}{p_1},\frac{1}{p_2},\frac{1}{p_3})\in {\bm D}$.

\medskip
\noindent \underline{Case $j_2=0$ and $j_1\geq 1$.} We have to deal with 
\[
\sum_{j_1 \geq 1} \frac{1}{{j_1}!} c(j_1,0)  \sum_{k \in \Z} \sum_{m_1, m_2<0} 2^{j_1 m_1} \,\Big(\frac{\xi_1-\xi_2}{2^{k+m_1}}\Big)^{j_1} \psi\left(\frac{\xi_1-\xi_2}{2^{k+m_1}}\right)    \psi\left(\frac{\xi_3-\xi_2}{2^{3k+m_2}}\right). 
\]
A similar reasoning with the one above, allows us to reduce matters to the symbol 
\[
\sum_{k \in \Z} \tilde \psi\left(\frac{\xi_1-\xi_2}{2^{k}}\right)\, \tilde\varphi\left(\frac{\xi_3-\xi_2}{2^{3k}}\right),
\]
whose boundedness is proved in Theorem \ref{thm:main-ani:biest}. This concludes the proof of Theorem \ref{thm:low}.
\end{proof}

\subsection{The high frequency non-stationary component $T^{\nonsta}$}

\label{subsec:nonsta}

\begin{proof}[Proof of Theorem \ref{thm:non-sta}]

By definition, $\mathfrak{m}_{k}^{\nonsta}$ is composed out of three pieces ($m_1\leq 0<m_2$, $m_2\leq 0<m_1$ and $m_1,m_2>0$ with $|m_1-m_2|\geq 3$), all of them corresponding to the situation when the phase function $\lambda(t)=-\frac{\xi_3-\xi_2}{2^{3k}} t^3 - \frac{\xi_1-\xi_2}{2^k}t$ has a large magnitude ($|\lambda(t)| \sim 2^{\max(m_1,m_2)}\geq 1$) and has no critical points ($|\lambda'(t)|\simeq 2^{\max(m_1,m_2)}$). As a consequence, we can apply the non-stationary phase principle in order to extract a suitable decay while preserving the structure of the phase. 

Indeed, recall that
\begin{align}
\label{eq:formula:symbol}
\mathfrak{m}_{k,m_1,m_2}(\xi_1,\xi_2,\xi_3) & = \big[\int_{\R} e^{- i \frac{\xi_1-\xi_2}{2^k}t  - i \frac{\xi_3-\xi_2}{2^{3k}} t^3}  \rho(t)\,dt\big] \, \psi\left(\frac{\xi_1-\xi_2}{2^{k+m_1}}\right)\, \psi\left(\frac{\xi_3-\xi_2}{2^{3k+m_2}}\right) \\
&:=  \big[\int_{\R} e^{ i \lambda(t)}  \rho(t)\,dt\big] \,  \psi\left(\frac{\xi_1-\xi_2}{2^{k+m_1}}\right)\, \psi\left(\frac{\xi_3-\xi_2}{2^{3k+m_2}}\right). \nonumber
\end{align}

In what follows, we treat separately each of the three cases:
\smallskip
\begin{itemize}[leftmargin=1.5em]
\item \underline{\textsf{Case 1.}} For $m_1\leq 0 < m_2$, one morally has
$$
\sum_{m_1\leq 0} \mathfrak{m}_{k,m_1,m_2}(\xi_1,\xi_2,\xi_3) \approx 2^{-N m_2} \tilde \varphi\left(\frac{\xi_1-\xi_2}{2^{k}}\right) \tilde \psi\left(\frac{\xi_3-\xi_2}{2^{3k+m_2}}\right)\,.
$$
Thus, if we fix $m_2 >0$, it suffices to study trilinear operators having as symbol
\begin{align*}
\sum_{k \in \Z} \varphi\left(\frac{\xi_1-\xi_2}{2^{k}}\right) \tilde \psi\left(\frac{\xi_3-\xi_2}{2^{3k+m_2}}\right),
\end{align*}
with the latter closely resembling the symbols discussed earlier in Part I.

To be more precise, summing in $m_1<0$, we first have
\begin{align*}
\sum_{m_2>0} \sum_{k \in \Z} \int_{\R}  \varphi\big(\frac{\xi_1-\xi_2}{2^{k}}\big)\, \psi\big(\frac{\xi_3-\xi_2}{2^{3k+m_2}}\big)  e^{- i \frac{\xi_1-\xi_2}{2^k}t  - i \frac{\xi_3-\xi_2}{2^{3k}} t^3}  \rho(t) \,dt := \sum_{m_2>0} \sum_{k \in \Z} \tilde{\mathfrak{m}}_{k,0,m_2}(\xi_1,\xi_2,\xi_3). 
\end{align*}
For each $\ds \tilde{\mathfrak{m}}_{k,0,m_2}(\xi_1,\xi_2,\xi_3)$, we perform a double Fourier series decomposition in the variables $\xi_1-\xi_2$ and $\xi_3-\xi_2$ on $[-2^k, 2^k]$ and $\big[ 2^{3k+m_2},  2^{3k+m_2+1}  \big]$, respectively:
\begin{align*}
\tilde{\mathfrak{m}}_{k,0,m_2}(\xi_1,\xi_2,\xi_3)& = \sum_{j_1, j_2 \in \Z} c_{j_1, j_2}^{0, m_2}  \varphi\big(\frac{\xi_1-\xi_2}{2^{k}}\big)\, e^{i \frac{\xi_1-\xi_2}{2^{k}} j_1} \, \psi\big(\frac{\xi_3-\xi_2}{2^{3k+m_2}}\big) e^{i \frac{\xi_3-\xi_2}{2^{3k+m_2}} j_2}
\end{align*}
where
\begin{align*}
c_{j_1,j_2}^{0, m_2} := & \frac{1}{2^{4k+m_2}} \iiint e^{-i\left( \frac{\eta_1}{2^k} t + \frac{\eta_2}{2^{3k}} t^3 \right)} e^{-ij_1 \frac{\eta_1}{2^{k}}} \varphi\left(\frac{\eta_1}{2^{k}}\right)\, e^{-ij_2 \frac{\eta_2}{2^{3k+m_2}}}\psi\left(\frac{\eta_2}{2^{3k+m_2}}\right) \rho(t) \, dt d\eta_1 d\eta_2 \\
= & \iiint  e^{-i (\eta_1 t + \eta_2 2^{m_2} t^3)} e^{-ij_1 \eta_1} \varphi\left(\eta_1 \right)\, e^{-ij_2 \eta_2} \psi\left(\eta_2\right) \rho(t) \, dt d\eta_1 d\eta_2.
\end{align*}
Due to the smoothness of the functions involved and since we are in a situation where the phase $\tilde \lambda(t):= \eta_1 t + \eta_2 2^{m_2} t^3$ has no critical points, one can see that 
\begin{equation}
\label{eq:decay:coeff:0:m_2}
|c_{j_1,j_2}^{0, m_2}|\lesssim_{M_1, M_2} (1+|j_1|+|j_2|)^{-M_1} 2^{-M_2 m_2}
\end{equation}
for any $M_1, M_2 \in \N$.

Once we fix $m_2>0$, we now deal with 
\begin{align*}
 \sum_{j_1, j_2 \in \Z} c_{j_1, j_2}^{0, m_2}  \sum_{k \in \Z}  \varphi\big(\frac{\xi_1-\xi_2}{2^{k}}\big)\, e^{i \frac{\xi_1-\xi_2}{2^{k}} j_1} \, \psi\big(\frac{\xi_3-\xi_2}{2^{3k+m_2}}\big) e^{i \frac{\xi_3-\xi_2}{2^{3k+m_2}} j_2},
\end{align*}
which can be rewritten as
\begin{align*}
 \sum_{j_1, j_2 \in \Z} c_{j_1, j_2}^{0, m_2} (1+|j_1|)^{N_1} (1+|j_2|)^{N_2}  2^{m_2 N_3}    \sum_{k \in \Z} (1+|j_1|)^{-N_1} (1+|j_2|)^{-N_2}  2^{- m_2 N_3}  \varphi\big(\frac{\xi_1-\xi_2}{2^{k}}\big)\, e^{i \frac{\xi_1-\xi_2}{2^{k}} j_1} \, \psi\big(\frac{\xi_3-\xi_2}{2^{3k+m_2}}\big) e^{i \frac{\xi_3-\xi_2}{2^{3k+m_2}} j_2}.
\end{align*}
Finally, one can see that for any $j_1, j_2$ fixed, the symbol
\begin{align*}
m_{j_1, j_2}^{0, m_2}(\xi_1-\xi_2, \xi_3-\xi_2) & :=   \sum_{k \in \Z} (1+|j_1|)^{-N_1} (1+|j_2|)^{-N_2}  2^{- m_2 N_3}  \varphi\Big(\frac{\xi_1-\xi_2}{2^{k}}\Big)\, e^{i \frac{\xi_1-\xi_2}{2^{k}} j_1} \, \psi\Big(\frac{\xi_3-\xi_2}{2^{3k+m_2}}\Big) e^{i \frac{\xi_3-\xi_2}{2^{3k+m_2}} j_2}\\
&=   \sum_{\tilde k \in \Z} (1+|j_1|)^{-N_1} (1+|j_2|)^{-N_2}  2^{- m_2 N_3} \varphi\Big(\frac{2^{m_2 \over 3}(\xi_1-\xi_2)}{2^{\tilde k}}\Big)\, e^{i \frac{2^{m_2 \over 3}(\xi_1-\xi_2)}{2^{\tilde k}} j_1} \, \psi\Big(\frac{\xi_3-\xi_2}{2^{3 \tilde k}}\Big) e^{i \frac{\xi_3-\xi_2}{2^{3 \tilde k}} j_2} \\
& = \sum_{\tilde k \in \Z}  \varphi_{(j_1, m_2)}\Big(\frac{\xi_1-\xi_2}{2^{\tilde k}}\Big)\, \psi_{(j_2)} \Big(\frac{\xi_3-\xi_2}{2^{3 \tilde k}}\Big) 
\end{align*}
with
\[
\varphi_{(j_1, m_2)}(\eta):=  (1+|j_1|)^{-N_1}  2^{- m_2 N_3} \varphi(2^{m_2 \over 3} \eta) e^{i 2^{m_2 \over 3} \, j_1 \eta } \quad \text{and } \quad \psi_{(j_2)}(\eta):=(1+|j_2|)^{-N_2}  \psi(\eta) e^{i j_2 \, \eta}.
\]
Now one observe that, by choosing $N_1,N_2,N_3$ large enough, the functions $\ds \varphi_{(j_1, m_2)}$ and $\ds\psi_{(j_2)}$ are of $\Phi$ and $\Psi$-type respectively, supported around the origin, and all of their first $100$ order derivatives are $O(1)$-uniformly bounded.

Hence the symbol $m_{j_1, j_2}^{0, m_2}$ is exactly of the form \eqref{eq:form:op:initial:gen}. By Theorem \ref{thm:ani}, we deduce that the corresponding trilinear operator is bounded, with a uniform bound with respect to $j_1,j_2,m_2$.
We now conclude our proof by summing over $j_1,j_2,m_2$ which is allowed due to \eqref{eq:decay:coeff:0:m_2}.
\smallskip

\item \underline{\textsf{Case 2.}} For $m_2\leq 0 < m_1$, one has to deal with 
\[
\sum_{m_2<0} \sum_{m_1>0} \sum_{k \in \Z} \mathfrak{m}_{k,m_1,m_2}.
\]

One proceeds as before:  we first sum in $m_2 \leq 0$, and then focus on understanding the symbol
\begin{align*}
\sum_{m_1>0} \sum_{k \in \Z} \Big( \int_{\R} e^{- i \frac{\xi_1-\xi_2}{2^k}t  - i \frac{\xi_3-\xi_2}{2^{3k}} t^3}  \rho(t) \,dt \Big) \psi\big(\frac{\xi_1-\xi_2}{2^{k+m_1}}\big) \varphi \big( \frac{\xi_3-\xi_2}{2^{3k}}\big).
\end{align*}
Then we perform a Fourier series decomposition, and, as before, reduce matters to the boundedness of the trilinear operator having as symbol
\begin{align*}
&\sum_{k \in \Z} (1+|j_1|)^{-N_1} (1+|j_2|)^{-N_2} 2^{-m_1 N_3} \psi\big(\frac{\xi_1-\xi_2}{2^{k+m_1}}\big) e^{i \frac{\xi_1-\xi_2}{2^{k+m_1}} j_1} \varphi\big(\frac{\xi_3-\xi_2}{2^{3k}}\big) e^{i \frac{\xi_3-\xi_2}{2^{3k}} j_2} \\
&:= \sum_{\tilde k \in \Z} (1+|j_1|)^{-N_1} (1+|j_2|)^{-N_2} 2^{-m_1 N_3} \psi\Big(\frac{\xi_1-\xi_2}{2^{\tilde k}}\Big) e^{i \frac{\xi_1-\xi_2}{2^{\tilde k}} j_1} \varphi\Big(2^{3 m_1}\,\frac{\xi_3-\xi_2}{2^{3 \tilde k}}\Big) e^{i 2^{3 m_1}\frac{(\xi_3-\xi_2)}{2^{3 \tilde k}} j_2}.
\end{align*}
Now we have two  possible ways of moving forward:
\begin{itemize}
\item[$\bullet$] We can see the above symbol as 
\[
\sum_{\tilde k \in \Z} \tilde \psi_{(j_1)}\Big(\frac{\xi_1-\xi_2}{2^{\tilde k}}\Big) \, \tilde\varphi_{(j_2, m_1)}\Big(\frac{\xi_3-\xi_2}{2^{3 \tilde k}}\Big)
\]
where $\ds  \tilde \psi_{(j_1)}(\eta):=(1+|j_1|)^{-N_1} \psi(\eta) \, e^{i \eta j_1}$ and $\ds \tilde\varphi_{(j_2, m_1)}(\eta):=(1+|j_2|)^{-N_2} 2^{-m_1 N_3}  \varphi(2^{3 m_1} \eta) \, e^{i 2^{3m_1} \eta j_2}$. Then we can invoke Theorem \ref{thm:ani}, since the above expression has the form in \eqref{eq:form:op:initial:gen} - the functions $ \tilde \psi_{(j_1)}$ and $\tilde\varphi_{(j_2, m_1)}$ are of $\Psi$ and  $\Phi$-type respectively, and their higher order derivatives are uniformly bounded provided $N_1,N_2,N_3$ are chosen large enough.

\item[$\bullet$] Alternatively, we can regard the above symbol as
\begin{equation}
\label{eq:anis:Mikhln:app}
\sum_{\tilde k \in \Z} (1+|j_1|)^{-N_1} (1+|j_2|)^{-N_2} 2^{-m_1 N_3} \psi_{(j_1)}\Big(\frac{(\xi_1-\xi_2)^3}{2^{3 \tilde k}}\Big)  \varphi_{(j_2, m_1)}\Big(\frac{\xi_3-\xi_2}{2^{3 \tilde k}}\Big),
\end{equation}
where
\[
\psi_{(j_1)}(\eta):= \psi(\eta^{1 \over 3}) e^{i \eta^{1 \over 3} j_1}, \quad \text{and} \quad \varphi_{(j_2, m_1)}(\eta):= \varphi(2^{3 m_1} \eta) e^{i 2^{3 m_1} \eta}.
\]
Now one can see that the expression in \eqref{eq:anis:Mikhln:app} satisfies condition \eqref{eq:Mikhlin} for which our Theorem \ref{thm:anisotropic:Mikhlin} applies, provided that $N_1,N_2,N_3$ are chosen large enough.
\end{itemize}
Obviously, the two approaches are morally equivalent since Theorems \ref{thm:ani} and \ref{thm:anisotropic:Mikhlin} are proved similarly.

\smallskip

\item \underline{\textsf{Case 3.}} For $m_1,m_2>0$ with $|m_1-m_2|\geq 3$, we now deal with 
\[
\sum_{\substack{m_1, m_2> 0   \\ |m_1-m_2|\geq 3}} \sum_{k \in \Z} \mathfrak{m}_{k,m_1,m_2}(\xi_1,\xi_2,\xi_3),
\]
with $\mathfrak{m}_{k,m_1,m_2}(\xi_1,\xi_2,\xi_3)$ given by \eqref{eq:formula:symbol}. Once at this point, we use a (windowed) Fourier series decomposition for each of these symbols:
\begin{align*}
\mathfrak{m}_{k,m_1,m_2}(\xi_1,\xi_2,\xi_3) & = \sum_{j_1,j_2 \in \Z} c_{j_1,j_2}^{m_1, m_2} e^{ij_1 \frac{\xi_1-\xi_2}{2^{k+m_1}}} \tilde \psi\left(\frac{\xi_1-\xi_2}{2^{k+m_1}}\right)\, e^{ij_2 \frac{\xi_3-\xi_2}{2^{3k+m_2}}}\tilde \psi\left(\frac{\xi_3-\xi_2}{2^{3k+m_2}}\right),
\end{align*}
with the Fourier coefficients equal to
\begin{align*}
c_{j_1,j_2}^{m_1, m_2}:&= \frac{1}{2^{4k+m_1+m_2}} \iiint e^{i\lambda(t)} e^{-ij_1 \frac{\eta_1}{2^{k+m_1}}} \psi\left(\frac{\eta_1}{2^{k+m_1}}\right)\, e^{-ij_2 \frac{\eta_2}{2^{3k+m_2}}}\psi\left(\frac{\eta_2}{2^{3k+m_2}}\right) \rho(t) \, dt d\eta_1 d\eta_2 \\
&= \int_{\R^2} \Big( \int_{\R} e^{-i \eta_1 2^{m_1} t - i \eta_2 2^{m_2} t^3} \rho(t) dt \Big)  \psi(\eta_1) \psi(\eta_2) e^{-i \eta_1 j_1} e^{-i \eta_2 j_2} d \eta_1 d \eta_2.
\end{align*}
The non-stationary phase principle (iterated integration by parts in $t$) yields that for any $N_1, N_2\in\N$
$$ |c_{j_1,j_2}^{m_1, m_2}|\lesssim_{N_1, N_2} (1+|j_1|+|j_2|)^{-N_1} 2^{-N_2 \max(m_1,m_2)}.$$
Since the Fourier coefficients carry enough decay in $|j_1|, |j_2|, m_1$ and $m_2$, we may consider these parameters fixed and focus on  
\begin{align*}
\sum_{k \in \Z} (1+|j_1|)^{-N_1} (1+|j_2|)^{-N_2} e^{ij_1 \frac{\xi_1-\xi_2}{2^{k+m_1}}} \tilde \psi\left(\frac{\xi_1-\xi_2}{2^{k+m_1}}\right)\, e^{ij_2 \frac{\xi_3-\xi_2}{2^{3k+m_2}}}\tilde \psi\left(\frac{\xi_3-\xi_2}{2^{3k+m_2}}\right)
\end{align*}
which can be controlled by invoking once again Proposition \ref{prop:facile} for the $\Psi$-type functions $\ds (1+|j_1|)^{-N_1} e^{i j_1 \eta} \tilde \psi(\eta)$ and $\ds (1+|j_2|)^{-N_2} e^{i j_2 \eta} \tilde \psi(\eta)$.
\end{itemize}
This concludes the proof of Theorem \ref{thm:non-sta}.
\end{proof}

\newpage

\part[II]{TREATMENT OF THE HIGH FREQUENCY STATIONARY COMPONENT \bf{$T^{\sta}$}}
\label{partII}
This second  part of our paper is devoted in its entirety to the proof of Theorem \ref{thm:sta}.

\section{Low and high resolution analysis of $T^{\sta}$: foundations}
\label{sec:discretization}
%

\subsection{Preliminaries and the basic strategy}
\label{subsec:prell}

We start by recalling the basic decomposition performed in Section \ref{preldiscret}, see in particular \eqref{op-de}---\eqref{dectht}. We thus deduce that our main object of study obeys
\begin{equation}\label{thid}
T^{\sta}\approx\ds \sum_{m \geq 0} \sum_{k \in\Z} T_{k,m}\,,
\end{equation}
where, applying a dualization procedure and writing
\begin{equation}\label{dualt}
\L_{k,m}(f_1,f_2,f_3,f_4):=\int_{\R} T_{k,m}(f_1,f_2,f_3)\,f_4, 
\end{equation}
we have 
\begin{equation}\label{fm4}
\L_{k,m}(f_1,f_2,f_3,f_4):=\iiint_{\R^3} \hat{f}_1(\xi_1)\,\hat{f}_2(\xi_2)\,\hat{f}_3(\xi_3)\,\hat{f}_4(-\xi_1-\xi_2-\xi_3)\, \mathfrak{m}_{k,m}(\xi_1,\xi_2,\xi_3)\,d\xi_1\,d\xi_2\,d\xi_3\,,
\end{equation}
with 
\beq\label{multkm}
\mathfrak{m}_{k,m}(\xi_1,\xi_2,\xi_3):= \left(\int_{\R}e^{-i\,\frac{\xi_1-\xi_2}{2^k}\,t}\,e^{-i\,\frac{\xi_3-\xi_2}{2^{3k}}\,t^3}\,\rho(t)\,dt\right)\,
\psi\left(\frac{\xi_1-\xi_2}{2^{k+m}}\right)\,\psi\left(\frac{\xi_3-\xi_2}{2^{3k+m}}\right),
\eeq
and $\psi\in C_{0}^{\infty}(\R)$ of $\Psi-$type (compactly supported away from the origin).

As noted in the introduction, the specific choice of the shape of HTHT in \eqref{hththyb} allows for an immediate symmetry manifestation between the cases $k\geq 0$ and $k\leq 0$---easy to notice\footnote{See also Section \ref{Seckneg}.} by changing the role of $f_1,f_3$.

Consequently, in what follows, we will only focus on the case $k\geq 0$, and hence, on bounding $\ds \sum_{m \geq 0} \sum_{k \geq 0} T_{k,m}$.

Our analysis will be split into three distinct regimes:
\smallskip
\begin{enumerate}[label=(\alph*)  , itemsep=10pt]
\item \label{case:c1}  $\sum_{m \geq 0} \sum_{k > {m \over 4}} T_{k, m} = \sum_{m \geq 0} \big(\sum_{k> {m \over 4}} T_{k,m}\big) :=\sum_{m \geq 0}T_m^{1}$;
\item \label{case:c2}  $\sum_{m \geq 0} \sum_{{m \over 4}\geq k \geq {m \over 8}} T_{k,m} = \sum_{\kappa \geq 0} \big(\sum\limits_{m= 4\kappa}^{8\kappa} T_{\frac{m}{4}-\frac{\kappa}{2},m}\big):=\sum_{\kappa \geq 0} T_\kappa^{2}$;
\item \label{case:c3}  $\sum_{m \geq 0} \sum_{{m \over 8}> k \geq {0}} T_{k, m} = \sum_{k \geq 0} \big(\sum_{m\geq 8k} T_{k,m} \big) :=\sum_{k \geq 0}T_k^{3}$.
\end{enumerate} 
\smallskip
With these said, we now present the following

\noindent

\smallskip
{\bf \underline{Strategy:}}
\smallskip

\noindent\underline{\textsf{Context:}} In the three cases above we deal with a sum $\sum_{n\geq 0} T_n$ (here $n$ corresponds to $m,\kappa$ or $k$) where $T_n$ itself is a multi-scale operator $T_n=\sum_{n'} T_{n,n'}$ (with $n'$ playing the role of $k$ or $m$). The main aim is then to prove the boundedness of $T_n$ with a decay of the form $2^{-\delta n}$ for some suitable $\delta>0$, thus allowing us to sum over $n\geq 0$. To accomplish this, for a  fixed $n \geq 0$, we proceed as follows:
\begin{itemize}[itemsep=6pt, leftmargin=2em]
 \item \underline{\textsf{Step 1: An $L^2$ based single-scale $n-$exponential decay estimate.}}\quad At this first step, we will show that there exists $\delta>0$, such that for every $n,n'$ and every spatial interval $I$ whose length depends on $n$ and $n'$, one has
 \begin{equation}
 |\langle T_{n,n'}(f_1,f_2,f_3), \chi_I f_4\rangle| \lesssim 2^{-\delta \, n} |I|^{-\frac{1}{2}} \prod_{j=1}^4 \|f_j\chi_I\|_{L^2}, \label{eq:*ff}
\end{equation}
after removing some suitable {\it sparse}\footnote{The precise meaning of this concept will be made explicit when discussing the proof of \eqref{eq:*ff}.} components from some of the functions $f_1,\cdots, f_4$. In order to prove this estimate we will employ the Rank II LGC method which further relies on a \emph{high-resolution single scale analysis}. This step will be completed in Sections \ref{sec:onescale_k_large}, \ref{Diagkpsmall124} and \ref{Diagkinterm}. 

\item \underline{\textsf{Step 2: A local-$L^2$ multi-scale $n-$exponential decay estimate.}} \quad At the second step, our aim is to propagate the $n-$exponential decay at \emph{one} scale obtained at Step 1 to \emph{all} scales. Indeed, by applying BHT and Biest-like time-frequency analysis, we will be able to sum over all the scales $n'$ and prove that for every exponents $p_1, p_2, p_3, p' \in(2,\infty]$ with $\frac{1}{p_1}+\frac{1}{p_2}+\frac{1}{p_3}+\frac{1}{p'}=1$ there exists $\delta>0$ such that 
\begin{equation} 
\|T_n(f_1, f_2, f_3)\|_p \lesssim 2^{- \delta \, n} \prod_{j=1}^3 \| f_j  \|_{p_j}.
\label{eq:decayk}
\end{equation}
This will require a \emph{low resolution multi-scale analysis} which will be presented in Sections \ref{sec:bilinear:analysis}, \ref{sec:bilinear:analysis2} and \ref{sec:bilinear:analysis3}.

\item \underline{\textsf{Step 3: A wider $L^p$ based multi-scale $n-$polynomial growth estimate.}} At this last, third step, we extend the boundedness range in exchange of turning the exponential decay into a tame polynomial growth in the parameter $n$. To be precise, in Sections \ref{extrange1} and \ref{extrange2}, we will prove that there exists $\alpha \in \N$ absolute constant such that, whenever $p_1, p_2,p_3\in (1,\infty)$ with $\frac{1}{p_1}+\frac{1}{p_2}+\frac{1}{p_3}=\frac{1}{p}<1$, we have
 \begin{equation}
\label{eq:est:pol:growth:k} 
 \|T_n(f_1, f_2, f_3)\|_p \lesssim (1+|n|)^\alpha \prod_{s=1}^3 \| f_s  \|_{p_s}\,.
 \end{equation}
\end{itemize}
Putting now together the three steps above and interpolating between \eqref{eq:decayk} and \eqref{eq:est:pol:growth:k}, we obtain the full boundendness range of $T^{\sta}$ stated in Theorem \ref{thm:sta}.

\begin{remark} Notice that the restriction in the boundedness range for the full operator $\THT$ in our main Theorem \ref{thm:main} is the byproduct of the comparatively weaker control on the low frequency component $T^{\low}$--- see Theorem \ref{thm:low}---which, in turn, is a consequence of the anisotropic Biest operator bounds in Theorem \ref{thm:main-ani:biest}.
\end{remark}

As a result of the above discussion, the main task for the remainder of the paper is to prove the following 

\begin{theorem} \label{thm:case1}
Fix exponents $1< p_1, p_2,p_3, p< \infty$ with $\frac{1}{p}:=\frac{1}{p_1}+\frac{1}{p_2}+\frac{1}{p_3}<1$. Then, there exists $\delta:=\delta(p_i,p)>0$ such that 
\begin{itemize}
\item for every $m \geq 0$ one has
$$ \| \sum_{k\geq \frac{m}{4}} T_{k,m}\|_{L^{p_1} \times L^{p_2} \times L^{p_3} \rightarrow L^p} \lesssim 2^{-\delta m}\,;$$
\item for every $\kappa\geq 0$ one has
$$ \| \sum_{m\geq 4 \kappa} T_{\frac{m}{4}-\frac{\kappa}{2},m}\|_{L^{p_1} \times L^{p_2} \times L^{p_3} \rightarrow L^p} \lesssim 2^{-\delta \kappa}\,;$$
\item for every $k\geq 0$ one has
$$ \| \sum_{m\geq 8k} T_{k,m} \|_{L^{p_1} \times L^{p_2} \times L^{p_3} \rightarrow L^p} \lesssim 2^{-\delta k}.$$
\end{itemize}
\end{theorem}

We end this section with the following

\begin{remark} A careful inspection of \eqref{dualt}--\eqref{multkm}, reveals the following:
\begin{itemize}
\item On the one hand, for each fixed pair $(k,m)$, the structure of the phase in \eqref{ph} combines a highly oscillatory behavior with the presence of the curvature in the $t-$parameter thus producing a cancelation phenomenon together with an entanglement of the variables $\xi_1-\xi_2$ and $\xi_3-\xi_2$. Consistent with Step 1 above, in order to be able to extract this cancelation, we employ the Rank II LGC method via a high resolution single scale analysis whose preliminary foundation is laid out in Section \ref{subsec:63}.
\smallskip
     
\item On the other hand, the same structure of \eqref{multkm} reveals the (double) translation invariance property of the multiplier relative to the $\xi_2-\xi_1$ or $\xi_3-\xi_2$ argument. This phenomenon propagates among different scales thus naturally bringing into the picture the recognizable patterns of the bilinear Hilbert transform as well as the \emph{biest} operator. Consistent with Step 2 above, in order to exploit this latter behavior we will need to develop a low resolution multiscale analysis  whose preliminary foundation is laid out in Section \ref{lres}.
\end{itemize}
\end{remark}

\subsection{Low-resolution discretization: macro-scale decomposition--$(k,m)$-discretization}\label{lres}

In what follows our analysis will focus on the form $\L_{+}$, as defined in \eqref{fmd}. Consequently, for a fixed curvature parameter $m\in\N$, and a fixed scale parameter $k\in\N$, we split the frequency band $|\xi_3-\xi_2|\sim 2^{3k+m}$ into frequency squares $(\xi_2,\xi_3)$ of type $\ds \big[ 2^{3k+m} \el,   2^{3k+m} (\el +1)  \big]  \times \big[ 2^{3k+m} (\el +2),   2^{3k+m} (\el+3)  \big]$ with $\el\in\Z$ while the frequency band $|\xi_1-\xi_2|\sim 2^{k+m}$ is split into frequency squares $(\xi_1,\xi_2)$ of type $\ds \big[ 2^{3k+m} \el + 2^{k+m} (\mathfrak{r}+2) ,   2^{3k+m} \el + 2^{k+m} (\r +3)  \big]  \times \big[2^{3k+m} \el + 2^{k+m} \r,  2^{3k+m} \el +  2^{k+m} (\r+1)  \big]$ where $0 \leq \r \leq 2^{2k}$.\footnote{Notice that, according to the notation introduced in the beginning of the paper, we can still write $\r \sim 2^{2k}$.} With these, it is immediate to see that $\mathfrak{m}_{k,m}(\xi_1,\xi_2,\xi_3)$ becomes a superposition of terms of the form 
\beq\label{multkms1}
\sum_{\el\in\Z} \sum_{\r\sim 2^{2k}}\mathfrak{m}^{\el,\r}_{k,m}(\xi_1,\xi_2,\xi_3),
\eeq
where
\begin{align*}
\mathfrak{m}^{\el,\r}_{k,m}(\xi_1,\xi_2,\xi_3) & := \left(\int_{\R}e^{-i\,\frac{\xi_1-\xi_2}{2^k}\,t}\,e^{i\,\frac{\xi_2-\xi_3}{2^{3k}}\,t^3}\,\rho(t)\,dt\right) \nonumber \\
& \cdot \phi_1\left(\frac{\xi_1}{2^{k+m}}-\el 2^{2k}-\r\right)\,\phi_2\left(\frac{\xi_2}{2^{k+m}}-\el 2^{2k}-\r\right)\,\phi_3\left(\frac{\xi_3}{2^{3k+m}}-\el\right) 
\end{align*}
for three suitable $\phi_1,\,\phi_2, \phi_3\in C_{0}^{\infty}(\R)$.

If we set, for $\ell, \kappa \in \Z$, 
\begin{equation}
\label{eq:notation:inv:greq}
W_{\ell}^{\kappa}:=\big[ \ell \, 2^\kappa,  (\ell+1) \, 2^\kappa  \big],
\end{equation}
we can see that $\mathfrak{m}^{\el,\r}_{k,m}$ is supported in frequency onto $\ds W_{(\mathfrak{r}+2)+2^{2k} \mathfrak l}^{k+m} \times  \big(W_{\mathfrak{r}+2^{2k} \mathfrak l}^{k+m} \cap  W_{\mathfrak l}^{3k+m} \big) \times W_{\mathfrak l+2}^{3k+m}$. As we consider the 4-linear form $\Lambda_{k,m}$ described in \eqref{fm4}, we need to localize all the functions to these frequency regions: for $s\in\{1,2\}$ we define
\begin{equation}\label{functsup}
\widehat{f_{s,\el,\r}}(\xi_s):= \hat{f}_s(\xi_s)\,\phi_s\left(\frac{\xi_s}{2^{k+m}}-\el2^{2k}-\r\right)\,,
\end{equation}
while for $f_3$ and $f_4$ we set\footnote{Again, here $\phi_4\in C_{0}^{\infty}(\R)$ is a suitable smooth function.}
\begin{equation}\label{functsup4}
\widehat{f_{3,\el}}(\xi_3):= \hat{f}_3(\xi_3)\,\phi_3\left(\frac{\xi_3}{2^{3k+m}}-\el \right)\qquad\textrm{and}\qquad \widehat{f_{4,\el}}(\xi_4):= \hat{f}_4(\xi_4)\,\phi_4\left(\frac{\xi_4}{2^{3k+m}}+3\el\right)\,.
\end{equation}

With these done, we have that
\begin{equation}\label{smd}
\L_{k,m}(f_1,f_2,f_3,f_4)\approx\sum_{\el\in\Z} \sum_{\r\sim 2^{2k}} \L_{k,m}^{\el,\r}(f_1,f_2,f_3,f_4)
\end{equation}
where
\begin{equation}\label{smd1}
\L_{k,m}^{\el,\r}(f_1,f_2,f_3,f_4):=\L_{k,m}(f_{1,\el,\r},f_{2,\el,\r},f_{3,\el},f_{4,\el})\,.
\end{equation}
Finally we denote by $I^k:=I_0^{0,k}=[2^{-k}, 2^{-k+1}]$ where, for $j\in\Z$, we set
\begin{equation}\label{intnj}
I_j^k:=I_j^{0,k}=\Big[\frac{j+1}{2^{k}},\,\frac{j+2}{2^{k}}\Big]\,,
\end{equation}
and define, for $\chi_{I_j^{k}}$ a suitable smooth function decaying fast away from $I_j^{k}$ and with compact frequency support,
\begin{equation}\label{smd2}
\L_{k,m}^{\el,\r,j}(f_1,f_2,f_3,f_4):=\L_{k,m}(f_{1,\el,\r},f_{2,\el,r},f_{3,\el},\chi_{I_j^{k}}\,f_{4,\el})\,.
\end{equation}
With this, the space variable $x$ is restricted to the interval $I_j^{k}$, so that we have localization in space as well.

Then, for $k\geq 0$ and $m\in \N$, $\L_{k,m}(f_1,f_2,f_3,f_4)$ becomes a superposition of terms of the form
\begin{equation}\label{Firststage}
\sum_{j \in \Z}\sum_{\el\in\Z} \sum_{\r\sim 2^{2k}} \L_{k,m}^{\el,\r, j}(f_1,f_2,f_3,f_4)
\end{equation}
where
\begin{equation}\label{fm11:2}
\L_{k,m}^{\el,\r,j}(f_1,f_2,f_3,f_4)=\int_{\R}\int_{\R} f_{1,\el,\r}(x-t)\,f_{2,\el,\r}(x+t+t^3)\,f_{3,\el}(x-t^3)\,f_{4,\el}(x)\, 2^{k}\,\rho(2^{k} t)\,\varphi(2^{k} x-j)\,dt\,dx\,.
\end{equation}

This completes the first stage of our low-resolution analysis which will be completed in Sections \ref{sec:bilinear:analysis}, \ref{sec:bilinear:analysis2} and \ref{sec:bilinear:analysis3}.

\subsection{High resolution discretization: micro-scale decomposition and time-frequency correlations via LGC methodology}
\label{subsec:63}

In this section we discuss a preliminary discretization of the single scale form $\L_{k,m}^{\el,\r,j}$ defined by \eqref{fm11:2}. Along the way, we will highlight certain time-frequency correlations emerging from the curvature features of $\L_{k,m}^{\el,\r,j}$. These time-frequency correlations are obtained via the Rank I LGC method, by decomposing the input functions into wave-packets adapted to Heisenberg tiles partitioning the so-called super Heisenberg tiles representing the time-frequency regions where the input functions are localized. Although we will not directly work with the model quadrilinear form obtained at the end of this section (see Proposition \ref{prop:GFF}), the plan presented here constitutes the foundation for the Rank II LGC method implemented in Sections \ref{sec:onescale_k_large}-\ref{Seckneg}.

\begin{obs}\label{Reduct}
i) [\textsf{Preview of the main goal}] As stated at Step 1 above,  our main focus for Sections \ref{sec:onescale_k_large}-\ref{Seckneg} is to obtain, (in the spatially uniform case which will be defined later),  an $L^2$ based single-scale exponential decay estimate of the form\footnote{By abusing notation, we sometimes disregard the bump functions in $x$ and $t$ and simply write
\begin{equation}\label{fm11}
\L_{k,m}^{\el,\r,j}(f_1,f_2,f_3,f_4)=|I_j^k|^{-1} \int_{I_j^k}\int_{I^k} f_{1,\el,\r}(x-t)\,f_{2,\el,\r}(x+t+t^3)\,f_{3,\el}(x-t^3)\,f_{4,\el}(x)\,\,dt\,dx\,.
\end{equation}}
\begin{equation}\label{F1}
|\L_{k,m}^{\el,\r,j}(f_1,f_2,f_3,f_4)|\lesssim 2^{-\ep \min\{2k,m\}}\,2^{k}\, \prod_{s=1}^2 \|f_{s, \mathfrak l, \mathfrak r} \cdot \chi_{I_j^k}\|_{L^2}\, \prod_{s=3}^4 \|f_{s, \mathfrak l} \cdot \chi_{I_j^k}\|_{L^2},
\end{equation}
for all $ k, m \in\N, \, \el,j\in\Z\,, 0 \leq \r \leq 2^{2k}$ and for some proper choice of $\ep>0$ (absolute constant). This estimate involving the space interval $I_j^k$ and the frequency region  $\ds W_{(\mathfrak{r}+2)+2^{2k} \mathfrak l}^{k+m} \times W_{\mathfrak{r}+2^{2k} \mathfrak l}^{k+m} \times W_{\mathfrak l+2}^{3k+m} \times W_{-(3 \mathfrak l+3)}^{3k+m}$,  will be used later on in the multi-scale analysis. 

\noindent ii) [\textsf{Reduction of $\L_{k,m}^{\el,\r,j}$ to $\L_{k,m}^{0,0,0}$}] Due to the translation and modulation invariance of the form $\L_{k,m}^{\el,\r,j}$, the proof of \eqref{F1} is independent of the parameters $\el,\r,j$, which is why in the subsequent sections concerning the single scale behavior we will only focus on $\L_{k,m}^{0,0,0}$.
\end{obs}

\medskip
\noindent\underline{\noindent\textbf{Stage 1.} \textsf{Multiplier analysis - phase linearization}}
\medskip

We start by restating for reader's convenience the multiplier expression in \eqref{multkm}:

\beq\label{multkmm1}
\mathfrak{m}_{k,m}(\xi_1,\xi_2,\xi_3)= \left(\int_{\R}e^{-i\,\frac{\xi_1-\xi_2}{2^k}\,t}\,e^{-i\,\frac{\xi_3-\xi_2}{2^{3k}}\,t^3}\,\rho(t)\,dt\right)\,
\psi\big(\frac{\xi_1-\xi_2}{2^{k+m}}\big)\,\psi\big(\frac{\xi_3-\xi_2}{2^{3k+m}}\big).
\eeq
Applying now the stationary phase analysis we have that the main term is given by\footnote{Here we assume we have stationary points.}

\beq\label{multkmmt}
\mathfrak{m}_{k,m}(\xi_1,\xi_2,\xi_3)\approx \frac{1}{2^{\frac{m}{2}}}\,e^{ic \frac{(\xi_1-\xi_2)^{\frac{3}{2}}}{(\xi_3-\xi_2)^{\frac{1}{2}}}}\,
\psi\big(\frac{\xi_1-\xi_2}{2^{k+m}}\big)\,\psi\big(\frac{\xi_3-\xi_2}{2^{3k+m}}\big)\,,
\eeq
with $c$ being a numerical constant.
Setting now
$$\Phi(\xi_1,\xi_2,\xi_3)=\frac{(\xi_1-\xi_2)^{\frac{3}{2}}}{(\xi_3-\xi_2)^{\frac{1}{2}}}$$
we have:
\beq\label{deriv1}
\left|\partial^2_{\xi_1}\Phi\right|,\:\left|\partial^2_{\xi_2}\Phi\right|, \:\left|\partial_{\xi_1}\partial_{\xi_2}\Phi\right|\lesssim \frac{1}{2^{2k+m}}\,,
\eeq
\beq\label{deriv2}
\left|\partial_{\xi_1}\partial_{\xi_3}\Phi \right|,\:\left|\partial_{\xi_1}\partial_{\xi_3} \Phi \right|\lesssim \frac{1}{2^{4k+m}}\,,
\eeq
and
\beq\label{deriv3}
\left|\partial^2_{\xi_3}\Phi\right|\lesssim \frac{1}{2^{6k+m}}\,.
\eeq
Thus, via a Taylor series argument, in order to have \emph{linear behavior of the phase} we need to impose the following parameter variation:

\beq\label{var}
|\triangle \xi_1|,\:|\triangle \xi_2|\sim 2^{\frac{m}{2}+k}\qquad\textrm{and}\qquad |\triangle \xi_3|\sim 2^{\frac{m}{2}+3k} \,.
\eeq

\medskip
\noindent \underline{\noindent\textbf{Stage 2.} \textsf{Gabor frame/wave packet discretization of the input functions}}
\medskip

Before proceeding, we introduce some notation that will be used in the context of wave packets. The frequency intervals will be denoted
\begin{equation}
\label{eq:intv:freq:small}
\omega_n^{m, k}:=\big[ n \, 2^{k+\frac{m}{2}}, (n+1) \, 2^{k+\frac{m}{2}} \big], \qquad \text{with  } n \in \Z,
\end{equation}
and the space intervals will be of the form $\ds I_{u}^{m,k} := \big[ u \, 2^{-(k+\frac{m}{2})} , (u+1) \, 2^{-(k+\frac{m}{2})} \big]$ for $u \in \Z$.
Let $\v\in C_0^{\infty}(\R)$ and set
\beq\label{gdec03}
\v_{u,n}^{k}(x):=2^{\frac{1}{2}(k+\frac{m}{2})}\,\v(2^{k+\frac{m}{2}}x-u)\,e^{i\,2^{k+\frac{m}{2}}\,n\,x}\,.
 \eeq
Notice that the wave packet is associated to the time-frequency region $\ds I_u^{m, k} \times \omega_n^{m, k}$.

We are now implementing the frequency support restrictions imposed by \eqref{var} via a Gabor frame decomposition of the input functions $f_1$, $f_2$ and $f_3$ under the assumptions $\el=\r=j=0$. Thus $f_1$ is decomposed into wave packets supported in frequency onto intervals of length $2^{\frac{m}{2}+k}$ included in $[2 \cdot 2^{k+m}, 3 \cdot 2^{k+m}]$, adapted in space onto intervals of length $2^{-(\frac{m}{2}+k)}$ :
\beq\label{gdecf1*}
f_1(x)=\sum_{u \in \Z}\sum_{n_1\sim 2^{\frac{m}{2}}} \langle f_1,\,\v_{u,n_1}^{k} \rangle \,\v_{u,n_1}^{k}(x)\:.
\eeq

The functions $f_2$ and $f_3$ are decomposed in a similar way, whilst we don't implement a wave packet decomposition for $f_4$ just yet. When we insert this decomposition into \eqref{fm11:2}, we realize that it suffices to consider the wave packets associated to spatial intervals contained in $I^k=I_0^k$, the general term being a superposition of similar expressions. 
Thus the functions $f_1, f_2, f_3$ in the expression $\Lambda_{k,m}^{0, 0, 0}$ will be replaced by the wave packet decompositions corresponding to the space interval $I^k$: 
\beq\label{gdecf1}
f_1(x)= \sum_{u\sim 2^{\frac{m}{2}}} \sum_{n_1\sim 2^{\frac{m}{2}}} \langle f_1,\,\v_{u,n_1}^{k} \rangle \,\v_{u,n_1}^{k}(x)\:,
\eeq
\beq\label{gdecf2}
f_2(x)= \sum_{v\sim 2^{\frac{m}{2}}} \sum_{n_2\sim 2^{\frac{m}{2}}} \langle f_2,\,\v_{v,n_2}^{k} \rangle \,\v_{v,n_2}^{k}(x)\:,
\eeq
and
\beq\label{gdecf3}
f_3(x)= \sum_{r\sim 2^{2k+\frac{m}{2}}} \sum_{n_3\sim 2^{\frac{m}{2}}} \langle f_3,\,\v_{r,n_3}^{3k} \rangle \,\v_{r,n_3}^{3k}(x)\:.
\eeq

With these $\L_{k,m}^{0,0,0}(f_1,f_2,f_3,f_4)$ becomes a superposition of terms of the form
\begin{align}\label{GF}
&\sum_{{u\sim 2^{\frac{m}{2}}}\atop{n_1\sim 2^{\frac{m}{2}}}}  \sum_{{v\sim 2^{\frac{m}{2}}}\atop{n_2\sim 2^{\frac{m}{2}}}}\sum_{{r\sim 2^{2k+\frac{m}{2}}}\atop{n_3\sim 2^{\frac{m}{2}}}} \langle f_1,\,\v_{u,n_1}^{k} \rangle \, \langle f_2,\,\v_{v,n_2}^{k} \rangle \, \langle f_3,\,\v_{r,n_3}^{3k}\rangle \,\I_{k,m}^{u,v,r,\vec{n}}[f_4]\,,
\end{align}
where
\begin{align}\label{GF1}
\I_{k,m}^{u,v,r,\vec{n}}[f_4]:&=\iint_{\R^2}\v_{u,n_1}^{k}(x-t)\,\v_{v,n_2}^{k}(x+t+t^3)\,\v_{r,n_3}^{3k}(x-t^3)\,
f_4(x)\,2^{k}\,\rho(2^{k} t)\,dt\,dx\\
&=\iint_{\R^2}2^{\frac{3m}{4}}\,2^{\frac{7k}{2}}\,\v(2^{\frac{m}{2}+k}(x-t)-u)\, \v(2^{\frac{m}{2}+k}(x+t+t^3)-v)\,\v(2^{\frac{m}{2}+3k}(x-t^3)-r)\,\nonumber\\
&\cdot\,e^{i 2^{\frac{m}{2}+k}\,(x-t)\,n_1}\,e^{i 2^{\frac{m}{2}+k}\,(x+t+t^3)\,n_2}\,e^{i 2^{\frac{m}{2}+3k}\,(x-t^3)\,n_3}\,f_4(x)\,dt\,dx\,.\nonumber
\end{align}
Deduce now that the $(x,t)$ variables in \eqref{GF1} obey
\begin{equation}
\begin{cases}
x-t=\frac{u}{2^{\frac{m}{2}+k}}\,+\,O(\frac{1}{2^{\frac{m}{2}+k}})\\\\
x+t+t^3=\frac{v}{2^{\frac{m}{2}+k}}\,+\,O(\frac{1}{2^{\frac{m}{2}+k}})\\\\
x-t^3=\frac{r}{2^{\frac{m}{2}+3k}}\,+\,O(\frac{1}{2^{\frac{m}{2}+3k}})
\end{cases}\,,
\end{equation}
and hence
\begin{equation}
\begin{cases}
x=\frac{1}{3}\cdot\frac{u+v}{2^{\frac{m}{2}+k}}\,+\,\frac{1}{3}\cdot\frac{r}{2^{\frac{m}{2}+3k}}\,+\,O(\frac{1}{2^{\frac{m}{2}+k}})\\\\
t=\frac{1}{3}\cdot\frac{-2u+v}{2^{\frac{m}{2}+k}}\,+\,\frac{1}{3}\cdot\frac{r}{2^{\frac{m}{2}+3k}}\,+\,O(\frac{1}{2^{\frac{m}{2}+k}})
    \end{cases}\,.\nonumber
\end{equation}
As a consequence of this discretization we remark that
\begin{align}\label{cube}
2^{\frac{m}{2}+3k}\,t^3&=-\frac{2^{\frac{m}{2}+1}}{27}\,\left(\frac{v-2u}{2^{\frac{m}{2}}}+\frac{r}{2^{\frac{m}{2}+2k}}\right)^3\,+\,\frac{t}{3}\,
2^{\frac{m}{2}+k}\,\left(\frac{v-2u}{2^{\frac{m}{2}}}+\frac{r}{2^{\frac{m}{2}+2k}}\right)^2\,+\,O(2^{-\frac{m}{2}})\\
&=\frac{2^{\frac{m}{2}}}{27}\,\left(\frac{v-2u}{2^{\frac{m}{2}}}+\frac{r}{2^{\frac{m}{2}+2k}}\right)^3\,+\,O(1)\,.
\end{align}
Inserting the above in \eqref{GF1} we get
\begin{align} 
\I_{k,m}^{u,v,r,\vec{n}}[f_4] & \approx 2^{\frac{3m}{4}}\,2^{\frac{7k}{2}}\,e^{i\,(n_3-\frac{n_2}{2^{2k}})\,\frac{2^{\frac{m}{2}+1}}{27}\,\left(\frac{v-2u}{2^{\frac{m}{2}}}+\frac{r}{2^{\frac{m}{2}+2k}}\right)^3}\, \nonumber  \\
& \quad \cdot\Big[\int_{\R} \v\Big(2^{\frac{m}{2}+k}t-\frac{v-2u}{3}-\frac{r}{3\cdot2^{2k}}\Big)\,e^{i 2^{\frac{m}{2}+k}t\,\left(-n_1+n_2+ (\frac{n_2}{2^{2k}}-n_3)\,\frac{1}{3}\,\left(\frac{v-2u}{2^{\frac{m}{2}}}+\frac{r}{2^{\frac{m}{2}+2k}}\right)^2\right)}\,dt\Big]\nonumber\\
& \quad \cdot\Big[\int_{\R}\v\Big(2^{\frac{m}{2}+k}x-\frac{v+u}{3}-\frac{r}{3\cdot2^{2k}}\Big)\, \v\Big(2^{\frac{m}{2}+3k}x-r-\frac{2^{\frac{m}{2}}}{27}\,\big(\frac{v-2u}{2^{\frac{m}{2}}}+\frac{r}{2^{\frac{m}{2}+2k}}\big)^3\Big)\, \cdots \nonumber \\
& \qquad \cdots \cdot f_4(x)\,e^{i 2^{\frac{m}{2}+k}\,x\,\left(n_1+n_2+2^{2k}\,n_3\right)}\,dx\Big]\,. \label{GF2}
\end{align}

\medskip
\noindent\underline{\noindent\textbf{Stage 3.} \textsf{Time-frequency correlations and the discrete model operator}}
\medskip

From \eqref{GF2} we deduce the following \emph{time-frequency correlations}:

\begin{equation}\label{TFCS}
\begin{cases}
\frac{2}{3}\cdot\frac{r}{2^{2k}}\,+\,\frac{2^{\frac{m}{2}-2k}}{27}\,\left(\frac{v-2u}{2^{\frac{m}{2}}}+\frac{r}{2^{\frac{m}{2}+2k}}\right)^3\,-\,\frac{u+v}{3}=O(1)\\\\
n_1-n_2+ (n_3-\frac{n_2}{2^{2k}})\,\frac{1}{3}\,\left(\frac{v-2u}{2^{\frac{m}{2}}}+\frac{r}{2^{\frac{m}{2}+2k}}\right)^2=O(1)
    \end{cases}\,,
\end{equation}

Applying now a final change of variable $\ds q=\frac{-2u+v}{3}\,+\,\frac{r}{3\cdot2^{2k}}$ and using \eqref{TFCS} we deduce that
\begin{equation}\label{TFCS1}
\begin{cases}
v=\frac{r}{2^{2k}}\,+\,2^{\frac{m}{2}-2k+1}\,\left(\frac{q}{2^{\frac{m}{2}}}\right)^3\,+\,q\,+O(1)\\\\
u=\frac{r}{2^{2k}}\,+\,2^{\frac{m}{2}-2k}\,\left(\frac{q}{2^{\frac{m}{2}}}\right)^3\,-\,q\,+O(1)\\\\
n_1-n_2+3(n_3-\frac{n_2}{2^{2k}})\,\left(\frac{q}{2^{\frac{m}{2}}}\right)^2=O(1).
    \end{cases}
\end{equation}

With these we conclude the final form of \eqref{GF} given by:

\begin{proposition} \label{prop:GFF} There exists an absolute constant $c$, such that for $k\geq 0$ and $\el=\r=j=0$, the form $\L_{k,m}^{0,0,0}$ is approximated (in the sense of being a finite sum - or fastly decreasing infinite sums - of similar terms) by
\begin{align}
\L_{k,m}^{0,0,0}(f_1,f_2,f_3,f_4) & \approx 2^{k}\,\sum_{{r\sim 2^{2k+\frac{m}{2}}}\atop{q\sim 2^{\frac{m}{2}}}} \,\sum_{{n_1,n_2,n_3\sim 2^{\frac{m}{2}}}\atop{n_1=n_2-3\left(n_3-\frac{n_2}{2^{2k}}\right)\,\left(\frac{q}{2^{\frac{m}{2}}}\right)^2}} e^{i\,c\,\left(n_3-\frac{n_2}{2^{2k}}\right)\,2^{\frac{m}{2}}\,\left(\frac{q}{2^{\frac{m}{2}}}\right)^3} \nonumber \\
& \qquad \cdot \Big\langle f_1,\,\v_{\frac{r}{2^{2k}}\,+\,2^{\frac{m}{2}-2k}\,\left(\frac{q}{2^{\frac{m}{2}}}\right)^3-q,\,n_1}^{k} \Big\rangle \,
\Big\langle f_2,\,\v_{\frac{r}{2^{2k}}\,+\,2^{\frac{m}{2}-2k+1}\,\left(\frac{q}{2^{\frac{m}{2}}}\right)^3+q,\,n_2}^{k} \Big \rangle
\nonumber \\ 
\nonumber \\
& \qquad \cdot \big\langle f_3,\,\v_{r,n_3}^{3k} \big \rangle \,
 \Big\langle f_4,\,\v_{\frac{r}{2^{2k}}\,+\,2^{\frac{m}{2}-2k}\,\big(\frac{q}{2^{\frac{m}{2}}}\big)^3, \frac{n_1}{2^{2k}}+\frac{n_2}{2^{2k}}+n_3}^{k} \Big\rangle \,.\label{GFF}
\end{align}
\end{proposition}


%
%
%
\section{The single scale estimate (I): Treatment of the case $k\geq \frac{m}{4} \geq 0$}  \label{sec:onescale_k_large}

This section is devoted to the task of obtaining an accurate $m-$decaying estimate of our four-linear form $\Lambda_{k,m}$ when the scale parameter $k$ is sufficiently large with respect to the oscillation parameter $m$, that is, when $k\geq \frac{m}{4}$. More precisely, we want to prove the estimate \eqref{F1} for the four-linear form $\ds \L_{k,m}^{\el,\r,j}(f_1,f_2,f_3,f_4)$ given by \eqref{fm11:2}, in the case $\el=\r=j=0$. In order to simplify the presentation, we temporarily forget the indices corresponding to $\el=\r=j=0$.

Before entering though into details, we record the following important fact:

\begin{remark}\label{summability} [\textsf{Rank I LGC versus Rank II LGC}] We have the following key contrasting behavior:
\begin{itemize}
\item  The Rank I LGC method employed in Section \ref{subsec:63} focuses on expressing the quadrilinear form $\ds \L_{k,m}^{\el,\r,j}(f_1,f_2,f_3,f_4)$ as a superposition of time-frequency correlated terms involving four-products of Gabor coefficients associated with the input functions---see \eqref{GFF}.
    However, the resulting sum is \emph{not absolutely summable} along spatial/frequency scales, which corresponds to the fact that if one places absolute values within the model sum \eqref{GFF}, no extra decaying/smoothing factor can be obtained. Thus this approach can not provide the desired decay estimate in \eqref{F1}; 
\smallskip

\item The solution to overcome the above mentioned obstacle is provided by the Rank II LGC method which employs an \emph{absolutely summable} model relying on the so-called joint Fourier coefficients which group together the information carried by two input functions with no direct appeal to their Gabor decomposition---see \emph{e.g.} \eqref{descompklarge}. It is thus this latter method that we embrace in the present section (as well as in Sections  \ref{Diagkpsmall124} and \ref{Diagkinterm}). 
 \end{itemize}
\end{remark}

With these being said, we now record several useful facts about the localized Fourier series decomposition (Gabor representation) of the four input functions:  for $f_1$ and $f_2$, the wave packets are supported in frequency on intervals of length $2^{k+\frac{m}{2}}$, and for $f_3, f_4$ they are supported onto frequency intervals of length $2^{3k+\frac{m}{2}}$. The main contribution comes from the wave packets which are adapted to intervals contained inside $I^k$, which is the interval of the space variable $x$, as one can see from \eqref{fm11:2}. For the wave packets adapted to intervals outside of $I^k$ we can invoke the fast decay of the wave packets and implement the same argument. Thus in formula \eqref{fm11:2} (and \eqref{fm11}) we can assume our functions to have the following form: 
\begin{align}
\label{eq:wave:packets:one:scale}
f_1(x)&=\sum_{{u\sim 2^{\frac{m}{2}}}\atop{n_1\sim 2^{\frac{m}{2}}}} \langle f_1,\,\v_{u,n_1}^{k} \rangle \,\v_{u,n_1}^{k}(x)\:, \quad f_2(x)=\sum_{{v\sim 2^{\frac{m}{2}}}\atop{n_2\sim 2^{\frac{m}{2}}}} \langle f_2,\,\v_{v,n_2}^{k} \rangle \,\v_{v,n_2}^{k}(x)\:, \\
f_3(x)&=\sum_{{r\sim 2^{2k+\frac{m}{2}}}\atop{n_3\sim 2^{\frac{m}{2}}}} \langle f_3,\,\v_{r,n_3}^{3k} \rangle \,\v_{r,n_3}^{3k}(x)\:, \quad f_4(x)=\sum_{{r'\sim 2^{2k+\frac{m}{2}}}\atop{n_4\sim 2^{\frac{m}{2}}}} \langle f_3,\,\v_{r',n_4}^{3k} \rangle \,\v_{r',n_4}^{3k}(x). \nonumber
\end{align}
Then for $1 \leq j \leq 4$, $\|f_j \|_{L^2(3 I^k)}$ represents the so-called time-frequency norm of these functions (as defined in Notation \ref{not:norms:TF:reg}), i.e. the $\ell^2$ norm of Gabor coefficients appearing above.



\subsection{A spatial sparse-uniform dichotomy}\label{spunifdichI}
Now we split the information carried by the input functions into spatially sparse and uniform components (the specifics of this analysis depends on the relative position of the parameters $k$ and $m$).

Indeed, given $k,m\in\N$ with $k\geq \frac{m}{4}$ and some fixed $\mu>0$ small (to be chosen later), for  $s\in\{1,2,3,4\}$ we set $\tilde{s}:=2 \left\lfloor\frac{s-1}{2}\right\rfloor+1\in\{1,3\}$ and introduce the space-sparse set of indices
\beq\label{sparseset}
\bbS_{\mu,1}^{0,k,\tilde{s}}:=\bbS^{1}_\mu(f_s):=\Big\{p\sim 2^{(\tilde{s}-1)k+\frac{m}{2}}\:\Big|\:\frac{1}{|I^{m,\tilde{s} k}_{p}|}\int_{3 I^{m,\tilde{s} k}_{p}}|f_s(x)|^2 dx\geq \frac{2^{(4 - \tilde s)\mu m}}{|I^k|}\,\|f_s\|^2_{L^2(I^k)}\Big\}\,.
\eeq
With these, we define the \emph{space-sparse component} of $f_1$ and $f_2$ as
    \begin{equation}\label{sparsefs1}
    f_s^{\bbS_\mu^{1}}:=\sum_{p \in \bbS_\mu^{1}(f_s)} \one_{I_{p}^{m, k}} f_s\,
    \end{equation}
and we retain that $\# \bbS^{1}_\mu(f_s) \leq 2^{{m \over 2} - 3 \mu m}$ for $s=1,2$.

For $f_3$ and $f_4$, the sparse components will be involved in the multiscale analysis in the Section \ref{sec:bilinear:analysis} bellow, so we need to preserve the time-frequency information: if $s=3,4$, then
\begin{equation}\label{sparsefs1:wavepackets}
f_s^{\bbS_\mu^{1}}:=\sum_{{r\sim 2^{2k+\frac{m}{2}}}\atop r \in \bbS^{1}_\mu(f_s)}  \sum_{n_s \sim 2^{\frac{m}{2}}} \langle f_s,\,\v_{r,n_s}^{3k} \rangle \,\v_{r,n_s}^{3k}(x)\:.
\end{equation}
Then for $1 \leq s \leq 4$, we define the \emph{space-uniform component} of $f_s$ by
    \begin{equation}\label{unifs1}
    f_s^{\bbU_\mu^{1}}:=f_s-f_s^{\bbS^{1}_\mu}\,.
    \end{equation}
With these done, we decompose the quadrilinear form $\Lambda_{k,m}(f_1,f_2,f_3,f_4)$ from \eqref{fm11:2} as
\begin{equation}\label{Ldecksmall}
\Lambda_{k,m}(f_1,f_2,f_3,f_4)=:\Lambda_{k,m}^{\bbU^1_\mu}(f_1,f_2,f_3,f_4)\,+\,\Lambda_{k,m}^{\bbS^1_\mu}(f_1,f_2,f_3,f_4),
\end{equation}
where
\begin{equation}\label{sparse1}
\Lambda_{k,m}^{\bbS^1_\mu}(f_1,f_2,f_3,f_4):=\Lambda_{k,m} \left(f_1, f_2, f_3^{\bbS_\mu}, f_4^{\bbS_\mu} \right)\,.
\end{equation}

Before moving on with the proof, we refer the reader to the conventions in Notation \ref{not:norms:TF:reg}. In particular, since the functions $f_s$ and their sparse/uniform components are given by \eqref{eq:wave:packets:one:scale} (for $s=1,2,3,4$) and \eqref{sparsefs1:wavepackets}(for $s=3,4$),
\[
\|f_s\|_2\simeq \|f\|_{\ii L^2(I^k)}, \qquad \|f_s^{\bbS_\mu^{1}}\|_2 \simeq \|f_s^{\bbS_\mu^{1}}\|_{\ii L^2(I^k)}= \Big( \sum_{{r\sim 2^{2k+\frac{m}{2}}}\atop r \in \bbS^{1}_\mu(f_s)}  \sum_{n_s \sim 2^{\frac{m}{2}}} |\langle f_s,\,\v_{r,n_s}^{3k} \rangle|^2 \Big)^{1 \over 2}.
\]
In general, we will not make explicit the collection of wave packets over which we are summing, but we will emphasize whether we are working with the sparse or uniform components.

\subsection{Treatment of the space-sparse component}
\begin{theorem} \label{sparsecont}[\textsf{$L^2$ control over the space-sparse component}] With the previous notations and conventions, we have for any $k, m \in \N$ with $k \geq \frac{m}{4}$ and for any $\mu>0$
\begin{equation}\label{sparsecomp}
|\Lambda_{k,m}^{\bbS^1_\mu}(f_1,f_2,f_3,f_4)|\lesssim 2^{k}\,\|f_1\|_{\ii L^2(I^k)}\,\|f_2\|_{\ii L^2(I^k)}\,\|f_3^{\bbS^1_\mu}\|_{\ii L^2(I^k)}\,\|f_4^{\bbS^1_\mu}\|_{\ii L^2(I^k)}\,.
\end{equation}
\end{theorem}

\begin{proof}
Using \eqref{fm11:2} with $\el=\r=0$ and $j=0$, we get
\begin{equation}\nonumber
|\Lambda_{k,m}^{\bbS^1_\mu}(f_1,f_2,f_3,f_4)|\lesssim  2^{k}\,\int_{I^k}\int_{I^k} \left| f_1(x-t)\,f_2(x+t+t^3)\,f_3^{\bbS^1_\mu}(x-t^3)\,f_4^{\bbS^1_\mu}(x)\right|\,dt\,dx\,,
\end{equation}
Applying first a Cauchy-Schwarz in $t$ followed by a Cauchy-Schwarz in $x$  we get
\begin{align}
& |\Lambda_{k,m}^{\bbS^1_\mu}(f_1,f_2,f_3,f_4)| \nonumber \\
& \qquad \qquad \lesssim 2^{k}\,\int_{I^k}\Big(\int_{I^k} \big| f_1(x-t)\,f_3^{\bbS^1_\mu}(x-t^3)\big|^2dt\Big)^{\frac{1}{2}} \Big(\int_{I^k} \big|f_2(x+t+t^3)\big|^2 \,dt\Big)^{\frac{1}{2}}\,\big|f_4^{\bbS^1_\mu}(x)\big|\,dx\nonumber\\
&\qquad \qquad \lesssim 2^{k}\,\|f_2\|_{L^2(3I^k)}\,\Big(\int_{I^k}\int_{I^k}\big| f_1(x-t)\,f_3^{\bbS^1_\mu}(x-t^3)\big|^2dt dx\Big)^{\frac{1}{2}}\,\Big(\int_{I^k}\big|f_4^{\bbS^1_\mu}(x)\big|^2\,dx\Big)^{\frac{1}{2}}\nonumber
\end{align}
which implies \eqref{sparsecomp}.
\end{proof}

\subsection{Treatment of the space-uniform component}~\\
Our main focus in the present section is proving Theorem \ref{sgscaleth} below.

\begin{theorem} \label{sgscaleth}[\textsf{Control over the space-uniform component}]. There exist two absolute constants $\d_0, \mathfrak c_0>0$ such that that for any $\ds 0<\mu< \frac{\delta_0}{\mathfrak c_0 +\frac{1}{2}}$ and for any $k, m\in\N$ with $k \geq \frac{m}{4}$, we have
\beq\label{unifcomp}
|\Lambda_{k,m}^{\bbU^1_\mu}(f_1,f_2,f_3,f_4)|\lesssim 2^{k}\, \min\big(  2^{- \mu \, m},   2^{-(\d_0- (\mathfrak c_0 +\frac{1}{2}) \mu)\, m} \big)\,\prod_{s=1}^4 \|f_s \|_{\ii L^2(I^k)}\,.
\eeq
\end{theorem}

\begin{remark}
If we set
\begin{equation}
\label{eq:decay:single:mu}
\mu_0:= \frac{\delta_0}{\mathfrak c_0 +{3 \over 2}},
\end{equation}
then the theorem above implies that for any $k, m\in\N$ with $k \geq \frac{m}{4}$,
\begin{equation*}
\label{eq:decay:single:mu:est}
|\Lambda_{k,m}^{\bbU^1_{\mu_0}}(f_1,f_2,f_3,f_4)|\lesssim   2^{- \mu_0 \, m} \, |I^k|^{-1}\,\prod_{s=1}^4 \|f_s \|_{\ii L^2(I^k)}.
\end{equation*}
\end{remark}

\subsubsection{The correlative time-frequency model}\label{Diagkpsmall}~\\
Our main goal in this section is to obtain an $m$-decaying single scale absolutely summable model for our form $\Lambda_m^{k, \bbU^1_\mu}(f_1,f_2,f_3,f_4)$. This will be achieved by applying the Rank II LGC model introduced in \cite{HuLie}.\footnote{For notational simplicity throughout this section we will drop the superindex and refer to $\Lambda_{k,m}^{\bbU^1_\mu}$ as simply $\Lambda_{k,m}$.}

Here we only employ the wave packet decomposition for the function $f_3$, as we try to keep together the information for $f_1$ and $f_2$. As we analyze the space-uniform component, at least one of the functions $f_3$ or $f_4$ are uniformly distributed in space, and so their wave packet decomposition is to be understood according to \eqref{sparsefs1:wavepackets} and \eqref{unifs1}. For simplicity, this will be omitted from the notation.

Next, we split the interval $I^k$ in the following way: the $x$ parameter varies inside intervals of length $2^{-(3k+\frac{m}{2})}$, and the $t$ parameter inside intervals of length $2^{-(k+\frac{m}{2})}$. Thus, we deduce that 
\begin{align}
&\Lambda_{k,m}(f_1,f_2,f_3,f_4)=\sum_{{p\sim 2^{2k+\frac{m}{2}}}\atop{q\sim 2^{\frac{m}{2}}}}\int_{I_p^{m,3k}}\int_{I_q^{m,k}}
f_1(x-t)\,f_2(x+t+t^3) \nonumber \\
&\Big(\sum_{{r\sim 2^{2k+\frac{m}{2}}}\atop{n\sim 2^{\frac{m}{2}}}}  \langle f_3,\,\v_{r,n}^{3k} \rangle \,\v_{r,n}^{3k}(x-t^3)\Big)\,f_4(x)\,2^{k}\,\rho(2^{k} t)\,\varphi(2^{k} x)\,dt\,dx\nonumber\\
&=\sum_{{p\sim 2^{2k+\frac{m}{2}}}\atop{q\sim 2^{\frac{m}{2}}}} \sum_{{r\sim 2^{2k+\frac{m}{2}}}\atop{n\sim 2^{\frac{m}{2}}}}
 \langle f_3,\,\v_{r,n}^{3k} \rangle \,\int_{I_p^{m,3k}} J_{q,n,r}^{k}(f_1,f_2)(x)\,f_4(x)\,2^{k}\,\varphi(2^{k} x)\,dx \,, \label{desc1}
\end{align}
where
\begin{equation}\label{Jdef}
J_{q,n,r}^{k}(f_1,f_2)(x):=2^{\frac{1}{2}(3k+\frac{m}{2})}\,\int_{I_q^{m,k}} f_1(x-t)\,f_2(x+t+t^3)\,\v(2^{3k+\frac{m}{2}}(x-t^3)-r)\,e^{i\,2^{3k+\frac{m}{2}}\,n\,(x-t^3)}\,dt\,.
\end{equation}

As a motivation/consequence of the LGC-methodology transferred to the input variable $t$, we linearize the argument of the function $\v$ as follows: since $t\in I_q^{m,k}$ we write
\begin{align}\label{linearize}
2^{3k+\frac{m}{2}}\,t^3=&2^{3k+\frac{m}{2}}\,\left(t-\frac{q}{2^{k+\frac{m}{2}}}+\frac{q}{2^{k+\frac{m}{2}}}\right)^3\\
=&2^{3k+\frac{m}{2}}\,\left(\frac{q}{2^{k+\frac{m}{2}}}\right)^3\,+
\,3\,2^{3k+\frac{m}{2}}\,\left(t-\frac{q}{2^{k+\frac{m}{2}}}\right)\cdot\left(\frac{q}{2^{k+\frac{m}{2}}}\right)^2\,+
\,O\left(\frac{1}{2^{\frac{m}{2}}}\right)\,.\nonumber
\end{align}

Applying now \eqref{linearize} to \eqref{Jdef} we get (up to a superposition argument)
\begin{align}\label{Jdef1}
&\qquad\qquad\qquad\qquad\qquad\qquad J_{q,n,r}^{k}(f_1,f_2)(x)\approx \\
&e^{i\,\frac{2 q^3}{2^{m}}\,n}\,2^{\frac{1}{2}(3k+\frac{m}{2})}\,\int_{I_q^{m,k}} f_1(x-t)\,f_2(x+t+t^3)\,\v\left(2^{3k+\frac{m}{2}}x-\frac{q^3}{2^m}-r\right)\,e^{i\,2^{3k+\frac{m}{2}}\,n\,
\left(x-3t\left(\frac{q}{2^{k+\frac{m}{2}}}\right)^2\right)}\,dt\,.\nonumber
\end{align}
Now from \eqref{desc1} we may assume $x\in I_p^{m,3k}$ and corroborating this with \eqref{Jdef1} we obtain\footnote{More precisely, the main term in the summation in $r$ comes from the value $r=p\,-\,\frac{q^3}{2^m}\,+\,O(1)$.} the \textit{time-frequency correlation} relation:
\begin{equation}\label{TFC}
r=p\,-\,\frac{q^3}{2^m}\,+\,O(1)\,.
\end{equation}
Putting together \eqref{Jdef1} and \eqref{TFC} -- for $r$ defined by \eqref{TFC} -- we get
\begin{align}\label{Jdef1bis}
&\qquad\qquad\qquad \sum_{r \sim 2^{2k+{m \over 2}}}J_{q,n,r}^{k}(f_1,f_2)(x)\approx \\
&\left(\int_{I_q^{m,k}} f_1(x-t)\,f_2(x+t+t^3)\,e^{-3\,i\,2^{k+\frac{m}{2}}\,t\,
\left(\frac{q}{2^{\frac{m}{2}}}\right)^2\,n}\,dt\right)
\,\v_{p,n}^{3k}(x)\,e^{i\,\frac{2 q^3}{2^{m}}\,n}\,.\nonumber
\end{align}
Thus inserting \eqref{Jdef1} in \eqref{desc1} we deduce
\begin{align}\label{descomp}
&\L_{k,m}(f_1,f_2,f_3,f_4)\approx2^{k}\,\sum_{p\sim 2^{2k+\frac{m}{2}}} \sum_{{q\sim 2^{\frac{m}{2}}}\atop{n\sim 2^{\frac{m}{2}}}}
e^{i\,\frac{2 q^3}{2^{m}}\,n}\, \langle f_3,\,\v_{p-\frac{q^3}{2^m},n}^{3k} \rangle \\
&\cdot\int_{I_p^{m,3k}} \Big(\int_{I_q^{m,k}} f_1(x-t)\,f_2(x+t+t^3)\,e^{-3\,i\,2^{k+\frac{m}{2}}\,t\,
\frac{q^2}{2^{m}}\,n}\,dt\Big)\,f_4(x)\,\v_{p,n}^{3k}(x)\,dx. \nonumber
\end{align}


We are now making the following key
\begin{observation} \label{Ob}
Taking into account that $f_1$ and $f_2$ are supported in frequency on intervals of length $2^{k+m}$ near the origin, we infer that $f_1,f_2$ \textit{are morally constant on intervals of length} $2^{-m-k}$. Deduce that since
\begin{equation}\label{krestr}
k\geq \frac{m}{4}\qquad \Rightarrow\qquad |I_p^{m,3k}|\leq 2^{-m-k}\,,
\end{equation}
and so we have that the function
\begin{equation}\label{krestrbis}
U(x):= \int_{I_q^{m,k}} f_1(x-t)\,f_2(x+t+t^3)\,e^{-3\,i\,2^{k+\frac{m}{2}}\,t\,
\frac{q^2}{2^{m}}\,n}\,dt\,
\end{equation}
is morally constant on intervals of the form $I_p^{m,3k}$.
\end{observation}

Based on the above observation we deduce
\begin{align}\label{descompklarge}
&\L_{k,m}(f_1,f_2,f_3,f_4)\approx 2^{k}\,\sum_{n\sim 2^{\frac{m}{2}}} \sum_{{q\sim 2^{\frac{m}{2}}}\atop{p\sim 2^{2k+\frac{m}{2}}}} e^{i\,\frac{2 q^3}{2^{m}}\,n}\,
 \langle f_3,\,\v_{p-\frac{q^3}{2^m},n}^{3k} \rangle \, \langle f_4,\,\v_{p,n}^{3k} \rangle \\
& \qquad \times \frac{1}{|I_p^{m,3k}|}\,\int_{I_p^{m,3k}} \int_{I_q^{m,k}} f_1(x-t)\,f_2(x+t+t^3)\,e^{-3\,i\,2^{k+\frac{m}{2}}\,t\,
\frac{q^2}{2^{m}}\,n}\,dt\,dx\nonumber\,.
\end{align}
Indeed, we consider the function $U$ as supported on $I_{p}^{m,3k}$ and then identify it to its periodic extension by writing its Fourier decomposition
$$ U(x)= \hat{U}(0) + \sum_{j\neq 0} \hat{U}(j) e^{i2^{\frac{m}{2}+3k}jx}.$$
We claim that the main contribution comes from the first coefficient $\hat{U}(0)$: for $j\neq 0$, the complex exponential can be moved on the wave packet with $e^{i2^{\frac{m}{2}+3k}jx}\v_{p,n}^{3k}(x) = \v_{p,n+j}^{3k}(x)$. So the argument used for $j=0$ can be reproduced in the general case $j \in \Z$ since the coefficient $\hat{U}(j)$ is fast decaying because of integration by parts in the integral formula for the coefficient (as explained in the previous observation, the derivatives of $U$ bring extra decay). This explains the fact that we can morally substitute $U$ by its average $\hat{U}(0)$ and deduce \eqref{descompklarge}.

Consequently, after applying absolute values, we further get

\begin{align}\label{descompklargeabs1}
&|\L_{k,m}(f_1,f_2,f_3,f_4)|\lesssim 2^{k}\,\sum_{n\sim 2^{\frac{m}{2}}} \sum_{{q\sim 2^{\frac{m}{2}}}\atop{p\sim 2^{2k+\frac{m}{2}}}}
| \langle f_3,\,\v_{p-\frac{q^3}{2^m},n}^{3k} \rangle |\,| \langle f_4,\,\v_{p,n}^{3k} \rangle |\\
&\qquad \times \frac{1}{|I_p^{m,3k}|}\,\Big|\int_{I_p^{m,3k}} \int_{I_q^{m,k}} f_1(x-t)\,f_2(x+t+t^3)\,e^{-3\,i\,2^{k+\frac{m}{2}}\,t\,
\frac{q^2}{2^{m}}\,n}\,dt\,dx\Big|.\nonumber
\end{align}
Applying now Cauchy-Schwarz successively, first in $q$ and then in $p$ (with the latter step preceded by a suitable partition of the $p$ range), this last quantity is bounded by
\begin{align*}
& 2^{k}\, \sum_{v\sim 2^{2k}}\sum_{n\sim 2^{\frac{m}{2}}} \Big(\sum_{q\sim 2^{\frac{m}{2}}} \sum_{p=v 2^{\frac{m}{2}}}^{(v+1) 2^{\frac{m}{2}}}| \langle f_3,\,\v_{p-\frac{q^3}{2^m},n}^{3k} \rangle |^2\, | \langle f_4,\,\v_{p,n}^{3k} \rangle |^2\Big)^{\frac{1}{2}}\\
& \quad \cdot \Big(\sum_{q\sim 2^{\frac{m}{2}}} \sum_{p=v 2^{\frac{m}{2}}}^{(v+1) 2^{\frac{m}{2}}} \frac{1}{|I_0^{m,3k}|^2}\, \Big|\int_{I_p^{m,3k}}
\int_{I_q^{m,k}} f_1(x-t)\,f_2(x+t+t^3)\,e^{-3\,i\,2^{k+\frac{m}{2}}\,t\,
\frac{q^2}{2^{m}}\,n}\,dt \,dx\Big|^2\Big)^{\frac{1}{2}}\nonumber\,.
\end{align*}
Applying now $\ell^2\times \ell^2\times \ell^{\infty}$ H\"older's inequality for the $n$ parameter followed by Parseval, we deduce
\begin{align}
\label{eqq:des:3}
&|\L_{k,m}(f_1,f_2,f_3,f_4)|\lesssim 2^{k}\,\sum_{v\sim 2^{2k}}
|f_3\|_{\ii L^2(10 I^{0,3k}_{v})} \, \|f_4\|_{\ii L^2(10 I^{0,3k}_{v})}\\
&\cdot  \Bigg(\sup_{n\sim 2^{\frac{m}{2}}} \sum_{q\sim 2^{\frac{m}{2}}} \sum_{p=v 2^{\frac{m}{2}}}^{(v+1) 2^{\frac{m}{2}}} \frac{1}{|I_0^{m,3k}|^2}\, \left|\int_{I_p^{m,3k}}
\int_{I_q^{m,k}} f_1(x-t)\,f_2(x+t+t^3)\,e^{-3\,i\,2^{k+\frac{m}{2}}\,t\,
\frac{q^2}{2^{m}}\,n}\,dt \,dx\right|^2\Bigg)^{\frac{1}{2}}\nonumber\,.
\end{align}


Since we are now in the space-uniform case, we know that at least one of $f_3$ or $f_4$ is space-uniform. For this function, averages over intervals of length $2^{-\big( 3k+\frac{m}{2}\big)}$ are under control, which implies that averages over intervals of length $2^{-3k}$ are also under control, by definitions \eqref{sparseset} and \eqref{sparsefs1:wavepackets}. Thus, by another Cauchy-Schwarz in $v$, we easily deduce that  
\begin{align}\label{prel}
&|\L_{k,m}(f_1,f_2,f_3,f_4)|\lesssim 2^{\frac{\mu}{2} m}\,\|f_3\|_{\ii L^2(I^k)}\,\|f_4\|_{\ii L^2(I^k)}\nonumber\\
&\cdot  \Bigg(\sum_{v\sim 2^{2k}}\sup_{n\sim 2^{\frac{m}{2}}} \Big(\sum_{p=v 2^{\frac{m}{2}}}^{(v+1) 2^{\frac{m}{2}}} \sum_{q\sim 2^{\frac{m}{2}}} \frac{1}{|I_0^{m,3k}|^2}\, \Big|\int_{I_p^{m,3k}}
\int_{I_q^{m,k}} f_1(x-t)\,f_2(x+t+t^3)\,e^{-3\,i\,2^{k+\frac{m}{2}}\,t\,
\frac{q^2}{2^{m}}\,n}\,dt \,dx\Big|^2\Big)\Bigg)^{\frac{1}{2}}\,.
\end{align}
Above, we use the convention introduced in Notation \ref{not:norms:TF:reg}. For example, if $f_3$ is space-uniform, then it can be written as a sum of wave packets, as one can see from \eqref{sparsefs1:wavepackets} and \eqref{unifs1}; then the $\ell^2$ norm of the wave packets represents a time-frequency norm which eventually is controlled by the $L^2$ norm of the  function $f_3$ localized to the present context, given by formula \eqref{eq:wave:packets:one:scale}.

Next, one further splits the analysis of \eqref{prel} in two subcases.
\subsubsection{At least one of $f_1$ or $f_2$ are \emph{space-sparse}} For simplicity, we assume that $f_1$ is space-sparse. In this case reasonings similar to those from \cite[Section 4.2.2]{HuLie} allow to conclude that
\begin{align}
\label{eq:f_1sparse:Decay}
|\Lambda_{k,m}^{\bbU^1_\mu}(f_1,f_2,f_3,f_4)|\lesssim 2^{k}\,2^{-\mu\, m} \,\prod_{s=1}^4 \|f_s \|_{\ii L^2(I^k)}\,.
\end{align}
Indeed, we can go back to \eqref{descompklargeabs1}, put the absolute values inside the integral in $f_1$ and $f_2$, do a Cauchy-Schwartz in $n$ and in the double integral and reorganize the summation in $p$ to get
\begin{align*}
|\Lambda_{k, m}(f_1, f_2, f_3, f_4)|\lesssim 2^k \sum_{\tilde p \sim 2^{m \over 2}} \sum_{q \sim 2^{m \over 2}} \sum_{p=2^{2k} \tilde p}^{2^{2k}(\tilde p+1)} \|f_3\|_{\ii L^2(I^{m, 3k}_{p - \frac{q^3}{2^m}})} \, \|f_4\|_{\ii L^2(I^{m, 3k}_p)}& \\
 \frac{|I_0^{m,k}|^{1 \over 2}}{|I_0^{m, 3k}|^{1 \over 2}}     \Big(\int_{I_p^{m,3k}} \int_{I_q^{m,k}} \big|f_1(x-t)\,f_2(x+t+t^3) \big|^2\,dt\,dx\Big)^{1 \over 2}.&
\end{align*}

Since $f_3$ is space-uniform, we have that 
\[
 \|f_3\|_{\ii L^2(I^{m, 3k}_{p - \frac{q^3}{2^m}})} \lesssim |I^{m, 3k}_0|^{1 \over 2} \,|I^k|^{ -{1 \over 2}} \, 2^{{\mu m} \over 2} \, \|f_3\|_{\ii L^2(I^k)}.
\]
On the other hand, since $f_1$ is space-sparse, the contribution in $p$ and $q$ is nonzero if and only if $\tilde p -q  = \lfloor 2^{-2k} p \rfloor - q \in \bbS^{1}_\mu(f_1)$, so that, for $\tilde p$ fixed the variable $q$ can take at most $2^{{m \over 2} - 3 \mu m}$ values. Thus, after extra Cauchy-Schwarz in $q$ and then in $p$ we have
\begin{align*}
|\Lambda_{k, m}(f_1, f_2, & f_3, f_4)|\lesssim 2^k  \frac{|I_0^{m,k}|^{1 \over 2}}{|I_0^{m, 3k}|^{1 \over 2}}  \frac{|I^{m, 3k}_0|^{1 \over 2}}{|I^k|^{1 \over 2} }   \, 2^{{\mu m} \over 2} \, 2^{{m \over 4} - {{3 \mu m} \over 2}}  \|f_3\|_{\ii L^2(I^k)}  \,  \|f_4\|_{\ii L^2(I^k)}  \\ &
 \Big(\sum_{p \sim 2^{2k+{m \over 2}}}\sum_{q \sim 2^{m \over 2}} \int_{I_p^{m,3k}} \int_{I_q^{m,k}} \big|f_1(x-t)\,f_2(x+t+t^3) \big|^2\,dt\,dx\Big)^{1 \over 2},
\end{align*}
which is exactly \eqref{eq:f_1sparse:Decay}.

\subsubsection{Both $f_1$ and $f_2$ are \emph{space-uniform}.} In this situation, applying Cauchy-Schwarz in \eqref{prel}, one gets
\begin{equation}\label{keyreduction}
|\L_{k,m}(f_1,f_2,f_3,f_4)|\lesssim 2^{k}\,2^{\frac{\mu}{2} m}\,\|f_3\|_{\ii L^2(I^k)}\,\|f_4\|_{\ii L^2(I^k)}\,\tilde{\Lambda}_{k,m}(f_1,f_2)\,,
\end{equation}
where here we set
\begin{itemize}
\item the joint (maximal) Fourier coefficient
$$J_q^{k,m}(f_1,f_2)(x):=\int_{I_q^{m,k}} f_1(x-t)\,f_2(x+t+t^3)\,e^{-3\,i\,2^{k+\frac{m}{2}}\,t\,
\frac{q^2}{2^{m}}\,\l(x)}\,dt\,;$$
\item the joint (maximal) Fourier coefficient $L_x^2 \ell^2_q-$average:
\begin{equation}\label{keyform}
\tilde{\Lambda}_{k,m}(f_1,f_2):=\left(\frac{1}{|I_0^{m,k}|}\,\int_{I^k}
\sum_{q\sim 2^{\frac{m}{2}}} \Big|\int_{I_q^{m,k}} f_1(x-t)\,f_2(x+t+t^3)\,e^{-3\,i\,2^{k+\frac{m}{2}}\,t\,
\frac{q^2}{2^{m}}\,\l(x)}\,dt\Big|^2\,dx \right)^{\frac{1}{2}}\,;
\end{equation}
\item the $\l$ function that achieves the supremum in $n$ in \eqref{prel}, having the key property that
\begin{equation}\label{lambd}
\l:\:[2^{-k},\,2^{-k+1}]\,\mapsto\,[2^{\frac{m}{2}-1},\,2^{\frac{m}{2}+1}] \textrm{  is measurable and constant on intervals of length}\:2^{-3k}\:.
\end{equation}
\end{itemize}

\subsubsection{Control over the joint (maximal) Fourier coefficient $L_x^2 \ell^2_q-$average}
From \eqref{keyreduction} and \eqref{keyform} it is now immediate to see that our fundamental single-scale Theorem \ref{sgscaleth}, is a direct consequence of the following

\begin{proposition} \label{sgscale} With the previous notations and assuming that\footnote{We remind here that $\tilde{\Lambda}_{k,m}$ and $\L_{k,m}(f_1,f_2,f_3,f_4)$ depend implicitly on the parameter $\mu>0$ used for defining the (space-)sparse and uniform components.}
\begin{itemize} 
\item both $f_1$ and $f_2$ are spatially uniform and supported in frequency on intervals of length $2^{k+m}$ near the origin, and
\item the linearizing function $\l$ verifies \eqref{lambd},
\end{itemize}
we have that there exist two absolute constants $\delta_0, \mathfrak{c}_0 >0$ such that, for any $0<\mu< \frac{\delta_0}{\mathfrak c_0}$ and any $k,\,m\in\N$ with $k\geq \frac{m}{4}$, we have
\beq\label{sgsc-bis}
\tilde{\Lambda}^{\mu}_{k,m}(f_1,f_2)\lesssim 2^{-(\delta_0 - \mathfrak{c}_0\, \mu) m} \,\|f_1\|_{\ii L^2(I^k)}\,\|f_2\|_{\ii L^2(I^k)}.
\eeq
\end{proposition}

\begin{proof}
As it turns out, the strategy here is to reduce the proof of Proposition \ref{sgscale} to that of \cite[Proposition 4.5]{BCa}. Indeed,\footnote{We skip here the technical details which are similar in nature with those in \cite{BCa}.} the outline of our approach is given by:

Firstly, using a similar reasoning as in \eqref{linearize} we notice that for $t\in I_q^{m,k}$ we have $t^3=\Big(\frac{q}{2^{k+\frac{m}{2}}}\Big)^3\,+\,O\Big(\frac{1}{2^{3k+\frac{m}{2}}}\Big)$ and hence since $2^{-3k-\frac{m}{2}}\leq 2^{-k-m}$ for $k\geq \frac{m}{4}$, due to the assumed compact support in frequency for $f_1$ and $f_2$ we have
\begin{equation}\label{cvp1}
\tilde{\Lambda}_{k,m}(f_1,f_2)\approx\left(\frac{1}{|I_0^{m,k}|}\,\int_{I^k}
\sum_{q\sim 2^{\frac{m}{2}}} \left|\int_{I_q^{m,k}} f_1(x-t)\,f_2\left(x+t+\left(\frac{q}{2^{k+\frac{m}{2}}}\right)^3\right)\,e^{-3\,i\,2^{k+\frac{m}{2}}\,t\,
\frac{q^2}{2^{m}}\,\l(x)}\,dt\right|^2\,dx \right)^{\frac{1}{2}}\,.
\end{equation}
Secondly, employing again that $2^{-3k}\leq 2^{-k-\frac{m}{2}}$ and using the change of variable $t\,\rightarrow\,t-\frac{1}{2}\,\left(\frac{q}{2^{k+\frac{m}{2}}}\right)^3$ we further have that \eqref{cvp1} is similar with
\begin{equation}\nonumber
\left(\sum_{q\sim 2^{\frac{m}{2}}} \frac{1}{|I_0^{m,k}|}\,\int_{I^k}
\left|\int_{I_q^{m,k}} f_1\left(x-t+\frac{1}{2}\,\left(\frac{q}{2^{k+\frac{m}{2}}}\right)^3\right)\,f_2\left(x+t+\frac{1}{2}\,\left(\frac{q}{2^{k+\frac{m}{2}}}\right)^3\right)\,e^{-3\,i\,2^{k+\frac{m}{2}}\,t\,
\frac{q^2}{2^{m}}\,\l(x)}\,dt\right|^2\,dx \right)^{\frac{1}{2}}\,.
\end{equation}
Thirdly, applying the change of variable $x\,\rightarrow\,x-\frac{1}{2}\,\left(\frac{q}{2^{k+\frac{m}{2}}}\right)^3$ and using \eqref{lambd} we deduce that 
\begin{equation}\label{cvp2}
\tilde{\Lambda}_{k,m}(f_1,f_2)\approx\underline{\Lambda}_{k,m}(f_1,f_2):=\left(\frac{1}{|I_0^{m,k}|}\,\int_{I^k}
\sum_{q\sim 2^{\frac{m}{2}}} \left|\int_{I_q^{m,k}} f_1(x-t)\,f_2(x+t)\,e^{-3\,i\,2^{k+\frac{m}{2}}\,t\,
\frac{q^2}{2^{m}}\,\l(x)}\,dt\right|^2\,dx \right)^{\frac{1}{2}}\,.
\end{equation}

With these done it remains now to control 
$$ \|\underline{\Lambda}_{k,m}\|_{L^2(3 I^k) \times L^2(3 I^k) \rightarrow {\mathbb R}}.$$
We then perform in \eqref{cvp2} the change of variables $t'=2^k t$ and $x'=2^k x$ (by a slight notational abuse let's still denote $t'$ by $t$ and $x'$ by $x$) and deduce that
$$ \|\underline{\Lambda}_{k,m}\|_{L^2(3 I^k) \times L^2(3 I^k) \rightarrow {\mathbb R}} = \| \tilde{\underline {\Lambda}}_{0,m}\|_{L^2(3 I^0) \times L^2(3 I^0) \rightarrow {\mathbb R}}$$
with
$$ \tilde{\underline{\Lambda}}_{0,m}(f_1,f_2)=\left(2^{\frac{m}{2}}\sum_{p\sim 2^{\frac{m}{2}}}\int_{I_p^{m,0}}
 \sum_{q\sim 2^{\frac{m}{2}}} \left|\int_{I_q^{m,0}} f_1(x-t)\,f_2(x+t)\,e^{-3\,i\,2^{\frac{m}{2}}\,t\,
\frac{q^2}{2^{m}}\,\tilde\l(x)}\,dt\right|^2\,dx \right)^{\frac{1}{2}}$$
and $\tilde \lambda(x)=\lambda(2^{-k} x)$ taking (integer) values within the set $\sim 2^{m/2}$.

Finally, defining $m'$ by $3m'=m$, one concludes that (up to conjugation)
$$ \tilde{\underline{\Lambda}}_{0,m}(f_1,f_2)=\Bigg(2^{\frac{3m'}{2}}\sum_{p\sim 2^{\frac{3m'}{2}}}\int_{I_p^{\frac{3m'}{2},0}}
\sum_{q\sim 2^{\frac{3m'}{2}}} \Big|\int_{\tilde I_q^{\frac{3m'}{2},0}} f_1(x-t)\,f_2(x+t)\,e^{3\,i\,2^{\frac{3m'}{2}}\,t\,
\frac{q^2}{2^{3m'}}\,\tilde\l(x)}\,dt\Big|^2\,dx \Bigg)^{\frac{1}{2}}\,.$$
However, this is precisely the form addressed by \cite[Proposition 4.5]{BCa} corresponding to the operator $BC^a$ with $a=3$ and having as a phase function $ \mu(x)=2^{3\frac{m'}{2}} \tilde \lambda(x)$ which has its order of magnitude $\simeq 2^{3m'}$.
\end{proof}
\section{The single scale estimate (II): Treatment of the case $0 \leq k \leq\frac{m}{8}$}  \label{Diagkpsmall124}
This section is devoted to the single scale estimate \eqref{eq:*ff} in the case $0\leq k \leq \frac{m}{8}$. In view if Observation \ref{Reduct}, we focus on the situation $j=\el=\r=0$. As before in Section \ref{sec:onescale_k_large}, we will exploit formula \eqref{fm11:2} and will make use of the time-frequency localization expressed in \eqref{eq:wave:packets:one:scale}. The sparse components will all be involved in the multi-scale analysis implemented in Section \ref{sec:bilinear:analysis2}.

\subsection{A spatial sparse-uniform dichotomy}\label{spunifdich}~\\
In this case, we define the \emph{space-sparse set} of indices with respect to the parameter $\mu>0$ as follows:
\begin{itemize}
\item for $s\in\{1, 2\}$
\begin{equation}\label{spatialsparsefs12}
\bbS^{1}_\mu \left(f_s \right):=\Big\{ q\sim 2^{2k} \:\Big|\: \frac{1}{|I_{q}^{0, 3k}|}\int_{3 I_{q}^{0, 3k}} \left|f_s\right|^2  \gtrsim \frac{2^{\mu k}}{|I^{k}|} \int_{3 I^{k}} \left|f_s \right|^2 \Big\},
\end{equation}
and define the \emph{space-sparse component} of $f_s$ as
    \begin{equation}\label{sparsefs3}
    f_s^{\bbS^1_\mu}(x):= \sum_{q \in \bbS^{1}_\mu \left(f_s \right)}\sum_{n_s\sim 2^{\frac{m}{2}}}\sum_{{u\sim 2^{\frac{m}{2}}}\atop I_u^{m, k} \subset I_{q}^{0, 3k}} \langle f_s,\,\v_{u,n_s}^{k} \rangle \,\v_{u,n_s}^{k}(x)\,.
    \end{equation}

\item for $s\in\{3, 4\}$
\begin{equation} \label{20230225eq23}
\bbS^{1}_\mu \left(f_s \right):=\Big\{ q\sim 2^{4k} \:\Big|\: \frac{1}{|I_{q}^{0, 5k}|}\int_{3 I_{q}^{0, 5k}} \left|f_s \right|^2  \gtrsim \frac{2^{\mu k}}{|I^{k}|} \int_{3 I^{k}} \left|f_s \right|^2  \Big\}\,,
\end{equation}
and define the \emph{space-sparse component} of $f_s$ as
    \begin{equation}\label{sparsefs31}
    f_s^{\bbS^1_\mu}:=\sum_{q \in \bbS^{1}_\mu \left(f_s \right)}\sum_{n_s\sim 2^{\frac{m}{2}}}\sum_{{r \sim 2^{2k+\frac{m}{2}}}\atop I_r^{m, 3k} \subset I_{q}^{0, 5k}} \langle f_s,\,\v_{r,n_s}^{3k} \rangle \,\v_{r,n_s}^{3k}(x)\,.
    \end{equation}
\end{itemize}

For each of the above cases we define the \emph{space-uniform sets} of indices $\bbU^{1}_\mu(f_s)$ as the complement of the corresponding sparse sets and define the \emph{space-uniform component} of $f_s$ from the equality
$f_s=f_s^{\bbU^1_\mu}+ f_s^{\bbS^1_\mu}$.

With these done, we decompose
\begin{equation}\label{LdecksmallII}
\Lambda_{k,m}(f_1,f_2,f_3,f_4)=:\Lambda_{k,m}^{\bbU^1_\mu}(f_1,f_2,f_3,f_4)\,+\,\Lambda_{k,m}^{\bbS^1_\mu}(f_1,f_2,f_3,f_4),
\end{equation}
where\footnote{The sparse component \eqref{sparse23} has three pieces, but in Theorem \ref{sparsecont1} it is dominated by only two quantities because $\|f_j^{\bbS^1_\mu}\|_{L^2(3I^k)} \leq \|f_j\|_{L^2(3I^k)}$. This is enough for the multi-scale analysis performed in Section \ref{sec:bilinear:analysis2}.}
\begin{align}
\Lambda_{k,m}^{\bbS^1_\mu}(f_1,f_2,f_3,f_4) & :=\Lambda_{k,m} \left(f_1^{\bbS^1_\mu}, f_2^{\bbS^1_\mu}, f_3, f_4\right) \nonumber \\
 & \quad +\, \Lambda_{k,m} \left(f_1, f_2, f_3^{\bbS^1_\mu}, f_4^{\bbS^1_\mu}\right)\,- \Lambda_{k,m} \left(f_1^{\bbS^1_\mu}, f_2^{\bbS^1_\mu}, f_3^{\bbS^1_\mu}, f_4^{\bbS^1_\mu}\right) \,. \label{sparse23}
\end{align}

One can now apply the same steps as in the proof of Theorem \ref{sparsecont} in order to obtain the analogue result:

\begin{theorem} \label{sparsecont1}[\textsf{$L^2$ control over the space-sparse component}]   With the previous notations and conventions, we have
\begin{align}
|\Lambda_{k,m}^{\bbS^1_\mu}(f_1,f_2,f_3,f_4)| & \lesssim 2^{k}\, \|f_1\|_{\ii L^2(I^k)}\,\|f_2\|_{\ii L^2(I^k)}\,\|f_3^{\bbS^1_\mu}\|_{\ii L^2(I^k)}\,\|f_4^{\bbS^1_\mu}\|_{\ii L^2(I^k)} \nonumber \\
& \quad +\, 2^{k} \|f_1^{\bbS^1_\mu}\|_{\ii L^2(I^k)}\,\|f_2^{\bbS^1_\mu}\|_{\ii L^2(I^k)}\,\|f_3\|_{\ii L^2(I^k)}\,\|f_4\|_{\ii L^2(I^k)}.
\label{sparsecomp1}
\end{align}
\end{theorem}

\subsection{Analysis of the space-uniform component}\label{1unifksmall}~\\

The main goal of this section is to provide the resolution of the following

\begin{theorem} \label{sgscale1unifksmal}[\textsf{Control over the space-uniform component}] There exists an absolute constant $\delta>0$ so that for any $\mu>0$ sufficiently\footnote{For example, for any $0<\mu<\frac{\delta}{500}$} small and for any $k,\,m\in\N$ with $0 \leq k\leq \frac{m}{8}$, we have 
\beq 
|\Lambda_{k,m}^{\bbU^1_\mu}(f_1,f_2,f_3,f_4)|\lesssim 2^{-\delta k}\, 2^k \,\prod_{s=1}^4 \|f_s\|_{\ii L^2(I^k)}\,.
\eeq
\end{theorem}

As one can track in the proof, the absolute constant $\delta$ is mainly given by the key estimate \eqref{gain}, together with some intermediate estimates such as \eqref{Uniform-I-Small} or \eqref{Uniform-Ifinal}.

\subsubsection{Preliminaries} Recalling \eqref{sparse23}, we may suppose without loss of generality  that
\begin{equation}\label{unif23}
\Lambda_{k,m}^{\bbU^1_\mu}(f_1,f_2,f_3,f_4):=\Lambda_{k,m} \Big(f_1^{\bbU^1_\mu}, f_2, f_3^{\bbU^1_\mu}, f_4\Big)\,.
\end{equation}

For notational simplicity we drop the superindex dependencies and simply refer to $f_1^{\bbU^1_\mu}$ as $f_1$, $ f_3^{\bbU^1_\mu}$ as $f_3$ and $\Lambda_{k,m}^{\bbU^1_\mu}$ as $\Lambda_{k,m}$. However, we keep in mind that $f_1$ and $f_3$ are \emph{space-uniform}.

\subsubsection{Preliminary frequency foliation}\label{Timefreqfoliation}

Following the strategy in \cite{HuLie}, we re-organize (simultaneously) each spatial fiber of $f_1$ and $f_2$ into frequency locations so that we achieve \emph{both} of the following key requirements:
\begin{itemize}[leftmargin=2em]
\item the length of the frequency subdivisions has to be dominated by the length of the frequency support of the wave-packet $\varphi_{p, n}^{3k}$,
\item the number of the subdivision corresponding to $f_1$ and $f_2$ is the same.
\end{itemize}

Concretely, we proceed as follows: for $j\in\{1, 2\}$ the functions $f_1$ and $f_2$ are supported in frequency onto intervals of length $2^{k+m}$, which we subdivide into intervals of length $2^{\frac{m}{2}+3k}$. We write
$$
\supp \ \widehat{f_j} \subseteq W_{4-2j}^{k+m}=\bigcup_{\ell_j \sim 2^{\frac{m}{2}-2k}} \omega_j^{\ell_j}:=\bigcup_{\ell_j \sim 2^{\frac{m}{2}-2k}} \left[\ell_j 2^{\frac{m}{2}+3 k} + (3-2j) 2^{k+m}, \left(\ell_j+1 \right) 2^{\frac{m}{2}+3 k} + (3-2j) 2^{k+m}\right]\,.$$
Further on, in view of the time-frequency correlation relation that will become one of the main pillars in our approach, we also subdivide the frequency supports of $f_3$ and $f_4$ into $2^{\frac{m}{2}-2k}$ subintervals, as follows: for $j\in\{3, 4\}$,
$$\supp \ \widehat{f_j}\subseteq W^{3k+m}_{-6j+20}:=\bigcup_{\ell_j \sim 2^{\frac{m}{2}-2k}} \omega_j^{\ell_j}
$$
where $\omega_j^{\ell_j}= \left[\ell_j 2^{\frac{m}{2}+5k}+ (-6j+19) 2^{m+3k}, \left(\ell_j+1 \right) 2^{\frac{m}{2}+5k} + (-6j+19) 2^{m+3k}\right]$.

\subsubsection{Action on the ``thickened'' frequency fibers}\label{Decksmall}~\\

In this section we utilize the frequency foliation explained above in order to obtain an analogue of the key heuristic of Observation \ref{Ob}. Thus, for $j\in\{1, 2\}$ we write 
\begin{equation} \label{20230505eq01}
f_j^{\omega_j^{\ell_j}}(x)=\sum_{u \sim 2^{\frac{m}{2}}} \sum_{{n_j\sim 2^{\frac{m}{2}}} \atop {\omega_{n_j}^{m, k} \subset \omega_j^{\ell_j}}} \langle f_j,\,\v_{u,n_j}^{k} \rangle \,\v_{u,n_j}^{k}(x)
\end{equation}
and set
\begin{equation} \label{20230408eq01}
\underline{f}_j^{\omega_j^{\ell_j}}(x):=f_j^{\omega_j^{\ell_j}}(x) e^{-i [\ell_j 2^{\frac{m}{2}+3 k} + (3-2j) 2^{k+m}]x}.
\end{equation}
Thus $\underline{f}_j^{\omega_j^{\ell_j}}$ is a modulation of $f_j^{\omega_j^{\ell_j}}$ which is supported in frequency around the origin. 

\begin{obs} \label{keyobs01} Fix $j\in\{1,2\}$.
Note that since $\supp\,\calF \Big(\underline{f}_j^{\omega_j^{\ell_j}} \Big) \subseteq (-2|\omega_j^{\ell_j}|,\,2|\omega_j^{\ell_j}|)$ with $\big|\omega_j^{\ell_j} \big|=2^{\frac{m}{2}+3k}$ we have that
$$
\underline{f}_1^{\omega_1^{\ell_1}}(\cdot),  \ \underline{f}_2^{\omega_2^{\ell_2}}(\cdot) \ \textnormal{are morally constant on} \ \Big(I_p^{m, 3k} \Big)_{p \sim 2^{\frac{m}{2}+2k}}.
$$
\end{obs}
Applying the above decomposition, we have
\begin{align}
{J}^{k}_{q, n}(f_1,f_2)(x)&=\int_{I_q^{m,k}} f_1(x-t)\,f_2(x+t+t^3)\,e^{-3\,i\,2^{k+\frac{m}{2}}\,t\,
\frac{q^2}{2^{m}}\,n}\,dt \nonumber \\
&=\sum_{\ell_1, \ell_2 \sim 2^{\frac{m}{2}-2k}} e^{i \left(\ell_1 2^{\frac{m}{2}+3k}+\ell_2 2^{\frac{m}{2}+3k} \right)x} \underline{J}^k_{q,n} \left( \underline{f}_1^{\omega_1^{\ell_1}}, \underline{f}_2^{\omega_2^{\ell_2}} \right)(x), \label{Dec}
\end{align}
where
\begin{equation*}
\underline{J}^k_{q,n} \left( \underline{f}_1^{\omega_1^{\ell_1}}, \underline{f}_2^{\omega_2^{\ell_2}} \right)(x):=\int_{I_q^{m, k}}  \underline{f}_1^{\omega_1^{\ell_1}}(x-t) \underline{f}_2^{\omega_2^{\ell_2}}(x+t+t^3) e^{i \left(-\ell_1 2^{\frac{m}{2}+3k}t+\ell_2 2^{\frac{m}{2}+3k} (t+t^3)-3 \cdot 2^{\frac{m}{2}+k} \cdot \frac{q^2}{2^m} n t \right)} dt.
\end{equation*}

\begin{obs} \label{keyobs-02}
As a consequence of the above decomposition, we obtain the desired analogue of Observation \ref{Ob}, that is,  we have that the mapping
$$
x \to \underline{J}^k_{q,n} \left( \underline{f}_1^{\omega_1^{\ell_1}}, \underline{f}_2^{\omega_2^{\ell_2}} \right)(x)
$$
is morally constant on intervals $\Big(I_p^{m, 3k} \Big)_{p \sim 2^{\frac{m}{2}+2k}}$.
\end{obs}

\subsubsection{Rescaled time-frequency correlation set}\label{RescTFC}~\\
In this section we focus on obtaining a rescaled version of the last equation in \eqref{TFCS1} that is adapted to the frequency foliation accomplished earlier and which will prove quintessential in controlling our quadrilinear form $\Lambda_{k,m}$ in the case $0\leq k\leq \frac{m}{4}$.

\begin{lem}\label{Rtfc}
Let $\ell_1, \ell_2, \ell_3\sim 2^{\frac{m}{2}-2k}$ and $q\sim 2^{\frac{m}{2}}$ be such that their contribution in \eqref{20230505eq01} is not zero, then one has
\begin{equation} \label{TFC-03}
\ell_1- \ell_2+\frac{3q^2}{2^m}\left( \ell_3-\frac{\ell_2}{2^{2k}}\right)=O(1).
\end{equation}
\end{lem}

\begin{proof}
Recall the decomposition \eqref{GFF} and the fact that $n_j \sim 2^{\frac{m}{2}}$. Fix $\ell_1,\,\ell_2,\,\ell_3$ and remark that both
\begin{itemize}
\item $n_j 2^{\frac{m}{2}+k} \in \omega_j^{\ell_j}$ for $j\in\{1, 2\}$,

\item $n_3 2^{\frac{m}{2}+3k} \in \omega_3^{\ell_3}$,
\end{itemize}
imply
\begin{equation} \label{20230203eq40}
n_j \in
\widetilde{\omega_j^{\ell_j}}:=\left[\ell_j 2^{2k}, \left(\ell_j+1 \right) 2^{2k} \right]\qquad j\in\{1,2,3\}.
\end{equation}

By \eqref{TFCS1}, that is, $n_1-n_2+3(n_3-\frac{n_2}{2^{2k}})\,\left(\frac{q}{2^{\frac{m}{2}}}\right)^2=O(1)$, we see that
$$
2^{2k} \ell_1-2^{2k} \ell_2+\frac{3q^2}{2^m} \cdot 2^{2k}\left( \ell_3-\frac{\ell_2}{2^{2k}}\right)=O(2^{2k}).
$$
The desired claim then follows by dividing by $2^{2k}$ on both sides of the above equation.
\end{proof}

As a result of the above lemma, it is now natural to introduce \emph{the time-frequency correlation set} as follows:
for each $q \sim 2^{\frac{m}{2}}$ we let
\begin{equation} \label{TFCq}
\TFC_q:=\Big\{ \substack{\left(\ell_1, \ell_2, \ell_3 \right)\\\ell_1, \ell_2, \ell_3 \sim 2^{\frac{m}{2}-2k}}\, :\, \ell_1-\ell_2+\frac{3q^2}{2^m}\Big( \ell_3-\frac{\ell_2}{2^{2k}}\Big)=O(1)\Big\}.
\end{equation}

\subsubsection{The correlative time-frequency model}\label{Model}~\\
With the above settled, we are now ready to refine the time-frequency decomposition \eqref{GFF} of our quadrilinear form  $\Lambda_{k,m}$. Indeed, using \eqref{Dec}, Observation \ref{keyobs-02}, Lemma \ref{Rtfc} and \eqref{TFCq}, we have
\begin{align*} 
& |\Lambda_{k,m}(f_1,f_2,f_3,f_4)| \\
& \lesssim 2^k \sum_{\substack{p \sim 2^{\frac{m}{2}+2k}\\ q \sim 2^{\frac{m}{2}}}} \sum_{(\ell_1, \ell_2, \ell_3) \in \TFC_q} \sum_{n \in \widetilde{\omega_3^{\ell_3}}} \Big| \Big\langle f_3^{\omega_3^{\ell_3}}, \varphi_{p-\frac{q^3}{2^m}, n}^{3k} \Big \rangle \Big|  \Big| \Big \langle f_4 e^{i \left(\ell_12^{\frac{m}{2}+3k}+\ell_22^{\frac{m}{2}+3k} \right) \left( \cdot \right) }, \varphi_{p, n}^{3k} \Big \rangle \Big|  \\
&\cdot \frac{1}{\left|I_p^{m, 3k} \right|} \bigg| \int_{I_p^{m, 3k}} \int_{I_q^{m, k}}  \underline{f}_1^{\omega_1^{\ell_1}}(x-t) \underline{f}_2^{\omega_2^{\ell_2}}(x+t+t^3) e^{i \left(-\ell_1 2^{\frac{m}{2}+3k}t+\ell_2 2^{\frac{m}{2}+3k} t\left(1+3 \left(\frac{q}{2^{\frac{m}{2}+k}}\right)^2\right)-3 \cdot 2^{\frac{m}{2}+k} \cdot \frac{q^2}{2^m}n t \right)} dt dx \bigg|\,.
\end{align*}
Finally, in order to decouple the information carried by $\ell_3$ and $n$, we apply the change variable $n \to n+\ell_3 2^{2k}$ in order to conclude the correlative time-frequency model:
\begin{align} 
&|\Lambda_{k,m}(f_1,f_2,f_3,f_4)| \nonumber \\
& \lesssim 2^k \sum_{\substack{p \sim 2^{\frac{m}{2}+2k}\\ q \sim 2^{\frac{m}{2}} \\ n \sim 2^{2k}}} \sum_{(\ell_1, \ell_2, \ell_3) \in \TFC_q}  \Big| \Big\langle f_3^{\omega_3^{\ell_3}}, \varphi_{p-\frac{q^3}{2^m}, \ell_3 2^{2k}+n}^{3k} \Big \rangle \Big|  \Big| \Big \langle f_4, \varphi_{p, \ell_3 2^{2k}+n+\ell_1+\ell_2}^{3k} \Big \rangle \Big| \nonumber \\
&\cdot \frac{1}{\left|I_p^{m, 3k} \right|} \bigg| \int_{I_p^{m, 3k}} \int_{I_q^{m, k}}  \underline{f}_1^{\omega_1^{\ell_1}}(x-t) \underline{f}_2^{\omega_2^{\ell_2}}(x+t+t^3)\,e^{-i 2^{\frac{m}{2}+k}t \cdot \frac{3q^2}{2^m}n}\,e^{i 2^{\frac{m}{2}+3k}t \cdot \left(-\ell_1 +\ell_2 -\left(\ell_3-\frac{\ell_2}{2^{2k}}\right) \frac{3q^2}{2^m}\right)} dt dx \bigg|. \label{20230224main01}
\end{align}

\subsection{A time-frequency sparse--uniform dichotomy}\label{SU1}~\\
 We start by reconsidering the foliation of the time-frequency space as follows:
\begin{itemize}
\item for $j\in\{1,2\}$, define
$$
\Omega_j^{L_j}:=\left[L_j 2^{m-k}, \left(L_j+1 \right) 2^{m-k} \right], \quad L_j \sim 2^{2k}\,;
$$
\item for $j\in\{3,4\}$, define
$$
\Omega_j^{L_j}:=\left[L_j 2^{m+k}, \left(L_j+1 \right) 2^{m+k} \right], \quad L_j \sim 2^{2k}\,;
$$
\end{itemize}
and notice that $\big|\Omega_j^{L_j} \big| \ge \big|\omega_j^{l_j} \big|$. In particular, since $f_4$ is the dualizing function, the interval $\Omega_4^{L_4}$ is fully determined by the remaining $\Omega_j^{L_j}$s, and so we will sometimes denote it $\Omega_4^{L_3}$.

Next, we perform a time-frequency sparse/uniform decomposition adapted to the above foliation. More precisely, for $\nu \in (0, 1)$ suitably small, we define the \emph{time-frequency sparse} ($2$-sparse) set of indices as follows:
\begin{itemize}
\item for $j\in\{1, 2\}$
$$
\bbS^2_\nu \left(f_j \right):=\left\{ \substack{(q_j, L_j)\\ q_j, L_j \sim 2^{2k}} : \int_{3 I_{q_j}^{0, 3k}} \Big|f_j^{\Omega_j^{L_j}} \Big|^2  \gtrsim 2^{-3\nu k} \int_{3 I_{q_j}^{0, 3 k}} \left|f_j \right|^2 \right\},
$$
\item for $j\in\{3, 4\}$
\begin{equation} \label{tfsparse34}
\bbS^2_\nu \left(f_j \right):=\left\{ \substack{(q_j, L_j)\\ q_j, L_j \sim 2^{2k}} : \int_{3 I_{q_j}^{0, 3k}} \Big|f_j^{3 \Omega_j^{L_j}} \Big|^2  \gtrsim 2^{-\nu k} \int_{3 I_{q_j}^{0, 3k}} \left|f_j \right|^2 \right\}.
\end{equation}
\end{itemize}

We also define the \emph{time-frequency uniform sets} of indices $\bbU^2_\nu(f_j)$ as the complement of the corresponding 2-sparse sets. With these, we decompose our functions into the time-frequency sparse and uniform components as follows: for $j\in\{1,\ldots,4\}$ we set
\begin{enumerate}[leftmargin=2em]
\item[$\bullet$] $f_j^{\Omega_j^{L_j}}$ is given by $\ds f_j^{\Omega_j^{L_j}}= \sum_{\substack{\ell_j \sim 2^{{m \over 2}-k} \\ \omega_j^{\ell_j} \subset \Omega_j^{L_j}}}f_j^{\omega_j^{\ell_j}}$; thus, using formula \eqref{20230505eq01} we deduce that $f_j^{\Omega_j^{L_j}}$ is the frequency projection of $f_j$ onto $\Omega_j^{L_j}$;
    \item [$\bullet$] recalling that $\tilde{j}=2 \left\lfloor\frac{j-1}{2}\right\rfloor+1\in\{1,3\}$, the \emph{time-frequency uniform component} of $f_j$ is given by
    \begin{equation}\label{unifj}
    f_j^{\bbU^2_\nu}:=\sum_{(q_j, L_j) \in \bbU_\nu^2(f_j)}  \sum_{{u\sim 2^{\frac{m}{2}}} \atop {I_{u}^{m, \tilde j k} \subset I_{q_j}^{0, 3k}}}  \sum_{{n_j\sim 2^{\frac{m}{2}}}\atop {\omega_{n_j}^{m, \tilde j k} \subset \Omega_j^{L_j}}} \langle f_j,\,\v_{u,n_j}^{\tilde j k} \rangle \,\v_{u,n_j}^{\tilde j k}(x);
    \end{equation}
    \item [$\bullet$] the \emph{time-frequency sparse component} of $f_j$ given by
    \begin{equation}\label{sparsj}
    f_j^{\bbS^2_\nu}:=\sum_{(q_j, L_j) \in \bbS_\nu^2(f_j)}  \sum_{{u\sim 2^{\frac{m}{2}}} \atop {I_{u}^{m, \tilde j k} \subset I_{q_j}^{0, 3k}}}  \sum_{{n_j\sim 2^{\frac{m}{2}}}\atop {\omega_{n_j}^{m, \tilde j k} \subset \Omega_j^{L_j}}} \langle f_j,\,\v_{u,n_j}^{\tilde j k} \rangle \,\v_{u,n_j}^{\tilde j k}(x).
    \end{equation}
\end{enumerate}

With these done, we decompose
\begin{equation*}
\Lambda_{k,m}(f_1,f_2,f_3,f_4)=:\Lambda_{k,m}^{\bbU_\nu^2}(f_1,f_2,f_3,f_4)\,+\,\Lambda_{k,m}^{\bbS_\nu^2}(f_1,f_2,f_3,f_4),
\end{equation*}
with $\Lambda_{k,m}^{\bbS_\nu^2}(f_1,f_2,f_3,f_4):=\Lambda_{k,m} \left(f_1^{\bbS^2_\nu}, f_2^{\bbS^2_\nu}, f_3^{\bbS^2_\nu}, f_4^{\bbS^2_\nu} \right)$.

\subsection{Treatment of the $2$-uniform component $\Lambda_{k,m}^{\bbU_\nu^2}$}

Applying Cauchy-Schwarz in \eqref{20230224main01}, we have
\begin{equation} \label{20230224eq01}
|\Lambda_{k,m}^{\bbU_\nu^2}(f_1,f_2,f_3,f_4)|\lesssim 2^k \sum_{p \sim 2^{\frac{m}{2}+2k}} \bbA_p \bbB_p,
\end{equation}
where
\begin{equation} \label{20230305eq01}
\bbA_p^2:=\sum_{\substack{n_3 \sim 2^{2k} \\ q \sim 2^{\frac{m}{2}}}} \sum_{\left(\ell_1, \ell_2, \ell_3 \right) \in \TFC_q} \Big| \Big\langle f_3^{\omega_3^{\ell_3}}, \varphi_{p-\frac{q^3}{2^m}, \ell_3 2^{2k}+n_3}^{3k} \Big \rangle \Big|^2  \Big| \Big \langle f_4, \varphi_{p, \ell_3 2^{2k}+n_3+\ell_1+\ell_2}^{3k} \Big \rangle \Big|^2,
\end{equation}
and
\begin{align} 
\bbB_p^2 & :=\sum_{\substack{n_3 \sim 2^{2k} \\ q \sim 2^{\frac{m}{2}}}} \sum_{\left(\ell_1, \ell_2, \ell_3 \right) \in \TFC_q}  \frac{1}{\big|I_p^{m, 3k} \big|}\,\int_{I_p^{m, 3k}} \bigg| \int_{I_q^{m, k}}  \underline{f}_1^{\omega_1^{\ell_1}}(x-t) \underline{f}_2^{\omega_2^{\ell_2}}(x+t+t^3) \nonumber \\
& \qquad \qquad \qquad \qquad \cdot e^{-i 2^{\frac{m}{2}+k}t \cdot \frac{3q^2}{2^m}n_3}   e^{i 2^{\frac{m}{2}+3k}t \cdot \left(-\ell_1 +\ell_2 -\left(\ell_3-\frac{\ell_2}{2^{2k}}\right) \frac{3q^2}{2^m}\right)}  dt  \bigg|^2 dx. \label{20230305eq02}
\end{align}

\medskip

\noindent\textit{Treatment of $\bbA_p$}. By making the change of variable $\ell_2+\ell_1 \to \ell_2$ we get that
\begin{equation}\label{ap}
\bbA_p^2\lesssim \sum_{\substack{n_3 \sim 2^{2k} \\ q \sim 2^{\frac{m}{2}}}} \sum_{\ell_2, \ell_3 \sim 2^{\frac{m}{2}-2k}}  \Big| \Big\langle f_3^{\omega_3^{\ell_3}}, \varphi_{p-\frac{q^3}{2^m}, \ell_3 2^{2k}+n_3}^{3k} \Big \rangle \Big|^2  \Big| \Big \langle f_4, \varphi_{p, \ell_3 2^{2k}+n_3+\ell_2}^{3k} \Big \rangle \Big|^2.
\end{equation}
Next, on each spatial fiber, summing up with respect to $\ell_2$ followed by $n_3$ and using then the fact that $2^{2k} \le 2^{\frac{m}{2}-2k}$, we see that the above term is further bounded from above by\footnote{Here we use the convention $\Big| \Big\langle f_3^{\omega_3^{\ell_3}}, \varphi_{p-\frac{q^3}{2^m}, \ell_3 2^{2k}+\left[0, 2^{2k+1} \right]}^{3k} \Big \rangle \Big|^2$ for the corresponding sum $\sum\limits_{n \in [0, 2^{2k+1}] \cap \Z} \Big| \Big\langle f_3^{\omega_3^{\ell_3}}, \varphi_{p-\frac{q^3}{2^m}, \ell_3 2^{2k}+n}^{3k} \Big \rangle \Big|^2.$}
$$
\lesssim \sum_{\substack{\ell_3 \sim 2^{\frac{m}{2}-2k} \\ q \sim 2^{\frac{m}{2}}}} \left| \left\langle f_3^{\omega_3^{\ell_3}}, \varphi_{p-\frac{q^3}{2^m}, \ell_3 2^{2k}+\left[0, 2^{2k+1} \right]}^{3k} \right \rangle \right|^2  \left| \left \langle f_4, \varphi_{p, \ell_3 2^{2k}+\left[0, 2^{\frac{m}{2}-2k+1} \right]}^{3k} \right \rangle \right|^2\,.
$$

By Parseval, we further have
\begin{align} 
\bbA_p^2 &\lesssim \sum_{\substack{\ell_3 \sim 2^{\frac{m}{2}-2k} \\ q \sim 2^{\frac{m}{2}}}} \left\|f_3^{3\omega_3^{\ell_3}} \right\|^2_{L^2 \left(I_{p-\frac{q^3}{2^m}}^{m, 3k} \right)} \left\|f_4^{2^{\frac{m}{2}+5k}\ell_3+\left[0, 2^{m+k+1} \right]} \right\|_{L^2 \left(I_p^{m, 3k} \right)}^2 \nonumber  \\
&\lesssim \sum_{\ell_3 \sim 2^{\frac{m}{2}-2k}} \left\|f_3^{3\omega_3^{\ell_3}} \right\|^2_{L^2 \left(3 I_{\frac{p}{2^{\frac{m}{2}}}}^{0, 3k} \right)} \left\|f_4^{2^{\frac{m}{2}+5k}\ell_3+\left[0, 2^{m+k+1} \right]} \right\|_{L^2 \left(I_p^{m, 3k} \right)}^2 \nonumber  \\
&\lesssim \sum_{L_3 \sim 2^{2k}} \left\|f_3^{3 \Omega_3^{L_3}} \right\|_{L^2\left(3I_{\frac{p}{2^{\frac{m}{2}}}}^{0, 3k} \right)}^2  \left\|f_4^{3\Omega_4^{L_3}} \right\|_{L^2 \left(I_p^{m, 3k} \right)}^2. \label{20230225eq01}
\end{align}

\medskip

\noindent\textit{Treatment of $\bbB_p$.} Observe that for any $h\in L^2(I^k)$, by Parseval, one has
\begin{equation} \label{20230225eq02}
\sum_{n_3 \sim 2^k} \Big|\int_{I_q^{m, k}} h(t) e^{i2^{\frac{m}{2}+k}t \cdot \frac{3q^2}{2^m} n_3} dt \Big|^2 \lesssim  \left|I_q^{m, k} \right| \int_{I_q^{m, k}} |h(s)|^2 ds.
\end{equation}
Applying \eqref{20230225eq02} for $\ds h(t)=\underline{f}_1^{\omega_1^{\ell_1}}(x-t) \underline{f}_2^{\omega_2^{\ell_2}}(x+t+t^3) e^{i 2^{\frac{m}{2}+3k}t \cdot \left(-\ell_1 +\ell_2 -\left( \ell_3 -\frac{\ell_2}{2^{2k}} \right) \frac{3q^2}{2^m}\right)}$, we have
\begin{equation} \label{20230225eq03}
\bbB_p^2 \lesssim 2^{2k}\,\sum_{\substack{\ell_1, \ell_2 \sim 2^{\frac{m}{2}-2k} \\q \sim 2^{\frac{m}{2}}}}  \int_{I_p^{m, 3k}} \int_{I_q^{m, k}} \left|\underline{f}_1^{\omega_1^{\ell_1}}(x-t) \underline{f}_2^{\omega_2^{\ell_2}}(x+t+t^3) \right|^2 dtdx.
\end{equation}

\medskip

Combining \eqref{20230224eq01}, \eqref{20230225eq01}, and \eqref{20230225eq03}, we conclude that
\begin{align} 
|\Lambda_{k,m}^{\bbU_\nu^2}(f_1,f_2,f_3,f_4)| &\lesssim 2^{2k} \sum_{p \sim 2^{\frac{m}{2}+2k}} \Bigg( \sum_{L_3 \sim 2^{2k}} \Big\|f_3^{3 \Omega_3^{L_3}} \Big\|_{L^2\big(3 I_{\frac{p}{2^{\frac{m}{2}}}}^{0, 3k} \big)}^2  \Big\|f_4^{3\Omega_3^{L_3}} \Big\|_{L^2 \big(I_p^{m, 3k} \big)}^2  \Bigg)^{\frac{1}{2}}  \nonumber \\
&\cdot \Bigg(\sum_{q \sim 2^{\frac{m}{2}}}\sum_{\ell_1, \ell_2 \sim 2^{\frac{m}{2}-2k}}  \int_{I_p^{m, 3k}} \int_{I_q^{m, k}} \left|\underline{f}_1^{\omega_1^{\ell_1}}(x-t) \underline{f}_2^{\omega_2^{\ell_2}}(x+t+t^3) \right|^2 dtdx \Bigg)^{\frac{1}{2}}\,. \label{20230305eq05}
\end{align}

Applying a partial Cauchy-Schwarz in $p$ and then grouping the intervals $\big(I_q^{m, k}\big)_{q \sim 2^{\frac{m}{2}}}$ into a collection of larger intervals $\big(I_{\widetilde{q}}^{0, 3k}\big)_{\widetilde{q} \sim 2^{2k}}$, we deduce that
\begin{align} 
|\Lambda_{k,m}^{\bbU_\nu^2}(f_1,f_2,f_3,f_4)| & \lesssim 2^{2k} \sum_{\widetilde{p} \sim  2^{2k}}  \Big( \sum_{L_3 \sim 2^{2k}} \big\|f_3^{3 \Omega_3^{L_3}} \big\|^2_{L^2 \left(3 I_{\widetilde{p}}^{0, 3k} \right)}\big\|f_4^{3\Omega_3^{L_3}} \big\|^2_{L^2 \left(3 I_{\widetilde{p}}^{0, 3k} \right)} \Big)^{\frac{1}{2}} \nonumber \\
& \cdot \Big(\sum_{\widetilde{q} \sim 2^{2k}}\sum_{\ell_1, \ell_2 \sim 2^{\frac{m}{2}-2k}} \int_{I_{\widetilde{p}}^{0, 3k}} \int_{I_{\widetilde{q}}^{0, 3k}} \left|\underline{f}_1^{\omega_1^{\ell_1}}(x-t) \underline{f}_2^{\omega_2^{\ell_2}}(x+t+t^3) \right|^2 dtdx \Big)^{\frac{1}{2}}.  \label{20230225eq04}
\end{align}
Finally, using Parseval, we conclude that
\begin{align} 
\left|\Lambda_{k,m}^{\bbU_\nu^2}(f_1,f_2,f_3,f_4)\right| &\lesssim 2^{2k} \sum_{\widetilde{p} \sim  2^{2k}}  \Big( \sum_{L_3 \sim 2^{2k}} \big\|f_3^{3 \Omega_3^{L_3}} \big\|^2_{L^2 \left(3 I_{\widetilde{p}}^{0, 3k} \right)}\big\|f_4^{3\Omega_3^{L_3}} \big\|^2_{L^2 \left(3 I_{\widetilde{p}}^{0, 3k} \right)} \Big)^{\frac{1}{2}}\nonumber \\
&\cdot \Big[\sum_{ \widetilde{q} \sim 2^{2k}}\sum_{L_1, L_2 \sim 2^{2k} }
\left(\int_{3 I^{0, 3k}_{\widetilde{p}-\widetilde{q}}} \left|\underline{f}_1^{\Omega_1^{L_1}}(u) \right|^2 du \right)  \left( \int_{3 I^{0, 3k}_{\widetilde{p}+\widetilde{q}}} \left| \underline{f}_2^{\Omega_2^{L_2}}(v) \right|^2 dv \right) \Big]^{\frac{1}{2}}.  \label{20230225main01}
\end{align}
Since we cannot simultaneously perfectly localize the functions in space and in frequency, the last application of Parseval should be interpreted in the following\footnote{This is a very standard approach and it has already been used in Part I.} way: due to a change of variables as in the proof of Proposition \ref{prop:Tk}, we have that 
\begin{align*}
 \int_{I_{\widetilde{p}}^{0, 3k}} \int_{I_{\widetilde{q}}^{0, 3k}} \left|\underline{f}_1^{\omega_1^{\ell_1}}(x-t) \underline{f}_2^{\omega_2^{\ell_2}}(x+t+t^3) \right|^2 dtdx  \lesssim \big\|  \underline{f}_1^{\omega_1^{\ell_1}} \big\|^2_{L^2(3 I^{0, 3k}_{\widetilde{p}-\widetilde{q}})} \, \big\|  \underline{f}_2^{\omega_2^{\ell_2}} \big\|^2_{L^2(3 I^{0, 3k}_{\widetilde{p}+\widetilde{q}})}.
\end{align*} 
Using the wave packet decomposition of $\underline{f}_j^{\omega_j^{\ell_j}}$ from \eqref{20230505eq01} and Notation \ref{not:norms:TF:reg}, we have that 
\begin{align}
\label{eq:tails:explicit}
\big\|  \underline{f}_1^{\omega_1^{\ell_1}} \big\|^2_{L^2(3 I^{0, 3k}_{\widetilde{p}-\widetilde{q}})} \lesssim \big\|  \underline{f}_1^{\omega_1^{\ell_1}} \big\|^2_{\tilde {\ii L}^2(3 I^{0, 3k}_{\widetilde{p}-\widetilde{q}})} := \sum_{\#_1 \in \Z} (1+|\#_1|)^{-M} \big\| f_1  \big\|^2_{L^2(3I^{0, 3k}_{\widetilde{p}-\widetilde{q}, \#_1} \times \omega_1^{\ell_1})},
\end{align}
where in the summation above we consider only wave packets with spatial interval contained inside a shifts $3 I^{0, 3k}_{\widetilde{p}-\widetilde{q}, \#_1}$ of $3 I^{0, 3k}_{\widetilde{p}-\widetilde{q}}$, which are still included in $I^k$. A similar inequality holds for $f_2$. Now we can sum over indices $\ell_j$ with $\omega_j^{\ell_j} \subset \Omega_j^{L_j}$ to obtain the expression in the last line of \eqref{20230225main01}, with tails. In order to simplify the presentation and since these tails do not affect the analysis to come neither in Subsection \ref{20230305subsec01} nor in Subsection \ref{subsec4.2.2}, we will not make them explicit in the choice of notation. Indeed, we will just need to sum contributions of the form \eqref{eq:tails:explicit} for intervals of length $2^{-3k}$ contained in $I^k$, which eventually are bounded above by $\|f_1\|_{L^2(I^k \times W_2^{k+m})}^2$.

Now we consider several different cases.
\subsubsection{Either $f_3$ or $f_4$ is time-frequency uniform} \label{20230305subsec01}~\\
 In this case, without loss of generality, we assume $f_3=f_3^{\bbU^2_\nu}$. Then from \eqref{20230225main01}, Parseval and Cauchy--Schwarz one has
 \begin{align} 
\left|\Lambda_{k,m}^{\bbU_\nu^2}(f_1,f_2,f_3,f_4)\right| & \lesssim 2^{2k} \sum_{\widetilde{p} \sim  2^{2k}}  \Big( \sum_{L_3 \sim 2^{2k}} 2^{-\nu k} \left\|f_3 \right\|^2_{L^2 \left(3I_{\widetilde{p}}^{0, 3k} \right)}\big\|f_4^{3\Omega_3^{L_3}} \big\|^2_{L^2 \left(3 I_{\widetilde{p}}^{0, 3k} \right)} \Big)^{\frac{1}{2}}\nonumber \\
&\quad\cdot \Bigg[\sum_{\substack{L_1, L_2 \sim 2^{2k} \\ \widetilde{q} \sim 2^{2k}}}
\Big(\int_{3 I^{0, 3k}_{\widetilde{p}-\widetilde{q}}} \big|\underline{f}_1^{\Omega_1^{L_1}}(u) \big|^2 du \Big)  \Big( \int_{3 I^{0, 3k}_{\widetilde{p}+\widetilde{q}}} \big| \underline{f}_2^{\Omega_2^{L_2}}(v) \big|^2 dv \Big) \Bigg]^{\frac{1}{2}}\nonumber \\
&\lesssim 2^{2k} \cdot 2^{-\frac{\nu k}{2}} \left\|f_1 \right\|_{L^2 \left(3 I^k \right)} \left\|f_2 \right\|_{L^2 \left(3 I^k \right)}\left( \sum_{\widetilde{p} \sim  2^{2k}}  \left\|f_3 \right\|^2_{L^2 \left(3 I_{\widetilde{p}}^{0, 3k} \right)}\left\|f_4 \right\|^2_{L^2 \left(3 I_{\widetilde{p}}^{0, 3k} \right)} \right)^{\frac{1}{2}}\nonumber\\
&\lesssim 2^{k} \cdot 2^{-\frac{(\nu-\mu)k}{2}} \prod_{s=1}^4 \left\|f_s \right\|_{L^2 \left(3 I^k \right)}\,, \label{Uniform-I-Small}
\end{align}
where here, recalling \eqref{20230225eq23}, we used the hypothesis that $f_3$ is also $1$-uniform. This estimate is acceptable in view of proving Theorem \ref{sgscale1unifksmal}, by setting $\delta \ll \nu$ and $\mu=\frac{\nu}{2}$.

\subsubsection{Both $f_3$ and $f_4$ are time-frequency sparse and at least one of the functions $f_1$ or $f_2$ is time-frequency uniform} \label{subsec4.2.2}~\\
In this case, by \eqref{20230225eq23}, we deduce that up to a factor of at most $2^{\nu k}$, there exists a unique measurable function $L_3: \left[2^{2k}, 2^{2k+1}\right]\cap \Z \to \left[2^{2k}, 2^{2k+1} \right]\cap\Z$ such that
\begin{equation} \label{20230225eq54}
 2^{-\nu k} \int_{3 I_{\widetilde{p}}^{0, 3k}} \left|f_j \right|^2\lesssim\int_{3 I_{\widetilde{p}}^{0, 3k}} \left| f_j^{\Omega_j^{3 L_3 \left(\widetilde{p} \right)}} \right|^2  \lesssim \int_{3 I_{\widetilde{p}}^{0, 3k}} |f_j|^2, \quad j\in\{3, 4\}.
\end{equation}
Recall that $\omega_3^{\ell_3}=\left[\ell_3 2^{\frac{m}{2}+5k}, \left(\ell_3+1 \right) 2^{\frac{m}{2}+5k} \right]$ and $\Omega_3^{L_3}=\left[L_3 2^{m+k}, \left(L_3+1 \right) 2^{m+k} \right]$, and hence if $\omega_3^{\ell_3} \subseteq \Omega_3^{L_3}$, one has
\begin{equation} \label{20230225eq32}
\ell_3 \in \left[L_3 2^{\frac{m}{2}-4k}, \left(L_3+1 \right) 2^{\frac{m}{2}-4k} \right].
\end{equation}
Now for each $p \sim 2^{\frac{m}{2}+2k}$ and $q \sim 2^{\frac{m}{2}}$, we set
\begin{equation*}
\TFC_{p, q}:=\left\{ \left(\ell_1, \ell_2, \ell_3 \right) \in \TFC_q: \ell_3 \in \left[L_3 \left(  \left\lfloor \frac{p}{2^{\frac{m}{2}}} \right\rfloor\right) 2^{\frac{m}{2}-4k}, \left(L_3 \left(\left\lfloor \frac{p}{2^{\frac{m}{2}}} \right\rfloor\right)+1 \right) 2^{\frac{m}{2}-4k} \right] \right\}.
\end{equation*}
Hence, by \eqref{20230224main01}, we have
\begin{eqnarray*}
|\Lambda_{k,m}^{\bbU_\nu^2}(f_1,f_2,f_3,f_4)|\lesssim 2^k \cdot 2^{\nu k} \sum_{\substack{p \sim 2^{\frac{m}{2}+2k}\\ q \sim 2^{\frac{m}{2}} \\ n_3 \sim 2^{2k}}} \quad\sum_{(\ell_1, \ell_2, \ell_3) \in \TFC_{p, q}}\Big| \Big\langle f_3^{\omega_3^{\ell_3}}, \phi_{p-\frac{q^3}{2^m}, \ell_3 2^{2k}+n_3}^{3k} \Big\rangle \Big|  \Big| \Big\langle f_4, \phi_{p, \ell_3 2^k+n_3+\ell_1+\ell_2}^{3k} \Big\rangle \Big| \nonumber \\
\cdot \frac{1}{\left|I_p^{m, 3k} \right|} \bigg| \int_{I_p^{m, 3k}} \int_{I_q^{m, k}}  \underline{f}_1^{\omega_1^{\ell_1}}(x-t) \underline{f}_2^{\omega_2^{\ell_2}}(x+t+t^3) e^{-i 2^{\frac{m}{2}+k}t \cdot \frac{3q^2}{2^m}n_3}\, e^{i 2^{\frac{m}{2}+3k}t \cdot \left(-\ell_1 +\ell_2 - \left(\ell_3-\frac{\ell_2}{2^{2k}} \right)\frac{3q^2}{2^m}\right)} dt dx \bigg|.
\end{eqnarray*}
Arguing as in \eqref{20230224eq01}, \eqref{20230225eq01}, \eqref{20230225eq03}, and \eqref{20230225eq04}, we have
\begin{equation} \label{20250225eq30}
|\Lambda_{k,m}^{\bbU_\nu^2}(f_1,f_2,f_3,f_4)| \lesssim 2^{2k} \cdot 2^{\nu k} \sum_{\widetilde{p} \sim 2^{2k}} \Big\|f_3^{3 \Omega_3^{L_3(\widetilde{p})}} \Big\|_{L^2 \left(3 I_{\widetilde{p}}^{0, 3k} \right)} \Big\|f_4^{3 \Omega_4^{L_3(\widetilde{p})}} \Big\|_{L^2 \left(3 I_{\widetilde{p}}^{0, 3k} \right)} \calR_{\widetilde{p}},
\end{equation}
where
$$
\calR_{\widetilde{p}}^2:=\sum_{q \sim 2^{\frac{m}{2}}} \sum_{\left(\ell_1, \ell_2, \ell_3 \right) \in \TFC_{2^{\frac{m}{2}} \widetilde{p}, q}}  \int_{I_{\widetilde{p}}^{0, 3k}} \int_{I_q^{m, k}} \left| \underline{f}_1^{\omega_1^{\ell_1}}(x-t) \underline{f}_2^{\omega_2^{\ell_2}}(x+t+t^3) \right|^2 dxdt.
$$
Now we estimate the term $\calR_{\widetilde{p}}$. Recall that $\left(\ell_1, \ell_2, \ell_3 \right) \in \TFC_{2^{\frac{m}{2}} \widetilde{p}, q}$ means
\begin{equation} \label{20230225eq31}
\begin{cases}
\ell_1-\ell_2+\frac{3q^2}{2^m} \left(\ell_3-\frac{\ell_2}{2^{2k}}\right)=O(1); \\
\\
\ell_3 \in \left[L_3\left(\widetilde{p} \right)2^{\frac{m}{2}-4k}, \left(L_3\left(\widetilde{p} \right)+1 \right) 2^{\frac{m}{2}-4k} \right].
\end{cases}
\end{equation}
Now if $\omega_j^{\ell_j} \subset \Omega_j^{L_j}$ or equivalently $\ell_j \in \left[L_j2^{\frac{m}{2}-4k}, \left(L_j+1 \right)2^{\frac{m}{2}-4k} \right]$ then by \eqref{20230225eq31}, we have
$$
L_1 - L_2 +3 \Big(\frac{\Big\lfloor\frac{q}{2^{\frac{m}{2}-2k}}\Big\rfloor}{2^{2k}}\Big)^2 \big(L_3 \left(\widetilde{p} \right)-\frac{L_2}{2^{2k}}\big)=O\left(1 \right).
$$
This motivates the introduction of the time-frequency correlations relative to $\widetilde{p}, \widetilde{q} \sim 2^{2k}$ as
$$
\overline{\TFC}_{\widetilde{p}, \widetilde{q}}:=\left\{ \substack{\left(L_1, L_2 \right)\\ L_1, L_2 \sim 2^{2k}} : L_1 - L_2 +\frac{3{\widetilde{q}}^2}{2^{4k}}L_3 \left(\widetilde{p} \right)=O\left(1 \right) \right\}.
$$
Then, by Parseval, we have that
\begin{eqnarray*}
&& \calR_{\widetilde{p}}^2  \lesssim \sum_{\substack{\widetilde{q}\sim 2^{2k} \\ (L_1, L_2) \in \overline{\TFC}_{\widetilde{p}, \widetilde{q}}}} \sum_{q=2^{{m \over 2}-2k} \widetilde q}^{2^{{m \over 2}-2k} (\widetilde q+1)} \sum_{\substack{ \ell_1 \in \left[L_12^{\frac{m}{2}-4k}, \left(L_1+1 \right)2^{\frac{m}{2}-4k} \right] \\ \ell_2 \in \left[L_22^{\frac{m}{2}-4k}, \left(L_2+1 \right)2^{\frac{m}{2}-4k} \right] \\  \ell_3 \in \left[L_3\left(\widetilde{p} \right)2^{\frac{m}{2}-4k}, \left(L_3\left(\widetilde{p} \right)+1 \right) 2^{\frac{m}{2}-4k} \right] \\ \ell_1-\ell_2+\frac{3q^2}{2^m} \left(\ell_3-\frac{\ell_2}{2^{2k}}\right)=O(1)}} \int_{I_{\widetilde{p}}^{0, 3k}} \int_{I_{q}^{m, k}} \left| \underline{f}_1^{\omega_1^{\ell_1}}(x-t) \underline{f}_2^{\omega_2^{\ell_2}}(x+t+t^3) \right|^2 dtdx \\
&& \quad \lesssim  \sum_{\widetilde{q}\sim 2^{2k}} \sum_{(L_1, L_2) \in \overline{\TFC}_{\widetilde{p}, \widetilde{q}}} \int_{I_{\widetilde{p}}^{0, 3k}}\int_{I_{\widetilde{q}}^{0, 3k}} \Big|f_1^{\Omega_1^{L_1}}(x-t) f_2^{\Omega_2^{L_2}} \left(x+t+t^3 \right) \Big|^2 dtdx,
\end{eqnarray*}
where in the above estimate, we used the fact that if $q, \ell_1$ and $\ell_2$ are chosen, then $\ell_3$ is (up to $O(1)$ factors) uniquely determined. It is now easy to see that the above relation can be rewritten as
\begin{equation} \label{20230225eq47}
\calR_{\widetilde{p}}^2 \lesssim \sum_{L_2,\widetilde{q} \sim 2^{2k}} \Bigg( \int_{3 I_{\widetilde{p}-\widetilde{q}}^{0, 3k}} \Big|f_1^{\Omega_1^{L_2-\frac{3\widetilde{q}^2}{2^{4k}} L_3\left(\widetilde{p}\right)}}\Big|^2 dy \Bigg) \Bigg(\int_{3 I_{\widetilde{p}+\widetilde{q}}^{0, 3k}} \Big| f_2^{\Omega_2^{L_2}} \Big|^2 dy \Bigg).
\end{equation}
Recalling now that we are in the regime in which at least one of the functions $f_1$, $f_2$ is $2$-uniform. Without loss of generality, we may assume $f_1=f_1^{\bbU^2_\nu}$. Then, by \eqref{20230225eq47} and Parseval, we have
\begin{equation}
\calR_{\widetilde{p}}^2 \lesssim 2^{-3\nu k} \sum_{\widetilde{q} \sim 2^k } \Big( \int_{3 I_{\widetilde{p}-\widetilde{q}}^{0, 3k}} \left |f_1 \right|^2 \Big) \Big(\int_{3 I_{\widetilde{p}+\widetilde{q}}^{0, 3k}} \left| f_2 \right|^2 \Big).
\end{equation}
This together with \eqref{20250225eq30} gives
\begin{align} 
&|\Lambda_{k,m}^{\bbU_\nu^2}(f_1,f_2,f_3,f_4)| \nonumber \\
& \quad \lesssim 2^{2k} \cdot 2^{-\frac{\nu k}{2}} \sum_{\widetilde{p} \sim 2^{2k}} \Big\|f_3^{3 \Omega_3^{L_3(\widetilde{p})}} \Big\|_{L^2 \left(3 I_{\widetilde{p}}^{0, 3k} \right)} \Big\|f_4^{3 \Omega_4^{L_3(\widetilde{p})}} \Big\|_{L^2 \left(3 I_{\widetilde{p}}^{0, 3k} \right)}\, \Bigg[  \sum_{\widetilde{q} \sim 2^{2k} } \Big( \int_{3 I_{\widetilde{p}-\widetilde{q}}^{0, 3k}} \left |f_1 \right|^2 \Big) \Big(\int_{3 I_{\widetilde{p}+\widetilde{q}}^{0, 3k}} \left| f_2 \right|^2 \Big)  \Bigg]^{\frac{1}{2}}. \label{20230305eq40}
\end{align}
Finally, applying a last Cauchy-Schwarz in $\widetilde{p}$ followed by information that $f_3$ is $1-$uniform, we have
\begin{align} 
\left|\Lambda_{k,m}^{\bbU_\nu^2}(f_1,f_2,f_3,f_4)\right|&\lesssim 2^{2k} \cdot 2^{-\frac{\nu k}{2}}  \left\|f_1 \right\|_{L^2 \left(3 I^k \right)} \left\|f_2 \right\|_{L^2 \left(3 I^k \right)} \Big( \sum_{\widetilde{p} \sim 2^{2k}} \left\|f_3 \right\|^2_{L^2 \left(3 I_{\widetilde{p}}^{0, 3k} \right)} \left\|f_4 \right\|^2_{L^2 \left(3 I_{\widetilde{p}}^{0, 3k} \right)} \Big)^{\frac{1}{2}} \nonumber \\
&\lesssim 2^{k} \cdot 2^{-\frac{(\nu-\mu) k}{2}}\,\prod_{s=1}^4 \|f_s\|_{L^2(3 I^k)}\,. \label{Uniform-II-Small}
\end{align}
This estimate is acceptable in view of proving Theorem \ref{sgscale1unifksmal}, by setting $\delta \ll \nu$ and $\mu=\frac{\nu}{2}$.

\subsection{Treatment of the time-frequency sparse component $\Lambda_{k,m}^{\bbS_\nu^2}$}~\\
Recall that in this case, we have that all $\left\{f_1, f_2, f_3, f_4\right\}$ are time-frequency sparse. Therefore, we see that
\begin{enumerate}
    \item [$\bullet$] up to a $2^{\nu k}$ factor, there is a unique measurable function $L_3: \left[2^{2k}, 2^{2k+1} \right]\cap\Z \mapsto \left[2^{2k}, 2^{2k+1}\right]\cap\Z$ such that\footnote{Note that the frequency localization of $f_3$ and $f_4$ are essentially determined by the same frequency interval $\left[L_3 2^{m+2k}, (L_3+1)2^{m+2k} \right]$ -- see \eqref{20230225eq01}.} for each $p \sim 2^{2k}$,
\begin{equation} \label{20230226eq01}
2^{-\nu k}\,\int_{3 I_{\widetilde{p}}^{0, 3k}} |f_j|^2\lesssim\int_{3 I_{\widetilde{p}}^{0, 3k}} \left| f_j^{3 \Omega_j^{L_3 \left(\widetilde{p} \right)}} \right|^2  \lesssim \int_{3 I_{\widetilde{p}}^{0, 3k}} |f_j|^2, \quad j\in\{3, 4\}.
\end{equation}

\item [$\bullet$] up to a $2^{3\nu k}$ factor and for given $j \in \{1, 2\}$, there exists a unique measurable function $L_j: \left[2^{2k}, 2^{2k+1} \right]\cap\Z \mapsto \left[2^{2k}, 2^{2k+1}\right]\cap\Z$ such that for each $q_j \sim 2^{2k}$,
 \begin{equation} \label{20230226eq02}
2^{-3\nu k}\,\int_{3 I_{q_j}^{0, 3k}} |f_j|^2\lesssim \int_{3 I_{q_j}^{0, 3 k}} \left|f_j^{\Omega_j^{L_j(q_j)}} \right|^2 \lesssim \int_{3 I_{q_j}^{0, 3 k}} |f_j|^2, \quad j\in\{1, 2\}.
 \end{equation}
\end{enumerate}

Notice from \eqref{20230226eq01} and \eqref{20230226eq02} that the frequency information of the functions $\{f_1, f_2, f_3, f_4\}$ is localized in a ``compressed'' interval within each spatial fiber. This will be later exploited via time-frequency correlations (that are encoding the curvature of our problem) in order to obtain the desired $k-$decaying estimate.

\subsubsection{Reductions via frequency sparsity of the spatial fibers} Starting from \eqref{20230224main01} and in the spirit provided by the approach in Subsection \ref{subsec4.2.2}, we first exploit the frequency localization as follows:
\begin{equation} \label{20230226eq10}
|\Lambda_{k,m}^{\bbS_\nu^2}(\vec{f})| \lesssim 2^k \cdot 2^{7\nu k} \sum_{\substack{ \widetilde{p},\,\widetilde{q} \sim 2^{2k}}} \bbA_{\widetilde{p}, \widetilde{q}}\:\bbB_{\widetilde{p}, \widetilde{q}},
\end{equation}
where
\begin{equation}
\bbA_{\widetilde{p}, \widetilde{q}}^2 :=\sum_{\substack{p \in \left[\widetilde{p}2^{\frac{m}{2}}, \left(\widetilde{p}+1 \right)2^{\frac{m}{2}} \right] \\ q \in \left[ \widetilde{q} 2^{\frac{m}{2}-2k}, \left( \widetilde{q}+1 \right)2^{\frac{m}{2}-2k} \right] \\ n_3 \sim 2^{2k}}} \sum_{\substack{(\ell_1, \ell_2, \ell_3) \in \TFC_{p, q} \\\eqref{20230226eq01}\,\&\, \eqref{20230226eq02}}} \Big| \Big\langle f_3^{\omega_3^{\ell_3}}, \varphi_{p-\frac{q^3}{2^m}, \ell_3 2^{2k}+n_3}^{3k} \Big \rangle \Big|^2  \Big| \Big \langle f_4,\, \varphi_{p, \ell_3 2^k+n_3+\ell_1+\ell_2}^{3k} \Big\rangle \Big|^2,
\end{equation}
and
\begin{eqnarray*}
&& \bbB_{\widetilde{p}, \widetilde{q}}^2:=\sum_{\substack{p \in \left[\widetilde{p}2^{\frac{m}{2}}, \left(\widetilde{p}+1 \right)2^{\frac{m}{2}} \right] \\ q \in \left[ \widetilde{q} 2^{\frac{m}{2}-2k}, \left( \widetilde{q}+1 \right)2^{\frac{m}{2}-2k} \right] \\ n_3 \sim 2^{2k}}} \sum_{\substack{(\ell_1, \ell_2, \ell_3) \in \TFC_{p, q} \\\eqref{20230226eq01}\,\&\, \eqref{20230226eq02}}}\\
&& \quad \quad \quad  \quad \frac{1}{\big|I_p^{m, 3k} \big|^2} \bigg| \int_{I_p^{m, 3k}} \int_{I_q^{m, k}}  \underline{f}_1^{\omega_1^{\ell_1}}(x-t) \underline{f}_2^{\omega_2^{\ell_2}}(x+t+t^3) e^{-i 2^{\frac{m}{2}+k}t \cdot \frac{3q^2}{2^m}n_3} \, e^{i 2^{\frac{m}{2}+3k}t \cdot \left(-\ell_1 +\ell_2 - \left(\ell_3-\frac{\ell_2}{2^{2k}} \right)\frac{3q^2}{2^m}\right)} dt dx \bigg|^2.
\end{eqnarray*}

\noindent\textit{Treatment of $\bbA_{\widetilde{p}, \widetilde{q}}$.} The estimate of this part is similar to the treatment of the term $\bbA_p$ in \eqref{20230225eq01}:
$$
\bbA_{\widetilde{p}, \widetilde{q}}^2 \lesssim \sum_{\substack{p \in \left[\widetilde{p}2^{\frac{m}{2}}, \left(\widetilde{p}+1 \right)2^{\frac{m}{2}} \right] \\ q \in \left[ \widetilde{q} 2^{\frac{m}{2}-2k}, \left( \widetilde{q}+1 \right)2^{\frac{m}{2}-2k} \right] \\ \ell_3 \in \left[L_3 \left(\widetilde{p}\right)2^{m-4k}, \left(L_3\left(\widetilde{p} \right)+1 \right)2^{m-4k}\right]}}\Big\|f_3^{3\omega_3^{\ell_3}} \Big\|^2_{L^2 \big(I_{p-\frac{q^3}{2^m}}^{m, 3k} \big)} \Big\|f_4^{2^{\frac{m}{2}+5k}\ell_3+\left[0, 2^{m+k+1} \right]} \Big\|_{L^2 \big(I_p^{m, 3k} \big)}^2.
$$
Applying Parseval in $\ell_3$ first, then in $q$, followed by the partition of $\left[\widetilde{p}2^{\frac{m}{2}}, \left(\widetilde{p}+1 \right)2^{\frac{m}{2}} \right]$ into smaller intervals of length $2^{\frac{m}{2}-2k}$, we have
\begin{eqnarray}  \label{20230306eq01}
\bbA_{\widetilde{p}, \widetilde{q}}^2 &\lesssim&   \sum_{\substack{p \in \left[\widetilde{p}2^{\frac{m}{2}}, \left(\widetilde{p}+1 \right)2^{\frac{m}{2}} \right] \\ q \in \left[ \widetilde{q} 2^{\frac{m}{2}-2k}, \left( \widetilde{q}+1 \right)2^{\frac{m}{2}-2k} \right]}} \Big\|f_3^{3 \Omega_3^{L_3 (\widetilde{p})}} \Big\|^2_{L^2 \big(I_{p-\frac{q^3}{2^m}}^{m, 3k} \big)} \Big\|f_4^{3\Omega_4^{L_3(\widetilde{p})}} \Big\|_{L^2 \big(I_p^{m, 3k} \big)}^2 \\
&\lesssim& \sum_{s=0}^{2^{2k}-1} \Big\|f_3^{3 \Omega_3^{L_3(\widetilde{p})}} \Big\|^2_{L^2 \big(3 I_{2^{2k} \widetilde{p}-s-\frac{\widetilde{q}^3}{2^{4k}}}^{0, 5k} \big)} \Big\|f_4^{3\Omega_4^{L_3(\widetilde{p})}} \Big\|_{L^2 \big(3 I_{2^{2k} \widetilde{p}+s}^{0, 5k} \big)}^2. \nonumber
\end{eqnarray}

\noindent\textit{Treatment of $\bbB_{\widetilde{p}, \widetilde{q}}$.} The estimate for $\bbB_{\widetilde{p}, \widetilde{q}}$ follows from similar arguments as in \eqref{20230225eq03}, \eqref{20230225eq04}, \eqref{20250225eq30} and \eqref{20230225eq47}:
$$
 \bbB_{\widetilde{p}, \widetilde{q}}^2 \lesssim 2^{2k} \sum_{
 \substack{L_2 \sim 2^{2k} \\ \eqref{20230226eq02}}}  \Bigg( \int_{3 I_{\widetilde{p}-\widetilde{q}}^{0, 3k}} \Big|f_1^{\Omega_1^{L_2+ \frac{3\widetilde{q}^2}{2^{4k}} L_3\left(\widetilde{p}\right)}}\Big|^2  \Bigg) \Bigg(\int_{3 I_{\widetilde{p}+\widetilde{q}}^{0, 3k}} \Big| f_2^{\Omega_2^{L_2}} \Big|^2 \Bigg).
$$

\medskip

Combining now the estimates of $\bbA_{\widetilde{p}, \widetilde{q}}$ and $\bbB_{\widetilde{p}, \widetilde{q}}$ with \eqref{20230226eq10}, we have
\begin{align} 
 |\Lambda^{k, \bbS^2_\nu}_m(f_1,f_2,f_3,f_4)|  & \lesssim 2^{2k} \cdot 2^{7\nu k} \sum_{\substack{\widetilde{p} \sim 2^{2k} \\ \widetilde{q} \sim 2^{2k}}}\Bigg[ \sum_{s=0}^{2^{2k}-1} \Big\|f_3^{3 \Omega_3^{L_3(\widetilde{p})}} \Big\|^2_{L^2 \big(3 I_{2^{2k} \widetilde{p}-s-\frac{\widetilde{q}^3}{2^{4k}}}^{0, 5k} \big)} \Big\|f_4^{3\Omega_4^{L_3(\widetilde{p})}} \Big\|_{L^2 \big(3 I_{2^{2k} \widetilde{p}+s}^{0, 5k} \big)}^2 \Bigg]^{\frac{1}{2}}  \nonumber  \\
 &  \quad \cdot \Bigg[ \sum_{
 \substack{L_2 \sim 2^{2k} \\ \eqref{20230226eq02}}}  \Bigg( \int_{3 I_{\widetilde{p}-\widetilde{q}}^{0, 3k}} \Big|f_1^{\Omega_1^{L_2+ \frac{3\widetilde{q}^2}{2^{4k}} L_3\left(\widetilde{p}\right)}}\Big|^2  \Big) \Big(\int_{3 I_{\widetilde{p}+\widetilde{q}}^{0, 3k}} \Big| f_2^{\Omega_2^{L_2}} \Big|^2 \Big)\Bigg]^{\frac{1}{2}}.  \label{20230226eq11}
\end{align}
Motivated by the time-frequency localization induced by the above estimate, we define the time-frequency correlation set
\begin{equation}
\sTFC:=\left\{\substack{\left( \widetilde{p}, \widetilde{q} \right)\\ \widetilde{p},\widetilde{q} \sim 2^{2k}}: L_2 \left(\widetilde{p}+\widetilde{q}\right)+\frac{3\widetilde{q}^2}{2^{4k}} L_3 \left( \widetilde{p} \right)=L_1 \left(\widetilde{p}-\widetilde{q} \right)+O(1) \right\}.
\label{eq:TFC-end}
\end{equation}
Therefore, by Cauchy-Schwarz, we deduce that
\begin{align}
|\Lambda_{k,m}^{\bbS_\nu^2}(f_1,f_2,f_3,f_4)| & \lesssim 2^{(2+7\nu) k} \cdot \left( \# \sTFC \right)^{\frac{1}{2}} \nonumber \\
& \quad \cdot \Bigg[ \sum_{\substack{ \widetilde{p} \sim 2^{2k} \\ \widetilde{q} \sim 2^{2k} }}\Big\|f_1 \Big\|^2_{L^2 \big(3 I_{\widetilde{p}-\widetilde{q}}^{0, 3k} \big)} \Big\|f_2 \Big\|^2_{L^2 \big(3 I_{\widetilde{p}+\widetilde{q}}^{0, 3k}\big)}\Bigg(\sum_{s=0}^{2^{2k}-1} \Big\|f_3 \Big\|^2_{L^2 \big(3 I_{2^{2k} \widetilde{p}-s-\frac{\widetilde{q}^3}{2^{4k}}}^{0, 5k} \big)} \Big\|f_4 \Big\|_{L^2 \big(3 I_{2^{2k} \widetilde{p}+s}^{0, 5k} \big)}^2 \Bigg) \Bigg]^{\frac{1}{2}}\,. \label{Uniform-I}
\end{align}
Recalling now that both $f_1$ and $f_3$ must be also $1$-uniform, from \eqref{spatialsparsefs12}--\eqref{sparse23} we conclude that
\begin{equation}  \label{Uniform-Ifinal}
|\Lambda_{k,m}^{\bbS_\nu^2}(f_1,f_2,f_3,f_4)|\lesssim 2^{-k}\,2^{(7\nu+\mu) k}\,\left( \# \sTFC \right)^{\frac{1}{2}} \prod_{s=1}^4 \|f_s\|_{L^2(3 I^k)}\,.
\end{equation}

Our desired estimate follows now if we can improve over the trivial estimate $\# \sTFC\lesssim 2^{4k}$. We are thus going to prove that $\# \sTFC\lesssim 2^{(4-\tilde \delta)k}$ (see below \eqref{gain}, for some numerical constant $\tilde \delta$), which then allows to conclude the proof of Theorem \ref{sgscale1unifksmal}, by setting $\delta \ll \tilde \delta -(7\nu+\mu)$ and choosing $\mu=\frac{\nu}{2} \lll \tilde \delta$.

\subsubsection{Improved control over the time-frequency correlation set: a level set analysis}\label{Levelsetanalysis}

As a matter of convenience we start by redenoting $a:=L_1$, $b:=L_2$, $c:=\frac{3}{4}L_3$, $u:=\widetilde{p}-\widetilde{q}$ and $v:=\widetilde{p}+\widetilde{q}$. 

Thus, in this new language, we have
$$
\sTFC:=\left\{\substack{\left(u, v \right)\\ u,v \sim 2^{2k}}: a(u)-b(v)=\frac{(u-v)^2}{2^{4k}}c(u+v)+O(1) \right\}.
$$
Assume now by contradiction, that for any ${\tilde \delta}>0$ small enough, we have that
\begin{equation} \label{contrad}
\#\sTFC\gtrsim_{\delta} 2^{(4-{\tilde \delta})k}\,.
\end{equation}
Then, this translates into\footnote{Here $\ds \lceil x \rceil :=(1+x^2)^{-1}$.}
$$
\sum_{u,v\sim 2^{2k}} \left\lceil a(u)-b(v)-\frac{(u-v)^2}{2^{4k}}c(u+v) \right\rceil \gtrsim 2^{(4-{\tilde \delta}) k}\,.
$$
Applying now Cauchy-Schwarz inequality we have
$$
\sum_{v\sim 2^{2k}} \sum_{u_1,u_2\sim 2^{2k}}\left\lceil a(u_1)-a(u_2)-\frac{(u_1-v)^2}{2^{4k}}c(u_1+v)+\frac{(u_2-v)^2}{2^{4k}}c(u_2+v) \right\rceil \gtrsim 2^{(6-2{\tilde \delta}) k}\,.
$$
We are next applying the change of variable $v\,\rightarrow\,v-u_2$ and $u_1\,\rightarrow\,u_1+u_2$ and Fubini in order to get
$$
\sum_{u_1,u_2\sim 2^{2k}}\sum_{v\sim 2^{2k}} \left\lceil a(u_1+u_2)-a(u_2)-\frac{(u_1+2 u_2-v)^2}{2^{4k}}c(u_1+v)+\frac{(2u_2-v)^2}{2^{4k}}c(v) \right\rceil \gtrsim 2^{(6-2{\tilde \delta}) k}\,.
$$
Applying now another Cauchy-Schwarz we further deduce
\begin{align}
\sum_{u_1,u_2\sim 2^{2k}}\sum_{v_1, v_2\sim 2^{2k}} &\left\lceil \frac{(u_1+2 u_2-v_1)^2}{2^{4k}}c(u_1+v_1)-\frac{(2u_2-v_1)^2}{2^{4k}}c(v_1)\right.\nonumber\\
&\left.-\frac{(u_1+2 u_2-v_2)^2}{2^{4k}}c(u_1+v_2)+\frac{(2u_2-v_2)^2}{2^{4k}}c(v_2)\right\rceil \gtrsim 2^{(8-4{\tilde \delta}) k}\,.\nonumber
\end{align}
Regrouping now the information above in terms of the ``smooth variable'' $u_2$ we obtain
\begin{align}
\sum_{u_1,v_1, v_2\sim 2^{2k}}\sum_{u_2\sim 2^{2k}} \Big\lceil \, &\frac{4 u_2^2}{2^{4k}}
\left(c(u_1+v_1)-c(v_1)-c(u_1+v_2)+c(v_2)\right)\nonumber\\
\quad + &\frac{4 u_2}{2^{4k}}\left((u_1-v_1)c(u_1+v_1)+v_1c(v_1)-(u_1-v_2)c(u_1+v_2)-v_2 c(v_2)\right)\nonumber\\
\quad + &\Big(\frac{(u_1-v_1)^2}{2^{4k}}c(u_1+v_1)-\frac{v_1^2}{2^{4k}}c(v_1)-
\frac{(u_1-v_2)^2}{2^{4k}}c(u_1+v_2)+\frac{v_2^2}{2^{4k}}c(v_2)\Big)\Big\rceil \gtrsim 2^{(8-4{\tilde \delta}) k}\,.\nonumber
\end{align}
Further, applying another Cauchy-Schwarz, we have
\begin{align}
\sum_{u_1,v_1, v_2\sim 2^{2k}}\sum_{u_2,u_3\sim 2^{2k}} &\left\lceil \frac{4 u_2^2}{2^{4k}}
\left(c(u_1+v_1)-c(v_1)-c(u_1+v_2)+c(v_2)\right)\right.\nonumber\\
&\frac{4 u_2}{2^{4k}}\left((u_1-v_1)c(u_1+v_1)+v_1c(v_1)-(u_1-v_2)c(u_1+v_2)-v_2 c(v_2)\right)\nonumber\\
&\left.+\left(\frac{(u_1-v_1)^2}{2^{4k}}c(u_1+v_1)-\frac{v_1^2}{2^{4k}}c(v_1)-
\frac{(u_1-v_2)^2}{2^{4k}}c(u_1+v_2)+\frac{v_2^2}{2^{4k}}c(v_2)\right)\right\rceil\nonumber\\
&\left\lceil \frac{u_3}{2^{2k}}\left(c(u_1+v_1)-c(v_1)-c(u_1+v_2)+c(v_2)\right)\right.\nonumber\\
&\left.+\frac{1}{2^{2k}}\left((u_1-v_1)c(u_1+v_1)+v_1c(v_1)-(u_1-v_2)c(u_1+v_2)-v_2 c(v_2)\right)\right\rceil \gtrsim 2^{(10-16{\tilde \delta}) k}\nonumber\,.
\end{align}
Above we have used that the main contribution comes from the case then $|u_2-u_3| \gtrsim  2^{(2-16{\tilde \delta}) k}$. After successive Cauchy-Schwarz applications, the above implies
\begin{align}
\sum_{u_1,u_2, v_1, v_2\sim 2^{2k}} \Big\lceil \, & \frac{4 u_2^2}{2^{4k}}
\left(c(u_1+v_1)-c(v_1)-c(u_1+v_2)+c(v_2)\right)\nonumber\\
+ & \frac{4 u_2}{2^{4k}}\left((u_1-v_1)c(u_1+v_1)+v_1c(v_1)-(u_1-v_2)c(u_1+v_2)-v_2 c(v_2)\right)\nonumber\\
+ & \Big(\frac{(u_1-v_1)^2}{2^{4k}}c(u_1+v_1)-\frac{v_1^2}{2^{4k}}c(v_1)-
\frac{(u_1-v_2)^2}{2^{4k}}c(u_1+v_2)+\frac{v_2^2}{2^{4k}}c(v_2)\Big)\Big\rceil^4\nonumber\\
&\Big\lceil c(u_1+v_1)-c(v_1)-c(u_1+v_2)+c(v_2)\Big\rceil\nonumber\\
&\Big\lceil\frac{1}{2^{2k}}\Big((u_1-v_1)c(u_1+v_1)+v_1c(v_1)-(u_1-v_2)c(u_1+v_2)-v_2 c(v_2)\Big)\Big\rceil \gtrsim 2^{(8-64{\tilde \delta}) k}\label{eq:*1}\,.
\end{align}
This further gives
\begin{align}
\sum_{u_1, v_1, v_2\sim 2^{2k}} &\left\lceil\frac{(u_1-v_1)^2}{2^{4k}}c(u_1+v_1)-\frac{v_1^2}{2^{4k}}c(v_1)-
\frac{(u_1-v_2)^2}{2^{4k}}c(u_1+v_2)+\frac{v_2^2}{2^{4k}}c(v_2)\right\rceil\nonumber\\
&\left\lceil c(u_1+v_1)-c(v_1)-c(u_1+v_2)+c(v_2)\right\rceil\nonumber\\
&\left\lceil\frac{1}{2^{2k}}\left((u_1-v_1)c(u_1+v_1)+v_1c(v_1)-(u_1-v_2)c(u_1+v_2)-v_2 c(v_2)\right)\right\rceil \gtrsim 2^{(6-128{\tilde \delta}) k}\nonumber\,.
\end{align}
After another Cauchy-Schwarz this implies
\begin{align}
\sum_{u_1, v_1, v_2\sim 2^{2k}} &\left\lceil\frac{(u_1-v_1)^2}{2^{4k}}(c(v_1)+c(u_1+v_2)-c(v_2))-\frac{v_1^2}{2^{4k}}c(v_1)-
\frac{(u_1-v_2)^2}{2^{4k}}c(u_1+v_2)+\frac{v_2^2}{2^{4k}}c(v_2)\right\rceil\nonumber\\
&\left\lceil\frac{1}{2^{2k}}\left((u_1-v_1)(c(v_1)+c(u_1+v_2)-c(v_2))+v_1c(v_1)-(u_1-v_2)c(u_1+v_2)-v_2 c(v_2)\right)\right\rceil \gtrsim 2^{(6-256{\tilde \delta}) k}\nonumber\,,
\end{align}
and hence, after regrouping the terms
\begin{align}
\sum_{u_1, v_1, v_2\sim 2^{2k}} &\left\lceil\frac{u_1(u_1-2v_1)}{2^{4k}}c(v_1)+\frac{(v_2-v_1)(2u_1-v_1-v_2)}{2^{4k}}c(u_1+v_2)+\frac{(u_1-v_1+v_2)
(v_2+v_1-u_1)}{2^{4k}}c(v_2)\right\rceil\nonumber\\
&\left\lceil\frac{1}{2^{2k}}\left(u_1 c(v_1)+(v_2-v_1)c(u_1+v_2)+(v_1-v_2-u_1)c(v_2)\right)\right\rceil \gtrsim 2^{(6-256{\tilde \delta}) k}\nonumber\,.
\end{align}
This gives us
\begin{align}
\sum_{u_1, v_1, v_2\sim 2^{2k}}
&\left\lceil\frac{u_1(u_1+v_1-v_2)}{2^{4k}}\,c(v_1)-\frac{u_1(u_1+v_2-v_1)}{2^{4k}}c(v_2)\right\rceil \gtrsim 2^{(6-256{\tilde \delta}) k}\label{eq:*2}\,
\end{align}
which can be rewritten as
\begin{align}
\sum_{u_1, v_1, v_2\sim 2^{2k}}
&\left\lceil\frac{u_1}{2^{2k}}\,(c(v_1)-c(v_2))+\frac{v_1-v_2}{2^{2k}}\,(c(v_1)+c(v_2))\right\rceil \gtrsim 2^{(6-256{\tilde \delta}) k}\nonumber\,.
\end{align}
After a final Cauchy-Schwarz we get
\begin{align}
\sum_{v_1, v_2\sim 2^{2k}}
\left\lceil c(v_1)-c(v_2)\right\rceil\,\left\lceil\frac{v_1-v_2}{2^{2k}}\,(c(v_1)+c(v_2))\right\rceil \gtrsim 2^{(4-1024{\tilde \delta}) k}\nonumber\,,
\end{align}
which implies
\begin{align}\label{Keydec}
k2^{2k}\gtrsim \sum_{v_1, v_2\sim 2^{2k}}
\left\lceil\frac{v_1-v_2}{2^{2k}}\,c(v_1)\right\rceil \gtrsim 2^{(4-1024{\tilde \delta}) k}\,,
\end{align}
where for the first inequality we used that $|c(v_1)|\sim 2^{2k}$.

Now it is obvious that the claim \eqref{contrad} contradicts \eqref{Keydec} if $\tilde \delta$ is chosen small enough (\textit{e.g.} ${\tilde \delta}<\frac{1}{600}$).

We thus conclude the desired exponential gain:
\begin{equation} \label{gain}
\#\sTFC\lesssim 2^{(4-{\tilde \delta})k}\,,
\end{equation}
for ${\tilde \delta}\approx\frac{1}{600}$.

\section{The single scale estimate (III): $0<\frac{m}{8}<k<\frac{m}{4}$}  \label{Diagkinterm}
This section is devoted to the intermediate case $\frac{m}{8} < k < \frac{m}{4}$. The strategy we employ is similar to the one used in the previous case (II), but now we are looking for a decay in the parameter $\kappa:=\frac{m}{2}-2k$ (instead of $k$ in case (II)).

\subsection{Single-scale estimate}~\\
The first step follows line by line the discussion in Section \ref{spunifdich} with the only modification that in the definitions \eqref{spatialsparsefs12} and \eqref{20230225eq23} one replaces $2^{\mu k}$ by $2^{\mu (\frac{m}{4}-k)}$, for $\mu>0$ a parameter whose value will be chosen later. Thus as a consequence of the analogous result for Theorem \ref{sparsecont1} it is enough to focus on the space-uniform component, which, as before, may be assumed without loss of generality to be
\begin{equation}\label{unif2311}
\Lambda_{k,m}^{\bbU^1_\mu}(f_1,f_2,f_3,f_4):=\Lambda_{k,m} \left(f_1^{\bbU^1_\mu}, f_2, f_3^{\bbU^1_\mu}, f_4\right)\,.
\end{equation}

Our goal in this subsection is to prove the following tamer version of our key result as given by

\begin{theorem} \label{sgscale1unifkintr}[\textsf{Control over the space-uniform component}]  There exists an absolute constant $\delta>0$ so that for any $\mu>0$ sufficiently small and for any $k,m \in \N$ with $0<\frac{m}{8}<k<\frac{m}{4}$, we have 
\beq
|\Lambda_{k,m}^{\bbU^1_\mu}(f_1,f_2,f_3,f_4)|\lesssim 2^{k}\,2^{-\delta(\frac{m}{2}-2k)}\,\prod_{s=1}^4 \|f_s\|_{L^2(3I^k)}\,.
\eeq
\end{theorem}

In what follows, for notational simplicity we will refer to $\Lambda_{k,m}^{\bbU^1_\mu}$ as $\Lambda_{k,m}$ and we will set $f_1:=f_1^{\bbU^1_\mu}$ and $f_3:=f_3^{\bbU^1_\mu}$.

Next, following the reasonings from Section \ref{1unifksmall} and exploiting the fact that in the case  $\frac{m}{8}< k < \frac{m}{4}$ we have that $\frac{\ell_2}{2^{2k}}=O(2^{\frac{m}{2}-4k})=O(1)$ we can reduce matters to the study of the analogue of \eqref{20230224main01} given by
\begin{align} 
& |\Lambda_{k,m}(f_1,f_2,f_3,f_4)| \nonumber \\
& \quad \lesssim 2^k \sum_{\substack{p \sim 2^{\frac{m}{2}+2k}\\ q \sim 2^{\frac{m}{2}} \\ n \sim 2^{2k}}} \sum_{(\ell_1, \ell_2, \ell_3) \in \TFC_q}  \left| \left\langle f_3^{\omega_3^{\ell_3}}, \varphi_{p-\frac{q^3}{2^m}, \ell_3 2^{2k}+n}^{3k} \right \rangle \right|  \left| \left \langle f_4, \varphi_{p, \ell_3 2^{2k}+n+\ell_1+\ell_2}^{3k} \right \rangle \right| \nonumber \\
& \quad \cdot \frac{1}{\left|I_p^{m, 3k} \right|} \bigg| \int_{I_p^{m, 3k}} \int_{I_q^{m, k}}  \underline{f}_1^{\omega_1^{\ell_1}}(x-t) \underline{f}_2^{\omega_2^{\ell_2}}(x+t+t^3) e^{-i 2^{\frac{m}{2}+k}t \cdot \frac{3q^2}{2^m}n}  e^{i 2^{\frac{m}{2}+3k}t \cdot \left(-\ell_1 +\ell_2 -\ell_3 \frac{3q^2}{2^m}\right)} dt dx \bigg|.
\label{corelmodksmall}
\end{align}

\subsection{A time-frequency sparse--uniform dichotomy}\label{SU1kint}~\\
We start by recalling that for $j\in\{1, 2\}$ and $\ell_j \sim 2^{\frac{m}{2}-2k}$  we defined
$$ \omega_j^{\ell_j}:=\left[\ell_j 2^{\frac{m}{2}+3 k}, \left(\ell_j+1 \right) 2^{\frac{m}{2}+3 k}\right]\,,$$
and
$$\Omega_j^{L_j}:=\left[L_j 2^{m-k}, \left(L_j+1 \right) 2^{m-k} \right], \quad L_j \sim 2^{2k}\,.$$
Notice that in our current regime for $k$ we have $\left|\Omega_j^{L_j} \right| \leq \left|\omega_j^{l_j} \right|$.
Also we recall that for $j\in\{3,4\}$ we defined
$$\Omega_j^{L_j}:=\left[L_j 2^{m+k}, \left(L_j+1 \right) 2^{m+k} \right], \quad L_j \sim 2^{2k}\,.$$

With these, we perform a time-frequency sparse/uniform decomposition adapted to the above foliation. More precisely, for $\nu \in (0, 1)$ suitably small, we define the \emph{time-frequency sparse} set ($2$-sparse) of indices as follows:
\begin{itemize}
\item for $j\in\{1, 2\}$
$$
\bbS^2_\nu \left(f_j \right):=\left\{ \substack{(q_j, \ell_j)\\ q_j\sim 2^{2k}\\\ell_j\sim2^{\frac{m}{2}-2k}} : \int_{3 I_{q_j}^{0, 3k}} \left|f_j^{\omega_j^{\ell_j}} \right|^2  \gtrsim 2^{-3\nu (\frac{m}{4}-k)} \int_{3 I_{q_j}^{0, 3 k}} \left|f_j \right|^2 \right\},
$$
\item for $j\in\{3, 4\}$
\begin{equation} \label{tfsparse34bis}
\bbS^2_\nu \left(f_j \right):=\left\{ \substack{(q_j, L_j)\\ q_j, L_j \sim 2^{2k}} : \int_{3 I_{q_j}^{0, 3k}} \left|f_j^{3 \Omega_j^{L_j}} \right|^2  \gtrsim 2^{-\nu (\frac{m}{4}-k)} \int_{3 I_{q_j}^{0, 3k}} \left|f_j \right|^2 \right\}.
\end{equation}
\end{itemize}

We also define the \emph{time-frequency uniform} sets of indices as the complement of the corresponding sparse sets and define in the usual fashion the \emph{time-frequency uniform} and \emph{time-frequency sparse} components of $f_j$.

With these done, we decompose
\begin{equation}\label{Ldecksmallin}
\Lambda_{k,m}(f_1,f_2,f_3,f_4)=:\Lambda_{k,m}^{\bbU_\nu^2}(f_1,f_2,f_3,f_4)\,+\,\Lambda_{k,m}^{\bbS_\nu^2}(f_1,f_2,f_3,f_4),
\end{equation}
with $\Lambda_{k,m}^{\bbS_\nu^2}(f_1,f_2,f_3,f_4):=\Lambda_{k,m} \left(f_1^{\bbS^2_\nu}, f_2^{\bbS^2_\nu}, f_3^{\bbS^2_\nu}, f_4^{\bbS^2_\nu} \right)$.

\subsection{Treatment of the time-frequency uniform component $\Lambda_{k,m}^{\bbU_\nu^2}$}~\\
Applying Cauchy-Schwarz in \eqref{corelmodksmall}, we can write exactly as in \eqref{20230224eq01}
\begin{equation} \label{20230224eq01int}
|\Lambda_{k,m}^{\bbU_\nu^2}(\vec{f})|\lesssim 2^k \sum_{p \sim 2^{\frac{m}{2}+2k}} \bbA_p \bbB_p,
\end{equation}
where $ \bbA_p$ and $\bbB_p$ are defined by \eqref{20230305eq01} and \eqref{20230305eq02}, respectively.

The treatment of $\bbB_p$ follows the same lines as before (see \eqref{20230225eq03}):

\begin{equation} \label{20230225eq03int}
\bbB_p^2 \lesssim 2^{2k}\,\sum_{\substack{\ell_1, \ell_2 \sim 2^{\frac{m}{2}-2k} \\q \sim 2^{\frac{m}{2}}}}  \int_{I_p^{m, 3k}} \int_{I_q^{m, k}} \left|\underline{f}_1^{\omega_1^{\ell_1}}(x-t) \underline{f}_2^{\omega_2^{\ell_2}}(x+t+t^3) \right|^2 dtdx.
\end{equation}

For the term $\bbA_p$, using the fact that $2^{2k}\geq 2^{\frac{m}{2}-2k}$, we proceed as follows
\begin{equation}\label{apbis}
\bbA_p^2\lesssim \sum_{\substack{\widetilde{n_3} \sim 2^{4k-\frac{m}{2}} \\ q \sim 2^{\frac{m}{2}}}} \: \sum_{n_3\in[\widetilde{n_3}2^{\frac{m}{2}-2k},\,(\widetilde{n_3}+1)2^{\frac{m}{2}-2k}]}\:\sum_{\ell_2, \ell_3 \sim 2^{\frac{m}{2}-2k}}  \left| \left\langle f_3^{\omega_3^{\ell_3}}, \varphi_{p-\frac{q^3}{2^m}, \ell_3 2^{2k}+n_3}^{3k} \right \rangle \right|^2  \left| \left \langle f_4, \varphi_{p, \ell_3 2^{2k}+n_3+\ell_2}^{3k} \right \rangle \right|^2.
\end{equation}
Next, summing up with respect to $\ell_2$ first and then $n_3$, we see that the above term is further bounded from above by
$$
\lesssim \sum_{\substack{\widetilde{n_3} \sim 2^{4k-\frac{m}{2}} \\ q \sim 2^{\frac{m}{2}}}}\sum_{\ell_3 \sim 2^{\frac{m}{2}-2k}} \left| \left\langle f_3^{\omega_3^{\ell_3}}, \varphi_{p-\frac{q^3}{2^m}, \ell_3 2^{2k}+\widetilde{n_3}2^{\frac{m}{2}-2k}+\left[0, 2^{\frac{m}{2}-2k+1} \right]}^{3k} \right \rangle \right|^2  \left| \left \langle f_4, \varphi_{p, \ell_3 2^{2k}+\widetilde{n_3}2^{\frac{m}{2}-2k}+3\left[0, 2^{\frac{m}{2}-2k+1} \right]}^{3k} \right \rangle \right|^2\,.
$$

By Parseval, we further have
\begin{align} 
\bbA_p^2 &\lesssim \sum_{\substack{\widetilde{n_3} \sim 2^{4k-\frac{m}{2}} \\ q \sim 2^{\frac{m}{2}}}}\sum_{\ell_3 \sim 2^{\frac{m}{2}-2k}}  \left\|f_3^{3\Omega_3^{2^{4k-\frac{m}{2}}\ell_3+\widetilde{n_3}}}\right\|^2_{L^2 \left(I_{p-\frac{q^3}{2^m}}^{m, 3k} \right)} \left\|f_4^{3\Omega_4^{2^{4k-\frac{m}{2}}\ell_3+\widetilde{n_3}}} \right\|_{L^2 \left(I_p^{m, 3k} \right)}^2 \nonumber  \\
&\lesssim \sum_{\widetilde{n_3} \sim 2^{4k-\frac{m}{2}}} \sum_{\ell_3 \sim 2^{\frac{m}{2}-2k}} \left\|f_3^{3\Omega_3^{2^{4k-\frac{m}{2}}\ell_3+\widetilde{n_3}}} \right\|^2_{L^2 \left(3 I_{\frac{p}{2^{\frac{m}{2}}}}^{0, 3k} \right)} \left\|f_4^{3\Omega_4^{2^{4k-\frac{m}{2}}\ell_3+\widetilde{n_3}}} \right\|_{L^2 \left(I_p^{m, 3k} \right)}^2 \nonumber  \\
&\lesssim \sum_{L_3 \sim 2^{2k}} \left\|f_3^{3 \Omega_3^{L_3}} \right\|_{L^2\left(3I_{\frac{p}{2^{\frac{m}{2}}}}^{0, 3k} \right)}^2  \left\|f_4^{3\Omega_4^{L_3}} \right\|_{L^2 \left(I_p^{m, 3k} \right)}^2. \label{20230225eq01ter}
\end{align}

Thus, we are able to recover exactly \eqref{20230305eq05} and therefore the case in which either $f_3$ or $f_4$ is $2$-uniform
can be treated exactly in the same fashion as in Section \ref{20230305subsec01} again, with just replacing $2^{\mu k}$ by $2^{\mu (\frac{m}{4}-k)}$ (and $2^{\nu k}$ by $2^{\nu (\frac{m}{4}-k)}$).

We move now to the case when \emph{both $f_3$ and $f_4$ are 2-sparse, and at least one of the functions $f_1$ or $f_2$ is 2-uniform}. Again, as in the case $0 \le k \le \frac{m}{4}$, we can pick -- up to a $2^{\nu \left(\frac{m}{4}-k \right)}$-factor -- a unique measurable function $L_3: \left[ 2^{2k}, 2^{2k+1} \right] \to \left[2^{2k}, 2^{2k+1} \right]$, such that, for $\widetilde{p} \sim 2^{2k}$,
$$
2^{-\nu \left(\frac{m}{4}-k \right)}\int_{3 I_{\widetilde{p}}^{0, 3k}} |f_j|^2\lesssim \int_{3 I_{\widetilde{p}}^{0, 3k}} \left| f_j^{3 \Omega_j^{L_3 \left(\widetilde{p} \right)}} \right|^2  \lesssim \int_{3 I_{\widetilde{p}}^{0, 3k}} |f_j|^2, \quad j\in\{3, 4\}\,.
$$
Note that in the regime $\frac{m}{8} \le k <\frac{m}{4}$, since $\omega_3^{\ell_3}=\left[\ell_3 2^{\frac{m}{2}+5k}, \left(\ell_3+1 \right) 2^{\frac{m}{2}+5k} \right]$ and $\Omega_3^{L_3}=\left[L_3 2^{m+k}, \left(L_3+1 \right) 2^{m+k} \right]$ we have $|\Omega_3^{L_3}|\leq |\omega_3^{\ell_3}|$  and hence if $\omega_3^{\ell_3} \cap\Omega_3^{L_3}\not=\emptyset$ then $\Omega_3^{L_3} \subseteq \omega_3^{\ell_3}$ case in which  $\ell_3=\left\lfloor \frac{L_3}{2^{4k-\frac{m}{2}}} \right\rfloor+ O(1)$. Consequently, it makes sense to define the map $\ell_3:\,[2^{2k}, 2^{2k+1}] \to [2^{\frac{m}{2}-2k}, 2^{\frac{m}{2}-2k+1}]$ by
$\ell_3 \left(\widetilde{p} \right):= \left\lfloor \frac{L_3(\widetilde{p})}{2^{4k-\frac{m}{2}}} \right\rfloor$.

Arguing as in \eqref{20230225eq03}--\eqref{20230225eq04} and making use of the analogue of  \eqref{20250225eq30}, we deduce that
\begin{equation} \label{20230305eq20}
\Lambda_{k,m}^{\bbU_\mu}(f_1,f_2,f_3,f_4) \lesssim 2^{2k} \cdot 2^{\nu \left(\frac{m}{2}-k \right)} \sum_{\widetilde{p} \sim 2^{2k}} \left\|f_3^{3 \Omega_3^{L_3(\widetilde{p})}} \right\|_{L^2 \left(3 I_{\widetilde{p}}^{0, 3k} \right)} \left\|f_4^{3 \Omega_4^{L_3(\widetilde{p})}} \right\|_{L^2 \left(I_{\widetilde{p}}^{0, 3k} \right)} \cdot \calR_{\widetilde{p}},
\end{equation}
where
$$
\calR_{\widetilde{p}}^2:=\sum_{q \sim 2^{\frac{m}{2}}} \sum_{\left(\ell_1, \ell_2, \ell_3 \left(\widetilde{p} \right) \right) \in \TFC_{2^{\frac{m}{2}} \widetilde{p}, q}} \int_{I_{\widetilde{p}}^{0, 3k}} \int_{I_q^{m, k}} \left| \underline{f}_1^{\omega_1^{\ell_1}}(x-t) \underline{f}_2^{\omega_2^{\ell_2}}(x+t+t^3) \right|^2 dtdx.
$$
Similar to the case $0 \le k \le \frac{m}{4}$, in order to estimate $\calR_{\widetilde{p}}$ we first split the regime of $q$ as follows $\left\{q \sim 2^{\frac{m}{2}} \right\}=\bigcup_{\widetilde{q} \sim 2^{2k}} \left[ \widetilde{q} 2^{\frac{m}{2}-2k}, \left( \widetilde{q}+1 \right) 2^{\frac{m}{2}-2k}\right]\,,$
and observe that if $q \in \left[ \widetilde{q} 2^{\frac{m}{2}-2k}, \left( \widetilde{q}+1 \right) 2^{\frac{m}{2}-2k}\right]$ and $\left(\ell_1, \ell_2, \ell_3( \widetilde{p}) \right) \in \TFC_q$, then $\left(\ell_1, \ell_2, \ell_3 (\widetilde{p}) \right) \in \TFC_{2^{\frac{m}{2}-2k} \widetilde{q}}$.  Indeed, this is an immediate consequence of
$$
\left| \left(  \ell_1-\ell_2+\frac{3q^2}{2^m} \ell_3 (\widetilde{p}) \right)- \left( \ell_1- \ell_2+\frac{3 \widetilde{q}^2  \cdot 2^{m-4k} }{2^m} \ell_3 (\widetilde{p}) \right)  \right| \lesssim 2^{\frac{m}{2}-4k} \le 1.
$$

Hence,  via Parseval and a similar argument with the one in Section \ref{subsec4.2.2}, one has
\begin{align} 
\calR_{\widetilde{p}}^2&\lesssim \sum_{\widetilde{q} \sim 2^{2k}} \sum_{\substack{ q \in \left[ \widetilde{q}2^{\frac{m}{2}-2k}, \left(\widetilde{q}+1 \right)2^{\frac{m}{2}-2k} \right] \\ \left(\ell_1, \ell_2, \ell_3 \left(\widetilde{p} \right) \right) \in \TFC_{2^{\frac{m}{2}-2k}\widetilde{q}}}} \int_{I_{\widetilde{p}}^{0, 3k}} \int_{I_q^{m, k}} \left| \underline{f}_1^{\omega_1^{\ell_1}}(x-t) \underline{f}_2^{\omega_2^{\ell_2}}(x+t+t^3) \right|^2 dtdx  \nonumber \\
&\lesssim \sum_{\substack{\widetilde{q} \sim 2^{2k} \\ \ell_2 \sim 2^{\frac{m}{2}-2k}}} \int_{I_{\widetilde{p}}^{0, 3k}} \int_{I_{\widetilde{q}}^{0, 3k}}\left|\underline{f}_1^{\omega_1^{\ell_2+\frac{3\widetilde{q}^2}{2^{2k}} \ell_3(\widetilde{p})}}(x-t)  \underline{f}_2^{\omega_2^{\ell_2}}(x+t+t^3) \right|^2 dtdx \nonumber  \\
&\lesssim \sum_{\substack{\widetilde{q} \sim 2^{2k} \\ \ell_2 \sim 2^{\frac{m}{2}-2k}}}  \left(\int_{3 I^{0, 3k}_{\widetilde{p}-\widetilde{q}}} \left| f_1^{\omega_1^{\ell_2+\frac{3\widetilde{q}^2}{2^{2k}} \ell_3(\widetilde{p})}}  \right|^2 \right) \left(\int_{3 I^{0, 3k}_{\widetilde{p}+\widetilde{q}}} \left|f_2^{\omega_2^{\ell_2}} \right|^2 \right). \label{20230306eq02}
\end{align}
Recall that from our assumption, at least one of the functions $f_1$, $f_2$ is 2-uniform. Without loss of generality, we may assume $f_1=f_1^{\bbU_\nu^2}$. Therefore, we have
\begin{equation*}
\calR_{\widetilde{p}}^2  \lesssim 2^{-3\nu \left(\frac{m}{4}-k \right)} \cdot  \sum_{\widetilde{q} \sim 2^{2k}}   \left(\int_{3 I^{0, 3k}_{\widetilde{p}-\widetilde{q}}} \left| f_1 \right|^2 \right) \left(\int_{3 I_{\widetilde{p}+\widetilde{q}}^{0, 3k}} \left|f_2 \right|^2 \right).
\end{equation*}
Inserting the above estimate back to \eqref{20230305eq20}, we see that the resulting estimate is precisely \eqref{20230305eq40} with $\nu k$ replaced by $\nu \left(\frac{m}{4}-k \right)$. The rest of the proof, therefore, follows now the same treatment as in \eqref{Uniform-II-Small} (with the obvious modifications).

 \subsection{Treatment of the time-frequency sparse component $\Lambda_m^{k, \bbS_\nu^2}$}\label{SparseII}~\\
In this case, since all the functions $f_1, f_2, f_3, f_4$ are $2$-sparse, we have
\begin{enumerate}[leftmargin=2 em]
    \item [$\bullet$] up to a $2^{\nu \left(\frac{m}{4}-k \right)}$ factor, there is a unique measurable $L_3: \left[2^{2k}, 2^{2k+1} \right] \to \left[2^{2k}, 2^{2k+1} \right]$ such that for $j\in\{3, 4\}$ and $\widetilde{p} \sim 2^{2k}$
   \begin{equation} \label{ctl34}
   2^{-\nu \left(\frac{m}{4}-k \right)}\,\int_{3 I_{\widetilde{p}}^{0, 3k}} |f_j|^2\lesssim \int_{3 I_{\widetilde{p}}^{0, 3k}} \left| f_j^{3 \Omega_j^{L_3 \left(\widetilde{p} \right)}} \right|^2 \lesssim \int_{3 I_{\widetilde{p}}^{0, 3k}} |f_j|^2\,.
   \end{equation}
    Moreover, as in the previous case, we denote $\ell_3 \left(\widetilde{p} \right):= \left\lfloor \frac{L_3(\widetilde{p})}{2^{4k-\frac{m}{2}}} \right\rfloor$ for $\widetilde{p} \sim 2^{2k}$.

    \item [$\bullet$] up to a $2^{3\nu \left(\frac{m}{4}-k\right)}$ factor, there exists unique measurable functions $\ell_j: \left[2^{2k}, 2^{2k+1} \right] \to \left[2^{\frac{m}{2}-2k}, 2^{\frac{m}{2}-2k+1} \right]$, $j\in\{1, 2\}$, such that for each $q_j \sim 2^{2k}$
    \begin{equation} \label{20230305eq60}
   2^{-3\nu \left(\frac{m}{4}-k\right)}\,\int_{3I_{q_j}^{0, 3k}} \left|f_j \right|^2\lesssim \int_{3I_{q_j}^{0, 3k}} \left|f_j^{\omega_j^{{\ell_j(q_j)}}} \right|^2 \lesssim  \int_{3I_{q_j}^{0, 3k}} \left|f_j \right|^2\,.
    \end{equation}
\end{enumerate}
Using these together with \eqref{corelmodksmall} and Cauchy--Schwarz, we deduce
\begin{equation} \label{20230305eq61}
\Lambda_{k,m}^{\bbS_\nu^2}  \lesssim  2^k \cdot 2^{7\nu  \left(\frac{m}{4}-k \right)} \sum_{\widetilde{p},\widetilde{q} \sim 2^{2k}} \bbA_{\widetilde{p}, \widetilde{q}}\bbB_{\widetilde{p}, \widetilde{q}},
\end{equation}
where
\begin{equation}
\bbA_{\widetilde{p}, \widetilde{q}}^2 :=\sum_{\substack{p \in \left[\widetilde{p}2^{\frac{m}{2}}, \left(\widetilde{p}+1 \right)2^{\frac{m}{2}} \right] \\ q \in \left[ \widetilde{q} 2^{\frac{m}{2}-2k}, \left( \widetilde{q}+1 \right)2^{\frac{m}{2}-2k} \right] \\ n_3 \sim 2^{2k}}} \sum_{\substack{\left(\ell_1, \ell_2, \ell_3\left( \frac{p}{2^{\frac{m}{2}}} \right) \right) \in \TFC_{p, q} \\ \eqref{ctl34}\&\eqref{20230305eq60}}} \left| \left\langle f_3^{\omega_3^{\ell_3}}, \varphi_{p-\frac{q^3}{2^m}, \ell_3 2^{2k}+n_3}^{3k} \right \rangle \right|^2  \left| \left \langle f_4,\, \varphi_{p, \ell_3 2^k+n_3+\ell_1+\ell_2}^{3k} \right \rangle \right|^2,
\end{equation}
and
\begin{eqnarray*}
&& \bbB_{\widetilde{p}, \widetilde{q}}^2:=\sum_{\substack{p \in \left[\widetilde{p}2^{\frac{m}{2}}, \left(\widetilde{p}+1 \right)2^{\frac{m}{2}} \right] \\ q \in \left[ \widetilde{q} 2^{\frac{m}{2}-2k}, \left( \widetilde{q}+1 \right)2^{\frac{m}{2}-2k} \right] \\ n_3 \sim 2^{2k}}} \sum_{\substack{\left(\ell_1, \ell_2, \ell_3\left( \frac{p}{2^{\frac{m}{2}}} \right) \right) \in \TFC_{p, q} \\ \eqref{ctl34}\&\eqref{20230305eq60}}}\\
&& \quad \quad \quad  \quad \frac{1}{\left|I_p^{m, 3k} \right|^2} \bigg| \int_{I_p^{m, 3k}} \int_{I_q^{m, k}}  \underline{f}_1^{\omega_1^{\ell_1}}(x-t) \underline{f}_2^{\omega_2^{\ell_2}}(x+t+t^3) e^{-i 2^{\frac{m}{2}+k}t \cdot \frac{3q^2}{2^m}n_3} \, e^{i 2^{\frac{m}{2}+3k}t \cdot \left(-\ell_1 +\ell_2 - \ell_3 \frac{3q^2}{2^m}\right)} dt dx \bigg|^2.
\end{eqnarray*}

We first estimate $\bbA_{\widetilde{p}, \widetilde{q}}$ by following the arguments in \eqref{20230225eq01}:
$$
\bbA_{\widetilde{p}, \widetilde{q}}^2  \lesssim \sum_{\substack{p \in \left[\widetilde{p}2^{\frac{m}{2}}, \left(\widetilde{p}+1 \right)2^{\frac{m}{2}} \right] \\ q \in \left[ \widetilde{q} 2^{\frac{m}{2}-k}, \left( \widetilde{q}+1 \right)2^{\frac{m}{2}-k} \right]}} \left\|f_3^{3\Omega_3^{L_3 (\widetilde{p})}} \right\|^2_{L^2 \left(3 I_{p-\frac{q^3}{2^m}}^{m, 3k} \right)} \left\|f_4^{3\Omega_4^{L_3(\widetilde{p})}} \right\|_{L^2 \left(3 I_p^{m, 3k} \right)}^2
$$
which matches the estimate \eqref{20230306eq01}, and therefore
$$
\bbA_{\widetilde{p}, \widetilde{q}}^2  \lesssim \sum_{s=0}^{2^{2k}-1} \left\|f_3^{3 \Omega_3^{L_3(\widetilde{p})}} \right\|^2_{L^2 \left(3 I_{2^{2k} \widetilde{p}-s-\frac{\widetilde{q}^3}{2^{4k}}}^{0, 5k} \right)} \left\|f_4^{3\Omega_4^{L_3(\widetilde{p})}} \right\|_{L^2 \left(3 I_{2^{2k} \widetilde{p}+s}^{0, 5k} \right)}^2.
$$
Next, we estimate $\bbB_{\widetilde{p}, \widetilde{q}}$. Arguing as in \eqref{20230225eq03} and \eqref{20230306eq02}, we have
$$
\bbB_{\widetilde{p}, \widetilde{q}}^2 \lesssim 2^{2k} \sum_{\substack{\ell_2 \sim 2^{\frac{m}{2}-2k} \\ \eqref{ctl34}\&\eqref{20230305eq60}}}  \left(\int_{3 I^{0, 3k}_{\widetilde{p}-\widetilde{q}}} \left| f_1^{\omega_1^{\ell_2-\frac{3\widetilde{q}^2}{2^{2k}} \ell_3(\widetilde{p})}}  \right|^2 \right) \left(\int_{3 I^{0, 3k}_{\widetilde{p}+\widetilde{q}}} \left|f_2^{\omega_2^{\ell_2}} \right|^2 \right).
$$
Combining the both estimates of $\bbA_{\widetilde{p}, \widetilde{q}}$ and $\bbB_{\widetilde{p}, \widetilde{q}}$ with \eqref{20230305eq61}, we have
\begin{eqnarray*}
\Lambda_{k,m}^{\bbS_\mu} &\lesssim& 2^{2k} \cdot 2^{7 \nu \left(\frac{m}{4}-k \right)} \sum_{\widetilde{p}, \widetilde{q} \sim 2^{2k}} \left[\sum_{s=0}^{2^{2k}-1} \left\|f_3^{3 \Omega_3^{L_3(\widetilde{p})}} \right\|^2_{L^2 \left(3 I_{2^{2k} \widetilde{p}-s-\frac{\widetilde{q}^3}{2^{4k}}}^{0, 5k} \right)} \left\|f_4^{3\Omega_4^{L_3(\widetilde{p})}} \right\|_{L^2 \left(3 I_{2^{2k} \widetilde{p}+s}^{0, 5k} \right)}^2 \right]^{\frac{1}{2}} \\
&\cdot& \left[ 2^{2k} \sum_{\substack{\ell_2 \sim 2^{\frac{m}{2}-2k} \\ \eqref{ctl34}\&\eqref{20230305eq60}}}  \left(\int_{3 I^{0, 3k}_{\widetilde{p}-\widetilde{q}}} \left| f_1^{\omega_1^{\ell_2-\frac{3\widetilde{q}^2}{2^{2k}} \ell_3(\widetilde{p})}}  \right|^2 \right) \left(\int_{3 I^{0, 3k}_{\widetilde{p}+\widetilde{q}}} \left|f_2^{\omega_2^{\ell_2}} \right|^2 \right) \right]^{\frac{1}{2}}.
\end{eqnarray*}
Again, using the inherited time-frequency localization in the above estimate, we define
$$
\sTFC:=\left\{ \substack{\left( \widetilde{p}, \widetilde{q} \right)\\\widetilde{p}, \widetilde{q} \sim 2^{2k}}\::\:  \ell_2 \left(\widetilde{p}+\widetilde{q}\right)-\frac{3\widetilde{q}^2}{2^{2k}} \ell_3 \left( \widetilde{p} \right)=\ell_1 \left(\widetilde{p}-\widetilde{q} \right)+O(1) \right\}.
$$

Following the argument in \eqref{Uniform-I}, we see that it suffices to show that there exists some $\epsilon_3>0$, such that
$$
\# \sTFC \lesssim  2^{4k-\epsilon_3 \left(\frac{m}{2}-2k \right)}.
$$

The proof of the above estimate follows essentially the argument in Subsection \ref{Levelsetanalysis} with only some simple modifications/adaptations among which we mention:
\begin{enumerate}
    \item [(1)] The three generic measurable functions $a(\cdot), b(\cdot)$, and $c(\cdot)$ have now the range
    $\newline$ $\left[2^{\frac{m}{2}-2k}, 2^{\frac{m}{2}-2k+1}\right] \cap \Z$, rather than the previous $[2^{2k}, 2^{2k+1} ]\cap \Z$;

    \item [(2)] In the current argument one replaces $\epsilon k$ by $\epsilon \left(\frac{m}{2}-2k \right)$;

     \item [(3)] Finally, the analogue of \eqref{Keydec} reads now
\begin{align}\label{Keydec1}
(1+|\frac{m}{2}-2k|) \cdot 2^{4k} 2^{- \left(\frac{m}{2}-2k \right)}\gtrsim \sum_{v_1, v_2\sim 2^{2k}}
\left\lceil\frac{v_1-v_2}{2^{2k}}\,c(v_1)\right\rceil \gtrsim 2^{4k-1024\epsilon_3\left(\frac{m}{2}-2k \right)}\,,
\end{align}
which gives  $\epsilon_3\approx\frac{1}{2000}$.
   \end{enumerate}


%
%

%
%
\section{The single scale estimate (IV): An outline of the case $k\leq 0$} \label{Seckneg}
In this section we only provide a very brief outline of the case $k< 0$, which, exploiting the versatility of the Rank II LGC method, follows very similar reasonings with the ones employed in the treatment of the case $k>0$.

Recalling the decomposition in Section \ref{subsec:prell} -- equation \eqref{fmd}, we now turn our attention on $\L_{-}$. We first notice that the equations \eqref{fm1}-\eqref{multkm} remain valid. However, in the regime $k<0$, the relative size ordering of the width of the frequency bands  $|\xi_1-\xi_2|\sim 2^{k+m}$ and $|\xi_3-\xi_2|\sim 2^{3k+m}$ is reversed, and thus the natural macroscale discretization is given by: the band $|\xi_1-\xi_2|\sim 2^{k+m}$ in squares $(\xi_1,\xi_2)$ of size $\ds \big[ 2^{k+m} (\el+2),   2^{k+m} (\el +3)  \big]  \times \big[ 2^{k+m} \el,   2^{k+m} (\el+1)  \big]$ with $\el\in\Z$ while the frequency band $|\xi_3-\xi_2|\sim 2^{3k+m}$ is split into frequency squares $(\xi_2,\xi_3)$ of size $\ds \big[ 2^{k+m} (\el+2) + 2^{3k+m} \mathfrak{r} ,   2^{k+m} (\el+2) + 2^{3k+m} (\r +1)  \big]  \times \big[2^{k+m} (\el+2) + 2^{3k+m} (\r +2),  2^{k+m} (\el+2) +  2^{3k+m} (\r+3)  \big]$ with $\r\sim 2^{-2k}$. With these, it is immediate to see that
\beq\label{multkms1n}
\mathfrak{m}_{k,m}(\xi_1,\xi_2,\xi_3)=\sum_{\el\in\Z} \sum_{\r\sim 2^{-2k}}\mathfrak{m}^{\el,\r}_{k,m}(\xi_1,\xi_2,\xi_3)
\eeq
where
\begin{align*} 
\mathfrak{m}^{\el,\r}_{k,m}(\xi_1,\xi_2,\xi_3) & := \left(\int_{\R}e^{-i\,\frac{\xi_1-\xi_2}{2^k}\,t}\,e^{-i\,\frac{\xi_3-\xi_2}{2^{3k}}\,t^3}\,\rho(t)\,dt\right) \\
& \qquad \cdot \phi_1\left(\frac{\xi_1}{2^{k+m}}-\el\right)\,\phi_2\left(\frac{\xi_2}{2^{3k+m}}-\el 2^{-2k}-\r\right)\,\phi_3\left(\frac{\xi_3}{2^{3k+m}}-\el 2^{-2k}-\r\right)
\end{align*}
for three suitable $\phi_1,\,\phi_2, \phi_3\in C_{0}^{\infty}(\R)$.

Further, for $s\in\{1,4\}$ we define
\begin{equation}\label{functsupn}
\widehat{f_{1,\el}}(\xi_1):= \hat{f}_1(\xi_1)\,\phi_1\left(\frac{\xi_1}{2^{k+m}}-\el\right)\qquad\textrm{and}\qquad\widehat{f_{4,\el}}(\xi_4):= \hat{f}_4(\xi_4)\,\phi_4\left(\frac{\xi_4}{2^{k+m}}+3\el\right)\,,
\end{equation}
while for $f_2$ and $f_3$ we set
\begin{equation}\label{functsup4n}
\widehat{f_{2,\el,\r}}(\xi_2):= \hat{f}_2(\xi_2)\,\phi_2\left(\frac{\xi_2}{2^{3k+m}}-\el 2^{-2k}-\r \right)\qquad\textrm{and}\qquad \widehat{f_{3,\el,\r}}(\xi_3):= \hat{f}_3(\xi_3)\,\phi_3\left(\frac{\xi_3}{2^{3k+m}}-\el 2^{-2k}-\r \right)\,.
\end{equation}
With these done, we have that
\begin{equation}\label{smdn}
\L_{k,m}(f_1,f_2,f_3,f_4)\approx\sum_{\el\in\Z} \sum_{\r\sim 2^{-2k}} \L_{k,m}^{\el,\r}(f_1,f_2,f_3,f_4)
\end{equation}
where
\begin{equation}\label{smd1n}
\L_{k,m}^{\el,\r}(f_1,f_2,f_3,f_4):=\L_{k,m}(f_{1,\el},f_{2,\el,\r},f_{3,\el,\r},f_{4,\el})\,.
\end{equation}
Finally, since in the current regime, $2^{-3k}\geq 2^{-k}$, the analogue of \eqref{smd2} becomes
\begin{equation}\label{smd2n}
\L_{k,m}^{\el,\r,j}(f_1,f_2,f_3,f_4):=\L_{k,m}(f_{1,\el},f_{2,\el,\r},f_{3,\el,\r},\chi_{I_j^{0,3k}}\,f_{4,\el})\,.
\end{equation}
or, equivalently
\begin{align}
& \L_{k,m}^{\el,\r,j}(f_1,f_2,f_3,f_4 )\approx \nonumber \\
& \qquad \int_{\R}\int_{\R} f_{1,\el}(x-t)\,f_{2,\el,\r}(x+t+t^3)\,f_{3,\el,\r}(x-t^3)\,f_{4,\el}(x)\, 2^{k}\,\rho(2^{k} t)\,\varphi(2^{3k} x-j)\,dt\,dx\,. \label{fm11n}
\end{align}
This completes the first stage of our discretization.  

However, making now the change of variable $t\,\rightarrow\,t^{\frac{1}{3}}$, $k\,\rightarrow\,\frac{k}{3}$ and switching the indices $1$ and $3$ we conclude that 
\begin{align}
& \L_{k,m}^{\el,\r,j}(f_1,f_2,f_3,f_4)\approx \nonumber \\
& \qquad \int_{\R}\int_{\R} f_{1,\el,\r}(x-t)\,f_{2,\el,\r}(x+t+t^{\frac{1}{3}})\,f_{3,\el}(x-t^{\frac{1}{3}})\,f_{4,\el}(x)\, 2^{k}\,\rho(2^{k} t)\,\varphi(2^{k} x-j)\,dt\,dx\,, \label{fm11ns}
\end{align}
which is essentially the same as \eqref{fm11:2} appearing in the case $k>0$. 

Thus, in a nutshell, if one compares the cases $k<0$ and $k>0$, the roles of the functions $f_1$ and $f_3$ are to be reversed. As a consequence, the analysis performed earlier in the case $k>0$ can now be applied line by line in the context $k<0$. We leave further details to the interested reader.
\section{The multi-scale analysis (I): Treatment of the case $k\geq \frac{m}{4} \geq 0$}
\label{sec:bilinear:analysis}

In this section, we present the multi-scale analysis for the case $k\geq m/4 $ , in order to prove Theorem \ref{thm:case1} for $k\geq \frac{m}{4}$. First, we need to re-interpret the single scale estimate and related to the discretized multi-linear form that will be treated in this section.

\subsection{Boundedness with decay in the local$-L^2$ range}


The purpose of this section is to provide bounds for a fixed $m\geq 0$ of
\[
\Lambda_m^3(f_1,f_2,f_3,f_4)= \sum_{k\geq {m \over 8}} \ \sum_{\el,j \in {\mathbb Z} \atop \r \sim 2^{2k} } \L_{k,m}^{\el,\r,j}(f_1,f_2,f_3,f_4).
\]

Recalling the definitions \eqref{smd2}, \eqref{functsupn}, \eqref{fm11:2} and the convention in Notation \ref{not:norms:TF:reg}, we have proved in Section \ref{sec:onescale_k_large} that, for $\mu_0 >0$ absolute constant given by \eqref{eq:decay:single:mu}, $\big| \L_{k,m}^{\el,\r,j}(f_1,f_2,f_3,f_4) \big|$ is bounded above by
\begin{align*}
2^{- \mu_0 \, m}  \,  |I_j^k|^{-1} \|f_1\|_{L^2(I_j^k \times W_{(\mathfrak{r}+2)+2^{2k} \mathfrak l}^{k+m})} \, \|f_2\|_{L^2(I_j^k \times W_{\mathfrak{r}+2^{2k} \mathfrak l}^{k+m})} \, \|f_3\|_{L^2(I_j^k \times W_{ \mathfrak l+2}^{3k+m})} \|f_4\|_{L^2(I_j^k \times W_{  - (3 \mathfrak l+3)}^{3k+m})} \\
+ |I_j^k|^{-1} \|f_1\|_{L^2(I_j^k \times W_{(\mathfrak{r}+2)+2^{2k} \mathfrak l}^{k+m})} \, \|f_2\|_{L^2(I_j^k \times W_{\mathfrak{r}+2^{2k} \mathfrak l}^{k+m})} \, \|f_3^{\bbS^1_{\mu_0}}\|_{L^2(I_j^k \times W_{ \mathfrak l+2}^{3k+m})} \|f_4^{\bbS^1_{\mu_0}}\|_{L^2(I_j^k \times W_{ - (3 \mathfrak l+3)}^{3k+m})}.
\end{align*}

Applying Cauchy-Schwarz in the $\r$ parameter, we are reduced to considering
\begin{align}
\label{eq:expression:multi:scale:parametrix}
\sum_{k\geq {m \over 8}}  \sum_{\el,j \in {\mathbb Z}} 2^{-\mu_0 \, m}  \,  |I_j^k|^{-1} \|f_1\|_{L^2(I_j^k \times W_{\mathfrak l}^{3k+m})} \, \|f_2\|_{L^2(I_j^k \times W_{ \mathfrak l}^{3k+m})} \, \|f_3\|_{L^2(I_j^k \times W_{ \mathfrak l+2}^{3k+m})} \|f_4\|_{L^2(I_j^k \times W_{ - (3 \mathfrak l+3)}^{3k+m})} \\
+ |I_j^k|^{-1} \|f_1\|_{L^2(I_j^k \times W_{\mathfrak l}^{3k+m})} \, \|f_2\|_{L^2(I_j^k \times W_{\mathfrak l}^{3k+m})} \, \|f_3^{\bbS^1_{\mu_0}}\|_{L^2(I_j^k \times W_{ \mathfrak l+2}^{3k+m})} \|f_4^{\bbS^1_{\mu_0}}\|_{L^2(I_j^k \times W_{ - (3 \mathfrak l+3)}^{3k+m})}. \nonumber
\end{align}
In the expression above, we recover a more classical model operator depending on three parameters: the scale $k$, the frequency parameter $\el$ and the space parameter $j$. It is now natural to switch to the \emph{tile formalism}, as in Section \ref{sec:proof:m:+}. Indeed, to each triple $(k, \el, j)$ we can uniquely associate a time frequency four-tile $P=(P_1, P_2, P_3, P_4)$ defined as follows: for all $1 \leq i \leq 4$, 
\begin{equation}
\label{eq:def:spatial:intv}
P_i= I_P \times \omega_{P_i}
\end{equation}
where $\ds I_P=I_j^k= \big[ 2^{-k} (j+1),  2^{-k} (j+2) \big]$, and 
\begin{equation}
\label{eq:def:freq:intvs}
\begin{gathered}
 \omega_{P_1}= \omega_{P_2}:=\big[ 2^{3k+m} \el ,   2^{3k+m} (\el+1) \big], \\
\omega_{P_3}:=\big[ 2^{3k+m} (\el +2),   2^{3k+m} (\el+3)  \big], \\ \omega_{P_4}:=-\big[2^{3k+m} (3\el+2)  , 2^{3k+m} (3\el+3) \big]. 
\end{gathered}
\end{equation}

Motivated by the expression \eqref{eq:expression:multi:scale:parametrix} and the definition of the time-frequency norm \eqref{eq:notation:time:freq}, we take a closer look at
\begin{align*}
\|f_1\|_{L^2(I_j^k \times W_{ \mathfrak l}^{3k+m})} = \|f_1\|_{L^2(I_P \times \omega_{P_1})} =  \Big(  \sum_{\substack{s=I \times \omega \\ I \subset I_P, |I|=2^{-{m \over 2}}|I_P| \\ \omega \subset \omega_{P_1}, |\omega|^3=2^{{m \over 2}}|\omega_{P_1}| }}  |\langle f_1, \phi_s  \rangle|   \Big)^{1 \over 2}.
\end{align*}
As in Section \ref{sec:proof:m:+} and in particular in Lemma \ref{lemma:wide:to:thin}, we prefer to work with wave packets associated to larger frequency intervals, of length $2^{-{m \over 2}}|\omega_{P_1}|$, in order to have a better space separation property for the tiles. So up to considering a superposition of similar terms, we replace $\|f_1\|_{L^2(I_j^k \times W_{ \mathfrak l}^{3k+m})}$ by the expression
\begin{equation}
\label{eq:def:coeff:multi:lin:form:f_1}
f_1(P_1):=\Big(  \sum_{\substack{I \subset I_P, |I|=2^{-{m \over 2}}|I_P|^3 \\ \omega \subset \omega_{P_1}, |\omega|= 2^{-{m \over 2}}|\omega_{P_1}|}}  |\langle f_1, \phi_{I \times \omega}  \rangle|^2   \Big)^{1 \over 2}.
\end{equation}
For any four-tile $P=(P_1, P_2, P_3, P_4)$ and any $i \in \{ 1, \ldots, 4 \}$, we set
\begin{equation*}
\label{eq:def:new:small:tiles}
\mathfrak{S}_i(P_i)=\big\lbrace    s=I \times \omega : I \subset I_{\pp}, |I|= |I_{\pp}|^3 \cdot 2^{-{m \over 2}}, \quad  \omega \subset \omega_{\pp_i}, |\omega|= 2^{-{m \over 2}} |\omega_{\pp_i}|  \big\rbrace
\end{equation*}
and we define
\begin{equation}
\label{def:quant:tile}
f_i(P_i):= \Big( \sum_{s_i \in \mathfrak{S}_i(P_i)} \big| \langle f_i, \phi_{s_i}  \rangle  \big|^2 \Big)^{1 \over 2},
\end{equation}
which matches the definition \eqref{eq:def:coeff:multi:lin:form:f_1} above.

For $f_3$ and $f_4$ we also need to take into account the spatial sparse and uniform components. For $\mu_0>0$ given by \eqref{eq:decay:single:mu} and $i \in \{ 3, 4\}$, we set, according to the definitions \eqref{sparsefs1}--\eqref{unifs1} 
\begin{equation*}
\mathfrak{S}^{\mathfrak{c}}_i(P_i) := \big\{ s=I \times \omega \in \mathfrak{S}_i(P_i) :  \aver {3I}|f_i(x)|^2 dx > 2^{\mu_0 m } \aver {I_P} |f_i(x)|^2 dx \big\}
\end{equation*}
and
\begin{equation*}
\mathfrak{S}^{\mathfrak{u}}_i(P_i) := \big\{ s=I \times \omega \in \mathfrak{S}_i(P_i) :  \aver {3I}|f_i(x)|^2 dx \leq  2^{\mu_0 m } \aver {I_P} |f_i(x)|^2 dx \big\},
\end{equation*}
respectively. In fact, we made use of the uniform component in Section \ref{sec:onescale_k_large} in order to deduce \eqref{eq:expression:multi:scale:parametrix}. So if $\ds i \in \{3, 4 \}$ the \emph{sparse} components are defined by:
\begin{equation}
\label{def:quant:tile:clustered}
f_i^{\mathfrak{c}}(P_i):= \Big( \sum_{s_i \in \mathfrak{S}^{\mathfrak{c}}_i(P_i)} \big| \langle f_i, \phi_{s_i}  \rangle  \big|^2 \Big)^{1 \over 2}.
\end{equation}

Now we consider the $4$-linear forms 
\begin{equation}
\label{def:disc:tile:4:form:c}
\Lambda^{\mathfrak{c}}_{\P}(f_1, f_2, f_3, f_4):=\sum_{P \in \P} |I_P|^{-1} f_1(P_1) \cdot f_2(P_2) \cdot f_3^{\mathfrak{c}}(P_3) \cdot f_4^{\mathfrak{c}}(P_4)
\end{equation}
and 
\begin{equation}
\label{def:disc:tile:4:form:u}
\Lambda^{\mathfrak{u}}_{\P}(f_1, f_2, f_3, f_4):=\sum_{P \in \P} |I_P|^{-1} f_1(P_1) \cdot f_2(P_2) \cdot f_3(P_3) \cdot f_4(P_4).
\end{equation}
We notice that estimating \eqref{eq:expression:multi:scale:parametrix} reduces to bounding these four-(sub)linear forms since we have
\begin{equation} \label{eq:decompL} \Big| \sum_{\substack{k \in \Z \\ k \geq {m \over 4} }} \ \sum_{\el \in \Z \atop \r \sim 2^{2k}} \ \sum_{j \in \Z} \ \L_{k,m}^{\el,\r,j}(f_1,f_2,f_3,f_4) \Big| \lesssim \Lambda^{\mathfrak{c}}_{\P}(f_1, f_2, f_3, f_4) + 2^{-\mu_0 m} \Lambda^{\mathfrak{u}}_{\P}(f_1, f_2, f_3, f_4).
\end{equation}
This will be done, at least for the uniform component, following the same strategy as in Section \ref{sec:proof:m:+}, which corresponds to the situation $m=0$.

More precisely, we prove the following result:
\begin{theorem}
\label{thm:multi:scale:aniso}
Let $m \in \N$ and let $\P$ be a collection of $4$-tiles as defined in \eqref{eq:def:spatial:intv} and \eqref{eq:def:freq:intvs}. Then we have for all $2 < p_1, p_2, p_3, p'_4 \leq \infty$ that 
\begin{equation}
\label{eq:wanted:res:u}
\big|\Lambda^{\mathfrak{u}}_{\P}(f_1, f_2, f_3, f_4)\big|  \lesssim \|f_1\|_{p_1}    \|f_2\|_{p_2} \|f_3\|_{p_3} \|f_4\|_{p'_4}        
\end{equation}
and that there exists $\delta:=\delta(p_1,p_2,p_3)$ such that
\begin{equation}
\label{eq:wanted:res:c}
\big|\Lambda^{\mathfrak{c}}_{\P}(f_1, f_2, f_3, f_4)\big|  \lesssim 2^{-\delta m} \|f_1\|_{p_1} \|f_2\|_{p_2} \|f_3\|_{p_3} \|f_4\|_{p'_4}.   
\end{equation}
\end{theorem}
It is easy to see that Theorem \ref{thm:multi:scale:aniso} and \eqref{eq:decompL} imply the boundedness of our operator (in the regime $k \geq {m \over 4}$), in the local $L^2$ range, and so prove the $(k\geq \frac{m}{4})$-part of Theorem \ref{thm:case1} in the local-$L^2$ range, which is
\begin{equation} \label{eq:T1-2}
\|T_m^1\|_{L^{p_1}\times L^{p_2} \times L^{p_3} \to L^p}:=\|\sum_{k\geq {m\over 4}}T_{k,m} \|_{L^{p_1}\times L^{p_2} \times L^{p_3} \to L^p} \lesssim 2^{-\delta m}
\end{equation}
for some $\delta>0$ depending on $p_1, p_2, p_3$ and $\mu_0$.

\subsection{Treatment of the \emph{uniform} component in Theorem \ref{thm:multi:scale:aniso}}~\\
Here we use the terminology that was introduced in Section \ref{sec:proof:m:+} (eventually with obvious adaptations to the current situation), although the tiles have $2^m$ larger frequency intervals. Before proceeding, the reader is invited to review this section. Since the ideas are similar, we only highlight the statements that allow to deduce \eqref{eq:wanted:res:u}.

As in Proposition \ref{prop:size1}, our sizes are only maximal local $L^2$ averages:
\begin{proposition}
\label{prop:size:est}
Given any subcollection $\P' \subseteq \P$ of tri-tiles, we have that 
\begin{equation}
\label{eq:size:f_3}
\size_{\P'}^j(f_j) \lesssim \sup_{P \in \P'} \frac{1}{|I_P|^{1 \over 2}} \left( \int_{\R} |f_j(x)|^2 \chi_{I_P}(x)dx \right)^\frac{1}{2}
\end{equation}
for $j \in \{ 1, 2, 3, 4\}$.
\end{proposition} 

We need the following \emph{decomposition lemma} which allows to organize our collection of tiles into subcollections having more structure, and over which we have control. 

\begin{lemma}
\label{lemma:decomp:lemma}
Let  $\PP' \subseteq \PP$ a collection of $4$-tiles,  $j \in \{1, 2,3, 4\}$ and $\lambda>0$. If $\ds \size_{\PP'}^j (f_j) \leq \lambda$, then there exists a decomposition $\ds \PP'= \PP'_{(1)} \cup \PP'_{(2)}$ so that 
\[
\size_{\PP'_{(1)}}^j (f_j) \leq {\lambda \over 2}
\]  
and $\ds \PP'_{(2)}$ can be written as a union of $j$-trees $\ds \PP'_{(2)}= \bigcup_{\tr{T} \in \ii F_{\lambda}} \tr{T}$ satisfying 
\begin{equation}
\label{eq:tree:top:cond}
\sum_{\tr{T} \in  \ii F_{\lambda}} |I_{\tr{T}}| \lesssim \lambda^{-2} \|f_j\|_2^2.
\end{equation}
\end{lemma}

The proof follows the same steps as that of Lemma \ref{lemma:decomposition:lemma:4:lin:form:no:curv}; in particular the key estimate \eqref{eq:est:trees:energy:wanted} still holds true, uniformly in $m$, although the frequency intervals have been dilated by a $2^m$ factor. 

We recall that a tree contains at least two lacunary directions: although $\omega_{P_1}=\omega_{P_2}$, the triple $(\omega_{P_1}, \omega_{P_3}, \omega_{P_4})$ has the same properties as a triple of intervals in the classical study of the bilinear Hilbert transform. As a consequence, we have the following \emph{tree estimate}:

\begin{lemma}
\label{lemma:tree:estimate}
If $\tr{T} \subseteq \P$ is a tree, then 
\begin{equation}
\Big| \Lambda_{\tr{T}} (f_1, f_2, f_3, f_4)  \Big| \lesssim \size_{\tr{T}}^1(f_1) \cdot \size_{\tr{T}}^2(f_2) \cdot \size_{\tr{T}}^3(f_3)  \cdot \size_{\tr{T}}^4(f_4) \cdot |I_{\tr{T}}|.
\end{equation}
\end{lemma}

By iterating the result in Lemma \ref{lemma:decomp:lemma}, we can decompose $\PP$ as in Proposition \ref{prop:decompo:full}:

\begin{proposition}
\label{prop:decompo:full:with:m}
Given $\PP$ a finite collection of $4$-tiles, we define $n_0$ so that
\[
2^{-n_0} \sim \max_{j \in \{ 1, 2, 3, 4 \} } \Big(  \frac{\size_{\PP}^j(f_j)}{\|  f_j \|_2}   \Big).
\]
Then one can decompose $\PP$ as 
\begin{equation}
\label{eq:coll:dec:a}
\PP= \bigcup_{n: 2^{-n} \leq 2^{-n_0}} \bigcup_{j \in \{ 1, 2, 3, 4  \}} \bigcup_{\tr{T} \in \ii T_n^{(j)}} \tr{T}:= \bigcup_{n: 2^{-n} \leq 2^{-n_0}}  \PP_n,
\end{equation} 
where for each $n \geq n_0$ we have, for $\PP_n \neq \emptyset$,
\[
\size_{\PP_n}^j(f_j) \leq \min \big( 2^{-n} \|f_j\|_2, \size_{\PP}^j(f_j) \big)
\]
and each $\ii T_n^{(j)}$, if non-empty, is a collection of trees such that 
\[
\sum_{\tr{T} \in \ii T_n^{(j)}}|I_{\tr{T}}| \lesssim 2^{2n}.
\] 
\end{proposition}

This allows us to obtain the equivalent of Proposition \ref{prop:final:est:0} used in the study of the anisotropic multilinear operator in Section \ref{sec:proof:m:+}:

\begin{proposition}\label{prop:final:est}
If $\PP$ is a collection of $4$-tiles, then 
\begin{align*}
\big| \Lambda_{\PP}(f_1, f_2, f_3, f_4)  \big| & \lesssim \big( \size_{\PP}^1(f_1)\big)^{\theta_1}  \big( \size_{\PP}^2(f_2)\big)^{\theta_2}  \big( \size_{\PP}^3(f_3)\big)^{\theta_3} \big( \size_{\PP}^4(f_4)\big)^{\theta_4} \\
&\cdot \|f_1\|_2^{1- \theta_1} \|f_2\|_2^{1- \theta_2}  \|f_3\|_2^{1- \theta_3} \|f_4\|_2^{1- \theta_4},
\end{align*}
for any $0 \leq \theta_1, \theta_2, \theta_3, \theta_4 <1$ with $\theta_1+\theta_2+\theta_3+\theta_4=1$.

Moreover, if $I_0$ is a fixed spatial interval and $\PP(I_0)$ denotes the collection of $4$-tiles in $\PP$ with spatial support in $I_0$, then 
\begin{align*}
\label{eq:local:est:decay}
\big| \Lambda_{\P(I_0)}(f_1, f_2, f_3, f_4)  \big|  \lesssim & \big( \size_{\P(I_0)}^1(f_1)\big)^{\theta_1}  \big( \size_{\P(I_0)}^2(f_2)\big)^{\theta_2}  \big( \size_{\P(I_0)}^3(f_3)\big)^{\theta_3} \big( \size_{\P(I_0)}^4(f_4)\big)^{\theta_4} \\
&\cdot \|f_1 \cdot \chi_{I_0}\|_2^{1- \theta_1} \|f_2 \cdot \chi_{I_0}\|_2^{1- \theta_2} \|f_3 \cdot \chi_{I_0}\|_2^{1- \theta_3} \|f_4 \cdot \chi_{I_0}\|_2^{1- \theta_4}, \nonumber
\end{align*}
for any $0 \leq \theta_1, \theta_2, \theta_3, \theta_4 <1$ with $\theta_1+\theta_2+\theta_3+\theta_4=1$.
\end{proposition}

Following \cite[Theorem 12]{BM-sparse}, and modulo a few details -- see the end of Section \ref{sec:proof:m:+} -- we have  that 
\begin{equation}
\label{eq:conc}
 \big|\Lambda^{\mathfrak{u}}_{\PP}(f_1, f_2, f_3, f_4)\big| \lesssim \|f_1\|_{p_1} \,  \|f_2\|_{p_2} \,  \|f_3\|_{p_3} \,  \|f_4\|_{p'_4},
\end{equation} 
whenever $\ds 2 < p_1, p_2, p_3, p'_4 \leq \infty$, with $\ds  \frac{1}{p_1}+\frac{1}{p_2}+\frac{1}{p_3}+\frac{1}{p'_4}=1$.

\subsection{Treatment of the \emph{sparse} component in Theorem \ref{thm:multi:scale:aniso}}~\\
\label{sebsec:clustered}
This subsection is devoted to the proof of \eqref{eq:wanted:res:c}, corresponding to the sparse component estimates and which concludes the proof of Theorem \ref{thm:multi:scale:aniso}. Once again, we follow the same line of ideas as in the previous section, and for the most part $\Lambda^{\mathfrak{c}}_{\PP}(f_1, f_2, f_3, f_4)$ can be treated in a similar way to $\Lambda^{\mathfrak{u}}_{\PP}(f_1, f_2, f_3, f_4)$.

The essential difference here consists in obtaining the extra exponential decay in the positive parameter $m$, as formulated in \eqref{eq:wanted:res:c}. This demands a more precise ``size'' estimate, where we recall that then $\mu_0$ is defined by \eqref{eq:decay:single:mu}. For that, the result below takes the place of Proposition \ref{prop:size:est}:

Recalling Definition \ref{def:size:bis}, we set 
\begin{equation}
\label{eq:def:size:sparse}
\size_{\PP'}^{j,\mathfrak{c}} (f_j) := \size_{\PP'}^{j} \big(f_j^{\mathfrak c}(P_j)\big)=\sup_{\substack{\tr T \subset \PP' \\ \tr T \, j\text{-lac tree} }} \frac{1}{|I_{\tr T}|^{1 \over 2}}\left(  \sum_{P \in \tr T} |f_j^{\mathfrak c}(P_j)| \right)^{1 \over 2}.
\end{equation}

\begin{proposition} \label{prop:sizesparse}
Given $\PP'\subseteq \PP$ a subcollection of $4$-tiles, $j \in \{ 3, 4  \}$ and $2<p<\infty$, we have that
\label{prop:size:est:clustered}
\begin{equation}
\label{eq:size:f_j:clustered}
\size_{\PP'}^{j,\mathfrak{c}} (f_j) \lesssim 2^{-{{\mu_0 m} \over 2}\left( \frac{1}{2}-\frac{1}{p}  \right)}  \sup_{P \in \PP'} \frac{1}{|I_P|^{1 \over p}} \left( \int_{\R} |f_j(x)|^p \chi_{I_P}(x)dx \right)^\frac{1}{p}.
\end{equation}
\end{proposition}

\begin{proof}
We prove \eqref{eq:size:f_j:clustered} for $j=4$, but the same proof works in the case $j=3$. To begin with, let $\tr{T}$ be a $4$-lacunary tree. It will be enough to prove, for $2<p < \infty$, that 
\begin{equation}
\frac{1}{|I_{\tr{T}}|^\frac{1}{2}} \left( \sum_{\pp \in \tr{T}} |f_j^{\mathfrak{c}}(\pp)|^2   \right)^{1 \over 2} \lesssim 2^{-{{\mu_0 m} \over 2}\left( \frac{1}{2}-\frac{1}{p}  \right)}  \frac{1}{|I_{\tr{T}}|^{1 \over p}} \left( \int_{\R} |f_j(x)|^p \chi_{I_{\tr{T}}}(x)dx \right)^\frac{1}{p}.
\end{equation}
This amounts to proving that
\begin{equation}
\label{eq:size:c:wavep}
\frac{1}{|I_{\tr{T}}|^\frac{1}{2}} \left( \sum_{\pp \in \tr{T}} \sum_{s \in \mathfrak{S}^{\mathfrak{c}}_4(\pp_4)} | \langle f_4, \phi_{s}  \rangle|^2   \right)^{1 \over 2} \lesssim 2^{-{{\mu_0 m} \over 2}\big( \frac{1}{2}-\frac{1}{p}  \big)}  \frac{1}{|I_{\tr{T}}|^{1 \over p}} \left( \int_{\R} |f_j(x)|^p \chi_{I_{\tr{T}}}(x)dx \right)^\frac{1}{p}.
\end{equation}

We decompose $\tr{T}$ into two parts $\tr{T}= {\tr{T}}^{(1)} \cup \tr{T}^{(2)} $ where 
\begin{equation}
\label{eq:def:tree:1}
 \tr{T}^{(1)} := \Bigl\{ \pp \in \tr{T} : \frac{1}{|I_{\pp}|} \int_{\R} |f_4(x)|^2 \cdot \chi_{I_{\pp}}(x) dx > 2^{-\frac{\mu_0 m}{2}} \frac{1}{|I_{\tr{T}}|} \int_{\R} |f_4(x)|^2 \chi_{I_{\tr{T}}}(x)dx  \Bigr\}
\end{equation}
and 
\[
 \tr{T}^{(2)} := \tr{T} \setminus {\tr{T}}^{(1)}.
\]
We consider both $\tr{T}^{(1)}$ and ${\tr{T}}^{(2)}$ to have the same top interval $I_{\tr{T}}$.

According to the definitions above, if $s \in \mathfrak{S}^{\mathfrak{c}}_4(\pp_4)$ with $\pp\in {\tr{T}}^{(1)}$, then the spatial information of $s$ is contained inside a set where a certain maximal function is large:
\begin{equation}
\label{eq:def:exceptional:set}
I_s \subset E(I_{\tr{T}}):= \bigl\{ x \in I_{\tr{T}} : M (|f_4|^2 \cdot \chi_{I_{\tr T}})(x) > 2^{\frac{\mu_0 m}{2}} \frac{1}{|I_{\tr{T}}|}  \int_{\R} |f_4(x)|^2 \chi_{I_{\tr{T}}}(x)dx \bigr\}.
\end{equation}

The weak-type boundedness of the maximal function implies that this set is small:
\begin{equation}
\label{eq:small:set}
\big| E(I_{\tr{T}}) \big| \lesssim 2^{- \frac{\mu_0 m}{2}} \big|  I_{\tr{T}} \big|.
\end{equation}

We define the following square function:
\begin{equation}
\label{eq:def:square:func}
\ii S_{\tr{T}}(f_4)(x):= \Big( \sum_{\pp \in \tr{T}} \sum_{s \in \mathfrak{S}^{\mathfrak{c}}_4(\pp_4)} | \langle f_4, \phi_{s}  \rangle|^2  \frac{\one_{I_s}(x)}{|I_s|}  \Big)^{1 \over 2} ,
\end{equation}
which allows us to interpret the left hand side of \eqref{eq:size:c:wavep} as $\ds \frac{1}{|I_{\tr{T}}|^\frac{1}{2}} \| \ii S_{\tr{T}}(f_4)\|_2$. By subadditivity, this is bounded by 
\[
\frac{1}{|I_{\tr{T}}|^\frac{1}{2}} \| \ii S_{\tr{T} ^{(1)}}(f_4)\|_2 + \frac{1}{|I_{\tr{T}}|^\frac{1}{2}} \| \ii S_{ \tr{T} ^{(2)}}(f_4)\|_2.
\]
Since $\ii S_{\tr{T}^{(1)}}(f_4)$ is supported on $E(I_{\tr{T}})$, we have for $2<p <\infty$
\begin{align*}
\frac{1}{|I_{\tr{T}}|^\frac{1}{2}} \| \ii S_{\tr{T}^{(1)}}(f_4)\|_2 &\lesssim \frac{1}{|I_{\tr{T}}|^\frac{1}{2}} \big| E(I_{\tr{T}}) \big|^{\frac{1}{2}-\frac{1}{p}} \| \ii S_{\tr{T} ^{(1)}}(f_4)\|_p \\
&\lesssim  2^{- \frac{\mu_0 m}{2} \big( \frac{1}{2}-\frac{1}{p}  \big)} \frac{1}{|I_{\tr{T}}|^\frac{1}{p}} \| \ii S_{\tr{T}^{(1)}}(f_4)\|_p.
\end{align*}

The wave packets in the definition of the square function \eqref{eq:def:square:func} are not associated to lacunary intervals only (indeed in order to apply the LGC method we cut each of them into $2^{m \over 2}$ parts), so the square function we are dealing with is not a classical (Littlewood-Paley type) square function, but a Rubio de Francia-type function \cite{RF}, which is still bounded from $L^q$ to $L^q$, provided that $2\leq q<\infty$. The local estimates that we employ here can be found in \cite{Cphdthesis}.

This proves precisely \eqref{eq:size:c:wavep} in the case $2<p<\infty$, where the summation runs over $\tr{T}^{(1)}$ instead of $\tr{T}$. In what concerns $\tr{T}^{(2)}$, we know that all the tiles $\pp \in \tr{T}^{(2)}$ satisfy 
\begin{equation}
\label{eq:small:aver:tree:clus}
\frac{1}{|I_{\pp}|} \int_{\R} |f_4(x)|^2 \cdot \chi_{I_{\pp}}(x) dx \leq 2^{- \frac{\mu_0 m}{2}} \frac{1}{|I_{\tr{T}}|} \int_{\R} |f_4(x)|^2 \chi_{I_{\tr{T}}}(x)dx. 
\end{equation}

We denote by $\ii I_{\text{max}}^{\tr{T}}$ the collection of such spatial intervals $I_{\pp}$ which are maximal with respect to inclusion. We note that for any fixed $J \in \ii I_{\text{max}}^{\tr{T}}$ the collection
\[
\{ \pp \in \tr{T}^{(2)} : I_{\pp} \subset J \}
\] 
is itself a tree. Then 
\begin{align*}
\| \ii S_{\tr{T} ^{(2)}}(f_4)\|_2^2 & = \sum_{J \in \ii I_{\text{max}}^{\tr{T}}} \sum_{\substack{\pp \in \tr{T}^{(2)} \\ I_{\pp} \subset J}} \sum_{s \in \mathfrak{S}^{\mathfrak{c}}_4(\pp_4)} | \langle f_4, \phi_{s}  \rangle|^2 \lesssim \sum_{J \in \ii I_{\text{max}}^{\tr{T}}} \| f_4 \cdot \chi_{J}  \|_2^2 \\
&\lesssim \sum_{J \in \ii I_{\text{max}}^{\tr{T}}}  |J| \cdot 2^{- \frac{\mu_0 m}{2}} \frac{1}{|I_{\tr{T}}|} \int_{\R} |f_4(x)|^2 \chi_{I_{\tr{T}}}(x)dx,
\end{align*}
where we used the $L^2$ boundedness of the square function and the fact that the intervals in $\ii I_{\text{max}}^{\tr{T}}$ satisfy \eqref{eq:small:aver:tree:clus}. Since the intervals in $\ii I_{\text{max}}^{\tr{T}}$ are maximal, 
\[
 \sum_{J \in \ii I_{\text{max}}^{\tr{T}}}  |J| \lesssim \big|I_{\tr{T}}   \big|
\] 
and thus we obtain that 
\[
\frac{1}{|I_{\tr{T}}|^{1 \over 2}}\| \ii S_{\tr{T} ^{(2)}}(f_4)\|_2 \lesssim 2^{- \frac{\mu_0 m}{4}} \frac{1}{|I_{\tr{T}}|^{1 \over 2}} \Big( \int_{\R} |f_4(x)|^2 \chi_{I_{\tr{T}}}(x)dx \Big)^{1 \over 2}
\]
which of course implies the desired \eqref{eq:size:c:wavep} (with $\tr{T}$ replaced by $\tr{T}^{(2)}$ on the left hand side).
\end{proof}

\subsection{Extension of the decay outside of the local$-L^2$ range}\label{extrange1}~\\
So far we have shown the boundedness of $T_m^{1}=\sum_{k\geq m/4} T_{k,m}$, which was introduced in Subsection \ref{subsec:prell}, \ref{case:c1}  with an exponential decay in the local-$L^2$ range -- see \eqref{eq:T1-2}. Now, in order to complete the proof of Theorem \ref{thm:case1} for the $k\geq m/4$ part, it remains to extend this decay in the Banach range. Using multilinear interpolation, it suffices to prove the following estimate which presents a logarithmic loss with respect to the decaying factor $2^{-\delta \, m}$:

\begin{proposition} \label{prop:outL2-I} For $p_1,p_2,p_3,p'\in (1,\infty)$ such that $\frac{1}{p}=\frac{1}{p_1}+\frac{1}{p_2}+\frac{1}{p_3}$ then
$$ \| T_m^1 \|_{L^{p_1}\times L^{p_2} \times L^{p_3} \to L^p} \lesssim (1+|m|)^{10}.$$
\end{proposition}

\begin{proof}
For this, we return to the $4$-linear form $\Lambda_{m}^{1}$ associated to $T_m^{1}$, which writes as
\begin{align*}
\sum_{k \geq {m \over 4}} \int_{\xi_1+\ldots+\xi_4=0} \prod_{j=1}^4 \hat f_j(\xi_j) \psi\Big(\frac{\xi_1-\xi_2}{2^{k+m}}\Big)\,\psi\Big(\frac{\xi_3-\xi_2}{2^{3k+m}}\Big)   \Big(\int_{\R}e^{-i\,\frac{\xi_1-\xi_2}{2^k}\,t}\,e^{-i\,\frac{\xi_3-\xi_2}{2^{3k}}\,t^3}\,\rho(t)\,dt\Big)\, d \xi.
\end{align*}

Unlike in Sections \ref{sec:onescale_k_large} and \ref{sec:bilinear:analysis}, here we do not exploit the potential decay produced by the oscillatory phase. Instead, we regard $\Lambda_{m}^{1}$ as
\begin{align*}
\int_{\R} \rho(t)\, \sum_{k \geq {m \over 4}} \left( \int_{\xi_1+\ldots+\xi_4=0} \prod_{j=1}^4 \hat f_j(\xi_j) \psi\Big(\frac{\xi_1-\xi_2}{2^{k+m}}\Big) e^{-i\,\frac{\xi_1-\xi_2}{2^{k+m}}\,2^m \,t}   \,\psi\Big(\frac{\xi_3-\xi_2}{2^{3k+m}}\Big) e^{-i\,\frac{\xi_3-\xi_2}{2^{3k+m}}\, 2^m t^3}\, d \xi \right) \, dt,
\end{align*}  
which allows to fix the parameter $t\simeq 1$.

Then for a fixed $t$, we can again use the approximation 
$$ \psi\Big(\frac{\xi_3-\xi_2}{2^{3k+m}}\Big) e^{-i\,\frac{\xi_3-\xi_2}{2^{3k+m}}\, 2^m t^3} \simeq \psi\Big(\frac{\xi_3-\frac{\xi_1+\xi_2}{2}}{2^{3k+m}}\Big) e^{-i\,\frac{\xi_3-\frac{\xi_1+\xi_2}{2}}{2^{3k+m}}\, 2^m t^3}$$
with a remainder term of order $|\xi_1-\xi_2|\cdot 2^{-3k} \simeq 2^{-2k-m}$ which is easily summable since we have decay both in $k,m$ with the uniform bounds for $T_{k,m}$ (see Proposition \ref{prop:Tk}).

So we are reduced to the study of
\begin{align*}
\int_{\R} \rho(t)\, \sum_{k \geq {m \over 4}} \left( \int_{\xi_1+\ldots+\xi_4=0} \prod_{j=1}^4 \hat f_j(\xi_j) \psi\Big(\frac{\xi_1-\xi_2}{2^{k+m}}\Big) e^{-i\,\frac{\xi_1-\xi_2}{2^{k+m}}\,2^m \,t}   \,\psi\Big(\frac{\xi_3-\frac{\xi_1+\xi_2}{2}}{2^{3k+m}}\right) e^{-i\,\frac{\xi_3-\frac{\xi_1+\xi_2}{2}}{2^{3k+m}}\, 2^m t^3}\, d \xi \Big) \, dt,
\end{align*}  
and such $4$-linear form can be decoupled (as explained in Proposition \ref{prop:facile}), making appear the product of two bilinear square functions (one acting on $f_1,f_2$ and another one on $f_3,f_4$). More precisely, in this situation, the square functions will be ``shifted'' with the shifts $\nu_1:=[2^m t]$ and $\nu_2:=[2^m t^3]$, both of order $2^m$. It is then known---see \emph{e.g.} \cite{Palsson}, \cite{BeckerDurcik} as well as related phenomena in \cite{Lie18} and \cite{HsuL26}---that we can incorporate the shifts in the BHT time-frequency analysis involving at most a logarithmic control of order $[1+\log(\nu_1)+\log(\nu_2)]^{10}$, which produces a loss of $(1+|m|)^{10}$. All the estimates are uniform in $t$ (because we are in the Banach range) and so we may conclude by integrating over $t\simeq 1$.
\end{proof}

\section{The multi-scale analysis (II) : Treatment of the case $0 \leq k \leq \frac{m}{8}$}
\label{sec:bilinear:analysis2}
The purpose of this section is to provide bounds for a fixed $k\geq 0$ of
\[
\Lambda_k^3(f_1,f_2,f_3,f_4)= \sum_{m\geq 8k} \ \sum_{\el,j \in {\mathbb Z} \atop \r \sim 2^{2k} } \L_{k,m}^{\el,\r,j}(f_1,f_2,f_3,f_4).
\]

As a consequence of Theorems \ref{sparsecont1} and \ref{sgscale1unifksmal}, we have that there exist some $\delta>0$ such that for any $0<\mu \ll \delta$, $\big| \L_{k,m}^{\el,\r,j}(f_1,f_2,f_3,f_4) \big|$ is bounded above by
\begin{align*}
 &  |I_j^k|^{-1} \|f_1^{\bbS^1_\mu}\|_{L^2(I_j^k \times W_{(\mathfrak{r}+2)+2^{2k} \mathfrak l}^{k+m})} \, \|f_2^{\bbS^1_\mu}\|_{L^2(I_j^k \times W_{\mathfrak{r}+2^{2k} \mathfrak l}^{k+m})} \, \|f_3\|_{L^2(I_j^k \times W_{ \mathfrak l+2}^{3k+m})} \|f_4\|_{L^2(I_j^k \times W_{  - (3 \mathfrak l+3)}^{3k+m})} \\
+ & |I_j^k|^{-1} \|f_1\|_{L^2(I_j^k \times W_{(\mathfrak{r}+2)+2^{2k} \mathfrak l}^{k+m})} \, \|f_2\|_{L^2(I_j^k \times W_{\mathfrak{r}+2^{2k} \mathfrak l}^{k+m})} \, \|f_3^{\bbS^1_\mu}\|_{L^2(I_j^k \times W_{ \mathfrak l+2}^{3k+m})} \|f_4^{\bbS^1_\mu}\|_{L^2(I_j^k \times W_{ - (3 \mathfrak l+3)}^{3k+m})} \\
+& 2^{-\delta k}  \,  |I_j^k|^{-1} \|f_1\|_{L^2(I_j^k \times W_{(\mathfrak{r}+2)+2^{2k} \mathfrak l}^{k+m})} \, \|f_2\|_{L^2(I_j^k \times W_{\mathfrak{r}+2^{2k} \mathfrak l}^{k+m})} \, \|f_3\|_{L^2(I_j^k \times W_{ \mathfrak l+2}^{3k+m})} \|f_4\|_{L^2(I_j^k \times W_{  - (3 \mathfrak l+3)}^{3k+m})}.
\end{align*}
We recall that the notion of ``{\it sparseness}'' has been defined in terms of parameters $\mu,\nu$ which have to be small enough and then $\delta$ can be fixed, small enough with respect to them.

As previously in Section \ref{sec:bilinear:analysis}, we apply Cauchy-Schwarz in the $\r$ parameter and then we move to the ``tiles'' formalism introduced there. Since $k \geq 0$ is fixed, all the tiles appearing have spatial interval of fixed length, and the frequency intervals are of length $2^{3k+m}$, as $m$ ranges over values $\geq 8k$. We still denote $\P$ the collection of multi-tiles $P=(P_1,P_2,P_3,P_4)$ with $P_i=I_P\times \omega_{P_i}$ where $I_P$ is a spatial interval of length $2^{-k}$ and the $\omega_{P_i}$ are described by \eqref{eq:def:spatial:intv} and \eqref{eq:def:freq:intvs}.

Then our study of $\Lambda_k^3(f_1,f_2,f_3,f_4)$ reduces to analyzing $\Lambda_{\P}^{\mathfrak{u}}(f_1,f_2,f_3,f_4)$ and  $\Lambda_{\P}^{\mathfrak{c}}(f_1,f_2,f_3,f_4)$ which resemble the corresponding quantities in the previous section and which are defined bellow (see \eqref{eq:lpu} and \eqref{eq:lpc}). Indeed Theorems \ref{sparsecont1} and \ref{sgscale1unifksmal} imply the existence of $\delta>0$ so that 
\begin{equation*} 
\Lambda_k^3(f_1,f_2,f_3,f_4) \lesssim 2^{-k\delta} \Lambda_{\P}^{\mathfrak{u}}(f_1,f_2,f_3,f_4) + \Lambda_{\P}^{\mathfrak{c}}(f_1,f_2,f_3,f_4),
\end{equation*}
where 
\begin{itemize}[leftmargin=1em]
\item the uniform component $\Lambda_{\P}^{\mathfrak{u}}$ is defined by
\begin{equation}
\Lambda_{\P}^{\mathfrak{u}}(f_1,f_2,f_3,f_4) := \sum_{P\in\P} |I_P|^{-1} f_1(P_1) \cdot f_2(P_2) \cdot f_3(P_3) \cdot f_4(\pp_4) \label{eq:lpu}
\end{equation}
with $f_i(P_i)$ defined by  $\displaystyle f_i(P_i):= \Big( \sum_{s_i \in \mathfrak{S}_i(P_i)} \big| \langle f_i, \phi_{s_i}  \rangle  \big|^2 \Big)^{1 \over 2}$, where for $i \in \{ 1, 2 \}$
\[
\mathfrak{S}_i(P_i):=\big\lbrace    s=I \times \omega : I \subset I_{\pp}, |I|= |I_{\pp}| \cdot 2^{-{m \over 2}}, \quad  \omega \subset \omega_{\pp_i}, |\omega|= 2^{{m \over 6}} |\omega_{\pp_i}|^{1 \over 3}  \big\rbrace
\]
and for $i \in \{ 3, 4 \}$
\[
\mathfrak{S}_i(P_i):=\big\lbrace    s=I \times \omega : I \subset I_{\pp}, |I|= |I_{\pp}|^3 \cdot 2^{-{m \over 2}}, \quad  \omega \subset \omega_{\pp_i}, |\omega|= 2^{-{m \over 2}} |\omega_{\pp_i}|  \big\rbrace.
\]
Notice that there sets of unitary tiles are different from the ones used in the previous section.

\item the sparse component $\Lambda_{\P}^{\mathfrak{c}}$ corresponds to the definition in \eqref{sparse23} when the functions $f_1$ and $f_2$ or the functions $f_3$ and $f_4$ are simultaneously made up of a ``sparse'' subcollection of unitary wave-packets. More precisely, for any multi-tile $P \in \P$ we consider
\begin{itemize}[leftmargin=1.5em]
\item for $i \in \{1,2 \}$,
\[ 
\ii I^{\mu}(P_i):= \Big\lbrace S \text{ dyadic interval} : |S|=2^{-3k}, S \subset I_P \text{  and  }  \frac{1}{|S|}\int_{3 S} \left|f_i\right|^2  \gtrsim \frac{2^{\mu k}}{|I_P|} \int_{3 I_P} \left|f_i \right|^2     \Big\rbrace;
\]
\item and for $i \in \{3,4\}$,
\[ 
\ii I^{\mu}(P_i):= \Big\lbrace S \text{ dyadic interval} : |S|=2^{-5k}, S \subset I_P \text{  and  }  \frac{1}{|S|}\int_{3 S} \left|f_i\right|^2  \gtrsim \frac{2^{\mu k}}{|I_P|} \int_{3 I_P} \left|f_i \right|^2     \Big\rbrace.
\]
\end{itemize}
Moreover we set for $i \in \{1, \ldots, 4 \}$
\[
E_i^\mu(I_P):=\cup\{ S: S \in  \ii I^{\mu}(P_i)\}
\]
and notice that $\displaystyle E_i^\mu(I_P) \subset \big\{ x\in I_P,\ M(|f_i|^2 \cdot \chi_{I_{\tr T}})(x) > 2^{\mu k} \aver{3I_P}|f_i|^2 dy \big\}$ and so $\ds| E_i^\mu(I_P) | \lesssim 2^{-\mu k} |I_P|$. We also  define
\[
\mathfrak{S}^{\mathfrak c}_i(P_i):=\big\lbrace    s=I \times \omega \in \mathfrak{S}_i(P_i) : \exists \, S \in  \ii I^{\mu}(P_i) \text{ such that } I \subset S \big\rbrace.
\]
and the corresponding quantities $f_i^{\mathfrak c}(P_i)$ by  $\displaystyle f_i^{\mathfrak c}(P_i):= \Big( \sum_{s_i \in \mathfrak{S}_i^{\mathfrak c}(P_i)} \big| \langle f_i, \phi_{s_i}  \rangle  \big|^2 \Big)^{1 \over 2}$.
Then the sparse component of the $4$-linear form from \eqref{sparse23} can be controlled by 
\begin{align}
\Lambda_{\P}^{\mathfrak{c}}(f_1,f_2,f_3,f_4) & := \sum_{P\in\P} |I_P|^{-1} \Big( f_1^{\mathfrak c}(P_1) \cdot f_2^{\mathfrak c}(P_2) \cdot f_3(P_3) \cdot f_4(P_4) \nonumber \\
& \qquad + f_1(P_1) \cdot f_2(P_2) \cdot f_3^{\mathfrak c}(P_3) \cdot f_4^{\mathfrak c}(P_4)\Big). \label{eq:lpc}
\end{align}
\end{itemize}

Our objective is to prove the following result:
\begin{theorem}
\label{thm:multi:scale:II}
Let $k\geq 0$ and let $\P$ be a collection of $4$-tiles as defined in \eqref{eq:def:spatial:intv} and \eqref{eq:def:freq:intvs}. Then for all $2 < p_1, p_2, p_3, p_4 < \infty$ with $\frac{1}{p_1}+\frac{1}{p_2}+\frac{1}{p_3}+\frac{1}{p_4}=1$ we have\footnote{For the \emph{uniform} component, some of the $p_j$ can be equal to $\infty$.} that
\begin{equation}
\label{eq:wanted:res:u:II}
\big|\Lambda^{\mathfrak{u}}_{\P}(f_1, f_2, f_3, f_4)\big|  \lesssim \|f_1\|_{p_1}    \|f_2\|_{p_2} \|f_3\|_{p_3} \|f_4\|_{p_4}        
\end{equation}
and there exists $\delta:=\delta(p_1,p_2,p_3,p_4)$ with
\begin{equation}
\label{eq:wanted:res:c:II}
\big|\Lambda^{\mathfrak{c}}_{\P}(f_1, f_2, f_3, f_4)\big|  \lesssim 2^{-\delta k} \|f_1\|_{p_1} \|f_2\|_{p_2} \|f_3\|_{p_3} \|f_4\|_{p_4}.   
\end{equation}
Moreover, these estimates are uniform with respect to the collection $\P$ and the parameter $k$.
\end{theorem}

\begin{proof}
We first focus on the uniform part.
By using H\"older inequality along the spatial intervals $I_0$ of fixed length $2^{-k}$, it suffices to prove that
\begin{equation}
\label{eq:llaprouver:II}
\Lambda^{\mathfrak{u}}_{\P,I_0}(f_1, f_2, f_3, f_4) \lesssim |I_0|^{-1} \|f_1 \cdot \chi_{I_0}\|_{2} \|f_2\cdot \chi_{I_0}\|_{2} \|f_3\cdot \chi_{I_0}\|_{2} \|f_4\cdot \chi_{I_0}\|_{2},
\end{equation}
where
$$ \Lambda^{\mathfrak{u}}_{\P,I_0}(f_1, f_2, f_3, f_4):=\sum_{P\in \P \atop I_P=I_0} |I_0|^{-1} f_1(P_1) \cdot f_2(P_2) \cdot f_3(P_3)\cdot f_4(P_4).$$

Before proving \eqref{eq:llaprouver:II}, we notice that the tiles contributing to the term $\Lambda^{\mathfrak{u}}_{\P,I_0}$ all have the same spatial interval $I_0$ of length $2^{-k}$, while in frequency the intervals $\omega_{P_i}$ are determined by two parameters, as one can see in the definition \eqref{eq:def:freq:intvs}: the length $2^{3k+m}$ and the position given by the parameter $\el$. We define for $m \geq 8k$
\[
\PP_{I_0}^m:=\{ P=(P_1, P_2, P_3, P_4) \in \PP : I_P=I_0 \text{  and  } |\omega_{P_i}|= 2^{3k+m} \}.
\]
As $P$ varies in $\PP_{I_0}^m$, the intervals $\omega_{P_i}$ are disjoint and of same length. One can also notice that, as $P$ ranges in $\bigcup_{m \geq 8k} \PP_{I_0}^m$, the squares $ \omega_{P_1} \times  \omega_{P_3}$ are mutually disjoint and they cover a subregion of $\{\xi_1<\xi_3\}$. Similarly, $\{  \omega_{P_2} \times  \omega_{P_4} : P \in \bigcup_{m \geq 8k} \PP_{I_0}^m \}$ represents a sub-covering of a similar region in frequency. This last observation is going to play a key role in our proof, as we will soon see.

Thus we estimate 
\begin{align*}
\Lambda^{\mathfrak{u}}_{\P,I_0}(f_1, f_2, f_3, f_4) & = |I_0|^{-1} \sum_{m \geq 8k} \sum_{P \in \PP_{I_0}^m} f_1(P_1) \cdot f_2(P_2) \cdot f_3(P_3)\cdot f_4(P_4) \\
&\lesssim |I_0|^{-1} \Big( \sum_{m \geq 8k} \sum_{P \in \PP_{I_0}^m} f_1(P_1)^2 f_3(P_3)^2  \Big)^{1 \over 2} \, \Big( \sum_{m \geq 8k} \sum_{P \in \PP_{I_0}^m} f_2(P_2)^2 f_4(P_4)^2  \Big)^{1 \over 2}.
\end{align*}

Invoking Remark \ref{remark:equiv:l2:tiles}, we notice that
\[
f_1(P_1)^2 \lesssim \sum_{\#_1 \in \Z} (1+|\#_1|)^{-M} \sum_{\substack{J_1 \subset I_{0, \#_1} \\ |J_1|=2^{-5k}}} \sum_{\substack{\omega_1 \subset 3 \omega_{P_1} \\ |\omega_1|=2^{5k}}} \big| \langle  f_1, \phi_{J_1 \times \omega_1}   \rangle  \big|^2. 
\]
Similarly, 
\[
f_3(P_3)^2 \lesssim \sum_{\#_3 \in \Z} (1+|\#_3|)^{-M} \sum_{\substack{J_3 \subset I_{0, \#_3} \\ |J_3|=2^{-5k}}} \sum_{\substack{\omega_3 \subset 3 \omega_{P_3} \\ |\omega_3|=2^{5k}}} \big| \langle  f_3, \phi_{J_3 \times \omega_3}   \rangle  \big|^2. 
\]
Thus 
\begin{align*}
\sum_{m \geq 8k} \sum_{P \in \PP_{I_0}^m} &f_1(P_1)^2 f_3(P_3)^2  \lesssim  \sum_{\#_1, \#_3 \in \Z} (1+|\#_1|)^{-M} (1+|\#_3|)^{-M} \cdots \\
& \cdots \sum_{\substack{J_1 \subset I_{0, \#_1} \\ |J_1|=2^{-5k}}} \sum_{\substack{J_3 \subset I_{0, \#_3} \\ |J_3|=2^{-5k}}}  \sum_{m \geq 8k} \sum_{P \in \PP_{I_0}^m} \sum_{\substack{\omega_1 \subset 3 \omega_{P_1} \\ |\omega_1|=2^{5k}}}  \sum_{\substack{\omega_3 \subset 3 \omega_{P_3} \\ |\omega_3|=2^{5k}}}\big| \langle  f_1, \phi_{J_1 \times \omega_1}   \rangle  \big|^2 \,  \big| \langle  f_3, \phi_{J_3 \times \omega_3}   \rangle  \big|^2.
\end{align*}
The collection $\{ \omega_1 \times \omega_3 \}$ appearing above is a collection of squares of fixed sidelength $2^{5k}$ contained inside the region $\{\xi_1<\xi_3\}$. Thus we can tensorize the summation in $\omega_1$ and $\omega_3$ to bound the term above by 
\begin{align*}
 \sum_{\#_1, \#_3 \in \Z} (1+|\#_1|)^{-M} (1+|\#_3|)^{-M}& \Big( \sum_{\substack{J_1 \subset I_{0, \#_1} \\ |J_1|=2^{-5k}}}  \sum_{\substack{\omega_1 : |\omega_1|=2^{5k}}} \big| \langle  f_1, \phi_{J_1 \times \omega_1}   \rangle  \big|^2  \Big) \cdot \Big( \sum_{\substack{J_3 \subset I_{0, \#_3} \\ |J_3|=2^{-5k}}}  \sum_{\substack{\omega_3 : |\omega_3|=2^{5k}}} \big| \langle  f_3, \phi_{J_2 \times \omega_2}   \rangle  \big|^2  \Big) \\
& \lesssim \|f_1 \cdot \chi_{I_0}\|_2^2 \cdot \|f_3 \cdot \chi_{I_0}\|_2^2.
\end{align*}

Treating the term grouping the $f_2$ and $f_4$ information in the same way, we deduce \eqref{eq:llaprouver:II}.

The analysis for the sparse component is similar, but we apply Remark \ref{remark:equiv:l2:tiles} to space intervals contained in $\ii I^{\mu}(P_i)$, after noticing that this collection only depends on the fixed space interval $I_0$. More precisely, if $f_1(P_1)$ is replaced by $f_1^{\mathfrak c}(P_1)$, then
\begin{align*}
f_1^{\mathfrak c}(P_1)^2=\sum_{S \in \ii I^{\mu}_1(I_0)} \sum_{\substack{s=I \times \omega \in \mathfrak{S}_1(P_1) \\ I \subset S}} |\langle f_1, \phi_s   \rangle|^2
\end{align*}
and by Remark \ref{remark:equiv:l2:tiles} we get that
\begin{align*}
f_1^{\mathfrak c}(P_1)^2 \lesssim \sum_{\#_1 \in \Z} (1+|\#_1|)^{-M} \sum_{S \in \ii I_1^{\mu}(I_0)} \sum_{\substack{J_1 \subset S_{\#_1} \\ |J_1|=2^{-5k}}} \sum_{\substack{\omega_1 \subset 3 \omega_{P_1} \\ |\omega_1|=2^{5k}}} \big| \langle  f_1, \phi_{J_1 \times \omega_1}   \rangle  \big|^2. 
\end{align*}
Together with the previous computations, we deduce that $\ds \sum_{m \geq 8k} \sum_{P \in \PP_{I_0}^m} f_1^{\mathfrak c}(P_1)^2 f_3(P_3)^2$ is bounded above by
\begin{align*}
\sum_{\#_1, \#_3 \in \Z} (1+|\#_1|)^{-M} (1+|\#_3|)^{-M} \Big( \sum_{S \in \ii I_1^{\mu}(I_0)} \sum_{\substack{J_1 \subset S_{\#_1} \\ |J_1|=2^{-5k}}}  \sum_{\substack{\omega_1 : |\omega_1|=2^{5k}}} \big| \langle  f_1, \phi_{J_1 \times \omega_1}   \rangle  \big|^2  \Big)& \\
 \cdot \Big( \sum_{\substack{J_3 \subset I_{0, \#_3} \\ |J_3|=2^{-5k}}}  \sum_{\substack{\omega_3 : |\omega_2|=2^{5k}}} \big| \langle  f_3, \phi_{J_3 \times \omega_3}   \rangle  \big|^2  \Big) &.
\end{align*}
While the term in $f_3$ can be bounded as before by $\| f_3 \cdot \chi_{I_0} \|_2^2$, we regard the term in $f_1$ as a Rubio de Francia square function localized to the set 
\[
E_{1, \#_1}^\mu(I_0):= \bigcup_{S \in \ii I_1^{\mu}(I_0)} S_{\#_1} ,
\]
which has the same length as $E_{1}^\mu(I_0)$. Indeed,
\begin{align*}
 \sum_{S \in \ii I_1^{\mu}(I_0)} \sum_{\substack{J_1 \subset S_{ \#_1} \\ |J_1|=2^{-5k}}}  \sum_{\substack{\omega_1 : |\omega_1|=2^{5k}}} \big| \langle  f_1, \phi_{J_1 \times \omega_1}   \rangle  \big|^2 = \Big\| \Big( \sum_{\substack{I \subset E_{1, \#_1}^\mu(I_0) \\ |I|=2^{-5k}  }} \sum_{|\omega_1|=2^{5k}} \frac{|\langle f_1, \phi_{I \times \omega_1}\rangle|^2}{|I|} \cdot \one_I   \Big)^{1 \over 2} \cdot \one_{E_{1, \#_1}^\mu(I_0)}   \Big\|_2^2.
\end{align*}
Now we use the fact that $|E_{1, \#_1}^\mu(I_0)|= |E_{1}^\mu(I_0)| \lesssim 2^{-\mu k} |I_0|$ and the $L^{p_1}$ boundedness (for $2 \leq p_1 <\infty$) of the Rubio de Francia-type square function in order to deduce that the above term is bounded by 
\[
\|f_1 \cdot \chi_{I_0, 2^{-2k}\#_1 }\|_{p_1}^2 \big( 2^{-\mu k} |I_0| \big)^{1-\frac{2}{p_1}}.
\]

The case when $f_3$ is spatially sparse can be treated in a similar way. In order to tensorize the information for $f_1$ and $f_3$, we express $f_3^{\mathfrak c}(P_3)$ using wave packets adapted to unitary tiles of dimension $2^{-5k}\times 2^{5k}$. We notice that the space intervals of these tiles have the same length as the intervals defining the sparse collection $\ii I^{\mu}(P_3)$, so the computations are somewhat simpler.

Finally, we can then run the same computations as before in the \emph{uniform} case, and obtain the global boundedness of the \emph{sparse} $4$-linear form $\Lambda_{\P}^{\mathfrak{c}}$  with an extra decay $2^{-\delta k}$ with $\delta:=\min_i\{\mu(\frac{1}{2}-\frac{1}{p_i})\}>0$, as desired. 
\end{proof}

\subsection{Extension of the decay outside of the local$-L^2$ range}\label{extrange2}~\\
In this section, we aim to conclude the proof of Theorem \ref{thm:case1}, for the regime $k\leq \frac{m}{8}$. We have just obtained the expected estimates for $T_k^{3}=\sum_{m \geq 8k} T_{k,m}$ in the local-$L^2$ range with an exponential decay $2^{-\nu k}$ and we now want to extend this result to the Banach range. As in the previous section, that can be done by interpolation with the following:

\begin{proposition} \label{prop:quasiBanach1} Let $k \geq 1$. Let $1<p_1,p_2,p_3 < \infty$ and $1 \leq p < \infty$ be Lebesgue exponents with $\frac{1}{p_1}+\frac{1}{p_2}+\frac{1}{p_3}=\frac{1}{p}$. Then one has, uniformly in $k$,
\begin{equation} \| \sum_{m\geq 8k} T_{k,m}\|_{L^{p_1} \times L^{p_2} \times L^{p_3} \rightarrow L^p} \lesssim 1. \label{eq:fullam1}
\end{equation}
\end{proposition}



\begin{proof} 
We recall that the operator $T_{k,m}$ is given by
\begin{align*}
 T_{k,m}(f_1,f_2,f_3)(x)= & \int_{\R} \rho(t)  \iiint_{{\mathbb R}^3} e^{ix(\xi_1+\xi_2+\xi_3)} \hat{f_1}(\xi_1) \hat{f_2}(\xi_2) \hat{f_3}(\xi_3)   \\
 & \qquad \psi\big( \frac{\xi_1-\xi_2}{2^{k+m}} \big) e^{-i \frac{\xi_1-\xi_2}{2^{k}} t}  \psi \big( \frac{\xi_3-\xi_2}{2^{3k+m}}\big)  e^{-i \frac{\xi_3-\xi_2}{2^{3k}} t^3} d \xi_1 d\xi_2 d\xi_3  dt.
\end{align*}
First, since we are in the Banach setting, it suffices to study the $3$-linear operator $T_{k,m,t}$ with a fixed parameter $t$ and to obtain uniform estimates in $t$. Then, since $k$ is also fixed, by considering the modulated functions 
$$ \tilde f_1(x)=f_1(x-2^{-k}t), \quad \tilde f_2(x)=f_2(x+2^{-k} t + 2^{-3k} t^3) \quad \textrm{and} \quad \tilde f_3(x)=f_3(x-2^{-3k} t^3)$$
which have the same Lebesgue norms as the initial functions, we see that it suffices to prove the boundedness of $\sum_{m\geq 8k} \tilde T_{k,m}$ with now
\begin{align*}
 \tilde T_{k,m}(f_1,f_2,f_3)(x)= & \iiint_{{\mathbb R}^3} e^{ix(\xi_1+\xi_2+\xi_3)} \hat{f_1}(\xi_1) \hat{f_2}(\xi_2) \hat{f_3}(\xi_3) \psi\big( \frac{\xi_1-\xi_2}{2^{k+m}} \big)   \psi \big( \frac{\xi_3-\xi_2}{2^{3k+m}}\big) d \xi_1 d\xi_2 d\xi_3.
\end{align*}
Notice that the operator $\sum_{m\geq 8k} \tilde T_{k,m}$ has symbol
\begin{equation}
\label{eq:symbol:decay:m>8k}
\sum_{m\geq 8k} \psi\big( \frac{\xi_1-\xi_2}{2^{k+m}} \big) \psi \big( \frac{\xi_3-\xi_2}{2^{3k+m}}\big), 
\end{equation}
which fails to satisfy the conditions in \cite{Mlinop}, for the case $k=1,n=4$, by a factor taking the form of a power of $2^{k}$.

Let $N_1$ be a sufficiently large integer. We perform a Taylor expansion of order $N_1$ for $ \psi \big( \frac{\xi_3-\xi_2}{2^{3k+m}}\big)$ in the spirit of \eqref{eq:tril:form:Taylor}
\begin{align*}
 \psi\big( \frac{\xi_3-\xi_2}{2^{3k+m}}\big) &=\sum_{\ell=0}^{N_1}\frac{1}{\ell! \, 2^\ell}  \big( \frac{\xi_1-\xi_2}{2^{3k+m}}  \big)^\ell  \psi\big( \frac{\xi_3-\frac{\xi_1+\xi_2}{2}}{2^{3k+m}}\big) \\
 & +\frac{1}{N_1!\, 2^{N_1}} \big(  \frac{\xi_1-\xi_2}{2^{3k+m}}   \big)^{N_1+1}\Big(\int_0^1 (1-\tau)^{N_1} \psi^{(N_1+1)} \big( \frac{\xi_3}{2^{3k+m}}- \frac{1-\tau}{2} \, \frac{\xi_1}{2^{3k+m}} - \frac{1+\tau}{2} \, \frac{\xi_2}{2^{3k+m}}  \big) d\, \tau \Big).
\end{align*}
With this, the symbol in \eqref{eq:symbol:decay:m>8k} splits into a main term given by 
\begin{equation}
\label{eq:symbol:decay:m>8k:main}
\sum_{\ell=0}^{N_1}\frac{1}{\ell! \, 2^\ell}\, \frac{1}{2^{2 k \ell}} \, \sum_{m\geq 8k} \big( \frac{\xi_1-\xi_2}{2^{k+m}}  \big)^\ell \psi\big( \frac{\xi_1-\xi_2}{2^{k+m}} \big)  \psi\big( \frac{\xi_3-\frac{\xi_1+\xi_2}{2}}{2^{3k+m}}\big)
\end{equation}
and a remainder given by 
\begin{equation}
\label{eq:symbol:decay:m>8k:remainder}
\frac{2^{-2 (N_1+1) k}}{N_1!\, 2^{N_1+1}} \sum_{m\geq 8k}  \big( \frac{\xi_1-\xi_2}{2^{k+m}}  \big)^{N_1+1} \psi\big( \frac{\xi_1-\xi_2}{2^{k+m}} \big)\Big( \int_0^1 (1-\tau)^{N_1} \psi^{(N_1+1)} \big( \frac{\xi_3}{2^{3k+m}}- \frac{1-\tau}{2} \, \frac{\xi_1}{2^{3k+m}} - \frac{1+\tau}{2} \, \frac{\xi_2}{2^{3k+m}}  \big) d\, \tau \Big).   
\end{equation}
The key observation is that, provided $N_1$ is sufficiently large, the remainder term is singular along the subspace $\Gamma':= (1,1,1){\mathbb R}$ and satisfies the conditions  thanks to the extra decay $2^{-2 N_1 k}$. Hence it gives rise to a trilinear operator similar to those studied in \cite{Mlinop}. Thus in particular the piece of the operator which corresponds to the symbol given by \eqref{eq:symbol:decay:m>8k:remainder} is bounded in the full Banach range.
On the other hand, each term\footnote{In fact, for $\ell \geq 1$, the symbol has some extra decay, thus behaving better those in Proposition \ref{prop:facile}.} in \eqref{eq:symbol:decay:m>8k:main} can be treated similarly to the operator studied in Proposition \ref{prop:facile}, as detailed in equation \eqref{eq:est:square:func:Prop:facile}. Thus the ``main'' part of the operator is bounded inside the Banach range. So overall the operator $\sum_{m\geq 8k} \tilde T_{k,m}$ and consequently the operator $\sum_{m\geq 8k} T_{k,m}$ are bounded from $L^{p_1}\times L^{p_2}\times L^{p_3}$ to $L^p$ whenever $1 < p_1, p_2, p_3, p' < \infty$, with an operator norm independent of $k$.

\end{proof}

\section{The multi-scale analysis (III) : Treatment of the case $0 \leq \frac{m}{8}< k \leq \frac{m}{4}$}
\label{sec:bilinear:analysis3}

In this section, we aim to prove Theorem \ref{thm:case1} for the regime $\frac{m}{8}<k<\frac{m}{4}$ about the boundedness of $T_{\kappa}^2=\sum_{m\geq 4\kappa} T_{{m \over 4}-{\kappa \over 2},m}$. In Section \ref{Diagkinterm}, it has been proved the single-scale estimates (see Theorem \ref{sgscale1unifkintr}) and it has been shown that to estimate $ \Lambda_{k,m}^{\el,\r,j}$, we have some exponential decay with respect to the quantity $\frac{m}{2}-2k$ that we name $\kappa:=\frac{m}{2}-2k$ (and so $k=\frac{m}{4}-\frac{\kappa}{2}$) and observe that
$$ \sum_{k\geq 0} \ \sum_{m=4k}^{8k} \cdots = \sum_{\kappa \geq 0} \ \sum_{m\geq 4\kappa} \cdots.$$

Having the single-scale estimates, the multi-scale analysis for $T_\kappa^2$ is exactly the one that we have developped for $T_m^1$ (see Section \ref{sec:bilinear:analysis}).

Indeed the case (III) is a mixture/intermediate case between case (I) and case (II): the single-scale estimate is ressembling to the one of case (II) and the multi-scale analysis is analogous to the one of case (I).

\section*{AI disclosure}   The authors declare that no generative artificial intelligence (AI) or AI-assisted tools were used in the research reported in this paper or in the preparation of its scientific content. This includes the conception and development of the research, the derivation and verification of the mathematical results, and the writing, proofreading, and editing of the body of the manuscript. The authors take full responsibility for the accuracy of this declaration and for the integrity of the paper.

\end{document}